\documentclass[moor]{informs3}              % for a regular run
\usepackage[caption=false]{subfig}
\usepackage[outdir=./]{epstopdf}
\usepackage{mathtools,dsfont,extarrows,mathrsfs}
\usepackage{enumitem,verbatim,algorithm}
\usepackage{algpseudocode,booktabs}
\usepackage{environ}

\providecommand{\cB}{\mathcal{B}}
\providecommand{\scrB}{\mathscr{B}}
\providecommand{\cP}{\mathcal{P}}
\providecommand{\R}{\mathbb{R}}

\providecommand{\E}{\mathbb{E}}
\providecommand{\calU}{\mathcal{U}}

\providecommand{\calL}{\mathcal{L}}

\providecommand{\frakM}{\mathfrak{M}}
\providecommand{\frakA}{\mathfrak{A}}
\providecommand{\frakB}{\mathfrak{B}}
\providecommand{\supp}{\mathrm{supp}\ }
\providecommand{\interior}{\mathrm{int}}

\providecommand{\hxi}{\widehat{\xi}}
\providecommand{\heta}{\widehat{\eta}}

\providecommand{\VaR}{\mathrm{VaR}}

\renewcommand{\Pr}{\mathbb{P}}

\providecommand{\Dom}{\mathrm{dom}}
\providecommand{\bp}{\mu}
\providecommand{\bq}{\nu}

\providecommand{\Dp}{\overline{D}}
\providecommand{\Dm}{\underline{D}}

\providecommand{\hfillqed}{\hfill\Halmos}

\DeclareMathOperator*{\esssupOp}{-ess\,sup}

\newcommand\esssup[2]{\underset{#2}{#1\!\esssupOp}}

\newcommand{\ignore}[1]{}
\def\defi{\vcentcolon=}%

\usepackage{natbib}
 \NatBibNumeric
 \def\bibfont{\small}%
 \bibpunct[, ]{[}{]}{,}{n}{}{,}%

\usepackage[colorlinks=true,breaklinks=true,bookmarks=true,urlcolor=blue,
     citecolor=blue,linkcolor=blue,bookmarksopen=false,draft=false]{hyperref}

\def\EMAIL#1{\href{mailto:#1}{#1}}% When hyperref is used, otherwise outcomment
\def\URL#1{\href{#1}{#1}}         % When hyperref is used, otherwise outcomment

\TheoremsNumberedThrough     % Preferred (Theorem 1, Lemma 1, Theorem 2)
\EquationsNumberedThrough    % Default: (1), (2), ...
\renewcommand{\theARTICLETOP}{}
\renewcommand{\theLRHFirstLine}{}
\renewcommand{\theRRHFirstLine}{}
\renewcommand{\theLRHSecondLine}{}
\renewcommand{\theRRHSecondLine}{}
\newcommand{\qedhere}{\tag*{$\square$}}

\makeatletter
\gdef\th@TH{%
  \def\theorem@headerfont{\normalfont\TheoremHeaderFont}%
  \theorempreskipamount=6pt plus 3pt minus 2pt \theorempostskipamount=6pt plus 3pt minus 2pt
  \if@MOOR\labelsep1em\else\labelsep0.5em\fi
  \def\@begintheorem##1##2{\normalfont\TheoremTextFont
        \item[\if@OPRE\else\hspace*{0em}\fi\hskip\labelsep
          \theorem@headerfont ##1\ ##2.]}%
  \def\@opargbegintheorem##1##2##3{\normalfont\TheoremTextFont
        \item[\if@OPRE\else\hspace*{0em}\fi\hskip\labelsep
          \theorem@headerfont ##1\ ##2\ {\bf(##3)}.]}}

\gdef\th@EX{%
  \def\theorem@headerfont{\normalfont\ExampleHeaderFont}%
  \if@OPRE
  \theorempreskipamount=6pt plus 3pt minus 2pt \theorempostskipamount=6pt plus 3pt minus 2pt
  \else
  \theorempreskipamount=6pt plus 3pt minus 2pt \theorempostskipamount=6pt plus 3pt minus 2pt
  \fi
  \if@MOOR\labelsep1em\else\labelsep0.5em\fi
  \def\@begintheorem##1##2{\normalfont\ExampleTextFont
        \item[\if@OPRE\else\hspace*{0em}\fi\hskip\labelsep \theorem@headerfont ##1\ ##2.]}%
  \def\@opargbegintheorem##1##2##3{\normalfont\ExampleTextFont
        \item[\if@OPRE\else\hspace*{0em}\fi\hskip\labelsep \theorem@headerfont ##1\ ##2\ {(##3)}.]}}

  \def\proof#1{\Trivlist\item[\hskip\labelsep{\it #1\enskip }]\ignorespaces}
  \def\endproof{\endTrivlist\addvspace{1.5em}}
\makeatother

\SingleSpacedXI
\begin{document}
%%%%%%%%%%%%%%%%

%\setcounter{page}{1}
%\DOI{10.1287/xxxx.0000.0000}%


%\RUNTITLE{}

% \TITLE{Bundling Information Goods of Decreasing Value}
\TITLE{Distributionally Robust Stochastic Optimization with Wasserstein Distance}

\ARTICLEAUTHORS{%
\AUTHOR{Rui Gao}
\AFF{Department of Information, Risk and Operations Management, University of Texas at Austin, Austin, TX 78705, \EMAIL{rui.gao@mccombs.utexas.edu}}
\AUTHOR{Anton J. Kleywegt}
\AFF{H.~Milton Stewart School of Industrial and Systems Engineering, Georgia Institute of Technology, Atlanta, GA 30332, \EMAIL{anton@isye.gatech.edu} \URL{}}
} % end of the block

\ABSTRACT{%
Distributionally robust stochastic optimization (DRSO) is an approach to optimization under uncertainty in which, instead of assuming that there is a known true underlying probability distribution, one hedges against a chosen set of distributions.
%It seeks a decision which hedges against the worst-case distribution in an ambiguity set, containing a family of distributions relevant to the considered problem.
In this paper we first point out that the set of distributions should be chosen to be appropriate for the application at hand, and that some of the choices that have been popular until recently are, for many applications, not good choices.
%Unfortunately, the worst-case distributions resulting from many widely used ambiguity sets are sometimes either unclear or unrealistic.
We next consider sets of distributions that are within a chosen Wasserstein distance from a nominal distribution.
Such a choice of sets has two advantages: (1)~The resulting distributions hedged against are more reasonable than those resulting from other popular choices of sets. (2)~The problem of determining the worst-case expectation over the resulting set of distributions has desirable tractability properties.
We derive a strong duality reformulation of the corresponding DRSO problem and construct approximate worst-case distributions (or an exact worst-case distribution if it exists) explicitly via the first-order optimality conditions of the dual problem.

Our contributions are four-fold.
(i) We identify necessary and sufficient conditions for the existence of a worst-case distribution, which are naturally related to the growth rate of the objective function.
(ii) We show that the worst-case distributions resulting from an appropriate Wasserstein distance have a concise structure and a clear interpretation.
(iii) Using this structure, we show that data-driven DRSO problems can be approximated to any accuracy by robust optimization problems, and thereby many DRSO problems become tractable by using tools from robust optimization.
(iv) Our strong duality result holds in a very general setting. As examples, we show that it can be applied to infinite dimensional process control  and intensity estimation for point processes.
}%

%\KEYWORDS{deterministic inventory theory; infinite linear programming duality;
%  existence of optimal policies; semi-Markov decision process; cyclic schedule}
%\MSCCLASS{Primary: 90B05; secondary: 90C40, 90C90}
%\ORMSCLASS{Primary: Inventory/production: deterministic multi-item;
%  secondary: dynamic programming/optimal control: deterministic
%  semi-Markov; programming: infinite dimensional}
%\HISTORY{Received November 20, 2003; revised March 8, 2004, and March 26, 2004.}

\MSCCLASS{Primary: 90C15; secondary: 90C46}
\ORMSCLASS{Primary: programming: stochastic}

\maketitle
%%%%%%%%%%%%%%%%%%%%%%%%%%%%%%%%%%%%%%%%%%%%%%%%%%%%%%%%%%%%%%%%%%%%%%
% Samples of sectioning (and labeling) in MOOR.
% NOTE: (1) all section levels end with a period,
%       (2) capitalization is as shown (sentence style, not title style).
%
%\section{Introduction.}\label{intro} %%1.
%\subsection{Duality and the classical EOQ problem.}\label{class-EOQ} %% 1.1.
%\subsection{Outline.}\label{outline1} %% 1.2.
%\subsubsection{Cyclic schedules for the general deterministic SMDP.}
%  \label{cyclic-schedules} %% 1.2.1
% Text of your paper here
\section{Introduction}
\label{sec:intro}

In decision making problems under uncertainty, a decision maker wants to choose a decision $x$ from a feasible region $X$.
The objective function $\Psi : X \times \Xi \mapsto \R$ depends on a quantity $\xi \in \Xi$ whose value is not known to the decision maker at the time that the decision has to be made.
In some settings it is reasonable to assume that $\xi$ is a random element with distribution $\mu$ supported on $\Xi$, for example, if multiple realizations of $\xi$ will be encountered.
In such settings, the decision making problems can be formulated as \textsl{stochastic optimization} problems as follows:
\[
\inf_{x \in X} \E_{\mu}[\Psi(x,\xi)].
\]
We refer to \citet{shapiro2014lectures} for a thorough study of stochastic optimization.
One major criticism of the formulation above for practical applications is the requirement that the underlying distribution $\mu$ be known to the decision maker.
Even if multiple realizations of $\xi$ are observed, $\mu$ still may not be known exactly, while use of a distribution different from $\mu$ may sometimes result in bad decisions.
Another major criticism is that in many applications there are not multiple realizations of $\xi$ that will be encountered, for example in problems involving events that may either happen once or not happen at all, and thus the notion of a ``true'' underlying distribution does not apply.
These criticisms motivate the notion of \textsl{distributionally robust stochastic optimization (DRSO)}, that does not rely on the notion of a known true underlying distribution.
One chooses a set $\frakM$ of probability distributions to hedge against,
and then finds a decision that provides the best hedge against the set $\frakM$ of distributions by solving the following minmax problem:
\begin{equation}
\label{eqn:DRSO}
\inf_{x \in X} \sup_{\mu \in \frakM} \E_{\mu}[\Psi(x,\xi)].
\tag{DRSO}
\end{equation}
Such a minmax approach has its roots in Von Neumann's game theory and has been used in many fields such as inventory management (\citet{scarf1958min,gallego1993distribution}), statistical decision analysis (\citet{10.2307/2723239}), as well as stochastic optimization (\citet{vzavckova1966minimax,dupavcova1987minimax,shapiro2002minimax}).
Recently it regained attention in the operations research literature, and sometimes is called data-driven stochastic optimization or ambiguous stochastic optimization.

A central question is: how to choose a good set of distributions $\frakM$ to hedge against?
A good choice of $\frakM$ should take into account the properties of the practical application as well as the tractability of problem \eqref{eqn:DRSO}.
Two typical ways of constructing $\frakM$ are moment-based and distance-based.
The moment-based approach considers distributions whose moments (such as mean and covariance) satisfy certain conditions (\citet{scarf1958min,delage2010distributionally,popescu2007robust,zymler2013distributionally}).
It has been shown that in many cases the resulting DRSO problem can be formulated as a conic quadratic or semi-definite program.
However, the moment-based approach is based on the curious assumption that certain conditions on the moments are known exactly but that nothing else about the relevant distribution is known.
More often in applications, either one has data from repeated observations of the quantity $\xi$, or one has no data, and in both cases the moment conditions do not describe exactly what is known about $\xi$.
In addition, the resulting worst-case distributions sometimes yield overly conservative decisions (\citet{Wang2015Likelihood,goh2010distributionally}).
For example, \citet{Wang2015Likelihood} shows that for the newsvendor problem, by hedging against all the distributions with fixed mean and variance, Scarf's moment approach yields a two-point worst-case distribution, and the resulting decision does not perform well under other more likely scenarios.

The distance-based approach considers distributions that are close, in the sense of a chosen statistical distance, to a \textit{nominal distribution} $\nu$, such as an empirical distribution or a Gaussian distribution (\citet{ghaoui2003worst,calafiore2006distributionally}).
Popular choices of the statistical distance are $\phi$-divergences (\citet{doi:10.1287/educ.2015.0134,ben2013robust}), which include Kullback-Leibler divergence (\citet{Jiang2015Data-driven}), Burg entropy (\citet{Wang2015Likelihood}), and Total Variation distance (\citet{sun2015convergence}) as special cases, Prokhorov metric (\citet{erdougan2006ambiguous}), and Wasserstein distance (\citet{wozabal2012framework,wozabal2014robustifying,esfahani2015data,zhao2018data}).

\subsection{Motivation: Potential Issues with \texorpdfstring{$\phi$}{phi}-divergence}
\label{sec:motivation}

Despite its widespread use, $\phi$-divergence has a number of shortcomings.
Here we highlight some of these shortcomings.
In a typical setup using $\phi$-divergence, $\Xi$ is partitioned into $\bar{B}+1$ bins represented by points $\xi^{0},\xi^{1},\ldots,\xi^{\bar{B}} \in \Xi$.
The nominal distribution $\bq$ associates $N_{i}$ observations with bin~$i$. That is, the nominal distribution is given by $\bq \defi (N_{0}/N, N_{1}/N, \ldots, N_{\bar{B}}/N)$, where $N \defi \sum_{i=0}^{\bar{B}} N_{i}$.
Let $\Delta_{\bar{B}} \defi \{(p_{0},p_{1},\ldots,p_{\bar{B}}) \in \R_{+}^{{\bar{B}}+1} \, : \, \sum_{j=0}^{\bar{B}} p_{j} = 1\}$ denote the set of probability distributions on the same set of bins.
Let $\phi : [0,\infty) \mapsto \R$ be a chosen convex function such that $\phi(1) = 0$, with the conventions that $0 \phi(a/0) \defi a \lim_{t \to \infty} \phi(t) / t$ for all $a > 0$, and $0 \phi(0/0) \defi 0$.
Then the $\phi$-divergence between $\bp = (p_{0},\ldots,p_{B}),\bq = (q_{0},\ldots,q_{B}) \in \Delta_{\bar{B}}$ is defined by
\[
\label{eqn:phi-divergence}
I_{\phi}(\bp,\bq) \ \ \defi \ \ \sum_{j=0}^{\bar{B}} q_{j} \phi\left(\frac{p_{j}}{q_{j}}\right).
\]
Let $\theta > 0$ denote a chosen radius.
Then $\frakM_{\phi} \defi \left\{\bp \in \Delta_{\bar{B}} \, : \, I_{\phi}(\bp,\bq) \leq \theta\right\}$
denotes the set of probability distributions given by the chosen $\phi$-divergence and radius $\theta$.
The DRSO problem corresponding to the $\phi$-divergence ball $\frakM_{\phi}$ is then given by
\[
\label{eqn:DRSO_phi}
\inf_{x \in X} \sup_{\bp \in \Delta_{\bar{B}}} \left\{\sum_{j=0}^{\bar{B}}  p_{j} \Psi(x,\xi^{j}) \; : \;  I_{\phi}(\bp,\bq) \leq \theta\right\}.
\]
It has been shown in \citet{ben2013robust} that the $\phi$-divergence ball $\frakM_{\phi}$ can be viewed as a statistical confidence region (\citet{pardo2005statistical}), and for several choices of $\phi$, the inner maximization of the problem above is tractable.

One well-known shortcoming of $\phi$-divergence balls is that, either they are not rich enough to contain distributions that are often relevant, or they they hedge against many distributions that are too extreme. For example, for some choices of $\phi$-divergence such as Kullback-Leibler divergence, if the nominal $q_{i} = 0$, then $p_{i} = 0$, that is, the $\phi$-divergence ball includes only distributions that are absolutely continuous with respect to the nominal distribution $\bq$, and thus does not include distributions with support on points where the nominal distribution $\bq$ is not supported.
As a result, if $\Xi = \R^{K}$ and $\bq$ is discrete, then there are no continuous distributions in the $\phi$-divergence ball $\frakM_{\phi}$.
Some other choices of $\phi$-divergence exhibit in some sense the opposite behavior.
For example, the Burg entropy ball includes distributions with some amount of probability allowed to be shifted from $\bq$ to any other bin, with the amount of probability allowed to be shifted depending only on $\theta$ and not on how extreme the bin is.
See Section~\ref{sec:newsvendor} for more details regarding this potential shortcoming.

Next we illustrate another shortcoming of $\phi$-divergence that will motivate the use of Wasserstein distance.

\begin{figure}
\subfloat[ Observed image with histogram $\bq$]{
\label{fig:bird_observe} %% label for first subfigure
\begin{minipage}[b]{0.3\linewidth}
\centering \includegraphics[width=\linewidth]{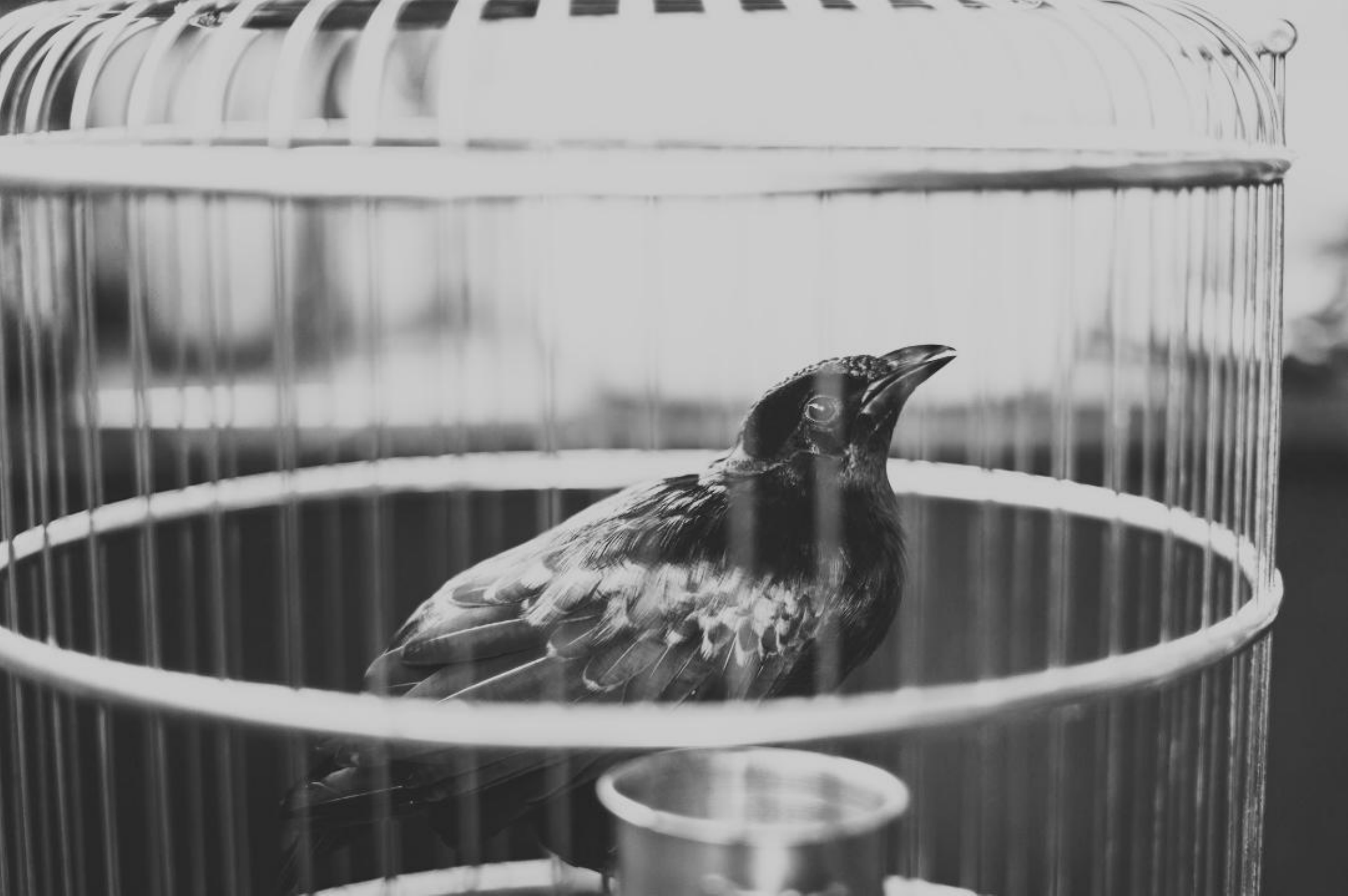}\\
\vspace{6pt}
\includegraphics[height=1.5in]{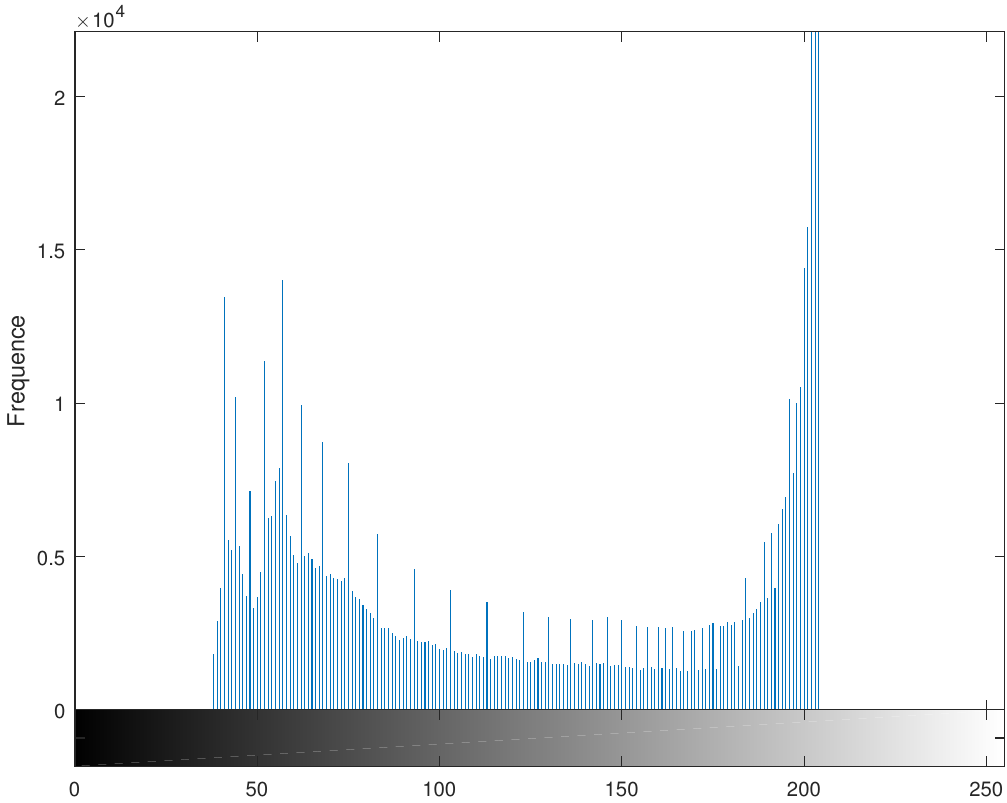}
\end{minipage}}%
\hfill
\subfloat[True image with histogram $\bp_{true}$]{
\label{fig:bird_true} %% label for second subfigure
\begin{minipage}[b]{0.3\linewidth}
\centering \includegraphics[width=\linewidth]{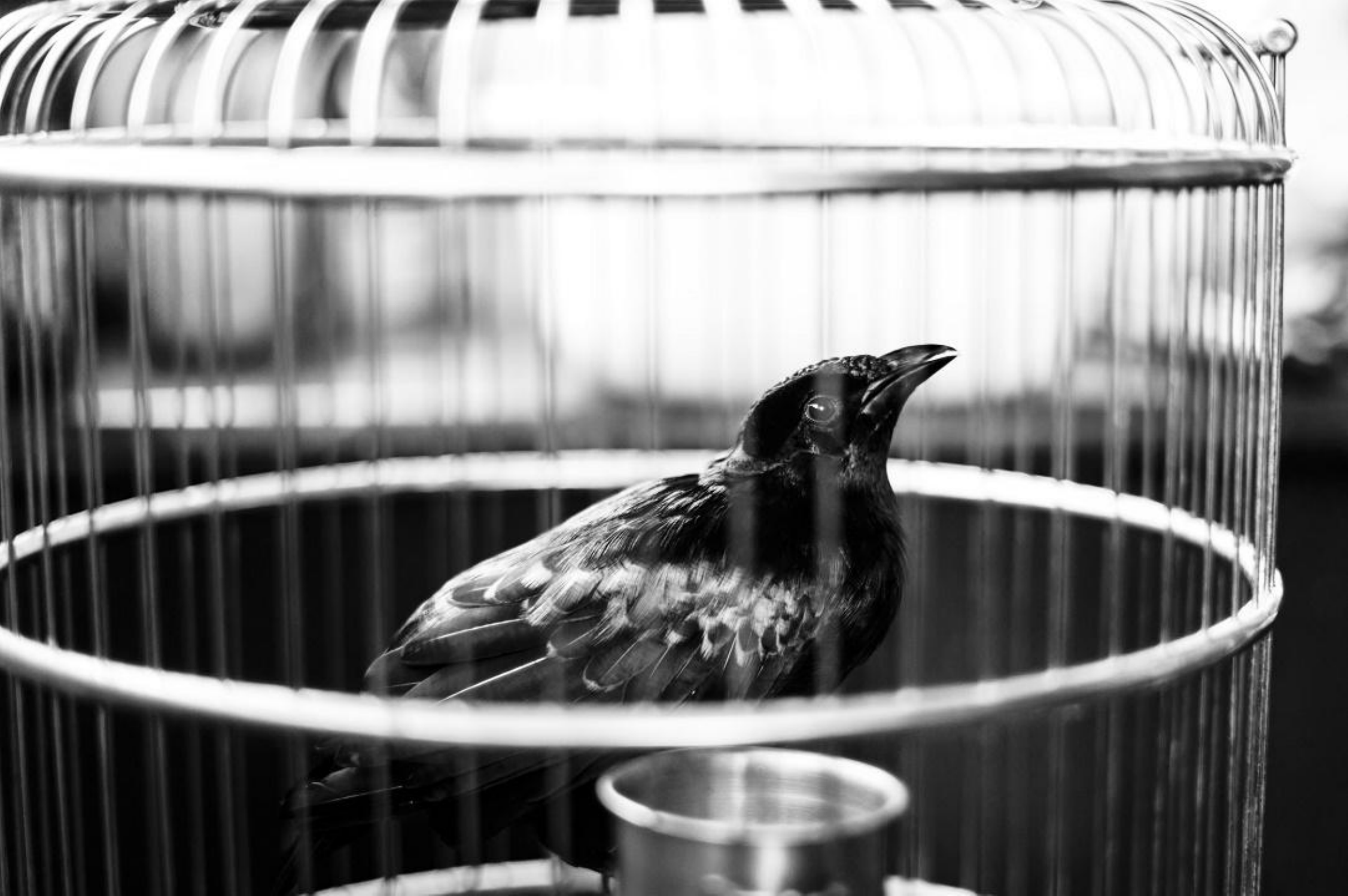}\\\vspace{6pt}
\includegraphics[height=1.5in]{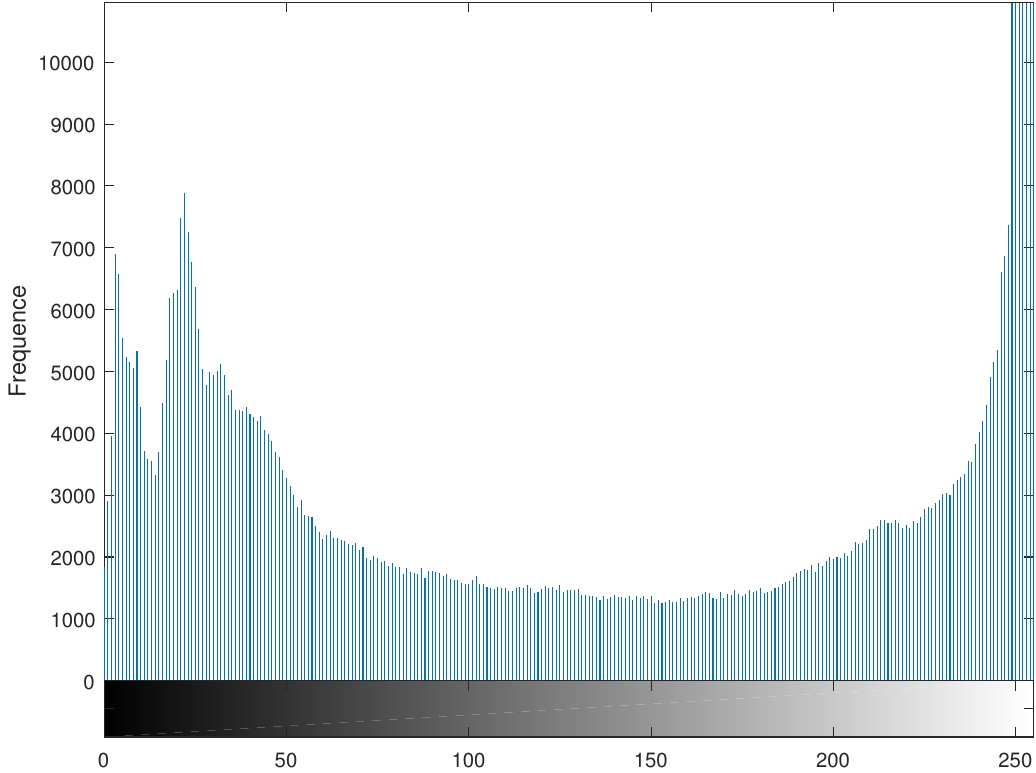}
\end{minipage}}
\hfill
\subfloat[Pathological image with histogram $\bp_{pathol}$]{
\label{fig:bird_pathol} %% label for second subfigure
\begin{minipage}[b]{0.3\linewidth}
\centering \includegraphics[width=\linewidth]{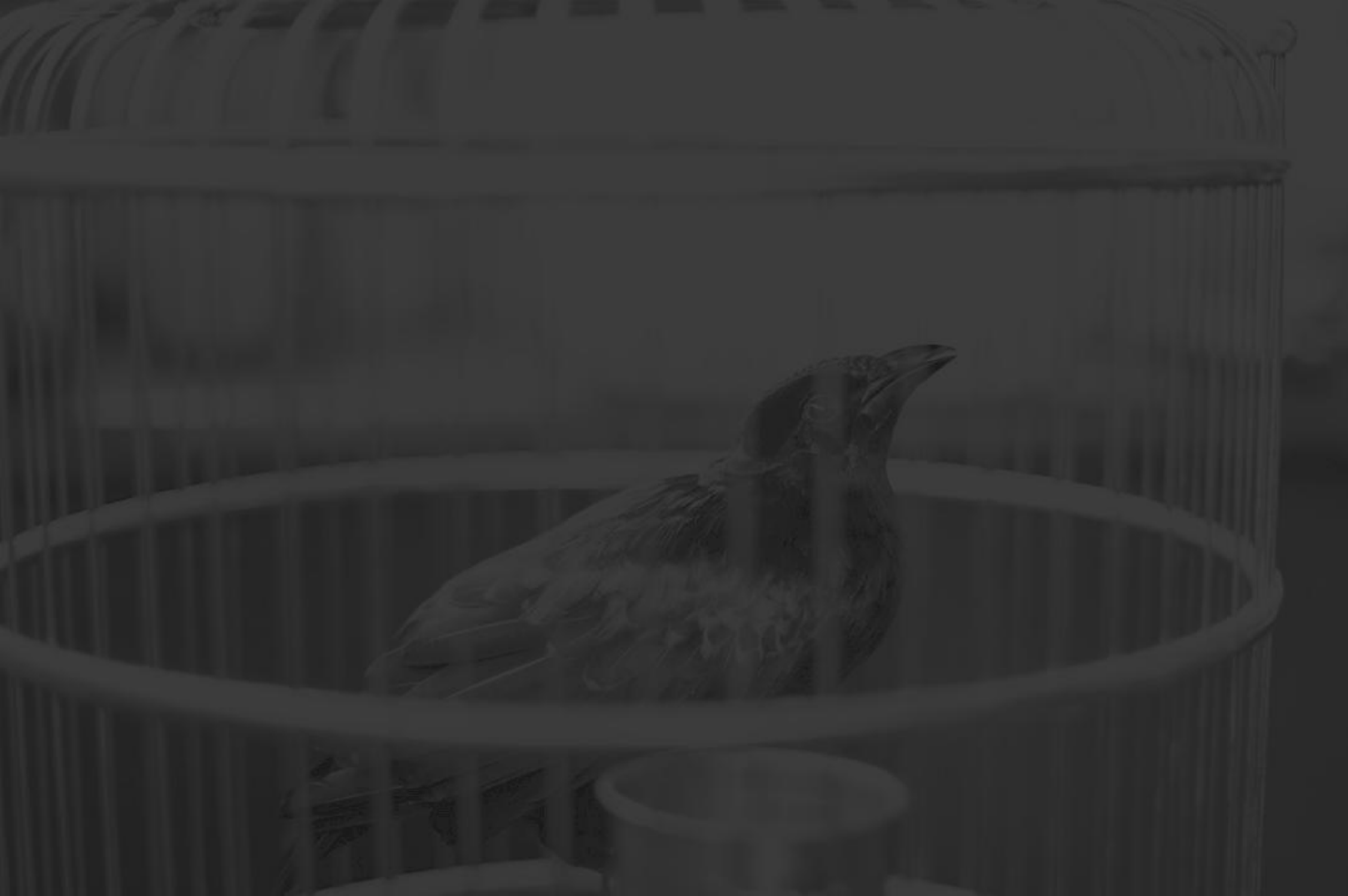}\\\vspace{6pt}
\includegraphics[height=1.5in]{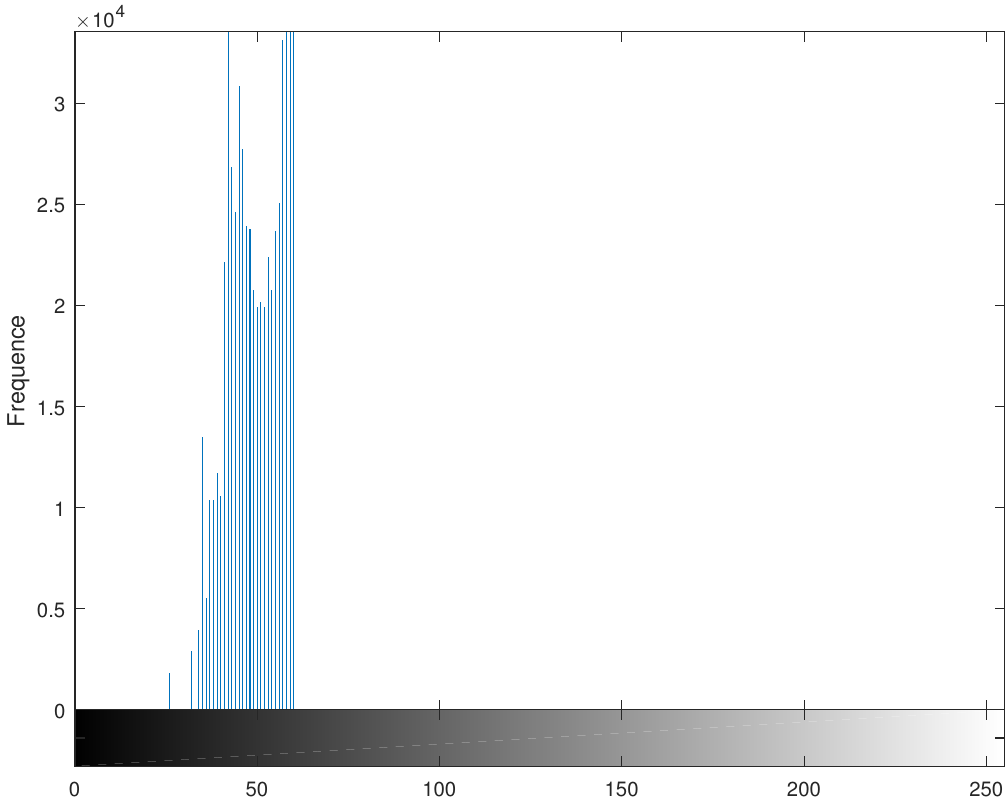}
\end{minipage}}
\caption{Three images and their gray-scale histograms. For KL divergence, it holds that $I_{\phi_{KL}}(\bp_{true},\bq) = 5.05 > I_{\phi_{KL}}(\bp_{pathol},\bq) = 2.33$, while in contrast, Wasserstein distance satisfies $W_{1}(\bp_{true},\bq) = 30.70 < W_{1}(\bp_{pathol},\bq) = 84.03$.}
\label{fig:flower} %% label for entire figure
\end{figure}

\begin{example}
\label{eg:bird}
Suppose that there is an underlying true image (\ref{fig:bird_true}), and a decision maker possesses, instead of the true image, an approximate image (\ref{fig:bird_observe}) obtained with a less than perfect device that loses some of the contrast.
The images are summarized by their gray-scale histograms.
(In fact, (\ref{fig:bird_observe}) was obtained from (\ref{fig:bird_true}) by a low-contrast intensity transformation (\citet{Gonzalez:2006:DIP:1076432}), by which the black pixels become somewhat whiter and the white pixels become somewhat blacker.
This type of transformation changes only the gray-scale value of a pixel and not the location of a gray-scale value, and therefore it can also be regarded as a transformation from one gray-scale histogram to another gray-scale histogram.)
As a result, roughly speaking, the observed histogram $\bq$ is obtained by shifting the true histogram $\bp_{true}$ inwards.
Also consider the pathological image (\ref{fig:bird_pathol}) that is too dark to see many details, with histogram $\bp_{pathol}$.
Suppose that the decision maker constructs a Kullback-Leibler (KL) divergence ball $\frakM_{\phi_{KL}} \defi \{\bp \in \Delta_{\bar{B}} \, : \, I_{\phi_{KL}}(\bp,\bq) \leq \theta\}$.
Note that $I_{\phi_{KL}}(\bp_{true},\bq) = 5.05 > I_{\phi_{KL}}(\bp_{pathol},\bq) = 2.33$.
Therefore, if $\theta$ is chosen small enough (less than $2.33$) for $\frakM_{\phi_{KL}}$ to exclude the pathological image~(\ref{fig:bird_pathol}), then $\frakM_{\phi_{KL}}$ will also exclude the true image~(\ref{fig:bird_true}).
If $\theta$ is chosen large enough (greater than $5.05$) for $\frakM_{\phi_{KL}}$ to include the true image~(\ref{fig:bird_true}), then $\frakM_{\phi_{KL}}$ also has to include the pathological image~(\ref{fig:bird_pathol}), and then the resulting decision may be overly conservative due to hedging against irrelevant distributions.
If an intermediate value is chosen for $\theta$ (between $2.33$ and $5.05$), then $\frakM_{\phi_{KL}}$ includes the pathological image~(\ref{fig:bird_pathol}) and excludes the true image~(\ref{fig:bird_true}).
In contrast, note that the Wasserstein distance $W_{1}$ satisfies $W_{1}(\bp_{true},\bq) = 30.7 < W_{1}(\bp_{pathol},\bq) = 84.0$, and thus Wasserstein distance does not exhibit the problem encountered with KL divergence (see also Example~\ref{eg:bird_w}).
\end{example}

The reason for such behavior is that $\phi$-divergence does not incorporate a notion of how close two points $\xi,\xi' \in \Xi$ are to each other, for example, how likely it is that observation is $\xi'$ given that the true value is $\xi$.
In Example~\ref{eg:bird}, $\Xi = \{0,1,\ldots,255\}$ represents 8-bit gray-scale values. In this case, we know that the likelihood that a pixel with gray-scale value $\xi \in \Xi$ is observed with gray-scale value $\xi' \in \Xi$ is decreasing in the absolute difference between $\xi$ and $\xi'$. However, in the definition of $\phi$-divergence, only the \textit{relative ratio} $p_{j}/q_{j}$ for the same gray-scale value $j$ is taken into account, while the distances between different gray-scale values is not taken into account. This phenomenon has been observed in studies of image retrieval (\citet{rubner2000earth,ling2007efficient}).

The drawbacks of $\phi$-divergence motivates us to consider sets $\frakM$ that incorporate a notion of how close two points $\xi,\xi' \in \Xi$ are to each other.
One such choice of $\frakM$ is based on Wasserstein distance.
Specifically, consider any underlying metric $d$ on $\Xi$ which measures the closeness of any two points in $\Xi$. Let $p \geq 1$, and let $\cP(\Xi)$ denote the set of Borel probability measures on $\Xi$.
Then the \textit{Wasserstein distance} of order $p$ between two distributions $\mu,\nu \in \cP(\Xi)$ is defined as
\[
W_{p}(\mu,\nu) \ \ \defi \ \ \min_{\gamma \in \cP(\Xi^2)} \; \Big\{\E_{(\xi,\zeta) \sim \gamma}^{1/p}[d^{p}(\xi,\zeta)] \; : \; \gamma \textrm{ has marginal distributions } \mu,\nu\Big\}.
\]
More detailed explanation and discussion on Wasserstein distance will be presented in Section~\ref{sec:wasserstein}.
Given a radius $\theta > 0$, the Wasserstein ball of probability distributions $\frakM$ is defined by
\[
\label{eqn:frakM}
\frakM \ \ \defi \ \ \{\mu \in \cP(\Xi) \; : \; W_{p}(\mu,\nu) \leq \theta\}.
\]

\subsection{Related Work}
\label{sec:related}

Wasserstein distance and the related field of \textit{optimal transport}, which is a generalization of the transportation problem, have been studied in depth.
In 1942, together with the linear programming problem (\citet{kantorovich1960mathematical}), \citet{kantorovich1942translocation} tackled Monge's problem originally brought up in the study of optimal transport.
In the stochastic optimization literature, Wasserstein distance has been used for single stage stochastic optimization (\citet{wozabal2012framework,wozabal2014robustifying}), and for multistage stochastic optimization (\citet{pflug2014multistage}).
The challenge for solving \eqref{eqn:DRSO} is that, the inner maximization involves a supremum over possibly an infinite dimensional space of distributions.
To tackle this problem, existing works focus on the setup when $\nu$ is the empirical distribution on a finite-dimensional space.
\citet{wozabal2012framework} transformed the inner maximization problem of \eqref{eqn:DRSO} into a finite-dimensional non-convex program, by using the fact that if $\nu$ is supported on at most $N$ points, then there are extreme distributions of $\frakM$ that are supported on at most $N+3$ points.
Recently, using duality theory of conic linear programming (\citet{shapiro2001duality}), \citet{esfahani2015data} and \citet{zhao2018data} showed that under certain conditions, the inner maximization problem of \eqref{eqn:DRSO} is actually equivalent to a finite-dimensional convex problem.

In this paper, we consider any arbitrary nominal distribution $\nu$ on a Polish space, and study the tractability of \eqref{eqn:DRSO} via strong duality.
By the time we completed the first version of this paper, we learned that \citet{blanchet2019quantifying} independently considered a similar problem with more general lower semi-continuous transport cost function and study its strong duality.
Our focus and our approach to this problem differ from theirs in several important ways.
First, we show that the strong duality holds for any measurable function $\Psi$, while \citet{blanchet2019quantifying} proves the result only when $\Psi$ is upper semi-continuous. The upper semi-continuity is used to ensure the existence of the worst-case distribution, but is not needed for the strong duality to hold.
Second, we prove the strong duality result for the inner maximization of \eqref{eqn:DRSO} using a novel, yet simple, \textit{constructive} approach, in contrast with the non-constructive approaches in their work and also in \citet{esfahani2015data} and \citet{zhao2018data}. This enables us to establish the structural characterization of the worst-case distributions of the data-driven DRSO (Corollary~\ref{cor:finite}\ref{itm:finiteDual_N+1}), which improves the result of \citet{wozabal2012framework} and the more recent result of \citet{owhadi1504extreme} on extremal distributions of Wasserstein balls (Remark~\ref{rmk:extreme}).  It also enables us to build a close connection between DRSO and robust optimization (Corollary~\ref{cor:finite}\ref{itm:finiteDual_robustapprox}).
Third, we focus on Wasserstein distance of order $p$ ($p \geq 1$), while they consider more general transport cost functions in which the distance between two points $\xi,\xi' \in \Xi$ is measured by a lower semi-continuous function rather than a metric $d^{p}(\xi,\xi')$ as in our case. Nevertheless, our proof remains valid for such transport cost functions, and in fact, for even more general cost functions that are not necessarily lower semi-continuous (Remark \ref{rmk:cost}). In the meantime, focusing on Wasserstein distance enables us to relate the condition for the existence of a worst-case distribution to the important notion of the ``growth rate'' of the objective function, and enables us to provide practical guidance for choosing the ambiguity set and controlling the degree of conservativeness based on the objective function (Remark~\ref{rmk:p}).

\subsection{Main Contributions}
\label{sec:contribution}

\begin{itemize}[itemsep=0.5em,topsep=0.5em]
\item
\textsl{General Setting.}
We prove a strong duality result for DRSO problems with Wasserstein distance in a very general setting.
We show that
\[
\sup_{\mu \in \cP(\Xi)} \big\{\E_{\mu}[\Psi(x,\xi)] \; : \; W_{p}(\mu,\nu) \leq \theta\big\}
\ \ = \ \ \min_{\lambda \geq 0} \left\{\lambda \theta^{p} - \int_{\Xi} \inf_{\xi \in \Xi} [\lambda d^{p}(\xi,\zeta) - \Psi(x,\xi)] \nu(d\zeta)\right\}
\]
holds for any Polish space $(\Xi,d)$ and measurable function $\Psi$ (Theorem~\ref{thm:strongDual}).
\begin{enumerate}
\item
Both \citet{esfahani2015data} and \citet{zhao2018data} assume that $\Xi$ is a convex subset of $\R^{K}$ with some associated norm.
The greater generality of our results enables one to consider interesting problems such as the process control problems in Sections~\ref{sec:process} and~\ref{sec:NHPP}, where $\Xi$ is the set of finite counting measures on $[0,1]$, which is infinite-dimensional and non-convex.

\item
Both \citet{esfahani2015data} and \citet{zhao2018data} assume that the nominal distribution $\nu$ is an empirical distribution, while we allow $\nu$ to be any Borel probability measure. %The greater generality enables one to study problems such as the worst-case Value-at-Risk analysis in Section~\ref{sec:wcVaR}.

\item
Both \citet{esfahani2015data} and \citet{zhao2018data} only consider Wasserstein distance of order $p=1$.
By considering a bigger family of Wasserstein distances, we establish the importance for DRSO problems of the notion of the ``growth rate'' of the objective function, which measures how fast the objective function grows compared to a polynomial of order $p$.
It turns out that the growth rate of the objective function determines the finiteness of the worst-case objective value (Proposition~\ref{prop:kappa_infty}), and it plays an important role in the existence conditions for the worst-case distribution (Corollary~\ref{cor:existence}).
This is of practical importance, since it provides guidance for choosing the proper Wasserstein distance and for controlling the degree of conservativeness based on the structure of the objective function.
\end{enumerate}

\item
\textsl{Constructive Proof of Duality.}
We prove the strong duality result using a novel, elementary, constructive approach.
The results of \citet{esfahani2015data} and \citet{zhao2018data} and other strong duality results in the literature are based on the established Hahn-Banach theorem for certain infinite dimensional vector spaces.
% For example, in \cite{esfahani2015data}, the strong duality result is obtained using a general result on conic linear programming (\citet{shapiro2001duality}).
In contrast, our proof idea is new and is relatively elementary and straightforward: we use the weak duality
result as well as the first-order optimality condition of the dual problem to construct a sequence of primal
feasible solutions whose objective values converge to the dual optimal value.
Our proof uses relatively elementary tools, without resorting to other ``big hammers''.

\item
\textsl{Existence Conditions and the Structure of Worst-case Distributions.}
As a by product of our constructive proof, we identify necessary and sufficient conditions for the existence of worst-case distributions, and a structural characterization of worst-case distributions (Corollary~\ref{cor:existence}).
Specifically, for data-driven DRSO problems where $\nu = \frac{1}{N} \sum_{i=1}^{N} \delta_{\hxi^{i}}$ (where $\delta_{\xi}$ denotes the unit mass on $\xi$), whenever a worst-case distribution exists, there is a worst-case distribution $\mu^{\ast}$ supported on at most $N+1$ points with the following concise structure:
\[
\mu^{\ast} \ \ = \ \ \frac{1}{N} \sum_{\underset{i \neq i_{0}}{i = 1}}^{N} \delta_{\xi_{\ast}^{i}}
+ \frac{p_{0}}{N} \delta_{\underline{\xi}_{\ast}^{i_{0}}} + \frac{1 - p_{0}}{N} \delta_{\overline{\xi}_{\ast}^{i_{0}}},
\]
for some $i_{0} \in \{1,\ldots,N\}$, $p_{0} \in [0,1]$ and
\[
\xi_{\ast}^{i} \ \ \in \ \ \argmin_{\xi \in \Xi} \left\{\lambda^{\ast} d^{p}(\xi,\hxi^{i}) - \Psi(x,\xi)\right\}, ~\forall \ i \neq i_{0},
\ \ \ \ \underline{\xi}_{\ast}^{i_{0}},\overline{\xi}_{\ast}^{i_{0}} \in \argmin_{\xi \in \Xi} \left\{\lambda^{\ast} d^{p}(\xi,\hxi^{i_{0}}) - \Psi(x,\xi)\right\},
\]
where $\lambda^{\ast}$ is the dual minimizer (Corollary~\ref{cor:finite}).
Thus $\mu^{\ast}$ can be viewed as a perturbation of $\nu$, where the mass on $\hxi^{i}$ is perturbed to $\xi_{\ast}^{i}$ for all $i \neq i_{0}$, a fraction $p_{0}$ of the mass on $\hxi^{i_{0}}$ is perturbed to $\underline{\xi}_{\ast}^{i_{0}}$, and the remaining fraction $1 - p_{0}$ of the mass on $\hxi^{i_{0}}$ is perturbed to $\overline{\xi}_{\ast}^{i_{0}}$.
In particular, uncertainty quantification problems have a worst-case distribution with this simple structure, and can be solved by a greedy procedure (Example~\ref{eg:UQ}).
Our result regarding the existence of a worst-case distribution with such a structure improves the result of \citet{wozabal2012framework} and the more recent result of \citet{owhadi1504extreme} regarding the extremal distributions of Wasserstein balls.

\item
\textsl{Connection with Robust Optimization.}
Using the structure of a worst-case distribution, we prove that data-driven DRSO problems can be approximated by robust optimization problems to any accuracy (Corollary~\ref{cor:finite}\ref{itm:finiteDual_robustapprox}). We use this result to show that two-stage linear DRSO problems with linear decision rules have a tractable semi-definite programming approximation (Section~\ref{sec:two-stage}).
Moreover, the robust optimization approximation becomes exact when the objective function $\Psi$ is concave in $\xi$.
In addition, if $\Psi$ is convex in $x$, then the corresponding DRSO problem can be formulated as a convex-concave saddle point problem.
\end{itemize}

The rest of this paper is organized as follows.
In Section~\ref{sec:wasserstein}, we review some results on the Wasserstein distance.
Next we prove strong duality for general nominal distributions in Section~\ref{sec:dual_general}, and in Section~\ref{sec:dual_finite} we derive additional results for finite-supported nominal distributions.
Then, in Sections~\ref{sec:application} and~\ref{sec:discussion}, we apply our results on strong duality and the structural description of the worst-case distributions to a variety of DRSO problems.
We conclude this paper in Section~\ref{sec:conclusion}.
Auxiliary results, as well as proofs of some Lemmas, Corollaries and Propositions, are provided in the Appendix.

\section{Notation and Preliminaries}
\label{sec:wasserstein}

In this section, we introduce notation and briefly outline some known results regarding Wasserstein distance. For a more detailed discussion we refer to \citet{villani2003topics,villani2008optimal}.

Let $\Xi$ be a Polish (separable complete metric) space with metric~$d$.
Let $\scrB(\Xi)$ denote the Borel $\sigma$-algebra on $\Xi$, and let $\scrB_{\nu}(\Xi)$ denote the \textit{completion} of $\scrB(\Xi)$ with respect to a measure $\nu$ on $\scrB(\Xi)$ such that the measure space $(\Xi,\scrB_{\nu}(\Xi),\nu)$ is complete (see, e.g., Definition~1.11 in \citet{ambrosio2000functions}).
Let $\cB(\Xi)$ denote the set of Borel measures on $\Xi$, let $\cP(\Xi)$ denote the set of Borel probability measures on $\Xi$, and let $\cP_{p}(\Xi)$ denote the subset of $\cP(\Xi)$ with finite $p$-th moment for $p \in [1,\infty)$:
\[
\cP_{p}(\Xi) \ \ \defi \ \ \left\{\mu \in \cP(\Xi) \; : \; \int_{\Xi} d^{p}(\xi,\zeta^{0}) \mu(d\xi) \, < \, \infty \textrm{ for some } \zeta^{0} \in \Xi\right\}.
\]
It follows from the triangle inequality that the definition above does not depend on the choice of $\zeta^{0}$.
A function $\Psi : \Xi \mapsto \R$ is called $\nu$-measurable if it is $(\scrB_{\nu}(\Xi),\scrB(\R))$-measurable, and a function $T : \Xi \mapsto \Xi$ is called $\nu$-measurable if it is $(\scrB_{\nu}(\Xi),\scrB(\Xi))$-measurable.
To facilitate later discussion, we introduce the push-forward operator on measures.

\begin{definition}[Push-forward Measure]
\label{def:pushforward}
Given measurable spaces $(\Xi,\scrB(\Xi))$ and $(\Xi',\scrB(\Xi'))$, a measurable function $T : \Xi \mapsto \Xi'$, and a measure $\nu \in \cB(\Xi)$, let $T_{\#}\nu \in \cB(\Xi')$ denote the \textsl{push-forward measure of $\nu$ through $T$}, defined by
\[
T_{\#}\nu(A) \ \ \defi \ \ \nu(T^{-1}(A))
\ \ = \ \ \nu\{\zeta \in \Xi \; : \; T(\zeta) \in A\}, \ \forall \textnormal{ measurable sets } A \subset \Xi'.
\]
\end{definition}

That is, $T_{\#}\nu$ is obtained by \textsl{transporting} (``pushing forward'') $\nu$ from $\Xi$ to $\Xi'$ using the function $T$.
For $i \in \{1,2\}$, let $\pi^{i} : \Xi \times \Xi \mapsto \Xi$ denote the canonical projections given by $\pi^{i}(\xi^{1},\xi^{2}) = \xi^{i}$.
Then for a measure $\gamma \in \cP(\Xi \times \Xi)$, $\pi^{i}_{\#}\gamma \in \cP(\Xi)$ is the $i$-th marginal of $\gamma$ given by $\pi^{1}_{\#}\gamma(A) = \gamma(A \times \Xi)$ and $\pi^{2}_{\#}\gamma(A) = \gamma(\Xi \times A)$.

\begin{definition}[Wasserstein distance]
\label{def:wasserstein}
The \textsl{Wasserstein distance} $W_{p}(\mu,\nu)$ between $\mu,\nu \in \cP_{p}(\Xi)$ is defined by
\begin{equation}
\label{eqn:def_wasserstein}
W_{p}^{p}(\mu,\nu) \ \ \defi \ \ \min_{\gamma \in \cP(\Xi \times \Xi)} \left\{\int_{\Xi \times \Xi} d^{p}(\xi,\zeta) \gamma(d\xi,d\zeta) \; : \; \pi^{1}_{\#}\gamma = \mu, \pi^{2}_{\#}\gamma = \nu\right\}.
\end{equation}
\end{definition}

That is, the Wasserstein distance between $\mu,\nu$ is the minimum cost (in terms of $d^{p}$) of redistributing mass from $\nu$ to $\mu$, which is why it is also called the ``earth mover's distance''.
Wasserstein distance is a natural way of comparing two distributions when one is obtained from the other by \textsl{perturbations}.
The minimum on the right side of~(\ref{eqn:def_wasserstein}) is attained, because $d$ is non-negative, continuous and thus lower semi-continuous (Theorem~1.3 of \cite{villani2003topics}).
The following example is a familiar special case of problem~(\ref{eqn:def_wasserstein}).

\begin{example}[Transportation problem]
\label{eg:transportationLP}
Consider $\mu = \sum_{i=1}^{M} p_{i} \delta_{\xi^{i}}$ and $\nu = \sum_{j=1}^{N} q_{j} \delta_{\hxi^{j}}$, where $M,N \geq 1$, $p_{i},q_{j} \geq 0$, $\xi^{i},\hxi^{j} \in \Xi$ for all $i,j$, and $\sum_{i=1}^{M} p_{i} = \sum_{j=1}^{N} q_{j} = 1$. Then problem~(\ref{eqn:def_wasserstein}) becomes the classical transportation problem in linear programming:
\[
\min_{\gamma_{ij} \geq 0} \left\{\sum_{i=1}^{M} \sum_{j=1}^{N} d^{p}(\xi^{i},\hxi^{j}) \gamma_{ij} \ : \ \sum_{j=1}^{N} \gamma_{ij} = p_{i}, \; \forall \; i, \ \sum_{i=1}^{M} \gamma_{ij} = q_{j}, \; \forall \; j\right\}.
\]
% In particular, when $M=N$ and $p_{i}=q_{j}$ for all $i,j$, it reduces to the classical assignment problem. In this case, by Birkhoff's Theorem, there exists an $N$-permutation, denoted as $\sigma$, such that $\gamma^{\ast}$ defined by $\gamma_{ij}^{\ast}\defi\frac{1}{N}\mathds{1}_{\{i=\sigma(j)\}}$ is the optimal solution. Setting $T(\hxi^{j})\defi\xi^{\sigma(j)}$, and $(T\times\mathrm{id})(\hxi^{j})\defi(T(\hxi^{j}),\hxi^{j})$ for all $j$, then $\gamma^{\ast}=(T\times\mathrm{id})_{\#}\nu$.
\end{example}

\begin{example}[Revisiting Example~\ref{eg:bird}]
\label{eg:bird_w}
Next we evaluate the Wasserstein distance between the histograms in Example~\ref{eg:bird}. To evaluate $W_{1}(\bp_{true},\bq)$, note that the least cost way of transporting mass from $\bq$ to $\bp_{true}$ is to move the mass outwards. In contrast, to evaluate $W_{1}(\bp_{pathol},\bq)$, one has to transport mass relatively long distances from right to left (changing the gray-scale values of pixels by large amounts), resulting in a larger cost than $W_{1}(\bp_{true},\bq)$.  Therefore $W_{1}(\bp_{pathol},\bq) > W_{1}(\bp_{true},\bq)$.
\end{example}

Wasserstein distance has a dual representation due to Kantorovich's duality (Theorem~5.10 in \cite{villani2008optimal}):
\begin{equation}
\label{eqn:Wp_dual}
W_{p}^{p}(\mu,\nu) \ \ = \ \ \sup_{u \in L^{1}(\mu), v \in L^{1}(\nu)} \left\{\int_{\Xi} u(\xi) \mu(d\xi) + \int_{\Xi} v(\zeta) \nu(d\zeta) \; : \; u(\xi) + v(\zeta) \leq d^{p}(\xi,\zeta), \; \forall \; \xi,\zeta \in \Xi\right\},
\end{equation}
where $L^{1}(\nu)$ represents the $L^{1}$ space of $\nu$-measurable functions.
In addition, $u \in L^{1}(\mu), v \in L^{1}(\nu)$ under the supremum above can be replaced by $u,v \in C_{b}(\Xi)$, where $C_{b}(\Xi)$ denotes the set of continuous and bounded real-valued functions on $\Xi$.
Particularly, when $p = 1$, by the Kantorovich-Rubinstein theorem, (\ref{eqn:Wp_dual}) can be simplified to (see, e.g., Equation (5.11) in \cite{villani2008optimal})
\[
W_{1}(\mu,\nu) \ \ = \ \ \sup_{u \in L^{1}(\mu)} \left\{\int_{\Xi} u(\xi) d(\mu-\nu)(\xi) \; : \; u \textrm{ is 1-Lipschitz}\right\}.
\]
So for an $L$-Lipschitz function $\Psi : \Xi \mapsto \R$, it holds that $\big|\E_{\mu}[\Psi(\xi)] - \E_{\nu}[\Psi(\xi)]\big| \leq L W_{1}(\mu,\nu) \leq L \theta$ for all $\mu \in \frakM$.
% The following lemma generalizes this statement.

% \begin{lemma}
%   \label{lemma:Emu-Enu}
%   Let $\Psi : \Xi \mapsto \R$.
%   Suppose that $\Psi$ satisfies $|\Psi(\xi) - \Psi(\zeta)| \leq L d^{p}(\xi,\zeta) + M$ for some $L,M \geq 0$ and all $\xi,\zeta \in \Xi$.
%   Then
%   \[
%   \big|\E_{\mu}[\Psi(\xi)] - \E_{\nu}[\Psi(\xi)]\big| ~\leq ~L \theta^{p} + M, \quad \forall \ \mu \in \frakM.
%   \]
% \end{lemma}

Definition~\ref{def:wasserstein} and the results above can be extended to finite Borel measures.
Moreover, we have the following result.

\begin{lemma}
\label{lemma:enlarge_{p}}
For any finite Borel measures $\mu,\nu \in \cB(\Xi)$ with $\mu(\Xi) \neq \nu(\Xi)$, it holds that $W_{p}(\mu,\nu) = \infty$.
\end{lemma}

Another important feature of Wasserstein distance is that $W_{p}$ metrizes weak convergence in $\cP_{p}(\Xi)$ (cf. Theorem~6.9 in \citet{villani2008optimal}). That is, for any sequence $\{\mu_{k}\}_{k=1}^{\infty}$ of measures in $\cP_{p}(\Xi)$ and $\mu \in \cP_{p}(\Xi)$, it holds that $\lim_{k \to \infty} W_{p}(\mu_{k},\mu) = 0$ if and only if $\mu_{k}$ converges weakly to $\mu$ and $\int_{\Xi} d^{p}(\xi,\zeta^{0}) \mu_{k}(d\xi) \to \int_{\Xi} d^{p}(\xi,\zeta^{0}) \mu(d\xi)$ as $k \to \infty$. Therefore, convergence in the Wasserstein distance of order $p$ implies convergence up to the $p$-th moment. \citet[chapter~6]{villani2008optimal} discusses the advantages of Wasserstein distance relative to other distances, such as the Prokhorov metric, that metrize weak convergence.

\section{Tractable Reformulation via Duality}
\label{sec:tractability}

In this section we develop a tractable reformulation by deriving its strong dual.
We suppress the variable $x$ of $\Psi$ in this section, and results are interpreted pointwise for each $x$.
Given $\nu \in \cP(\Xi)$ and $\Psi \in L^{1}(\nu)$, for any $\theta > 0$ and $p \in [1,\infty)$,
the inner maximization problem of \eqref{eqn:DRSO} is written as
\begin{equation}
\label{eqn:problem_primal}
v_{P} \ \ \defi \ \ \sup_{\mu \in \frakM} \; \int_{\Xi} \Psi(\xi) \mu(d\xi)
\ \ = \ \ \sup_{\mu \in \cP(\Xi)} \left\{\int_{\Xi} \Psi(\xi) \mu(d\xi) \; : \; W_{p}(\mu,\nu) \leq \theta\right\}.
\tag{Primal}
\end{equation}
Our main goal is to derive its strong dual
\begin{equation}
\label{eqn:problem_dual}
v_{D} \ \ \defi \ \ \inf_{\lambda \geq 0} \left\{\lambda \theta^{p} - \int_{\Xi} \inf_{\xi \in \Xi} \big[\lambda d^{p}(\xi,\zeta) - \Psi(\xi)\big] \nu(d\zeta)\right\}.
\tag{Dual}
\end{equation}
The dual problem is a one-dimensional convex minimization problem with respect to $\lambda$, the Lagrangian multiplier of the Wasserstein constraint in the primal problem.
The term $\inf_{\xi \in \Xi} [\lambda d^{p}(\xi,\zeta) - \Psi(\xi)]$ is called the \textsl{Moreau-Yosida regularization} of $-\Psi$ with parameter $1/\lambda$ in the literature (cf. \citet{parikh2014proximal}).
Its measurability with respect to $\nu$ is established in Lemma~\ref{lemma:measurable}\ref{itm:measurable_phi} in Section~\ref{sec:dual_general}.

\subsection{General Nominal Distribution}
\label{sec:dual_general}

In this section, we prove the strong duality result for a general nominal distribution~$\nu$ on a Polish space~$\Xi$. Such generality broadens the applicability of the result for \eqref{eqn:DRSO}. For example, the result is useful when the nominal distribution is some distribution such as a Gaussian distribution on $\R^{K}$ (Section~\ref{sec:wcVaR}), or even some stochastic process (Sections~\ref{sec:process} and~\ref{sec:NHPP}).
We begin with a weak duality result, which is an application of Lagrangian weak duality.

\begin{proposition}[Weak duality]
\label{prop:weakDuality}
Consider any $\nu \in \cP(\Xi)$ and $\Psi \in L^{1}(\nu)$. Then for any $p \in [1,\infty)$ and $\theta > 0$, it holds that $v_{P} \leq v_{D}$.
\end{proposition}

To prove the strong duality result, we consider two separate cases: $v_{D} = \infty$ and $v_{D} < \infty$.
As can be seen from \eqref{eqn:problem_dual}, if the term $-\int_{\Xi} \inf_{\xi \in \Xi} \big[\lambda d^{p}(\xi,\zeta) - \Psi(\xi)\big] \nu(d\zeta)$ is infinite for all $\lambda \geq 0$, then $v_{D} = \infty$.
Thus, to facilitate our analysis, we introduce the following definitions.

\begin{definition}[Regularization Operator $\Phi$]
\label{def:phi}
Let $\Phi : \R \times \Xi \mapsto \R \cup \{-\infty\}$ be given by
\[
\label{eqn:Phi}
\Phi(\lambda,\zeta) \ \ \defi \ \ \inf_{\xi \in \Xi} \big\{\lambda d^{p}(\xi,\zeta) - \Psi(\xi)\big\}.
\]
\end{definition}

%For any function $f : \R_{+} \mapsto \R \cup\{-\infty\}$, let $\Dom(f)$ denote its effective domain
%\[
%\Dom(f) \ \ \defi \ \ \{\lambda \geq 0 \; : \; f(\lambda) > -\infty\},
%\]
%and let $\interior(\Dom(f))$ denote the interior of $\Dom(f)$.

For any $\lambda \ge 0$ and any $\zeta \in \Xi$ such that $\Phi(\lambda,\zeta) > -\infty$, let
\begin{equation}
\label{eqn:Dpm}
\begin{aligned}
\Dp(\lambda,\zeta) &&& \defi &&& \limsup_{\varepsilon \downarrow 0} \Big\{\sup_{\xi \in \Xi} \big\{d^{p}(\xi,\zeta) \; : \; \lambda d^{p}(\xi,\zeta) - \Psi(\xi) \, \leq \, \Phi(\lambda,\zeta) + \varepsilon\big\}\Big\}, \\
\Dm(\lambda,\zeta) &&& \defi &&& \liminf_{\varepsilon \downarrow 0} \Big\{\inf_{\xi \in \Xi} \big\{d^{p}(\xi,\zeta) \; : \; \lambda d^{p}(\xi,\zeta) - \Psi(\xi) \, \leq \, \Phi(\lambda,\zeta) + \varepsilon\big\}\Big\}.
\end{aligned}
\end{equation}
For any $\lambda \ge 0$ and any $\zeta \in \Xi$ such that $\argmin_{\xi \in \Xi} \{\lambda d^{p}(\xi,\zeta) - \Psi(\xi)\}$ is nonempty, let
\begin{equation}
\label{eqn:D0}
\begin{aligned}
\overline{D}_{0}(\lambda,\zeta) &&& \defi &&& \sup_{\xi \in \Xi} \{d^{p}(\xi,\zeta) \; : \; \lambda d^{p}(\xi,\zeta) - \Psi(\xi) = \Phi(\lambda,\zeta)\} \\
\underline{D}_{0}(\lambda,\zeta) &&& \defi &&& \inf_{\xi \in \Xi} \{d^{p}(\xi,\zeta) \; : \; \lambda d^{p}(\xi,\zeta) - \Psi(\xi) = \Phi(\lambda,\zeta)\}
\end{aligned}
\end{equation}
Then $\underline{D}_{0}(\lambda,\zeta)$ and $\overline{D}_{0}(\lambda,\zeta)$ represent respectively the closest and furthest distances between $\zeta$ and any point in $\argmin_{\xi \in \Xi} \{\lambda d^{p}(\xi,\zeta) - \Psi(\xi)\}$.
Note that $\overline{D}_{0}(\lambda,\zeta)$ (resp. $\underline{D}_{0}(\lambda,\zeta)$) may not be equal to $\Dp(\lambda,\zeta)$ (resp. $\Dm(\lambda,\zeta)$).

\begin{definition}[Growth rate]
\label{def:kappa}
Define the \textit{growth rate} $\kappa$ of $\Psi$ as
\[
\label{eqn:kappa_function}
\kappa \ \ \defi \ \ \inf\left\{\lambda \geq 0 \; : \; \int_{\Xi} \Phi(\lambda,\zeta) \nu(d\zeta) > -\infty\right\}.
% \begin{cases} 0, & \textrm{ if } \sup_{\xi,\zeta\in\Xi} d^{p}(\xi,\zeta) <\infty,\\
% % \displaystyle
% % \esssup{\nu}{\zeta\in\Xi} \limsup_{\xi\in\Xi:\; d^{p}(\xi,\zeta) \to \infty} \frac{\max\{0,\Psi(\xi) - \Psi(\zeta)\}}{d^{p}(\xi,\zeta)},
% & \textrm{ otherwise}.
% \end{cases}
\]
Particularly, if $\int_{\Xi} \Phi(\lambda,\zeta) \nu(d\zeta) = -\infty$ for all $\lambda \geq 0$, then $\kappa = \infty$.
\end{definition}

If $\Xi$ is bounded and $\Psi$ is bounded above, then $\kappa = 0$, and if $\Xi$ is bounded and $\Psi$ is not bounded above, then $\kappa = \infty$.
The possibilities when $\Xi$ is not bounded are more interesting.
Next, Lemma~\ref{lemma:kappa} establishes some additional properties of $\kappa$, including the property that if $\Xi$ is not bounded and $\kappa < \infty$, then
\[
\kappa \ \ = \ \ \limsup_{\xi \in \Xi \; : \; d^{p}(\xi,\zeta) \to \infty} \frac{\max\{0, \Psi(\xi) - \Psi(\zeta)\}}{d^{p}(\xi,\zeta)}
\]
for any $\zeta \in \Xi$, which motivates why we call $\kappa$ the growth rate of $\Psi$.

\begin{lemma}[Properties of the growth rate $\kappa$]
\label{lemma:kappa}
$\mbox{}$
\begin{enumerate}[label=(\roman*)]
\item
\label{lemma:kappa0}
Suppose that $\Xi$ is unbounded.
Then the quantity
\[
\limsup_{\xi \in \Xi \, : \, d^{p}(\xi,\zeta) \to \infty} \frac{\max\{0, \Psi(\xi) - \Psi(\zeta)\}}{d^{p}(\xi,\zeta)}
\]
is independent of $\zeta$.
\item
\label{lemma:kappa_condition}
Suppose that
%\[
%\nu \ \ \in \ \ cP_{p}(\Xi) \ \ \defi \ \ \left\{\mu \in \cP(\Xi) \; : \; \int_{\Xi} d^{p}(\zeta,\zeta^{0}) \nu(d\zeta) < \infty \textrm{ for some } \zeta^{0} \in \Xi\right\}.
%\]
$\nu \in \cP_{p}(\Xi)$.
Then the growth rate $\kappa$ is finite if and only if there exists $\zeta^{0} \in \Xi$ and $L, M > 0$ such that
\begin{equation}
\label{eqn:kappa_condition}
\Psi(\xi) - \Psi(\zeta^{0}) \ \ \leq \ \ L d^{p}(\xi,\zeta^{0}) + M \ \ \ \forall \ \xi \in \Xi.
\end{equation}
\item
\label{lemma:kappa_expression}
Suppose that $\nu \in \cP_{p}(\Xi)$.
If $\Xi$ is unbounded and $\kappa < \infty$, then
\[
\kappa \ \ = \ \ \limsup_{\xi \in \Xi \, : \, d^{p}(\xi,\zeta) \to \infty} \frac{\max\{0, \Psi(\xi) - \Psi(\zeta)\}}{d^{p}(\xi,\zeta)}
\]
for any $\zeta \in \Xi$.
\end{enumerate}
\end{lemma}

\begin{lemma}[Measurability]
\label{lemma:measurable}
For any $p \in [1,\infty)$, $\nu \in \cP(\Xi)$, and $\Psi \in L^{1}(\nu)$, the following hold:
\begin{enumerate}[label=(\roman*)]
\item
\label{itm:measurable_phi}
$\Phi(\lambda,\cdot)$, $\Dp(\lambda,\cdot)$, $\Dm(\lambda,\cdot)$, $\overline{D}_{0}(\lambda,\cdot)$, and $\underline{D}_{0}(\lambda,\cdot)$ are $\nu$-measurable.

\item
\label{itm:measurable_selection}
Suppose that $\kappa < \infty$.
%Let $\lambda \in \Dom\big(\int_{\Xi} \Phi(\cdot,\zeta) \nu(d\zeta)\big)$ be such that $\Dp(\lambda,\zeta) < \infty$.
Then for any $\lambda,\delta,\varepsilon \geq 0$ such that the sets
\begin{eqnarray*}
\overline{F}^{\varepsilon}_{\delta}(\lambda,\zeta) & \defi & \Big\{\xi \in \Xi  \; : \; \lambda d^{p}(\xi,\zeta) - \Psi(\xi) \, \leq \, \Phi(\lambda,\zeta) + \varepsilon, \ d^{p}(\xi,\zeta) \, \geq \, \Dp(\lambda,\zeta) - \delta\Big\}, \\
\underline{F}^{\varepsilon}_{\delta}(\lambda,\zeta) & \defi & \Big\{\xi \in \Xi \; : \; \lambda d^{p}(\xi,\zeta) - \Psi(\xi) \, \leq \, \Phi(\lambda,\zeta) + \varepsilon, \ d^{p}(\xi,\zeta) \, \leq \, \Dm(\lambda,\zeta) + \delta\Big\}
\end{eqnarray*}
are non-empty for $\nu$-almost all $\zeta \in \Xi$, there exists $\nu$-measurable mappings $\overline{T}^{\varepsilon}_{\delta}(\lambda,\cdot),\underline{T}^{\varepsilon}_{\delta}(\lambda,\cdot) : \Xi \mapsto \Xi$ such that $\overline{T}^{\varepsilon}_{\delta}(\lambda,\zeta) \in \overline{F}^{\varepsilon}_{\delta}(\lambda,\zeta)$ and $\underline{T}^{\varepsilon}_{\delta}(\lambda,\zeta) \in \underline{F}^{\varepsilon}_{\delta}(\lambda,\zeta)$ for $\nu$-almost all $\zeta \in \Xi$.

\item
\label{itm:measurable_selection0}
Suppose that $\kappa < \infty$.
Then for any $\lambda,\delta \geq 0$ such that the sets
\begin{eqnarray*}
\overline{F}(\lambda,\zeta) & \defi & \Big\{\xi \in \Xi  \; : \; \lambda d^{p}(\xi,\zeta) - \Psi(\xi) \, = \, \Phi(\lambda,\zeta), \ d^{p}(\xi,\zeta) \, \geq \, \overline{D}_{0}(\lambda,\zeta) - \delta\Big\}, \\
\underline{F}(\lambda,\zeta) & \defi & \Big\{\xi \in \Xi \; : \; \lambda d^{p}(\xi,\zeta) - \Psi(\xi) \, = \, \Phi(\lambda,\zeta), \ d^{p}(\xi,\zeta) \, \leq \, \underline{D}_{0}(\lambda,\zeta) + \delta\Big\}
\end{eqnarray*}
are non-empty for $\nu$-almost all $\zeta \in \Xi$, there exists $\nu$-measurable mappings $\overline{T}(\lambda,\cdot),\underline{T}(\lambda,\cdot) : \Xi \mapsto \Xi$ such that $\overline{T}(\lambda,\zeta) \in \overline{F}(\lambda,\zeta)$ and $\underline{T}(\lambda,\zeta) \in \underline{F}(\lambda,\zeta)$ for $\nu$-almost all $\zeta \in \Xi$.

\item
\label{itm:measurable_n}
For any $E \in \scrB_{\nu}(\Xi)$ and $a,b > 0$ such that
\[
F(\zeta) \ \ \defi \ \ \{\xi \in \Xi \; : \; \Psi(\xi) - \Psi(\zeta) > a d^{p}(\xi,\zeta) + b\}
\]
is non-empty for $\nu$-almost all $\zeta \in E$, there exists a $\nu$-measurable mapping $T : E \mapsto \Xi$ such that $T(\zeta) \in F(\zeta)$ for $\nu$-almost all $\zeta \in E$.

\item
\label{itm:measurable_M}
Suppose that $\kappa < \infty$.
For any $\kappa' \in (0,\kappa)$, and any $\nu$-measurable function $M : \Xi \mapsto \R$ such that the set
\[
F(\zeta) \ \ \defi \ \ \{\xi \in \Xi \; : \; \Psi(\xi) - \Psi(\zeta) \geq \kappa' d^{p}(\xi,\zeta), \; d^{p}(\xi,\zeta) \geq M(\zeta)\}
\]
is non-empty for $\nu$-almost all $\zeta \in \Xi$, there exists a $\nu$-measurable mapping $T : \Xi \mapsto \Xi$ such that $T(\zeta) \in F(\zeta)$ for $\nu$-almost all $\zeta \in \Xi$.
\end{enumerate}

% \[\label{eqn:selection overline}
%   \overline{T}(\lambda,\zeta) \ \in \ \Big\{\xi \in \Xi \; : \; \lambda d^{p}(\xi,\zeta) - \Psi(\xi) \, \leq \, \Phi(\lambda,\zeta) + \delta, \ d^{p}(\xi,\zeta) \, \leq \, \Dp(\lambda,\zeta) - \varepsilon\Big\},
% \]
% \[ \label{eqn:selection underline}
%   \underline{T}(\lambda,\zeta) \ \in \ \Big\{\xi \in \Xi \; : \; \lambda d^{p}(\xi,\zeta) - \Psi(\xi) \, \leq \, \Phi(\lambda,\zeta) + \delta, \ d^{p}(\xi,\zeta) \, \geq \, \Dm(\lambda,\zeta) + \varepsilon\Big\}.
% \]
% If $\overline{D}(\kappa,\zeta)=\infty$, and the set $\{\xi \in \Xi \; : \; \kappa d^{p}(\xi,\zeta) - \Psi(\xi) \, \leq \, \Phi(\kappa,\zeta) + \delta \}$
% is non-empty for some $\delta\geq0$.
% Then for any $\nu$-measurable  $R > 0$, there exists a $\nu$-measurable mapping $\overline{T} : \Xi \mapsto \Xi$ such that for $\nu$-almost all $\zeta$,
% \[ \label{eqn:selection_R}
%   \overline{T}(\zeta) \  \in  \ \Big\{\xi \in \Xi \; : \; \kappa d^{p}(\xi,\zeta) - \Psi(\xi) \, \leq \, \Phi(\kappa,\zeta) + \delta,\ d^{p}(\xi,\zeta) >  R\Big\}.
% \]
\end{lemma}

\begin{proposition}[Strong duality with infinite optimal value]
\label{prop:kappa_infty}
Consider any $p \in [1,\infty)$, $\nu \in \cP(\Xi)$, and $\Psi \in L^{1}(\nu)$.
Suppose that $\theta > 0$ and $\kappa = \infty$.  Then $v_{P} = v_{D} = \infty$.
\end{proposition}
% We note that when $\theta=0$ and $\kappa=\infty$, the strong duality does not hold.
% Indeed, $W_p(\mu,\nu)=0$ implies that $\mu=\nu$ and thus $v_{P}=\E_{\nu}[\Psi]$.
% On the other hand, for the dual problem, $\int_{\Xi} \inf_{\xi\in\Xi}[\lambda d^{p}(\xi,\zeta) - \Psi(\xi)] \nu(d\zeta)=-\infty$ for any $\lambda\geq0$ and thus $v_{D} = \infty$.

% The idea of the proof is straightforward, though we have to take care of some technical details, such as the measurability of the inner infimum involved in the dual problem, and the difficulty resulting from the unboundedness of $\Xi$.

\begin{remark}[Choosing Wasserstein order $p$]
\label{rmk:p}
Let
\[
\underline{p} \ \ \defi \ \ \inf\Big\{p \geq 1 \; : \; \limsup_{d(\zeta,\zeta^{0}) \to \infty} \frac{\Psi(\zeta) - \Psi(\zeta^{0})}{d^{p}(\zeta,\zeta^{0})} < \infty\Big\}.
\]
Proposition~\ref{prop:kappa_infty} suggests that a meaningful formulation of \eqref{eqn:DRSO} should be such that the Wasserstein order $p$ is greater than or equal to $\underline{p}$.
% Noticing that in Definition \ref{def:phi}, parameter $p$ controls the extent of regularization. When $p$ is much greater than $\underline{p}$, only small perturbation is allowed.
In both \citet{esfahani2015data} and \citet{zhao2018data} only $p = 1$ is considered.  By considering higher orders $p$ in our analysis, we can accommodate a greater set of functions $\Psi$, and we also have more flexibility to choose the ambiguity set and to control the degree of conservativeness.
\end{remark}

\begin{lemma}[Properties of the regularization operator $\Phi$]
\label{lemma:phi} %\leavevmode
Let $(\Xi,d)$ be a Polish space.
Consider any $p \in [1,\infty)$, $\nu \in \cP(\Xi)$, and $\Psi \in L^{1}(\nu)$ such that $\kappa < \infty$.
Then there is a set $B \in \scrB_{\nu}(\Xi)$ such that $\nu(B) = 1$, and the following holds.
\begin{enumerate}[label=(\roman*),itemsep=1ex,topsep=1ex]
%\item
%\label{lemma:phi finite} [Finite values]
%For all $\lambda \ge 0$ and all $\zeta \in \Xi$, it holds that $\Phi(\lambda,\zeta) \leq -\Psi(\zeta)$.
%Also, if $\nu \in cP_{p}(\Xi)$ and $\kappa < \infty$, then for all $\lambda > \kappa$ and all $\zeta \in \Xi$, it holds that $\Phi(\lambda,\zeta) > -\infty$, and for all $\lambda < \kappa$ and all $\zeta \in \Xi$, it holds that $\Phi(\lambda,\zeta) = -\infty$.
%That is, $(\kappa,\infty) \subseteq \Dom(\Phi(\cdot,\zeta)) \subseteq [\kappa,\infty)$ and $\interior(\Dom(\Phi(\cdot,\zeta))) = (\kappa,\infty)$.
%\item
%\label{lemma:phi_bound}[Minimizers]
%For all $\lambda > \kappa$ and all $\zeta \in \Xi$, it holds that $\argmin_{\xi \in \Xi}  \big\{\lambda d^{p}(\xi,\zeta) - \Psi(\xi)\big\}$ is nonempty and compact.

\item
\label{lemma:phi_monotone}
[Monotonicity]
$\Phi(\cdot,\zeta)$ is nondecreasing and upper-semi-continuous for all $\zeta \in \Xi$.
$\Phi(\lambda,\zeta) > -\infty$ for all $\lambda > \kappa$ and all $\zeta \in B$.
$\Phi(\cdot,\zeta)$ is concave for all $\zeta \in B$.
For any $\varepsilon \ge 0$, any $\lambda_{2} > \lambda_{1}$, and any $\zeta \in \Xi$ such that $\Phi(\lambda_{1},\zeta) > -\infty$, it holds that \[\sup_{\xi \in \Xi} \left\{d^{p}(\xi,\zeta) \, : \, \lambda_{2} d^{p}(\xi,\zeta) - \Psi(\xi) \leq \Phi(\lambda_{2}, \zeta) + \varepsilon\right\} \le \sup_{\xi \in \Xi} \left\{d^{p}(\xi,\zeta) \, : \, \lambda_{1} d^{p}(\xi,\zeta) - \Psi(\xi) \leq \Phi(\lambda_{1}, \zeta) + \varepsilon\right\}.
\]
For any $\lambda_{2} > \lambda_{1}$ and any $\zeta \in \Xi$ such that $\Phi(\lambda_{1},\zeta) > -\infty$, it holds that $\Dp(\lambda_{2},\zeta) \leq \Dm(\lambda_{1},\zeta) \leq \Dp(\lambda_{1},\zeta)$.

%\item
%\label{lemma:phi limit}[Limit]
%$\lim_{\lambda \to \infty} \Phi(\lambda,\zeta) >-\infty$.

\item
\label{lemma:phi_bound}
[Bounds]
For any $\lambda_{2} > \lambda_{1}$ such that $\Phi(\lambda_{1},\zeta) > -\infty$, it holds that
\[
(\lambda_{2} - \lambda_{1}) \Dp(\lambda_{2},\zeta) \ \ \leq \ \ - \Psi(\zeta) - \Phi(\lambda_{1},\zeta).
\]
% Then there is a constant $C > 0$ such that
% \[
%   \frac{\lambda - \lambda_{1}}{2} \Dp(\lambda,\zeta) ~\leq ~\Phi(\lambda,\zeta) - \Phi(\lambda_{1},\zeta^{0}) + C d^{p}(\zeta,\zeta^{0})
% \]
% for all $\zeta,\zeta^{0} \in \Xi$.

%\item
%\label{lemma:phi_continuous} [Continuity]
%$\Phi$ is continuous in both $\lambda$ and $\zeta$ on its domain.
%In addition, $\lim_{\lambda\downarrow\kappa}\Phi(\lambda,\zeta)=\Phi(\kappa,\zeta)$ provided that $\Phi(\kappa,\zeta)>-\infty$.
%If $\Xi=\R^{K}$, then $\frac{\partial\Phi(\lambda,\zeta)}{\partial\lambda+}= \Dm(\lambda,\zeta)\leq\Dp(\lambda,\zeta)= \frac{\partial\Phi(\lambda,\zeta)}{\partial\lambda-}$.

% %If it is unbounded,
% When $\Xi=\R^{K}$ and $\{\xi\in\Xi:\lambda d^{p}(\xi,\zeta)-\Psi(\xi)=\Phi(\lambda,\zeta)\}$ is non-empty, (\ref{eqn:selection_underline}) holds also for $\varepsilon=\delta=0$.
% If the set $\{\xi\in\Xi: \lambda d^{p}(\xi',\zeta)-\Psi(\xi')\leq\Phi(\lambda,\zeta)\}$ is bounded, then (\ref{eqn:selection_overline}) holds also for $\varepsilon=\delta=0$, otherwise (\ref{eqn:selection_R}) holds also for $\delta=0$.

%\item
%\label{lemma:phi integrable} [Integrable]
%For all $\lambda \ge 0$, it holds that $\int_{\Xi} \Phi(\lambda,\zeta) \nu(d\zeta) \leq - \int_{\Xi} \Psi(\zeta) \nu(d\zeta) < \infty$.
%Also, for all $\lambda > \kappa$, it holds that $\int_{\Xi} \Phi(\lambda,\zeta) \nu(d\zeta) > -\infty$.

\item
\label{lemma:phi_derivative}
[Derivative]
For all $\lambda > \kappa$ and all $\zeta \in B$, the left partial derivative $\partial \Phi(\lambda,\zeta) / \partial \lambda-$ exists and satisfies
\[
\Dp(\lambda,\zeta) \ \ \leq \ \ \frac{\partial \Phi(\lambda,\zeta)}{\partial \lambda-} \ \ \leq \ \ \lim_{\lambda_{1} \uparrow \lambda} \Dm(\lambda_{1},\zeta).
\]
For all $\lambda \ge 0$ and $\zeta \in \Xi$ such that $\Phi(\lambda,\zeta) > -\infty$, the right partial derivative $\partial \Phi(\lambda,\zeta) / \partial \lambda+$ exists and satisfies
\[
\lim_{\lambda_{2} \downarrow \lambda} \Dp(\lambda_{2},\zeta) \ \ \leq \ \ \frac{\partial \Phi(\lambda,\zeta)}{\partial \lambda+} \ \ \leq \ \ \Dm(\lambda,\zeta).
\]

%For $\lambda = \kappa > 0$, and any $R, \delta > 0$, there exists a $\nu$-measurable map $\overline{T} : \Xi \mapsto \Xi$ such that
%\begin{equation}
%\label{eqn:selection_R}
%\overline{T}(\zeta) ~\in ~\Big\{\xi \in \Xi \; : \; d^{p}(\xi,\zeta) \, > \, R, \, \kappa d^{p}(\xi,\zeta) - \Psi(\xi) < \lim_{\lambda \downarrow \kappa} \Phi(\lambda,\zeta) + \delta\Big\}.
%\end{equation}
%If $\Phi(\kappa,\zeta) > -\infty$ and $\kappa > 0$, then for any $R, \delta > 0$, there exists a $\nu$-measurable map $\overline{T}(\kappa,\cdot) : \Xi \mapsto \Xi$ such that
%\begin{equation}
%\label{eqn:selection_R}
%\overline{T}(\kappa,\zeta) ~\in ~\Big\{\xi \in \Xi \; : \; \kappa d^{p}(\xi,\zeta) - \Psi(\xi) \, \leq \, \Phi(\kappa,\zeta) + \delta, \ d^{p}(\xi,\zeta) \, > \, R\Big\}.
%\end{equation}
%If $\Phi(\kappa,\zeta) > -\infty$, then for any $\delta,\varepsilon > 0$, there exists a $\nu$-measurable map $\underline{T}(\kappa,\cdot) : \Xi \mapsto \Xi$ such that
%\[
%\underline{T}(\kappa,\zeta) ~\in ~\Big\{\xi \in \Xi \; : \; \kappa d^{p}(\xi,\zeta) - \Psi(\xi) \, \leq \, \Phi(\kappa,\zeta) + \delta, \ d^{p}(\xi,\zeta) \, \leq \, \Dm(\kappa,\zeta) + \varepsilon\Big\}.
%\]
\end{enumerate}
\end{lemma}

% Properties~\ref{lemma:phi_monotone}, \ref{lemma:phi_derivative}, and~\ref{lemma:measurable} will be used in the proof of Theorem~\ref{thm:strongDual}.
% Property~\ref{lemma:phi_bound} shows that for any $\lambda > \kappa$ and any $\zeta$, the set of $\delta$-minimizers in~\eqref{eqn:Phi} is bounded, which will be used in the proof of Corollary~\ref{cor:existence}.
% \begin{remark}\label{rmk:phi}
%   Lemma \ref{lemma:phi} remains to hold if we replace $d^{p}(\xi,\zeta)$ in~\eqref{eqn:Phi} with any $(\scrB_{\nu} \otimes \scrB_{\nu}, \scrB(\R))$-measurable function $c(\xi,\zeta)$ such that $c(\xi,\zeta)=0$ if and only if $\xi=\zeta$.
%   In this case, we need to replace $\kappa$ with $\kappa(\zeta)$ as defined in \eqref{eqn:kappa}.
% \end{remark}

% More specifically, consider the following condition: there exists $L,M : \Xi \mapsto \R_+$ with $\esssup{\nu}{\zeta\in\Xi}\; L(\zeta)<\infty$ and $M\in L^{1}(\nu)$, such that
% \begin{equation} \label{eqn:growth_condition}
%   \Psi(\xi)-\Psi(\zeta) ~\leq~ L(\zeta)\cdot d^{p}(\xi,\zeta) + M(\zeta),~\forall \xi,\zeta\in\Xi.
%   \tag{Growth condition}
% \end{equation}
% When \eqref{eqn:growth_condition} does not hold, the optimal value is infinity, as shown below.

Let $h : \R_{+} \mapsto \R \cup \{\infty\}$ denote the dual objective function given by
\[
h(\lambda) \ \ \defi \ \ \lambda \theta^{p} - \int_{\Xi} \Phi(\lambda,\zeta) \nu(d\zeta).
\]

\begin{lemma}[Dual objective function]
\label{lemma:dual objective}
The dual objective function $h$ has the following properties:
\begin{enumerate}[label=(\roman*)]
\item
\label{itm:dual objective finite}
$h(\lambda) = \infty$ for all $\lambda \in [0,\kappa)$ and $h(\lambda) < \infty$ for all $\lambda \in (\kappa,\infty)$.
\item
\label{itm:dual objective convex}
$h$ is a convex function.
\item
\label{itm:dual objective lower semi-continuous}
$h$ is a lower semi-continuous function.
\item
\label{itm:dual objective infinite}
$h(\lambda) \to \infty$ as $\lambda \to \infty$.
\item
\label{itm:dual minimizer}
$h$ has a minimizer $\lambda^{\ast} \in [\kappa,\infty)$.
\end{enumerate}
\end{lemma}

\begin{lemma}[Structure of $\varepsilon$-optimal primal solution with $\lambda^{\ast} > \kappa$]
\label{lem:e-optimal lambda > kappa}
Consider any $p \in [1,\infty)$, any $\nu \in \cP(\Xi)$, any $\theta > 0$, and any $\Psi \in L^{1}(\nu)$ such that $\kappa < \infty$.
Suppose that $h$ has a minimizer $\lambda^{\ast} > \kappa$.
Then, for any $\varepsilon > 0$, there are a $\lambda_{1}^{\varepsilon} \in \left(\max\{\kappa, \lambda^{\ast} - \varepsilon\}, \, \lambda^{\ast}\right)$, a $\lambda_{2}^{\varepsilon} \in \left(\lambda^{\ast}, \, \lambda^{\ast} + \varepsilon\right)$, a $\delta^{\varepsilon} \in (0,\varepsilon)$, $\nu$-measurable mappings $T_{1}^{\varepsilon},T_{2}^{\varepsilon} : \Xi \mapsto \Xi$, and $p_{1}^{\varepsilon},p_{2}^{\varepsilon},p_{3}^{\varepsilon} \ge 0$, such that $p_{1}^{\varepsilon} + p_{2}^{\varepsilon} + p_{3}^{\varepsilon} = 1$,
\begin{align*}
T_{1}^{\varepsilon}(\zeta) \ \ & \in \ \ \Big\{\xi \in \Xi  \; : \; \lambda_{1}^{\varepsilon} d^{p}(\xi,\zeta) - \Psi(\xi) \, \leq \, \Phi(\lambda_{1}^{\varepsilon},\zeta) + \varepsilon, \ d^{p}(\xi,\zeta) \, \geq \, \Dp(\lambda_{1}^{\varepsilon},\zeta) - \delta^{\varepsilon}\Big\}, \\
T_{2}^{\varepsilon}(\zeta) \ \ & \in \ \ \Big\{\xi \in \Xi \; : \; \lambda_{2}^{\varepsilon} d^{p}(\xi,\zeta) - \Psi(\xi) \, \leq \, \Phi(\lambda_{2}^{\varepsilon},\zeta) + \varepsilon, \ d^{p}(\xi,\zeta) \, \leq \, \Dm(\lambda_{2}^{\varepsilon},\zeta) + \delta^{\varepsilon}\Big\}
\end{align*}
for $\nu$-almost all $\zeta \in \Xi$, and
\[
\mu^{\varepsilon} \ \ \defi \ \ p_{1}^{\varepsilon} 
{T_{1}^{\varepsilon}}_{\#} \nu + p_{2}^{\varepsilon} {T_{2}^{\varepsilon}}_{\#} \nu + p_{3}^{\varepsilon} \nu
\]
satisfies
\begin{align}
\int_{\Xi} \Psi(\xi) \mu^{\varepsilon}(d\xi) \ \ \ge \ \ \lambda^{\ast} \theta^{p} - \int_{\Xi} \Phi(\lambda^{\ast},\zeta) \nu(d\zeta) - \varepsilon,
\label{eqn:primal e-optimal 1}
\end{align}
and
\begin{align}
W_{p}^{p}(\mu^{\varepsilon},\nu)
\ \ \leq \ \ p_{1} \int_{\Xi} d^{p}(T_{1}^{\varepsilon}(\zeta),\zeta) \nu(d\zeta) + p_{2} \int_{\Xi} d^{p}(T_{2}^{\varepsilon}(\zeta),\zeta) \nu(d\zeta)
\ \ \leq \ \ \theta^{p}.
\label{eqn:primal feasible 1}
\end{align}
\end{lemma}

\proof{Proof of Lemma~\ref{lem:e-optimal lambda > kappa}.}
Note that for any $\lambda > \kappa$ and $\delta,\varepsilon > 0$, the sets $\overline{F}(\lambda,\zeta),\underline{F}(\lambda,\zeta)$ in Lemma~\ref{lemma:measurable}\ref{itm:measurable_selection} are non-empty for $\nu$-almost all $\zeta \in \Xi$.
Hence there exists $\nu$-measurable mappings $\overline{T}^{\varepsilon}_{\delta}(\lambda,\cdot), \underline{T}^{\varepsilon}_{\delta}(\lambda,\cdot) : \Xi \mapsto \Xi$ such that
\begin{align*}
\overline{T}^{\varepsilon}_{\delta}(\lambda,\zeta) \ \ & \in \ \ \Big\{\xi \in \Xi  \; : \; \lambda d^{p}(\xi,\zeta) - \Psi(\xi) \, \leq \, \Phi(\lambda,\zeta) + \varepsilon, \ d^{p}(\xi,\zeta) \, \geq \, \Dp(\lambda,\zeta) - \delta\Big\}, \\
\underline{T}^{\varepsilon}_{\delta}(\lambda,\zeta) \ \ & \in \ \ \Big\{\xi \in \Xi \; : \; \lambda d^{p}(\xi,\zeta) - \Psi(\xi) \, \leq \, \Phi(\lambda,\zeta) + \varepsilon, \ d^{p}(\xi,\zeta) \, \leq \, \Dm(\lambda,\zeta) + \delta\Big\}
\end{align*}
for $\nu$-almost all $\zeta \in \Xi$.

The first-order optimality conditions $\frac{\partial}{\partial \lambda-} h(\lambda^{\ast}) \leq 0$ and $\frac{\partial}{\partial \lambda+} h(\lambda^{\ast}) \geq 0$ imply that
\begin{equation}
\label{eqn:derivative bound}
\frac{\partial}{\partial \lambda+} \left(\int_{\Xi} \Phi(\lambda^{\ast},\zeta) \nu(d\zeta)\right) \ \ \leq \ \ \theta^{p}
\ \ \leq \ \ \frac{\partial}{\partial \lambda-} \left(\int_{\Xi} \Phi(\lambda^{\ast},\zeta) \nu(d\zeta)\right).
\end{equation}
Next we verify that we can interchange the partial derivative and integration in~\eqref{eqn:derivative bound}.
Recall from Lemma~\ref{lemma:phi}\ref{lemma:phi_monotone} that there is a set $B \in \scrB_{\nu}(\Xi)$ such that $\nu(B) = 1$, and $\Phi(\lambda,\zeta) > -\infty$ for all $\lambda > \kappa$ and all $\zeta \in B$, and $\Phi(\cdot,\zeta)$ is concave for all $\zeta \in B$.
To show it for the right derivative in~\eqref{eqn:derivative bound}, consider any $\zeta \in B$ and any decreasing sequence $\lambda^{n} \downarrow \lambda^{\ast}$.
Let
\[
f_{n}(\zeta) \ \ \defi \ \ \frac{\Phi(\lambda^{n},\zeta) - \Phi(\lambda^{\ast},\zeta)}{\lambda^{n} - \lambda^{\ast}}
\]
Since $\Phi(\cdot,\zeta)$ is nondecreasing for all $\zeta$, $f_{n}(\zeta) \ge 0$ for all $\zeta$.
In addition, since $\Phi(\cdot,\zeta)$ is concave, $f_{n} \leq f_{n+1}$ for all $n$.
Note that $\lim_{n \to \infty} f_{n}(\zeta) = \frac{\partial}{\partial \lambda+} \Phi(\lambda^{\ast},\zeta)$.
Thus it follows from the monotone convergence theorem that
\begin{equation}
\label{eqn:derivative interchange +}
\frac{\partial}{\partial \lambda+} \left(\int_{\Xi} \Phi(\lambda^{\ast},\zeta) \nu(d\zeta)\right)
\ \ = \ \ \int_{\Xi} \frac{\partial}{\partial \lambda+} \Phi(\lambda^{\ast},\zeta) \nu(d\zeta).
\end{equation}
To show it for the left derivative in~\eqref{eqn:derivative bound}, consider any $\zeta \in B$ and any increasing sequence $\lambda^{n} \uparrow \lambda^{\ast}$ with $\lambda^{1} > \kappa$.
Let $f_{n}$ be defined as before, and thus $f_{n}(\zeta) \ge 0$ for all $\zeta$.
In addition, since $\Phi(\cdot,\zeta)$ is concave, $f_{n} \geq f_{n+1}$ for all $n$.
That is, for all $n$ it holds that
\[
\left|f_{n}(\zeta)\right| \ \ = \ \ f_{n}(\zeta)
\ \ \le \ \ f_{1}(\zeta)
\ \ \le \ \ \frac{\left|\Phi(\lambda^{1},\zeta)\right| + \left|\Phi(\lambda^{\ast},\zeta)\right|}{\lambda^{\ast} - \lambda^{1}}.
\]
It follows from $\lambda^{1} > \kappa$ that
\[
\int_{\Xi} f_{1}(\zeta) \nu(d\zeta)
\ \ \le \ \ \frac{\int_{\Xi} \left|\Phi(\lambda^{1},\zeta)\right| \nu(d\zeta) + \int_{\Xi} \left|\Phi(\lambda^{\ast},\zeta)\right| \nu(d\zeta)}{\lambda^{\ast} - \lambda^{1}}
\ \ < \ \ \infty
\]
Also, $\lim_{n \to \infty} f_{n}(\zeta) = \frac{\partial}{\partial \lambda-} \Phi(\lambda^{\ast},\zeta)$.
Thus it follows from the dominated convergence theorem that
\begin{equation}
\label{eqn:derivative interchange -}
\frac{\partial}{\partial \lambda-} \left(\int_{\Xi} \Phi(\lambda^{\ast},\zeta) \nu(d\zeta)\right) \ \ = \ \ \int_{\Xi} \frac{\partial}{\partial \lambda-} \Phi(\lambda^{\ast},\zeta) \nu(d\zeta).
\end{equation}
%Then it follows from monotonicity and concavity of $\Phi(\cdot,\zeta)$, and from Lemma~\ref{lemma:phi}\ref{lemma:phi_derivative}, that $\Dp(\lambda_{1},\zeta) \geq f_{n}(\zeta) \geq f_{n+1}(\zeta) \geq 0$ for all $n$.
%Consider any $\lambda_{0} \in (\kappa,\lambda_{1})$.
%It follows from Lemma~\ref{lemma:phi}\ref{lemma:phi_bound} that
%\[
%f_{1}(\zeta) \ \ \leq \ \ \Dp(\lambda_{1},\zeta) \ \ \leq \ \ - \frac{1}{\lambda_{1} - \lambda_{0}} (\Psi(\zeta) + \Phi(\lambda_{0},\zeta)) \ \ \in \ \ L^{1}(\nu).
%\]
% It follows from  that
% \[
% \int_{\Xi} f_{1}(\zeta) \nu(d\zeta)
% ~\le ~\frac{\int_{\Xi} \left|\Phi(\lambda_{1},\zeta)\right| \nu(d\zeta) + \int_{\Xi} \left|\Phi(\lambda^{\ast},\zeta)\right| \nu(d\zeta)}{\lambda^{\ast} - \lambda_{1}}
% ~< ~\infty
% \]
% Also, $\lim_{n \to \infty} f_{n}(\zeta) = \frac{\partial}{\partial \lambda-} \Phi(\lambda^{\ast},\zeta)$.
%Thus, applying the reverse Fatou's lemma to the sequence $\{f_{n}\}_{n}$ gives
%\begin{equation}
%\label{eqn:derivative interchange -}
%\frac{\partial}{\partial \lambda-} \left(\int_{\Xi} \Phi(\lambda^{\ast},\zeta) \nu(d\zeta)\right) \ \ = \ \ \int_{\Xi} \frac{\partial}{\partial \lambda-} \Phi(\lambda^{\ast},\zeta) \nu(d\zeta).
%\end{equation}
Therefore it follows from~(\ref{eqn:derivative bound}), (\ref{eqn:derivative interchange +}), (\ref{eqn:derivative interchange -}), and Lemma~\ref{lemma:phi}\ref{lemma:phi_derivative} that
\begin{equation}
\label{eqn:first-order}
\begin{aligned}
\theta^{p} & \ \ \geq \ \ \frac{\partial}{\partial \lambda+} \left(\int_{\Xi} \Phi(\lambda^{\ast},\zeta) \nu(d\zeta)\right)
\ \ = \ \ \int_{\Xi} \frac{\partial}{\partial \lambda+} \Phi(\lambda^{\ast},\zeta) \nu(d\zeta)
\ \ \geq \ \ \int_{\Xi} \lim_{\lambda \downarrow \lambda^{\ast}} \Dp(\lambda,\zeta) \nu(d\zeta), \\
\theta^{p} & \ \ \leq \ \ \frac{\partial}{\partial \lambda-} \left(\int_{\Xi} \Phi(\lambda^{\ast},\zeta) \nu(d\zeta)\right)
\ \ = \ \ \int_{\Xi} \frac{\partial}{\partial \lambda-} \Phi(\lambda^{\ast},\zeta) \nu(d\zeta)
\ \ \leq \ \ \int_{\Xi} \lim_{\lambda \uparrow \lambda^{\ast}} \Dm(\lambda,\zeta) \nu(d\zeta).
\end{aligned}
\end{equation}
In particular, for any $\lambda_{1},\lambda_{2}$ with $\kappa < \lambda_{1} < \lambda^{\ast} < \lambda_{2}$, it follows from~(\ref{eqn:first-order}) and Lemma~\ref{lemma:phi}\ref{lemma:phi_monotone} that
\begin{equation}
\label{eqn:first-order2}
\begin{aligned}
\theta^{p} & \ \ \geq \ \ \int_{\Xi} \lim_{\lambda \downarrow \lambda^{\ast}} \Dp(\lambda,\zeta) \nu(d\zeta)
\ \ \geq \ \ \int_{\Xi} \Dp(\lambda_{2},\zeta) \nu(d\zeta)
\ \ \geq \ \ \int_{\Xi} \Dm(\lambda_{2},\zeta) \nu(d\zeta)
\ \ \geq \ \ \int_{\Xi} d^{p}(\underline{T}^{\varepsilon}_{\delta}(\lambda_{2},\zeta),\zeta) \nu(d\zeta) - \delta, \\
\theta^{p} & \ \ \leq \ \ \int_{\Xi} \lim_{\lambda \uparrow \lambda^{\ast}} \Dm(\lambda,\zeta) \nu(d\zeta)
\ \ \leq \ \ \int_{\Xi} \Dm(\lambda_{1},\zeta) \nu(d\zeta)
\ \ \leq \ \ \int_{\Xi} \Dp(\lambda_{1},\zeta) \nu(d\zeta)
\ \ \leq \ \ \int_{\Xi} d^{p}(\overline{T}^{\varepsilon}_{\delta}(\lambda_{1},\zeta),\zeta) \nu(d\zeta) + \delta.
\end{aligned}
\end{equation}
% Also,
% \begin{align*}
%   \lambda_{1} d^{p}(\overline{T}_{\lambda_{1}}(\zeta),\zeta) - \Psi(\overline{T}_{\lambda_{1}}(\zeta)) & \leq & \Phi(\lambda_{1},\zeta) + \delta, \\
%   \lambda_{2} d^{p}(\underline{T}_{\lambda_{2}}(\zeta),\zeta) - \Psi(\underline{T}_{\lambda_{2}}(\zeta)) & \leq & \Phi(\lambda_{2},\zeta) + \delta.
% \end{align*}
Based on \eqref{eqn:first-order2}, we now construct a feasible primal solution.
Note that there is a $q^{\varepsilon}_{\delta}(\lambda_{1},\lambda_{2}) \in [0,1]$ such that
\begin{equation}
\label{eqn:q1_theta}
q^{\varepsilon}_{\delta}(\lambda_{1},\lambda_{2}) \left[\int_{\Xi} d^{p}(\overline{T}^{\varepsilon}_{\delta}(\lambda_{1},\zeta),\zeta) \nu(d\zeta) + \delta\right] + \big(1 - q^{\varepsilon}_{\delta}(\lambda_{1},\lambda_{2})\big) \left[\int_{\Xi} d^{p}(\underline{T}^{\varepsilon}_{\delta}(\lambda_{2},\zeta),\zeta) \nu(d\zeta) - \delta\right] \ \ = \ \ \theta^{p}.
\end{equation}
Let $q^{\delta} \defi \theta^{p} / \left(\theta^{p} + \delta\right)$.
Define a distribution $\mu^{\varepsilon}_{\delta}(\lambda_{1},\lambda_{2})$ by
\begin{equation}
\label{eqn:mu_lambda>kappa}
\mu^{\varepsilon}_{\delta}(\lambda_{1},\lambda_{2}) \ \ \defi \ \ q^{\delta} q^{\varepsilon}_{\delta}(\lambda_{1},\lambda_{2}) \overline{T}^{\varepsilon}_{\delta}(\lambda_{1},\cdot)_{\#} \nu + q^{\delta} \big(1 - q^{\varepsilon}_{\delta}(\lambda_{1},\lambda_{2})\big) \underline{T}^{\varepsilon}_{\delta}(\lambda_{2},\cdot)_{\#} \nu + (1 - q^{\delta}) \nu.
\end{equation}
Then $\mu^{\varepsilon}_{\delta}(\lambda_{1},\lambda_{2})$ is primal feasible, because
\begin{align}
W_{p}^{p}(\mu^{\varepsilon}_{\delta}(\lambda_{1},\lambda_{2}),\nu) \ \ & \leq \ \ q^{\delta} q^{\varepsilon}_{\delta}(\lambda_{1},\lambda_{2}) \int_{\Xi} d^{p}(\overline{T}^{\varepsilon}_{\delta}(\lambda_{1},\zeta),\zeta) \nu(d\zeta) + q^{\delta} \big(1 - q^{\varepsilon}_{\delta}(\lambda_{1},\lambda_{2})\big) \int_{\Xi} d^{p}(\underline{T}^{\varepsilon}_{\delta}(\lambda_{2},\zeta),\zeta) \nu(d\zeta) \nonumber \\
& = \ \ q^{\delta} \bigg(\theta^{p} + \big[1 - 2 q^{\varepsilon}_{\delta}(\lambda_{1},\lambda_{2})\big] \delta\bigg)
\ \ \leq \ \ \theta^{p}.
\label{eqn:primal feasible 2}
\end{align}
Furthermore, recall that
\begin{eqnarray*}
& \lambda_{1} d^{p}(\overline{T}^{\varepsilon}_{\delta}(\lambda_{1},\zeta),\zeta) - \Phi(\lambda_{1},\zeta) - \varepsilon
\ \ \leq \ \ \Psi(\overline{T}^{\varepsilon}_{\delta}(\lambda_{1},\zeta))
\ \ \leq \ \ \lambda_{1} d^{p}(\overline{T}^{\varepsilon}_{\delta}(\lambda_{1},\zeta),\zeta) - \Phi(\lambda_{1},\zeta), \\
& \lambda_{2} d^{p}(\underline{T}^{\varepsilon}_{\delta}(\lambda_{2},\zeta),\zeta) - \Phi(\lambda_{2},\zeta) - \varepsilon
\ \ \leq \ \ \Psi(\underline{T}^{\varepsilon}_{\delta}(\lambda_{2},\zeta))
\ \ \leq \ \ \lambda_{2} d^{p}(\underline{T}^{\varepsilon}_{\delta}(\lambda_{2},\zeta),\zeta) - \Phi(\lambda_{2},\zeta).
\end{eqnarray*}
for $\nu$-almost all $\zeta \in \Xi$.
This, together with~(\ref{eqn:q1_theta}), implies that
\begin{align}
& \int_{\Xi} \Psi(\xi) \mu^{\varepsilon}_{\delta}(\lambda_{1},\lambda_{2})(d\xi) \nonumber \\
& = \ \ q^{\delta} q^{\varepsilon}_{\delta}(\lambda_{1},\lambda_{2}) \int_{\Xi} \Psi(\overline{T}^{\varepsilon}_{\delta}(\lambda_{1},\zeta)) \nu(d\zeta) \nonumber \\
& \qquad + q^{\delta} \big(1 - q^{\varepsilon}_{\delta}(\lambda_{1},\lambda_{2})\big) \int_{\Xi} \Psi(\underline{T}^{\varepsilon}_{\delta}(\lambda_{2},\zeta)) \nu(d\zeta) + (1 - q^{\delta}) \int_{\Xi} \Psi(\zeta) \nu(d\zeta) \nonumber \\
& \geq \ \ q^{\delta} q^{\varepsilon}_{\delta}(\lambda_{1},\lambda_{2}) \int_{\Xi} \big[\lambda_{1} d^{p}(\overline{T}^{\varepsilon}_{\delta}(\lambda_{1},\zeta),\zeta) - \Phi(\lambda_{1},\zeta) - \varepsilon\big] \nu(d\zeta) \nonumber \\
& \qquad + q^{\delta} \big(1 - q^{\varepsilon}_{\delta}(\lambda_{1},\lambda_{2})\big) \int_{\Xi} \big[\lambda_{2} d^{p}(\underline{T}^{\varepsilon}_{\delta}(\lambda_{2},\zeta),\zeta) - \Phi(\lambda_{2},\zeta) - \varepsilon\big] \nu(d\zeta) + (1 - q^{\delta}) \int_{\Xi} \Psi(\zeta) \nu(d\zeta) \nonumber \\
& \geq \ \ q^{\delta} \lambda_{1} \bigg[\theta^{p} + \big(1 - 2 q^{\varepsilon}_{\delta}(\lambda_{1},\lambda_{2})\big) \delta\bigg] - q^{\delta} q^{\varepsilon}_{\delta}(\lambda_{1},\lambda_{2}) \int_{\Xi} \Phi(\lambda_{1},\zeta) \nu(d\zeta) \nonumber \\
& \qquad - q^{\delta} \big(1 - q^{\varepsilon}_{\delta}(\lambda_{1},\lambda_{2})\big) \int_{\Xi} \Phi(\lambda_{2},\zeta) \nu(d\zeta) - q^{\delta} \varepsilon + (1 - q^{\delta}) \int_{\Xi} \Psi(\zeta) \nu(d\zeta).
\label{eqn:primal e-optimal 2}
\end{align}
Recall that $\Phi(\lambda,\zeta) \leq -\Psi(\zeta)$ for all $\zeta \in \Xi$.
Also, consider any $\lambda_{0} \in (\kappa,\lambda_{1})$.
Recall that $\Phi(\cdot,\zeta)$ is non-decreasing, and thus $\left|\Phi(\lambda,\zeta)\right| \ \ \le \ \ \left|\Psi(\zeta)\right| + \left|\Phi(\lambda_{0},\zeta)\right|$ for all $\lambda \geq \lambda_{0}$ and all $\zeta \in \Xi$.
Also, it follows from $\lambda_{0} > \kappa$ that $\int_{\Xi} \Phi(\lambda_{0},\zeta) \nu(d\zeta) > -\infty$.
Hence it follows from the dominated convergence theorem that
$\lim_{\lambda_{1} \uparrow \lambda^{\ast}} \int_{\Xi} \Phi(\lambda_{1},\zeta) \nu(d\zeta) = \int_{\Xi} \Phi(\lambda^{\ast},\zeta) \nu(d\zeta)$ and
$\lim_{\lambda_{2} \downarrow \lambda^{\ast}} \int_{\Xi} \Phi(\lambda_{2},\zeta) \nu(d\zeta) = \int_{\Xi} \Phi(\lambda^{\ast},\zeta) \nu(d\zeta)$.
Thus, given any $\varepsilon > 0$, choose $\lambda_{1}^{\varepsilon} \in \left(\max\{\kappa, \lambda^{\ast} - \varepsilon\}, \, \lambda^{\ast}\right)$ such that $\int_{\Xi} \Phi(\lambda_{1}^{\varepsilon},\zeta) \nu(d\zeta) \le \int_{\Xi} \Phi(\lambda^{\ast},\zeta) \nu(d\zeta) + \varepsilon$ and $\left(\lambda^{\ast} - \lambda_{1}^{\varepsilon}\right) \theta^{p} \le \varepsilon$, choose $\lambda_{2}^{\varepsilon} \in \left(\lambda^{\ast}, \, \lambda^{\ast} + \varepsilon\right)$ such that $\int_{\Xi} \Phi(\lambda_{2}^{\varepsilon},\zeta) \nu(d\zeta) \le \int_{\Xi} \Phi(\lambda^{\ast},\zeta) \nu(d\zeta) + \varepsilon$, and choose $\delta^{\varepsilon} \in (0,\varepsilon)$ such that $2 \lambda_{1}^{\varepsilon} \delta^{\varepsilon} \le \varepsilon$, $- \delta^{\varepsilon} \int_{\Xi} \Phi(\lambda^{\ast},\zeta) \nu(d\zeta) \le \theta^{p} \varepsilon$, and $- \delta^{\varepsilon} \int_{\Xi} \Psi(\zeta) \nu(d\zeta) \le \theta^{p} \varepsilon$.
Set
\[
T_{1}^{\varepsilon} \defi \overline{T}^{\varepsilon}_{\delta^{\varepsilon}}(\lambda_{1}^{\varepsilon},\cdot), \ \
T_{2}^{\varepsilon} \defi \underline{T}^{\varepsilon}_{\delta^{\varepsilon}}(\lambda_{2}^{\varepsilon},\cdot), \ \
p_{1}^{\varepsilon} \defi q^{\delta^{\varepsilon}} q^{\varepsilon}_{\delta^{\varepsilon}}(\lambda_{1}^{\varepsilon},\lambda_{2}^{\varepsilon}), \ \
p_{2}^{\varepsilon} \defi q^{\delta^{\varepsilon}} \big(1 - q^{\varepsilon}_{\delta^{\varepsilon}}(\lambda_{1}^{\varepsilon},\lambda_{2}^{\varepsilon})\big), \ \
p_{3}^{\varepsilon} \defi 1 - q^{\delta^{\varepsilon}}.
\]
Thus \eqref{eqn:primal feasible 1} follows from~\eqref{eqn:primal feasible 2}.
Also, it follows from~\eqref{eqn:primal e-optimal 2} that
\begin{align*}
& \int_{\Xi} \Psi(\xi) \mu^{\varepsilon}(d\xi) \nonumber \\
& \geq \ \ \frac{\theta^{p}}{\left(\theta^{p} + \delta^{\varepsilon}\right)} \lambda_{1}^{\varepsilon} \left(\theta^{p} - \delta^{\varepsilon}\right) - \frac{\theta^{p}}{\left(\theta^{p} + \delta^{\varepsilon}\right)} \left(\int_{\Xi} \Phi(\lambda^{\ast},\zeta) \nu(d\zeta) + \varepsilon\right) - \frac{\theta^{p}}{\left(\theta^{p} + \delta^{\varepsilon}\right)} \varepsilon + \frac{\delta^{\varepsilon}}{\left(\theta^{p} + \delta^{\varepsilon}\right)} \int_{\Xi} \Psi(\zeta) \nu(d\zeta) \\
& = \ \ \lambda^{\ast} \theta^{p} - \left(\lambda^{\ast} - \lambda_{1}^{\varepsilon}\right) \theta^{p} - 2 \frac{\delta^{\varepsilon}}{\left(\theta^{p} + \delta^{\varepsilon}\right)} \lambda_{1}^{\varepsilon} \theta^{p} - \int_{\Xi} \Phi(\lambda^{\ast},\zeta) \nu(d\zeta) + \frac{\delta^{\varepsilon}}{\left(\theta^{p} + \delta^{\varepsilon}\right)} \int_{\Xi} \Phi(\lambda^{\ast},\zeta) \nu(d\zeta) \\
& \qquad - 2 \frac{\theta^{p}}{\left(\theta^{p} + \delta^{\varepsilon}\right)} \varepsilon + \frac{\delta^{\varepsilon}}{\left(\theta^{p} + \delta^{\varepsilon}\right)} \int_{\Xi} \Psi(\zeta) \nu(d\zeta) \\
& \geq \ \ \lambda^{\ast} \theta^{p} - \int_{\Xi} \Phi(\lambda^{\ast},\zeta) \nu(d\zeta) - 6 \varepsilon.
\end{align*}
\hfillqed
\endproof

\begin{lemma}[Structure of $\varepsilon$-optimal primal solution with $\lambda^{\ast} = \kappa$]
\label{lem:e-optimal lambda = kappa}
Consider any $p \in [1,\infty)$, any $\nu \in \cP(\Xi)$, any $\theta > 0$, and any $\Psi \in L^{1}(\nu)$ such that $\kappa < \infty$.
Suppose that $\kappa$ is the unique minimizer of $h$.
If $\kappa = 0$, then, for any $\varepsilon > 0$, there are a $\lambda^{\varepsilon} \in \left(0, \varepsilon\right)$ and a $\nu$-measurable mapping $T^{\varepsilon} : \Xi \mapsto \Xi$, such that $\lambda^{\varepsilon} d^{p}(T^{\varepsilon}(\zeta),\zeta) - \Psi(T^{\varepsilon}(\zeta)) \leq \Phi(\lambda^{\varepsilon},\zeta) + \varepsilon$ for $\nu$-almost all $\zeta \in \Xi$, and
\[
\mu^{\varepsilon} \ \ \defi \ \ T^{\varepsilon}_{\#}\nu
\]
satisfies
\begin{align}
\int_{\Xi} \Psi(\xi) \mu^{\varepsilon}(d\xi) \ \ \ge \ \ \kappa \theta^{p} - \int_{\Xi} \Phi(\kappa,\zeta) \nu(d\zeta) - \varepsilon,
\label{eqn:primal e-optimal 3}
\end{align}
and
\begin{align}
W_{p}^{p}(\mu^{\varepsilon},\nu)
\ \ \leq \ \ \int_{\Xi} d^{p}(T^{\varepsilon}(\zeta),\zeta) \nu(d\zeta)
\ \ \leq \ \ \theta^{p}.
\label{eqn:primal feasible 3}
\end{align}
If $\kappa > 0$, then for any $\varepsilon > 0$, there is a $\lambda_{1}^{\varepsilon} \in \left(\kappa - \varepsilon, \, \kappa\right)$, a $\lambda_{2}^{\varepsilon} \in \left(\kappa, \, \kappa + \varepsilon\right)$, $\nu$-measurable mappings $T_{1}^{\varepsilon},T_{2}^{\varepsilon} : \Xi \mapsto \Xi$, and $p^{\varepsilon} \in (0,1)$, such that $\Psi(T_{1}^{\varepsilon}(\zeta)) - \Psi(\zeta) \geq \lambda_{1}^{\varepsilon} d^{p}(T_{1}^{\varepsilon}(\zeta),\zeta)$ and $\lambda_{2}^{\varepsilon} d^{p}(T_{2}^{\varepsilon}(\zeta),\zeta) - \Psi(T_{2}^{\varepsilon}(\zeta)) \leq \Phi(\lambda_{2}^{\varepsilon},\zeta) + \varepsilon$ for $\nu$-almost all $\zeta \in \Xi$, $\int_{\Xi} d^{p}(T_{1}^{\varepsilon}(\zeta),\zeta) \nu(d\zeta) > \theta^{p} + \theta^{p} / \varepsilon$, and
\[
\mu^{\varepsilon} \ \ \defi \ \ p^{\varepsilon} {T_{1}^{\varepsilon}}_{\#}\nu + \left(1 - p^{\varepsilon}\right) {T_{2}^{\varepsilon}}_{\#}\nu
\]
satisfies~\eqref{eqn:primal e-optimal 3} and
\begin{align}
W_{p}^{p}(\mu^{\varepsilon},\nu)
\ \ \leq \ \ p^{\varepsilon} \int_{\Xi} d^{p}(T_{1}^{\varepsilon}(\zeta),\zeta) \nu(d\zeta) + \left(1 - p^{\varepsilon}\right) \int_{\Xi} d^{p}(T_{2}^{\varepsilon}(\zeta),\zeta) \nu(d\zeta)
\ \ = \ \ \theta^{p}.
\label{eqn:primal feasible 4}
\end{align}
\end{lemma}

\proof{Proof of Lemma~\ref{lem:e-optimal lambda = kappa}.}
If $\kappa$ is the unique minimizer of $h$, then $h$ is increasing and convex on $[\kappa,\infty)$.
For any $\lambda > \kappa$, it follows from $h$ being increasing that
\begin{equation}
\label{eqn:lambda=kappa_ast}
%h(\lambda) \ \ > \ \ h(\lambda_{1})
%\ \ \ \Leftrightarrow \ \ \
\int_{\Xi} \left[\Phi(\lambda,\zeta) - \Phi(\kappa,\zeta)\right] \nu(d\zeta) \ \ < \ \ (\lambda - \kappa) \theta^{p}.
\end{equation}
Consider any $\varepsilon \in \left(0, (\lambda - \kappa) \theta^{p} - \int_{\Xi} \left[\Phi(\lambda,\zeta) - \Phi(\kappa,\zeta)\right] \nu(d\zeta)\right)$.
It follows from Lemma~\ref{lemma:measurable} that there exists a $\nu$-measurable map $\underline{T}_{\varepsilon}(\lambda,\cdot) : \Xi \mapsto \Xi$ such that $\lambda d^{p}(\underline{T}_{\varepsilon}(\lambda,\zeta),\zeta) - \Psi(\underline{T}_{\varepsilon}(\lambda,\zeta)) \leq \Phi(\lambda,\zeta) + \varepsilon$ for $\nu$-almost all $\zeta \in \Xi$.
Also, note that $\Phi(\kappa,\zeta) \leq \kappa d^{p}(\underline{T}_{\varepsilon}(\lambda,\zeta),\zeta) - \Psi(\underline{T}_{\varepsilon}(\lambda,\zeta))$.
Thus,
\[
\Phi(\lambda,\zeta) - \Phi(\kappa,\zeta) \ \ \geq \ \ (\lambda - \kappa) d^{p}(\underline{T}_{\varepsilon}(\lambda,\zeta),\zeta) - \varepsilon
\]
for $\nu$-almost all $\zeta \in \Xi$.
This together with \eqref{eqn:lambda=kappa_ast} yield that
\begin{align}
(\lambda - \kappa) \int_{\Xi} d^{p}(\underline{T}_{\varepsilon}(\lambda,\zeta),\zeta) \nu(d\zeta)
\ \ & \leq \ \ \int_{\Xi} \left[\Phi(\lambda,\zeta) - \Phi(\kappa,\zeta)\right] \nu(d\zeta) + \varepsilon
\ \ < \ \ (\lambda - \kappa) \theta^{p} \nonumber \\
\Rightarrow \ \ \ \int_{\Xi} d^{p}(\underline{T}_{\varepsilon}(\lambda,\zeta),\zeta) \nu(d\zeta)
\ \ & < \ \ \theta^{p}.
\label{eqn:primal feasible 5}
\end{align}
Hence, the distribution $\underline{T}_{\varepsilon}(\lambda,\cdot)_{\#}\nu$ is primal feasible.

Next, we separately consider the cases $\kappa = 0$ and $\kappa > 0$.
If $\kappa = 0$, then for the distribution
\begin{equation}
\label{eqn:mu_lambda=0}
\mu_{\varepsilon}(\lambda) \ \ = \ \ \underline{T}_{\varepsilon}(\lambda,\cdot)_{\#}\nu,
\end{equation}
it holds that
\begin{align}
\int_{\Xi} \Psi(\xi) \mu_{\varepsilon}(\lambda)(d\xi)
\ \ & = \ \ \int_{\Xi} \Psi(\underline{T}_{\varepsilon}(\lambda,\zeta)) \nu(d\zeta)
\ \ \geq \ \ \int_{\Xi} [\lambda d^{p}(\underline{T}_{\varepsilon}(\lambda,\zeta),\zeta) - \Phi(\lambda,\zeta) - \varepsilon] \nu(d\zeta) \nonumber \\
& \geq \ \ - \int_{\Xi} \Phi(\lambda,\zeta) \nu(d\zeta) - \varepsilon.
\label{eqn:primal e-optimal 4}
\end{align}
Thus, given any $\varepsilon > 0$, choose $\lambda^{\varepsilon} \in \left(0, \varepsilon\right)$ such that $\lambda^{\varepsilon} \theta^{p} \le \varepsilon$, and choose $\varepsilon' \in \left(0, \lambda^{\varepsilon} \theta^{p} - \int_{\Xi} \left[\Phi(\lambda^{\varepsilon},\zeta) - \Phi(0,\zeta)\right] \nu(d\zeta)\right)$ (note that $\varepsilon' \le \varepsilon$).
Set $T^{\varepsilon} \defi \underline{T}_{\varepsilon'}(\lambda^{\varepsilon},\cdot)$.
Thus~\eqref{eqn:primal feasible 3} follows from~\eqref{eqn:primal feasible 5}.
Also, it follows from~\eqref{eqn:primal e-optimal 4} that
\begin{align*}
\int_{\Xi} \Psi(\xi) \mu^{\varepsilon}(d\xi)
\ \ & \geq \ \ - \int_{\Xi} \Phi(\lambda^{\varepsilon},\zeta) \nu(d\zeta) - \varepsilon'
\ \ \geq \ \ \lambda^{\varepsilon} \theta^{p} - \int_{\Xi} \Phi(\lambda^{\varepsilon},\zeta) \nu(d\zeta) - 2 \varepsilon \\
& \geq \ \ \kappa \theta^{p} - \int_{\Xi} \Phi(\kappa,\zeta) \nu(d\zeta) - 2 \varepsilon.
\end{align*}

Otherwise, if $\kappa > 0$, then consider any $\kappa' \in (0,\kappa)$.
First, note that $\zeta \in \{\xi \in \Xi \, : \, \Psi(\xi) - \Psi(\zeta) \geq \kappa' d^{p}(\xi,\zeta)\}$, and hence $\{\xi \in \Xi \, : \, \Psi(\xi) - \Psi(\zeta) \geq \kappa' d^{p}(\xi,\zeta)\} \neq \varnothing$ for all $\zeta \in \Xi$.
Let
\[
\hat{D}(\kappa',\zeta) \ \ \defi \ \ \sup_{\xi \in \Xi} \{d^{p}(\xi,\zeta) \; : \; \Psi(\xi) - \Psi(\zeta) \geq \kappa' d^{p}(\xi,\zeta)\}.
\]
Next we show that $\int_{\Xi} \hat{D}(\kappa',\zeta) \nu(d\zeta) = \infty$.
Note that $\Psi(\xi) \leq \lambda d^{p}(\xi,\zeta) - \Phi(\lambda,\zeta)$ for all $\xi \in \Xi$.
Thus
\begin{eqnarray*}
\int_{\Xi} \Phi(\kappa',\zeta) \nu(d\zeta) & = & \int_{\Xi} \inf_{\xi \in \Xi} \Big\{\kappa' d^{p}(\xi,\zeta) - \Psi(\xi) \; : \; \Psi(\xi) - \Psi(\zeta) \geq \kappa' d^{p}(\xi,\zeta)\Big\} \nu(d\zeta) \\
& \geq & \int_{\Xi} \inf_{\xi \in \Xi} \Big\{-\Psi(\xi) \; : \; \Psi(\xi) - \Psi(\zeta) \geq \kappa' d^{p}(\xi,\zeta)\Big\} \nu(d\zeta) \\
& \geq & \int_{\Xi} \inf_{\xi \in \Xi} \Big\{-\lambda d^{p}(\xi,\zeta) + \Phi(\lambda,\zeta) \; : \; \Psi(\xi) - \Psi(\zeta) \geq \kappa' d^{p}(\xi,\zeta)\Big\} \nu(d\zeta) \\
& = & -\lambda \int_{\Xi} \hat{D}(\kappa',\zeta) \nu(d\zeta) + \int_{\Xi} \Phi(\lambda,\zeta) \nu(d\zeta).
\end{eqnarray*}
It follows from the definition of $\kappa$ that $\int_{\Xi} \Phi(\kappa',\zeta) \nu(d\zeta) = -\infty$, from $\lambda > \kappa$ that $\int_{\Xi} \Phi(\lambda,\zeta) \nu(d\zeta) > -\infty$, and thus $\int_{\Xi} \hat{D}(\kappa',\zeta) \nu(d\zeta) = \infty$.
Consider any $R > \theta^{p}$.
It follows that there exists $M \in L^{1}(\nu)$ such that $\int_{\Xi} M(\zeta) \nu(d\zeta) > R$ and
\[
\{\xi \in \Xi \; : \; \Psi(\xi) - \Psi(\zeta) \geq \kappa' d^{p}(\xi,\zeta),\ d^{p}(\xi,\zeta) \geq M(\zeta)\} \ \ \neq \ \ \varnothing
\]
for $\nu$-almost all $\zeta \in \Xi$.
Thus it follows from Lemma~\ref{lemma:measurable}\ref{itm:measurable_M} that for any $\kappa' \in (0,\kappa)$ and $R > \theta^{p}$, there exists a $\nu$-measurable mapping $T^{R}(\kappa',\cdot) : \Xi \mapsto \Xi$ such that
\[
\Psi\big(T^{R}(\kappa',\zeta)\big) - \Psi(\zeta) \ \ \geq \ \ \kappa' d^{p}(T^{R}(\kappa',\zeta),\zeta)
\]
for $\nu$-almost all $\zeta \in \Xi$, and
\[
\int_{\Xi} d^{p}(T^{R}(\kappa',\zeta),\zeta) \nu(d\zeta) \ \ \ge \ \ \int_{\Xi} M(\zeta) \nu(d\zeta) \ \ > \ \ R.
\]
If $\int_{\Xi} d^{p}(T^{R}(\kappa',\zeta),\zeta) \nu(d\zeta) = \infty$, let $\Xi_{r} \defi \{\zeta \in \Xi \, : \, d^{p}(T^{R}(\kappa',\zeta),\zeta) \leq r\}$.
Note that $\lim_{r \to \infty} \Xi_{r} = \Xi$, and thus there exists a $\bar{r} > R$ such that $R < \int_{\Xi_{\bar{r}}} d^{p}(T^{R}(\kappa',\zeta),\zeta) \nu(d\zeta) \leq \bar{r} < \infty$.
Then, let $\overline{T}^{R}(\kappa',\cdot) : \Xi \mapsto \Xi$ be given by $\overline{T}^{R}(\kappa',\zeta) \defi T^{R}(\kappa',\zeta)$ for all $\zeta \in \Xi_{\bar{r}}$, and $\overline{T}^{R}(\kappa',\zeta) \defi \zeta$ for all $\zeta \in \Xi \setminus \Xi_{\bar{r}}$.
Note that $\overline{T}^{R}(\kappa',\cdot)$ is $\nu$-measurable, and that $\int_{\Xi} d^{p}(\overline{T}^{R}(\kappa',\zeta),\zeta) \nu(d\zeta) = \int_{\Xi_{\bar{r}}} d^{p}(T^{R}(\kappa',\zeta),\zeta) \nu(d\zeta) \in (R,\infty)$.
Otherwise, if $\int_{\Xi} d^{p}(T^{R}(\kappa',\zeta),\zeta) \nu(d\zeta) < \infty$, then let $\Xi_{\bar{r}} = \Xi$, and $\overline{T}^{R}(\kappa',\cdot) \defi T^{R}(\kappa',\cdot)$.
Also in this case $\overline{T}^{R}(\kappa',\cdot)$ is $\nu$-measurable, and $\int_{\Xi} d^{p}(\overline{T}^{R}(\kappa',\zeta),\zeta) \nu(d\zeta) = \int_{\Xi_{\bar{r}}} d^{p}(T^{R}(\kappa',\zeta),\zeta) \nu(d\zeta) \in (R,\infty)$.
Let $q_{\varepsilon}^{R}(\kappa',\lambda) \in (0,1)$ be such that
\begin{equation}
q_{\varepsilon}^{R}(\kappa',\lambda) \int_{\Xi} d^{p}\big(\underline{T}_{\varepsilon}(\lambda,\zeta),\zeta\big) \nu(d\zeta) + \left(1 - q_{\varepsilon}^{R}(\kappa',\lambda)\right) \int_{\Xi} d^{p}\left(\overline{T}^{R}(\kappa',\zeta),\zeta\right) \nu(d\zeta) \ \ = \ \ \theta^{p}.
\label{eqn:convex combination = radius}
\end{equation}
Define a distribution $\mu_{\varepsilon}^{R}(\kappa',\lambda)$ by
\begin{equation}
\label{eqn:mu_lambda=kappa}
\mu_{\varepsilon}^{R}(\kappa',\lambda) \ \ \defi \ \ q_{\varepsilon}^{R}(\kappa',\lambda) \underline{T}_{\varepsilon}(\lambda,\cdot)_{\#}\nu + \left(1 - q_{\varepsilon}^{R}(\kappa',\lambda)\right) \overline{T}^{R}(\kappa',\cdot)_{\#}\nu.
\end{equation}
By construction, $\mu_{\varepsilon}^{R}(\kappa',\lambda)$ is primal feasible.
Also,
\begin{align}
& \int_{\Xi} \Psi(\xi) \mu_{\varepsilon}^{R}(\kappa',\lambda)(d\xi) \nonumber \\
& = \ \ q_{\varepsilon}^{R}(\kappa',\lambda) \int_{\Xi} \Psi(\underline{T}_{\varepsilon}(\lambda,\zeta)) \nu(d\zeta) + \left(1 - q_{\varepsilon}^{R}(\kappa',\lambda)\right) \int_{\Xi} \Psi\left(\overline{T}^{R}(\kappa',\zeta)\right) \nu(d\zeta) \nonumber \\
& \geq \ \ q_{\varepsilon}^{R}(\kappa',\lambda) \int_{\Xi} \big[\lambda d^{p}(\underline{T}_{\varepsilon}(\lambda,\zeta),\zeta) - \Phi(\lambda,\zeta) - \varepsilon\big] \nu(d\zeta) \nonumber \\
& \qquad + \left(1 - q_{\varepsilon}^{R}(\kappa',\lambda)\right) \int_{\Xi_{\bar{r}}} \left[\kappa' d^{p}\left(T^{R}(\kappa',\zeta),\zeta\right) + \Psi(\zeta)\right] \nu(d\zeta) + \left(1 - q_{\varepsilon}^{R}(\kappa',\lambda)\right) \int_{\Xi \setminus \Xi_{\bar{r}}} \Psi(\zeta) \nu(d\zeta) \nonumber \\
& \geq \ \ \kappa' \left(q_{\varepsilon}^{R}(\kappa',\lambda) \int_{\Xi} d^{p}(\underline{T}_{\varepsilon}(\lambda,\zeta),\zeta) \nu(d\zeta) + \left(1 - q_{\varepsilon}^{R}(\kappa',\lambda)\right) \int_{\Xi_{\bar{r}}} d^{p}(T^{R}(\kappa',\zeta),\zeta) \nu(d\zeta)\right) \nonumber \\
& \qquad - q_{\varepsilon}^{R}(\kappa',\lambda) \int_{\Xi} \Phi(\lambda,\zeta) \nu(d\zeta) - q_{\varepsilon}^{R}(\kappa',\lambda) \varepsilon + \left(1 - q_{\varepsilon}^{R}(\kappa',\lambda)\right) \int_{\Xi} \Psi(\zeta) \nu(d\zeta) \nonumber \\
& = \ \ \kappa' \theta^{p} - q_{\varepsilon}^{R}(\kappa',\lambda) \int_{\Xi} \Phi(\lambda,\zeta) \nu(d\zeta) - q_{\varepsilon}^{R}(\kappa',\lambda) \varepsilon + \left(1 - q_{\varepsilon}^{R}(\kappa',\lambda)\right) \int_{\Xi} \Psi(\zeta) \nu(d\zeta).
\label{eqn:primal e-optimal 5}
\end{align}
Thus, given any $\varepsilon > 0$, choose $\lambda_{1}^{\varepsilon} \in \left(\kappa - \varepsilon, \, \kappa\right)$ such that $(\kappa - \lambda_{1}^{\varepsilon}) \theta^{p} \le \varepsilon$, choose $\lambda_{2}^{\varepsilon} \in \left(\kappa, \, \kappa + \varepsilon\right)$ such that $(\lambda_{2}^{\varepsilon} - \kappa) \theta^{p} \le \varepsilon$, choose $\varepsilon' \in \left(0, (\lambda_{2}^{\varepsilon} - \kappa) \theta^{p} - \int_{\Xi} \left[\Phi(\lambda_{2}^{\varepsilon},\zeta) - \Phi(\kappa,\zeta)\right] \nu(d\zeta)\right)$ such that $\varepsilon' \left|\int_{\Xi} \Phi(\lambda_{2}^{\varepsilon},\zeta) \nu(d\zeta)\right| \le \varepsilon$ and $\varepsilon' \left|\int_{\Xi} \Psi(\zeta) \nu(d\zeta)\right| \le \varepsilon$ (note that $\varepsilon' \le \varepsilon$), and choose $R \ge \theta^{p} + \theta^{p} / \varepsilon'$.
Set
\[
T_{1}^{\varepsilon} \ \defi \ \overline{T}^{R}(\lambda_{1}^{\varepsilon},\cdot), \quad
T_{2}^{\varepsilon} \ \defi \ \underline{T}_{\varepsilon'}(\lambda_{2}^{\varepsilon},\cdot), \quad
p^{\varepsilon} \ \defi \ 1 - q_{\varepsilon'}^{R}(\lambda_{1}^{\varepsilon},\lambda_{2}^{\varepsilon}).
\]
Thus~\eqref{eqn:primal feasible 4} follows from~\eqref{eqn:convex combination = radius}.
It also follows from~\eqref{eqn:convex combination = radius} that
\[
p^{\varepsilon} \ \ = \ \ \frac{\theta^{p} - \int_{\Xi} d^{p}\big(\underline{T}_{\varepsilon'}(\lambda_{2}^{\varepsilon},\zeta),\zeta\big) \nu(d\zeta)}{\int_{\Xi} d^{p}\left(\overline{T}^{R}(\lambda_{1}^{\varepsilon},\zeta),\zeta\right) \nu(d\zeta) - \int_{\Xi} d^{p}\big(\underline{T}_{\varepsilon'}(\lambda_{2}^{\varepsilon},\zeta),\zeta\big) \nu(d\zeta)}
\ \ \le \ \ \dfrac{\theta^{p}}{\theta^{p} + \dfrac{\theta^{p}}{\varepsilon'} - \theta^{p}}
\ \ = \ \ \varepsilon'.
\]
Thus it follows from~\eqref{eqn:primal e-optimal 5} that
\begin{align*}
& \int_{\Xi} \Psi(\xi) \mu^{\varepsilon}(d\xi) \nonumber \\
& \ge \ \ \lambda_{1}^{\varepsilon} \theta^{p} - \left(1 - p^{\varepsilon}\right) \int_{\Xi} \Phi(\lambda_{2}^{\varepsilon},\zeta) \nu(d\zeta) - \left(1 - p^{\varepsilon}\right) \varepsilon + p^{\varepsilon} \int_{\Xi} \Psi(\zeta) \nu(d\zeta) \\
& = \ \ \lambda_{2}^{\varepsilon} \theta^{p} - \left(\lambda_{2}^{\varepsilon} - \kappa\right) \theta^{p} - \left(\kappa - \lambda_{1}^{\varepsilon}\right) \theta^{p} - \int_{\Xi} \Phi(\lambda_{2}^{\varepsilon},\zeta) \nu(d\zeta) + p^{\varepsilon} \int_{\Xi} \Phi(\lambda_{2}^{\varepsilon},\zeta) \nu(d\zeta) - \left(1 - p^{\varepsilon}\right) \varepsilon + p^{\varepsilon} \int_{\Xi} \Psi(\zeta) \nu(d\zeta) \\
& \ge \ \ \lambda_{2}^{\varepsilon} \theta^{p} - \int_{\Xi} \Phi(\lambda_{2}^{\varepsilon},\zeta) \nu(d\zeta) - \varepsilon' \left|\int_{\Xi} \Phi(\lambda_{2}^{\varepsilon},\zeta) \nu(d\zeta)\right| - \varepsilon' \left|\int_{\Xi} \Psi(\zeta) \nu(d\zeta)\right| - 3 \varepsilon \\
& \ge \ \ \kappa \theta^{p} - \int_{\Xi} \Phi(\kappa,\zeta) \nu(d\zeta) - 5 \varepsilon.
\end{align*}
\hfillqed
\endproof

The next theorem establishes a strong duality result when the growth rate $\kappa$ is finite.
The result follows from Lemmas~\ref{lem:e-optimal lambda > kappa} and~\ref{lem:e-optimal lambda = kappa} and by noting that
\[
v_{P} \ \ \ge \ \ \int_{\Xi} \Psi(\xi) \mu^{\varepsilon}(d\xi)
\ \ \ge \ \ \lambda^{\ast} \theta^{p} - \int_{\Xi} \Phi(\lambda^{\ast},\zeta) \nu(d\zeta) - \varepsilon
\ \ = \ \ v_{D} - \varepsilon
\]
for any $\varepsilon > 0$.

\begin{theorem}[Strong duality with finite optimal value]
\label{thm:strongDual}
Consider any $p \in [1,\infty)$, any $\nu \in \cP(\Xi)$, any $\theta > 0$, and any $\Psi \in L^{1}(\nu)$ such that $\kappa < \infty$.
Then $v_{P} = v_{D} < \infty$.
\end{theorem}

\begin{remark}
\label{rmk:cost}
All the above results and proofs in Section~\ref{sec:dual_general} (except for Lemma~\ref{lemma:kappa}) continue to hold if we replace the transportation cost $d^{p}(\cdot,\cdot)$ with any measurable, non-negative cost function $c(\cdot,\cdot)$ that satisfies $c(\xi,\zeta) = 0$ if $\xi = \zeta$.
\end{remark}

%\begin{remark}
%\label{rmk:near-optimal distribution}
%It follows from the proof of Theorem~\ref{thm:strongDual} (see equations \eqref{eqn:mu_lambda>kappa}, \eqref{eqn:mu_lambda=0}, and~\eqref{eqn:mu_lambda=kappa}), that if $\kappa < \infty$,
% if $\lambda^{\ast} > \kappa$,
%then for any $\varepsilon > 0$,
% by choosing $\lambda_{1}, \lambda_{2} \in (\kappa, \lambda^{\ast})$ sufficiently close to $\lambda^{\ast}$, 
% $\lambda_{2} > \lambda^{\ast}$ sufficiently close to $\lambda^{\ast}$, 
% and $\delta,\varepsilon > 0$ sufficiently small, there are $\nu$-measurable mappings $\overline{T}(\lambda_{1},\cdot),\underline{T}(\lambda_{2},\cdot) : \Xi \mapsto \Xi$ such that $\mu^{\varepsilon}_{\delta}(\lambda_{1},\lambda_{2})$ given in~(\ref{eqn:mu_lambda>kappa}) is primal feasible and $\epsilon$-optimal.
% That is, for any $\epsilon > 0$, 
%there is an approximately optimal primal solution $\mu^{\varepsilon}$ of the form
%\[
%\mu^{\varepsilon} \ \ = \ \ p_{1} \underline{T}_{\#} \nu + p_{2} \overline{T}_{\#} \nu + (1 - p_{1} - p_{2}) \nu,
%\]
%where $p_{1},p_{2} \ge 0$, $p_{1} + p_{2} \le 1$, and $\underline{T},\overline{T} : \Xi \mapsto \Xi$ are $\nu$-measurable mappings, such that
%\[
%p_{1} \int_{\Xi} d^{p}\left(\underline{T}(\zeta),\zeta\right) \nu(d\zeta) + p_{2} \int_{\Xi} d^{p}\left(\overline{T}(\zeta),\zeta\right) \nu(d\zeta) \ \ \le \ \ \theta^{p}.
%\]
%\end{remark}

Next, we investigate existence conditions for worst-case distributions and their structure.
In the remainder of this section, we assume that $\Psi$ is upper-semi-continuous, and every bounded subset in $(\Xi,d)$ is totally bounded (see, for example, Section 45 in \cite{munkres2000topology}), which is satisfied by, for example, any finite-dimensional normed space.
Under this assumption, if $\lambda > \kappa$ and $\zeta \in B$, then Lemma~\ref{lemma:phi}\ref{lemma:phi_bound} and the upper semi-continuity of $\Psi$ imply that the set $\argmin_{\xi \in \Xi}\{\lambda d^{p}(\xi,\zeta) - \Psi(\xi)\}$ is nonempty, and that $\min/\max_{\xi \in \Xi} \{d^{p}(\xi,\zeta) \, : \, \lambda d^{p}(\xi,\zeta) - \Psi(\xi) = \Phi(\lambda,\zeta)\}$ can be attained.
If $\lambda = \kappa$ and $\nu\big(\{\zeta \in \Xi \, : \, \argmin_{\xi \in \Xi} \{\kappa d^{p}(\xi,\zeta) - \Psi(\xi)\} = \varnothing\}\big) = 0$, then the upper semi-continuity of $\Psi$ imply that $\min_{\xi \in \Xi} \{d^{p}(\xi,\zeta) \, : \, \kappa d^{p}(\xi,\zeta) - \Psi(\xi) = \Phi(\kappa,\zeta)\}$ can be attained for $\nu$-almost all $\zeta \in \Xi$, but $\sup_{\xi \in \Xi} \{d^{p}(\xi,\zeta) \, : \, \kappa d^{p}(\xi,\zeta) - \Psi(\xi) = \Phi(\kappa,\zeta)\}$ can be infinite.
Thus, if (i) $\lambda > \kappa$, or (ii) $\lambda = \kappa$ and $\nu\big(\{\zeta \in \Xi \, : \, \argmin_{\xi \in \Xi} \{\kappa d^{p}(\xi,\zeta) - \Psi(\xi)\} = \varnothing\}\big) = 0$, then the quantities
$\underline{D}_{0}(\lambda,\zeta)$ and $\overline{D}_{0}(\lambda,\zeta)$ in~(\ref{eqn:D0}) are well-defined for $\nu$-almost all $\zeta \in \Xi$ (where $\overline{D}_{0}(\lambda,\zeta)$ can be infinite if $\lambda = \kappa$).
%Then $\underline{D}_{0}(\lambda,\zeta)$ and $\overline{D}_{0}(\lambda,\zeta)$ represent respectively the closest and furthest distances between $\zeta$ and any point in $\argmin_{\xi \in \Xi} \{\lambda d^{p}(\xi,\zeta) - \Psi(\xi)\}$.
% , and are finite when $\lambda > \kappa$.  In addition, if $\Phi(\kappa,\zeta)$ is finite, then $\underline{D}_{0}(\lambda,\zeta)$ is also finite (but $\overline{D}_{0}(\lambda,\zeta)$ can be infinite).
%Note that $\overline{D}_{0}(\lambda,\zeta)$ (resp. $\underline{D}_{0}(\lambda,\zeta)$) may not be equal to $\overline{D}(\lambda,\zeta)$ (resp. $\underline{D}(\lambda,\zeta)$) as defined in~(\ref{eqn:Dpm}).

\begin{corollary}[Worst-case distribution]
\label{cor:existence}
Consider any $p \in [1,\infty)$, $\nu \in \cP(\Xi)$, $\theta > 0$, and $\Psi \in L^{1}(\nu)$ such that $\kappa < \infty$.
Assume that $\Psi$ is upper-semi-continuous, and that bounded subsets of $(\Xi,d)$ are totally bounded. Then the following holds:
\begin{enumerate}[label=(\roman*)]
\item
\label{itm:iff}
[Existence condition]
A worst-case distribution exists if and only if any of the following conditions hold:
\begin{enumerate}[leftmargin=*,label=(\alph*)]
\item
\label{itm:existenceCond1}
There exists a dual minimizer $\lambda^{\ast} > \kappa$.
\item
\label{itm:existenceCond2}
$\lambda^{\ast} = \kappa > 0$ is the unique dual minimizer, $\nu\big(\{\zeta \in \Xi \, : \, \argmin_{\xi \in \Xi} \{\kappa d^{p}(\xi,\zeta) - \Psi(\xi)\} = \varnothing\}\big) = 0$, and
\[
\int_{\Xi} \underline{D}_{0}(\kappa,\zeta) \nu(d\zeta) \ \ \leq \ \ \theta^{p} \ \ \leq \ \ \int_{\Xi} \overline{D}_{0}(\kappa,\zeta) \nu(d\zeta).
\]
\item
\label{itm:existenceCond3}
$\lambda^{\ast} = \kappa = 0$ is the unique dual minimizer, $\argmax_{\xi \in \Xi} \{\Psi(\xi)\}$ is nonempty, and
\[
\int_{\Xi} \underline{D}_{0}(0,\zeta) \nu(d\zeta) \ \ \leq \ \ \theta^{p}.
\]
\end{enumerate}

%\item
%\label{itm:small}
In addition, if $\nu\big(\zeta \in \Xi \, : \, -\Psi(\zeta) > \inf_{\xi \in \Xi} \left\{\kappa d^{p}(\xi,\zeta) - \Psi(\xi)\right\}\big) = 0$, then $\lambda^{\ast} = \kappa$ for any $\theta > 0$.
Otherwise, there is $\theta_{0} > 0$ such that $\lambda^{\ast} > \kappa$ for any $\theta < \theta_{0}$.

\item
\label{itm:form}
[Structure]
Whenever a worst-case distribution exists, there exists a worst-case distribution $\mu^{\ast}$ which can be represented as a convex combination of two distributions $\overline{T}^{\ast}_{\#} \nu$ and $\underline{T}^{\ast}_{\#} \nu$, each of which is a perturbation of $\nu$, as follows:
\[
\mu^{\ast} \ \ = \ \ p^{\ast} \overline{T}^{\ast}_{\#} \nu + (1 - p^{\ast}) \underline{T}^{\ast}_{\#} \nu,
\]
where $p^{\ast} \in [0,1]$, and $\overline{T}^{\ast}, \underline{T}^{\ast} : \Xi \mapsto \Xi$ satisfy
\begin{equation}
\label{eqn:transportMap}
\nu\left(\left\{\zeta \in \Xi \; : \; \overline{T}^{\ast}(\zeta), \underline{T}^{\ast}(\zeta) \notin \argmin_{\xi \in \Xi} \{\lambda^{\ast} d^{p}(\xi,\zeta) - \Psi(\xi)\}\right\}\right) \ \ = \ \ 0.
\end{equation}
Let $\gamma^{T} \in \cP(\Xi \times \Xi)$ be the joint distribution given by $\overline{T}^{\ast}, \underline{T}^{\ast}$ and $p^{\ast}$.
For any measurable set $A \subset \Xi \times \Xi$, let
\[
\gamma^{T}(A) \ \ \defi \ \ p^{\ast} \nu\left(\left\{\zeta \; : \; \big(\zeta,\overline{T}^{\ast}(\zeta)\big) \in A\right\}\right) + \left(1 - p^{\ast}\right) \nu\left(\Big\{\zeta \; : \; \big(\zeta,\underline{T}^{\ast}(\zeta)\big) \in A\Big\}\right),
\]
If $\lambda^{\ast} > 0$, then $\gamma^{T}$ is an optimal joint distribution in the definition~(\ref{eqn:def_wasserstein}) of $W_{p}^{p}(\mu^{\ast},\nu)$.

\item
\label{itm:strongDual_concave}
Suppose $\Xi$ is convex, $\Psi$ is a concave function, and $d^{p}(\cdot,\zeta)$ is a convex function for $\nu$-almost all $\zeta \in \Xi$, then
\[
v_{P} \ \ = \ \ v_{D} \ \ = \ \ \sup_{T : \Xi \mapsto \Xi} \left\{\E_{T_{\#}\nu}[\Psi(\xi)] \; : \; W_{p}(T_{\#}\nu,\nu) \leq \theta, \; T \textrm{ is } \nu\textrm{-measurable}\right\}.
\]
Moreover, whenever a worst-case distribution exists, there exists $T^{\ast} : \Xi \mapsto \Xi$ such that $T^{\ast}_{\#} \nu$ is primal optimal and
\[
\nu\left(\left\{\zeta \in \Xi \; : \; T^{\ast}(\zeta) \notin \argmin_{\xi \in \Xi} \{\lambda^{\ast} d^p(\xi,\zeta) - \Psi(\xi)\}\right\}\right) \ \ = \ \ 0.
\]
% \begin{equation}
%   \label{eqn:dual_frakM'}
%   v_{P} \ = \ v_{D} \ = \ \sup_{\mu \in \frakM'} \E_{\mu}[\Psi(\xi)],
% \end{equation}
% where
% \begin{equation}
%   \label{eqn:frakM'}
%   \frakM' \ \defi \ \left\{\mu = T_{\#}\nu \ \Big| \ T : \Xi \mapsto \Xi, \; \int_{\Xi} d^{p}(T(\zeta),\zeta) \nu(d\zeta) \leq \theta^{p}\right\}.
% \end{equation}
\end{enumerate}
\end{corollary}

\begin{remark}
Compared with Corollary~4.7 in \citet{esfahani2015data}, Corollary~\ref{cor:existence}\ref{itm:iff} provides a complete description of the necessary and sufficient conditions for the existence of a worst-case distribution. Note that Example~1 in \citet{esfahani2015data} corresponds to $\lambda^{\ast} = \kappa = 1$ and $p = 1$.
We also remark that \citet{yue2021linear} derives a sufficient condition on the existence of the worst-case distribution that does not require the knowledge of $\lambda^\ast$.
\end{remark}

% \noindent
% \begin{remark}
% When $\lambda^{\ast} > \kappa$, if $\Psi$ is upper semi-continuous and every bounded subset of $\Xi$ is \textit{totally bounded}, then the set $\argmin_{\xi \in \Xi}\{\lambda^{\ast} d^{p}(\xi,\zeta) - \Psi(\xi)\}$ is non-empty.
% Indeed, recall that Lemma \ref{lemma:phi}\ref{lemma:phi_bound} suggests that the set $\argmin_{\xi\in\Xi}\{ \lambda d^{p}(\xi,\zeta)-\Psi(\xi)\}$ is bounded for all $\lambda>\kappa$. If it is also totally bounded (which holds, for example, $\Xi=\R^{K}$ or a subset of $\R^{K}$), then together with the completeness of $\Xi$ would imply the compactness of the set (Theorem 45.1 in \citet{munkres2000topology}).
% 	  This, together with upper semicontinuity of $\Psi$, implies that the set is non-empty.
% \end{remark}
% \vspace{0.5em}

% \begin{definition}[Optimal transport maps $\underline{T}^{\ast}$ and $\overline{T}^{\ast}$]\label{def:T_{p}m}
%   The optimal transport map $\underline{T}^{\ast},\overline{T}^{\ast} : \Xi \mapsto \Xi$ are defined by \eqref{eqn:transportMap}.
% \end{definition}

\begin{example}
We present several examples that correspond to different cases in Corollary~\ref{cor:existence}\ref{itm:iff}.
In all these examples, $\Xi = [0,\infty)$, $d(\xi,\zeta) = |\xi - \zeta|$ for all $\xi,\zeta \in \Xi$, $p = 1$, $\theta > 0$, and $\nu = \delta_{0}$.

\begin{figure}[h]
\subfloat[$\Psi_{a}(\xi) = \max\{0, \xi - a\}$]{
\label{fig:psi_abs} %% label for first subfigure
\begin{minipage}[b]{0.3\linewidth}
\centering
\includegraphics[height=1.2in]{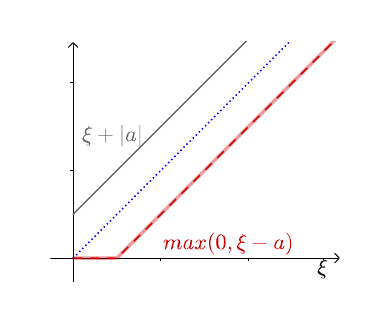}
\end{minipage}}%
\subfloat[$\Psi(\xi) = \max\{0, 1 - \xi^{2}\}$]{
\label{fig:psi_parabola} %% label for second subfigure
\begin{minipage}[b]{0.3\linewidth}
\centering
\includegraphics[height=1.2in]{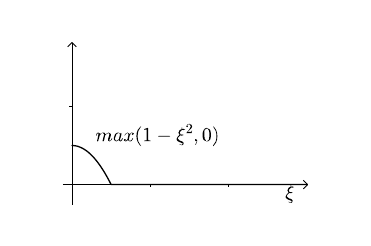}
\end{minipage}}
\subfloat[$\Psi_{\pm}(\xi) = 1 + \xi \pm \frac{1}{\xi+1}$]{
\label{fig:psi_asymptote} %% label for second subfigure
\begin{minipage}[b]{0.4\linewidth}
\centering
\includegraphics[height=1.2in]{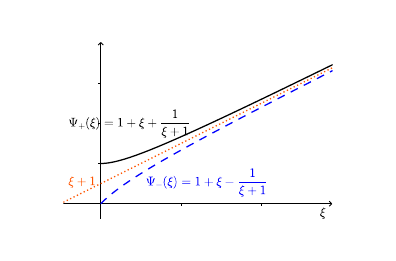}
\end{minipage}}
\caption{Examples for existence and non-existence of the worst-case distribution}
\label{fig:existence}
\end{figure}

\begin{enumerate}[label=\arabic*.,itemsep=1ex,topsep=1ex]
\item
$\Psi_{a}(\xi) \defi \max\{0, \xi - a\}$ for some $a \in \R$.  It follows that $\lambda^{\ast} = \kappa = 1$.
\begin{itemize}[leftmargin=*,label=--]
\item
If $a \leq 0$, then $\argmin_{\xi \in \Xi} \{d^{p}(\xi,0) - \Psi_{a}(\xi)\} = [0,\infty)$, hence $\underline{D}_{0}(\kappa,\zeta) = 0$ and $\overline{D}_{0}(\kappa,\zeta) = \infty$.
Thus condition~\ref{itm:existenceCond2} is satisfied. One of the worst-case distributions is $\mu^{\ast} = \delta_{\theta}$ with $v_{P} = v_{D} = \theta - a$.
\item
If $a > 0$, then $\argmin_{\xi \in \Xi} \{d^{p}(\xi,0) - \Psi_{a}(\xi)\} = \{0\}$, hence $\underline{D}_{0}(\kappa,\zeta) = \overline{D}_{0}(\kappa,\zeta) = 0 < \theta$.
Thus condition~\ref{itm:existenceCond2} is violated. There is no worst-case distribution, but the objective value of $\mu_{\varepsilon} = (1 - \varepsilon) \delta_{0} + \varepsilon \delta_{\theta/\varepsilon}$ converges to $v_{P} = v_{D} = \theta$ as $\varepsilon \to 0$.
\end{itemize}

\item
$\Psi(\xi) = \max\{0, 1 - \xi^{2}\}$. It follows that $\lambda^{\ast} = \kappa = 0$, and $\argmax_{\xi \in \Xi} \Psi(\xi) = \{0\}$.
Thus condition~\ref{itm:existenceCond3} is satisfied, and the worst-case distribution is $\mu^{\ast} = \delta_{0} = \nu$.

\item
$\Psi_{\pm}(\xi) = 1 + \xi \pm \frac{1}{\xi+1}$.  It follows that $\kappa = 1$.  Note that $\Psi_{\pm}'(\xi) = 1 \mp \frac{1}{(\xi+1)^{2}}$.
\begin{itemize}[leftmargin=*,label=--]
\item
Note that $\Psi_{+}'(\xi) < \kappa = 1$ on $\Xi$.
Also, $\Psi_{+}$ satisfies the condition in \ref{itm:small}, thus for all $\theta > 0$ it holds that $\lambda_{+}^{\ast} = \kappa = 1$ and $\argmin_{\xi \in \Xi} \{\lambda_{+}^{\ast} d^{p}(\xi,0) - \Psi_{+}(\xi)\} = \{0\}$.  There is no worst-case distribution, but the objective value of $\mu_{\varepsilon} = (1 - \varepsilon) \delta_{0} + \varepsilon \delta_{\theta/\varepsilon}$ converges to $v_{P} = v_{D} = 2 + \theta$ as $\varepsilon \to 0$.
\item
Note that $\Psi_{-}'(\xi) > \kappa = 1$ on $\Xi$.
Also, $\argmin_{\lambda \geq 0} \left\{\lambda \theta - \inf_{\xi \in \Xi} \left\{\lambda \xi - \left(1 + \xi - \frac{1}{\xi+1}\right)\right\}\right\} = \argmin_{\lambda \geq 1} \left\{\lambda (\theta + 1) - 2 \sqrt{\lambda - 1}\right\} = \left\{1 + \frac{1}{(\theta+1)^{2}}\right\}$.  Thus $\lambda_{-}^{\ast} > 1 = \kappa$.
\end{itemize}
\end{enumerate}
\end{example}

\subsection{Finite-Supported Nominal Distribution}
\label{sec:dual_finite}

In this section, we consider the setting in which the nominal distribution has finite support.  Let $\nu = \frac{1}{N} \sum_{i=1}^{N} \delta_{\hxi^{i}}$ for some $\hxi^{i} \in \Xi$, $i = 1,\ldots,N$.  This occurs, for example, in data-driven settings in which the nominal distribution is given by an empirical distribution constructed with $N$ observations.

\begin{corollary}[Data-Driven DRSO]
\label{cor:finite}
Consider any $\nu = \frac{1}{N} \sum_{i=1}^{N} \delta_{\hxi^{i}}$, $p \in [1,\infty)$, and $\theta > 0$.
The following hold:
\begin{enumerate}[label=(\roman*),itemsep=6pt]
\item{[Strong duality]}
\label{itm:finiteDual_lambda}
The primal problem~(\ref{eqn:problem_primal}) has a strong dual problem
\begin{equation}
\label{eqn:dual_finite}
v_{P} \ \ = \ \ v_{D} \ \ = \ \ \inf_{\lambda \geq 0} \left\{\lambda \theta^{p} - \frac{1}{N} \sum_{i=1}^{N} \inf_{\xi \in \Xi} \big[\lambda d^{p}(\xi,\hxi^{i}) - \Psi(\xi)\big]\right\}.
\end{equation}
% Moreover, $v_{P} = v_{D}$ is also equal to
% \begin{equation}
%   \label{eqn:2Nprogram}
%   \sup_{\substack{\underline{\xi}^{i},\overline{\xi}^{i} \in \Xi, i=1,\ldots,N \\ q \in [0,1]}} \bigg\{\frac{1}{N} \sum_{i=1}^{N} \Big[q \Psi(\underline{\xi}^{i}) + (1 - q) \Psi(\overline{\xi}^{i})\big] \; : \; \frac{1}{N} \sum_{i=1}^{N} \big[q d^{p}(\underline{\xi}^{i},\hxi^{i}) + (1 - q) d^{p}(\overline{\xi}^{i},\hxi^{i}) \big] \, \leq \, \theta^{p}\bigg\}.
% \end{equation}
% \item
%   \label{itm:finiteDual_N}
%   If $\kappa < \infty$, $\Xi$ is convex, and $\Psi$ is concave, then (\ref{eqn:2Nprogram}) is further simplified to
%   \begin{equation}
%     \label{eqn:dual_robust}
%     \sup_{\xi^{i} \in \Xi, i=1,\ldots,N} \left\{\frac{1}{N} \sum_{i=1}^{N} \Psi(\xi^{i}) \; : \; \frac{1}{N} \sum_{i=1}^{N} d^{p}(\xi^{i},\hxi^{i}) \, \leq \, \theta\right\}.
%   \end{equation}
\item{[Structure of the worst-case distribution]}
\label{itm:finiteDual_N+1}
Whenever a worst-case distribution exists, there exists one which is supported on at most $N+1$ points and has the form
\begin{equation}
\label{eqn:finiteDual_N+1}
\mu^{\ast} \ \ = \ \ \frac{1}{N} \sum_{i \neq i_{0}} \delta_{\xi^{i}_{\ast}} + \frac{p_{0}}{N} \delta_{\underline{\xi}^{i_{0}}_{\ast}} + \frac{1 - p_{0}}{N} \delta_{\overline{\xi}^{i_{0}}_{\ast}},
\end{equation}
where $i_{0} \in \{1,\ldots,N\}$, $p_{0} \in [0,1]$, $\underline{\xi}^{i_{0}}_{\ast},\overline{\xi}^{i_{0}}_{\ast} \in \argmin_{\xi \in \Xi} \{\lambda^{\ast} d^{p}(\xi,\hxi^{i_{0}}) - \Psi(\xi)\}$, and $\xi^{i}_{\ast} \in \argmin_{\xi \in \Xi} \{\lambda^{\ast} d^{p}(\xi,\hxi^{i}) - \Psi(\xi)\}$ for all $i \neq i_{0}$.
\item{[Approximation by robust optimization]}
\label{itm:finiteDual_robustapprox}
Suppose that there exists $\zeta^{0} \in \Xi$, $L, M \geq 0$ such that $|\Psi(\xi) - \Psi(\zeta^{0})| < L d^{p}(\xi,\zeta^{0}) + M$ for all $\xi \in \Xi$.
For any positive integer $K$, consider the robust optimization problem
\[
v_{K} \ \ \defi \ \ \sup_{(\xi^{ik})_{i,k} \in \frakM_{K}} \frac{1}{NK} \sum_{i=1}^{N} \sum_{k=1}^{K} \Psi(\xi^{ik}),
\]
with uncertainty set
\[
\label{eqn:frakM_K}
\frakM_{K} \ \ \defi \ \ \bigg\{(\xi^{ik})_{i,k} \; : \; \frac{1}{NK} \sum_{i=1}^{N} \sum_{k=1}^{K} d^{p}(\xi^{ik},\hxi^{i}) \leq \theta^{p}, \, \xi^{ik} \in \Xi \; \forall \; i,k\bigg\}.
\]
% Then $v_{K} \uparrow \sup_{\mu \in \frakM} \E_{\mu}[\Psi(\xi)]$ as $K \to \infty$.  In particular,
% if $\lambda^{\ast}=\kappa=0$, $v_{1}=\sup_{\mu\in\frakM}\E_{\mu}[\Psi(\xi)]$, and
If $\lambda^{\ast} > \kappa$, then there exists a constant $D$ independent of $K$ such that
\[
v_{K} \ \ \leq \ \ \sup_{\mu \in \frakM} \E_{\mu}[\Psi(\xi)] \ \ \leq \ \ v_{K} + \frac{L D + M}{NK},
\]
In addition, if $\Xi$ is convex and $\Psi$ is concave, then $v_{1} = v_{P} = v_{D}$.
\end{enumerate}
\end{corollary}

Statement~\ref{itm:finiteDual_N+1} shows that the worst-case distribution $\mu^{\ast}$ is a \textit{perturbation} of $\nu = \frac{1}{N} \sum_{i=1}^{N} \delta_{\hxi^{i}}$, where $N-1$ out of the $N$ points, $\{\hxi^{i}\}_{i \neq i_{0}}$, are perturbed with all their probability mass to an associated maximizer $\xi^{i}_{\ast}$, while at most one point $\hxi^{i_{0}}$ is split and perturbed to two maximizers $\underline{\xi}^{i_{0}}_{\ast}$ and $\overline{\xi}^{i_{0}}_{\ast}$.  If for each $\hxi^{i}$ the set of maximizers is a singleton, then there is no need to split, and the worst-case distribution $\mu^{\ast}$ has support on $N$ points as well.
%It is easy to see that the same result applies when $\nu(\hxi^{i})$ has values other than $1/N$.
Using this structure, we obtain result~\ref{itm:finiteDual_robustapprox}, which states that the primal problem can be approximated by a robust optimization problem with uncertainty set $\frakM_{K}$, which is a subset of $\frakM$ that contains all distributions supported on $NK$ points (not necessarily distinct) with probability $1 / (NK)$ each.
Particularly, when $\Psi$ is concave, such approximation is exact; and when $\Psi$ is Lipschitz and $p = 1$, then $v_{1}$ is an $O(1/N)$-approximation of $v_{P} = v_{D}$.

\begin{remark}
As can easily be seen in the proof, the result of statement~\ref{itm:finiteDual_N+1} can be generalized as follows.
Suppose that $\nu = \sum_{i=1}^{N} \nu_{i} \delta_{\hxi^{i}}$, then whenever a worst-case distribution exists, there exists one of the form
\[
\mu^{\ast} \ \ = \ \ \sum_{i \neq i_{0}} \nu_{i} \delta_{\xi^{i}_{\ast}} + p_{0} \nu_{i_{0}} \delta_{\underline{\xi}^{i_{0}}_{\ast}} + (1 - p_{0}) \nu_{i_{0}} \delta_{\overline{\xi}^{i_{0}}_{\ast}}.
\]
\end{remark}

\begin{remark}
\label{rmk:metric_space}
The results in Corollary~\ref{cor:finite} hold for arbitrary metric spaces $(\Xi,d)$.  The Polish space assumption on $(\Xi,d)$ is used only for the measurability results in Lemma~\ref{lemma:measurable}, and in the finite-supported setting such an assumption is not needed.
\end{remark}

\begin{remark}%[Extreme distribution of a Wasserstein ball]
\label{rmk:extreme}
Under a compactness assumption on $(\Xi,d)$, \citet{wozabal2012framework} pointed out that to solve~\eqref{eqn:problem_primal}, it suffices to consider the extreme points of the Wasserstein ball~$\frakM$, which are distributions that are supported on at most $N+3$ points.
Later, in~\citet{owhadi1504extreme}, this result was improved for Polish spaces or Borel subsets of Polish spaces to distributions that are supported on at most $N+2$ points.
Statement~\ref{itm:finiteDual_N+1} further strengthens these results --- for arbitrary metric spaces (see Remark~\ref{rmk:metric_space}), it suffices to consider distributions that are supported on at most $N+1$ points, and this bound is tight as shown by Example~\ref{eg:UQ} below.
Moreover, $N-1$ out of the at most $N+1$ points in the support of the extreme distribution has the same probability masses as the associated points in the support of the nominal distribution.
We also remark that after this work, \citet{yue2021linear} shows that the worst-case distribution is supported on at most $N+1$ points by employing the Richter-Rogosinski theorem.
\end{remark}

\begin{remark}[Total Variation metric]
By choosing the discrete metric $d(\xi,\zeta) = \mathds{1}_{\{\xi \neq \zeta\}}$ on $\Xi$, the Wasserstein distance is equal to the Total Variation distance (\citet{gibbs2002choosing}), which can be used for settings in which it matters whether two points are the same or different, but no other notion of difference is relevant.
In this case, if $\theta$ is chosen such that $N \theta$ is an integer, then there is no point indexed by $i_{0}$ in~(\ref{eqn:finiteDual_N+1}) such that fractions of its probability are transported to different points to obtain worst-case distribution $\mu^{\ast}$, and the primal problem~\eqref{eqn:problem_primal} is reduced to the robust optimization problem with uncertainty set $\frakM_{1}$, whether $\Xi$ ($\Psi$) is convex (concave) or not.
\end{remark}

\begin{remark}[Choosing radius $\theta$]
\label{rmk:theta}
There are multiple ways to choose the radius $\theta$ of the Wasserstein ball.
Practically we found that the cross-validation method seems to work well for many problems.
We used this method for the numerical experiments on intensity estimation in Section~\ref{sec:NHPP}, and obtained good performance.
It has also been used in \cite{esfahani2015data} for various numerical examples including portfolio optimization and uncertainty quantification.
On the other hand, one can use concentration inequalities to determine an upper bound on the radius with a statistical guarantee.
For example, in the context of a model with an underlying true distribution, one can estimate a theoretical upper bound on the distance between the empirical distribution and the true distribution.
The classical concentration inequalities (see, e.g., \cite{esfahani2015data,fournier2014rate}) provides a radius $\theta$ of the order $O(N^{-1/s})$, where $s$ is the dimension of the random variable.
We use this technique for a one-dimensional newsvendor problem in Section~\ref{sec:newsvendor}.
However, note that such a bound is pessimistic for problems with high-dimensional random variables, since it suggests that the number of data points needed grows exponentially (in the dimension $s$) to make the radius~$\theta$ smaller while maintaining the statistical guarantee.
Fortunately, for many decision problems it is \textit{unnecessary} to choose a radius~$\theta$ such that the ball $\mathfrak{M}$ contains a true distribution with guaranteed probability --- our goal is not to approximate the true distribution, but only the objective value.
This also means that even if we would know the exact distance between the empirical distribution and the true distribution, we should not choose the distance to be the radius.
Instead, we can choose the radius such that with high probability, for every $x$ the worst-case expectation in the DRSO problem dominates the true expectation.
This idea has been exploited recently in \citet{gao2020finite}, which shows that the radius can be in the order of $N^{-1/2}$ in many settings, independent of the dimension of the random variable.
Another principle is to choose the radius $\theta$ such that the uncertainty set yields a confidence region for the true minimizer, and it has been shown in \citet{blanchet2019robust,blanchet2019confidence} that such choice would lead to a $N^{-1/2}$-bound asymptotically under proper conditions.
% for distributionally robust stochastic optimization with Burg entropy \cite{lam2016recovering}, divergence measures \cite{namkoong2017variance}, and Wasserstein distance \cite{blanchet2016robust,gao2017Wasserstein}.
\end{remark}

\begin{remark}[Algorithmic Complexity]
According to the dual reformulation, fixing the dual variable $\lambda$, the dual objective function is the sum of the optimal objective values of $N$~separate maximization problems.
Suppose that we have an oracle for solving the inner maximization problem $\sup_{\xi \in \Xi} \{\Psi(x,\xi) - \lambda \|\xi - \hxi^{i}\|\}$ for any $\hxi^{i}$.
Thereby the number of oracle calls to evaluate the worst-case expectation scales linearly with the sample size~$N$.
To optimize the finite sum problem over $x$, there are optimal algorithms whose complexity match the theoretical lower bound (see, e.g., \cite{agarwal2014lower,lan2018optimal}).

For an important class of problems, the inner maximization problem is easy.
Examples include the following problem types:
\begin{enumerate}
\item
Problems in which the inner maximization problem has a closed-form solution, for example, when $\Psi(x,\xi)$ is linear in $\xi$, or more generally, $\Psi(x,\xi)=\max_{1\le k \le K} a_k(x)^\top \xi + b_k(x)$ and $\Xi = \R^{s}$ (\citet[Remark 6.6]{esfahani2015data} and \citet{gao2017Wasserstein}).
In this case, the minimax DRSO problem is equivalent to a finite-sum optimization problem with regularization, so the complexity is linear in the sample size~$N$, and can be (nearly) independent of the dimension of $\Xi$ under suitable conditions.
\item
Problems in which $\Psi(x,\xi)$ is (piecewise) concave in $\xi$, for which the overall DRSO is equivalent to a saddle-point problem (Examples \ref{eg:saddle} and \ref{eg:piecewiseConcave}).
In this case, the complexity is more involved since it may also depend on the size $\theta$ of the Wasserstein ball and the geometry of the sample space.
For example, consider use of the Mirror-Prox algorithm \cite{nemirovski2004prox} to solve the saddle-point problem in which $\Psi(x,\xi)$ is convex in $x$.
The complexity is proportional to the square of the radius $N \theta$ of the Wasserstein ball, and is nearly independent of the dimension of $\Xi$ for uncertainty sets with nice geometry.
As indicated in Remark~\ref{rmk:theta}, a good way to choose $\theta$ is roughly $O(N^{-1/2})$, and similar to the first case above this leads to a complexity that is linear in $N$ and nearly independent of the dimension of the random variable.
\end{enumerate}
\end{remark}

\proof{Proof of Corollary~\ref{cor:finite}.}\leavevmode

\noindent
\ref{itm:finiteDual_lambda}
The result follows directly from Theorem~\ref{thm:strongDual} and Proposition~\ref{prop:kappa_infty}.

\noindent
\ref{itm:finiteDual_N+1}
By Corollary~\ref{cor:existence}\ref{itm:form}, whenever a worst-case distribution exists, there is one supported on at most $2N$ points with the form
\begin{equation}
\label{eqn:finiteDual_2N}
\hat{\mu} \ \ = \ \ \hat{p} \frac{1}{N} \sum_{i=1}^{N} \delta_{\underline{\xi}^{i}_{\ast}} + \left(1 - \hat{p}\right) \frac{1}{N} \sum_{i=1}^{N} \delta_{\overline{\xi}^{i}_{\ast}},
\end{equation}
where $\hat{p} \in [0,1]$, and $\underline{\xi}^{i}_{\ast},\overline{\xi}^{i}_{\ast} \in \argmin_{\xi \in \Xi} \{\lambda^{\ast} d^{p}(\xi,\hxi^{i}) - \Psi(\xi)\}$.
%(In fact, Corollary~\ref{cor:existence}\ref{itm:form} proves the stronger statement that there exists a worst-case distribution of the form~\eqref{eqn:finiteDual_2N} such that all $p^{i}$ are equal, but here we allow them to vary in order to show that a worst-case distribution has a more specific form.)
% Therefore, we obtain the following (possibly non-convex) reformulation of (\ref{eqn:problem_primal}):
% \[\begin{aligned}
%   \max_{\underline{\xi}^{i},\overline{\xi}^{i}\in\Xi,0\leq p^{i}\leq1} \Bigg\{\frac{1}{N}\sum_{i=1}^{N}\Psi(\hxi^{i}) + \frac{1}{N}\sum_{i=1}^{N} \Big[p^{i}(\Psi(\underline{\xi}^{i})&-\Psi(\hxi^{i}))  +(1-p^{i})(\Psi(\overline{\xi}^{i})-\Psi(\hxi^{i}))\Big]: \\
%   &\frac{1}{N}\sum_{i=1}^{N} p^{i}d^{p}(\underline{\xi}^{i},\hxi^{i})+(1-p^{i})d^{p}(\overline{\xi}^{i},\hxi^{i})\leq \theta^{p}\Bigg\},
%   \end{aligned}
% \]
Given $\underline{\xi}^{i}_{\ast},\overline{\xi}^{i}_{\ast}$ for all $i$, note that
\[
\E_{\mu^{\ast}}[\Psi(\xi)] \ \ \le \ \
\max_{0 \leq p_{i} \leq 1} \left\{\frac{1}{N} \sum_{i=1}^{N} \left[p_{i} \Psi(\underline{\xi}_{\ast}^{i}) + (1 - p_{i}) \Psi(\overline{\xi}_{\ast}^{i})\right] \; : \; \frac{1}{N} \sum_{i=1}^{N} \left[p_{i} d^{p}(\underline{\xi}_{\ast}^{i},\hxi^{i}) + (1 - p_{i}) d^{p}(\overline{\xi}_{\ast}^{i},\hxi^{i})\right] \leq \theta^{p}\right\}.
\]
The problem above is a linear program with $N$ variables, one perturbation distance constraint, and $2N$ constraints $0 \leq p_{i} \leq 1$, $i = 1,\ldots,N$.
At an extreme point optimal solution $p^{\ast}$ of the linear program, if the perturbation distance constraint holds as an equality, then at least $N-1$ of the constraints $0 \leq p_{i} \leq 1$ hold as equalities, or equivalently, for at most one index $i_{0}$ it holds that $p_{i_{0}}$ is fractional.
Otherwise, if the perturbation distance constraint holds as a strict inequality, then exactly $N$ of the constraints $0 \leq p_{i} \leq 1$ hold as equalities, or equivalently, none of the $p_{i}$ is fractional.
Such an extreme point optimal solution $p^{\ast}$ gives a worst-case distribution
\[
\mu^{\ast} \ \ = \ \ \frac{1}{N} \sum_{i=1}^{N} \left[p^{\ast}_{i} \delta_{\underline{\xi}^{i}_{\ast}} + \left(1 - p^{\ast}_{i}\right) \delta_{\overline{\xi}^{i}_{\ast}}\right]
\ \ = \ \ \frac{1}{N} \sum_{i \neq i_{0}} \delta_{\xi^{i}_{\ast}} + \frac{p_{0}}{N} \delta_{\underline{\xi}^{i_{0}}_{\ast}} + \frac{1 - p_{0}}{N} \delta_{\overline{\xi}^{i_{0}}_{\ast}}
\]
where $\xi^{i}_{\ast} = \underline{\xi}^{i}_{\ast}$ if $p^{\ast}_{i} = 1$, $\xi^{i}_{\ast} = \overline{\xi}^{i}_{\ast}$ if $p^{\ast}_{i} = 0$, and $p_{0} = p^{\ast}_{i_{0}}$.
% \[
%   \mu^{\ast}=\frac{1}{N}\sum_{i\neq i_{0}} \delta_{\xi^{i}_{\ast}}+ \frac{p^{\ast}}{N}\delta_{\underline{\xi}^{i_{0}}_{\ast}}+\frac{1-p^{\ast}}{N}\delta_{\overline{\xi}^{i_{0}}_{\ast}} ,
% \]
% where $1\leq i_{0}\leq N$, $p^{\ast}\in[0,1]$, $\underline{\xi}_{\ast}^{i_{0}}\defi\underline{T}^{\ast}(\hxi^{i_{0}})$, $\overline{\xi}_{\ast}^{i_{0}}\defi\overline{T}^{\ast}(\hxi^{i_{0}})$, and $\xi^{i}_{\ast}\in\{\underline{T}^{\ast}(\hxi^{i}),\overline{T}^{\ast}(\hxi^{i})\}$ for all $i\neq i_{0}$.

\noindent
\ref{itm:finiteDual_robustapprox}
Note that $\nu \in \cP_{p}(\Xi)$.
If $\Xi$ is bounded, then $\kappa = 0$, and if $\Xi$ is unbounded, then it follows from the assumption and Lemma~\ref{lemma:kappa}\ref{lemma:kappa_expression} that $\kappa \leq L < \infty$.

It follows from Lemma~\ref{lemma:phi}\ref{lemma:phi_bound} that $\overline{D}(\lambda,\hxi^{i}) < \infty$ for all $\lambda > \kappa$ and all~$i$.
Let $D \defi \max\left\{\overline{D}([\lambda^{\ast} + \kappa] / 2, \hxi^{i}) \, : \, i \in \{1,\ldots,N\}\right\} < \infty$.
Consider any $\varepsilon' > 0$.
For each~$i$, let $\varepsilon^{i} > 0$ be such that $\sup_{\xi \in \Xi} \left\{d^{p}(\xi,\hxi^{i}) \, : \, [\lambda^{\ast} / 2 + \kappa / 2] d^{p}(\xi,\hxi^{i}) - \Psi(\xi) \, \leq \, \Phi([\lambda^{\ast} + \kappa] / 2, \hxi^{i}) + \varepsilon\right\} \le \overline{D}\left([\lambda^{\ast} + \kappa] / 2, \hxi^{i}\right) + \varepsilon'$ for all $\varepsilon \in (0, \varepsilon^{i})$.
Consider any $\varepsilon \in (0, [\lambda^{\ast} - \kappa] / 2)$ such that $\varepsilon < \varepsilon^{i}$ for all~$i$.
By Lemma~\ref{lem:e-optimal lambda > kappa}, there are a $\lambda_{1}^{\varepsilon} \in \left(\lambda^{\ast} - \varepsilon, \, \lambda^{\ast}\right)$, a $\lambda_{2}^{\varepsilon} \in \left(\lambda^{\ast}, \, \lambda^{\ast} + \varepsilon\right)$,  $\hat{p}_{1}^{\varepsilon}, \hat{p}_{2}^{\varepsilon}, \hat{p}_{3}^{\varepsilon} \ge 0$, and $\hat{\xi}_{i1}^{\varepsilon}, \hat{\xi}_{i2}^{\varepsilon}, \hat{\xi}_{i3}^{\varepsilon} \in \Xi$ for $i = 1,\ldots,N$, such that $\hat{p}_{1}^{\varepsilon} + \hat{p}_{2}^{\varepsilon} + \hat{p}_{3}^{\varepsilon} = 1$,
\begin{align*}
\lambda_{1}^{\varepsilon} d^{p}(\hat{\xi}_{i1}^{\varepsilon},\hxi^{i}) - \Psi(\hat{\xi}_{i1}^{\varepsilon}) \ \ & \leq \ \ \Phi(\lambda_{1}^{\varepsilon},\hxi^{i}) + \varepsilon \\
\lambda_{2}^{\varepsilon} d^{p}(\hat{\xi}_{i2}^{\varepsilon},\hxi^{i}) - \Psi(\hat{\xi}_{i2}^{\varepsilon}) \ \ & \leq \ \ \Phi(\lambda_{2}^{\varepsilon},\hxi^{i}) + \varepsilon \\
\hat{\xi}_{i3}^{\varepsilon} \ \ & = \ \ \hxi^{i}
\end{align*}
for all $i = 1,\ldots,N$, $\frac{1}{N} \sum_{i=1}^{N} \sum_{j=1}^{3} \hat{p}_{j}^{\varepsilon} d^{p}(\hat{\xi}_{ij}^{\varepsilon},\hxi^{i}) \leq \theta^{p}$, and
\[
\hat{\mu}^{\varepsilon} \ \ \defi \ \ \hat{p}_{1}^{\varepsilon} \frac{1}{N} \sum_{i=1}^{N} \delta_{\hat{\xi}_{i1}^{\varepsilon}} + \hat{p}_{2}^{\varepsilon} \frac{1}{N} \sum_{i=1}^{N} \delta_{\hat{\xi}_{i2}^{\varepsilon}} + \hat{p}_{3}^{\varepsilon} \frac{1}{N} \sum_{i=1}^{N} \delta_{\hat{\xi}_{i3}^{\varepsilon}}
\]
satisfies $\sup_{\mu \in \frakM} \E_{\mu}[\Psi(\xi)] \le \E_{\hat{\mu}^{\varepsilon}}[\Psi(\xi)] + \varepsilon$.
If $\Psi(\hat{\xi}_{ij}^{\varepsilon}) < \Psi(\hxi^{i})$ for any $(i,j)$, then $\hat{\xi}_{ij}^{\varepsilon}$ can be replaced with $\hxi^{i}$, and the statements above continue to hold.
Therefore, without loss of generality, we assume that $\Psi(\hat{\xi}_{ij}^{\varepsilon}) \ge \Psi(\hxi^{i})$.
Given $\{\hat{\xi}_{ij}^{\varepsilon} \, : \, 1 \leq i \leq N, 1 \leq j \leq 3\}$, note that
\[
\E_{\hat{\mu}^{\varepsilon}}[\Psi(\xi)] \ \ \le \ \
\max_{p_{ij} \ge 0} \left\{\frac{1}{N} \sum_{i=1}^{N} \sum_{j=1}^{3} p_{ij} \Psi(\hat{\xi}_{ij}^{\varepsilon}) \; : \; \frac{1}{N} \sum_{i=1}^{N} \sum_{j=1}^{3} p_{ij} d^{p}(\hat{\xi}_{ij}^{\varepsilon},\hxi^{i}) \leq \theta^{p}, \; \sum_{j=1}^{3} p_{ij} = 1 \; \forall \; i = 1,\ldots,N\right\}.
\]
The problem above is a linear program with $3N$ variables and $N+1$ constraints in addition to the variable bounds $p_{ij} \geq 0$.
At an extreme point optimal solution of the linear program, at most $N+1$ of the variables are nonzero.
In addition, for each $i = 1,\ldots,N$, the constraint $p_{i1} + p_{i2} + p_{i3} = 1$ requires at least one of the variables $p_{i1},p_{i2},p_{i3}$ to be nonzero.
Thus, for at least $N-1$ indices~$i$, exactly one of $p_{i1},p_{i2},p_{i3}$ is equal to~$1$, and for at most one index $i_{0}$, two of $p_{i_{0},1},p_{i_{0},2},p_{i_{0},3}$ are positive.
Hence, an extreme point optimal solution of the linear program gives a solution of the form
\[
\mu^{\varepsilon} \ \ = \ \ \frac{1}{N} \sum_{i \neq i_{0}} \delta_{\xi^{i}_{\varepsilon}} + \frac{p_{\varepsilon}}{N} \delta_{\underline{\xi}^{i_{0}}_{\varepsilon}} + \frac{1 - p_{\varepsilon}}{N} \delta_{\overline{\xi}^{i_{0}}_{\varepsilon}},
\]
where $\xi^{i}_{\varepsilon} \in \left\{\hat{\xi}_{i1}^{\varepsilon}, \hat{\xi}_{i2}^{\varepsilon}, \hat{\xi}_{i3}^{\varepsilon}\right\}$ for all $i \in \{1,\ldots,N\} \setminus \{i_{0}\}$, $\underline{\xi}^{i_{0}}_{\varepsilon}, \overline{\xi}^{i_{0}}_{\varepsilon} \in \left\{\hat{\xi}_{i_{0},1}^{\varepsilon}, \hat{\xi}_{i_{0},2}^{\varepsilon}, \hat{\xi}_{i_{0},3}^{\varepsilon}\right\}$, such that $d^{p}\left(\underline{\xi}^{i_{0}}_{\varepsilon},\hxi^{i_{0}}\right) \le d^{p}\left(\overline{\xi}^{i_{0}}_{\varepsilon},\hxi^{i_{0}}\right)$,
\[
\frac{1}{N} \sum_{i \neq i_{0}} d^{p}\left(\xi^{i}_{\varepsilon},\hxi^{i}\right) + \frac{p_{\varepsilon}}{N} d^{p}\left(\underline{\xi}^{i_{0}}_{\varepsilon},\hxi^{i_{0}}\right) + \frac{1 - p_{\varepsilon}}{N} d^{p}\left(\overline{\xi}^{i_{0}}_{\varepsilon},\hxi^{i_{0}}\right) \ \ \leq \ \ \theta^{p},
\]
and $\E_{\hat{\mu}^{\varepsilon}}[\Psi(\xi)] \le \E_{\mu^{\varepsilon}}[\Psi(\xi)]$.
Consider
\[
\xi^{ik}_{\varepsilon} \ \ \defi \ \ \left\{\begin{array}{ll}
\xi^{i}_{\varepsilon} & \mbox{ for } i \neq i_{0}, k = 1, \ldots, K, \\
\underline{\xi}^{i_{0}}_{\varepsilon} & \mbox{ for } k = 1, \ldots, \lceil{K p_{\varepsilon}}\rceil, \\
\overline{\xi}^{i_{0}}_{\varepsilon} & \mbox{ for } k = \lceil{K p_{\varepsilon}}\rceil + 1, \ldots, K.
\end{array}\right.
\]
%Note that
%\begin{align*}
%& \frac{1}{NK} \sum_{i=1}^{N} \sum_{k=1}^{K} d^{p}(\xi^{ik}_{\varepsilon},\hxi^{i}) \\
%& = \ \ \frac{1}{N} \sum_{i \neq i_{0}} d^{p}(\xi^{i}_{\varepsilon},\hxi^{i}) + \frac{p_{\varepsilon}}{N} d^{p}(\underline{\xi}^{i_{0}}_{\varepsilon},\hxi^{i_{0}}) + \frac{1 - p_{\varepsilon}}{N} d^{p}(\overline{\xi}^{i_{0}}_{\varepsilon},\hxi^{i_{0}}) + \frac{\lceil{K p_{\varepsilon}}\rceil / K - p_{\varepsilon}}{N} \left[d^{p}(\underline{\xi}^{i_{0}}_{\varepsilon},\hxi^{i_{0}}) - d^{p}(\overline{\xi}^{i_{0}}_{\varepsilon},\hxi^{i_{0}})\right] \\
%& \le \ \ \frac{1}{N} \sum_{i \neq i_{0}} d^{p}(\xi^{i}_{\varepsilon},\hxi^{i}) + \frac{p_{\varepsilon}}{N} d^{p}(\underline{\xi}^{i_{0}}_{\varepsilon},\hxi^{i_{0}}) + \frac{1 - p_{\varepsilon}}{N} d^{p}(\overline{\xi}^{i_{0}}_{\varepsilon},\hxi^{i_{0}})
%\ \ \leq \ \ \theta^{p}
%\end{align*}
Then, since $d^{p}\left(\underline{\xi}^{i_{0}}_{\varepsilon},\hxi^{i_{0}}\right) \le d^{p}\left(\overline{\xi}^{i_{0}}_{\varepsilon},\hxi^{i_{0}}\right)$, it follows that
\begin{align*}
\frac{1}{NK} \sum_{i=1}^{N} \sum_{k=1}^{K} d^{p}\left(\xi^{ik}_{\varepsilon},\hxi^{i}\right)
\ \ & = \ \ \frac{1}{N} \sum_{i \neq i_{0}} d^{p}\left(\xi^{i}_{\varepsilon},\hxi^{i}\right) + \frac{\lceil{K p_{\varepsilon}}\rceil / K}{N} \; d^{p}\left(\underline{\xi}^{i_{0}}_{\varepsilon},\hxi^{i_{0}}\right) + \frac{1 - \lceil{K p_{\varepsilon}}\rceil / K}{N} \; d^{p}\left(\overline{\xi}^{i_{0}}_{\varepsilon},\hxi^{i_{0}}\right) \\
& \le \ \ \frac{1}{N} \sum_{i \neq i_{0}} d^{p}\left(\xi^{i}_{\varepsilon},\hxi^{i}\right) + \frac{p_{\varepsilon}}{N} \; d^{p}\left(\underline{\xi}^{i_{0}}_{\varepsilon},\hxi^{i_{0}}\right) + \frac{1 - p_{\varepsilon}}{N} \; d^{p}\left(\overline{\xi}^{i_{0}}_{\varepsilon},\hxi^{i_{0}}\right)
\ \ \leq \ \ \theta^{p},
\end{align*}
and thus $\{\xi^{ik}_{\varepsilon}\}_{i,k} \in \frakM_{K}$.
Since $0 \le \lceil{K p_{\varepsilon}}\rceil / K - p_{\varepsilon} < 1/K$, it follows that
\begin{align*}
\E_{\mu^{\varepsilon}}[\Psi(\xi)]
\ \ & = \ \ \frac{1}{NK} \sum_{i=1}^{N} \sum_{k=1}^{K} \Psi(\xi^{ik}_{\varepsilon}) + \frac{1}{N} \left(\frac{\lceil{K p_{\varepsilon}}\rceil}{K} - p_{\varepsilon}\right) \big[\Psi(\overline{\xi}^{i_{0}}_{\varepsilon}) - \Psi(\underline{\xi}^{i_{0}}_{\varepsilon})\big] \\
\ \ & \le \ \ \frac{1}{NK} \sum_{i=1}^{N} \sum_{k=1}^{K} \Psi(\xi^{ik}_{\varepsilon}) + \frac{1}{N} \left(\frac{\lceil{K p_{\varepsilon}}\rceil}{K} - p_{\varepsilon}\right) \big[\Psi(\overline{\xi}^{i_{0}}_{\varepsilon}) - \Psi(\hxi^{i_{0}})\big] \\
\ \ & \le \ \ v_{K} + \frac{1}{NK} \left[L d^{p}\left(\overline{\xi}^{i_{0}}_{\varepsilon},\hxi^{i_{0}}\right) + M\right].
\end{align*}
Since $\overline{\xi}^{i_{0}}_{\varepsilon} \in \left\{\hat{\xi}_{i_{0},1}^{\varepsilon}, \hat{\xi}_{i_{0},2}^{\varepsilon}, \hat{\xi}_{i_{0},3}^{\varepsilon}\right\}$, it follows that
\begin{align*}
d^{p}\left(\overline{\xi}^{i_{0}}_{\varepsilon},\hxi^{i_{0}}\right) \ \ & \le \ \ \sup_{\xi \in \Xi} \left\{d^{p}(\xi,\hxi^{i_{0}}) \, : \, \lambda_{1}^{\varepsilon} d^{p}(\xi,\hxi^{i_{0}}) - \Psi(\xi) \leq \Phi(\lambda_{1}^{\varepsilon}, \hxi^{i_{0}}) + \varepsilon\right\} \\
& \le \ \ \sup_{\xi \in \Xi} \left\{d^{p}(\xi,\hxi^{i_{0}}) \, : \, \frac{\lambda^{\ast} + \kappa}{2} d^{p}(\xi,\hxi^{i_{0}}) - \Psi(\xi) \leq \Phi([\lambda^{\ast} + \kappa] / 2, \hxi^{i_{0}}) + \varepsilon\right\} \\
& \le \ \ \overline{D}\left([\lambda^{\ast} + \kappa] / 2, \hxi^{i_{0}}\right) + \varepsilon',
\end{align*}
where the second inequality follows from Lemma~\ref{lemma:phi}\ref{lemma:phi_monotone} and $\lambda_{1}^{\varepsilon} > [\lambda^{\ast} + \kappa] / 2$.
Therefore
\begin{align*}
\sup_{\mu \in \frakM} \E_{\mu}[\Psi(\xi)] \ \ & \le \ \ \E_{\mu^{\varepsilon}}[\Psi(\xi)] + \varepsilon
\ \ \le \ \ v_{K} + \frac{1}{NK} \left[L \left(\overline{D}\left([\lambda^{\ast} + \kappa] / 2, \hxi^{i_{0}}\right) + \varepsilon'\right) + M\right] + \varepsilon \\
\Rightarrow \ \ \ \sup_{\mu \in \frakM} \E_{\mu}[\Psi(\xi)] \ \ & \le \ \ v_{K} + \frac{1}{NK} \left[L D + M\right].
\end{align*}
\hfillqed
\endproof

\begin{example}[Saddle-point Problem]
\label{eg:saddle}
When $\Psi(x,\xi)$ is convex in $x$ and concave in $\xi$, $p = 1$, and $d(\xi,\zeta) = \|\xi - \zeta\|_{2}$, then Corollary~\ref{cor:finite}\ref{itm:finiteDual_robustapprox} shows that~\eqref{eqn:DRSO} is equivalent to a convex-concave saddle point problem
\[
\label{eqn:saddle}
\min_{x \in X} \max_{(\xi^{1},\ldots,\xi^{N}) \in Y} \frac{1}{N} \sum_{i=1}^{N} \Psi(x,\xi^{i}),
\]
with $\ell_{1} / \ell_{2}$-norm uncertainty set
\[
Y \ \ = \ \ \left\{(\xi^{1},\ldots,\xi^{N}) \in \Xi^{N} \; : \; \sum_{i=1}^{N} \left\|\xi^{i} - \hxi^{i}\right\|_{2} \leq N \theta\right\}.
\]
\end{example}

\begin{example}[Piecewise concave objective]
\label{eg:piecewiseConcave}
\citet{esfahani2015data} showed that if $\Xi$ is a closed convex subset of $\R^{K}$ with norm $\|\cdot\|$, $\nu$ has finite support, $p = 1$, and $\Psi(\xi) = \max_{1 \leq j \leq J} \Psi^{j}(\xi)$, where each $\Psi^{j}$ is concave, then the DRSO can be formulated as a convex optimization problem.
Here we show that the result can be obtained as a corollary from the structure of a worst-case distribution.
If $\lambda^{\ast} > \kappa$, then
%the set
% \[
% \argmin_{\xi \in \Xi} \{\lambda^\ast \|\xi-\hxi^{i}\| - \Psi(\xi)\} \subset \bigcup_{j=1}^{J} \argmin_{\xi \in \Xi} \{\lambda^\ast \|\xi-\hxi^{i}\| - \Psi^{j}(\xi)\}
% \]
% is non-empty, thus 
by Corollary~\ref{cor:existence}\ref{itm:iff}\ref{itm:existenceCond1}, a worst-case distribution exists and, by Corollary~\ref{cor:finite}, has the form $\frac{1}{N} \sum_{i=1}^{N} \sum_{k=1}^{2} p_{ik} \delta_{\xi^{ik}}$, where for each $i$ it holds that $\sum_{k=1}^{2} p_{ik} = 1$ and $\xi^{ik} \in \argmin_{\xi \in \Xi} \{\lambda^{\ast} \|\xi - \hxi^{i}\| - \Psi(\xi)\}$.
Moreover, due to the concavity of $\Psi^{j}$, if $\Psi(\xi^{i1}) = \Psi^{j}(\xi^{i1})$ and $\Psi(\xi^{i2}) = \Psi^{j}(\xi^{i2})$ for some~$j$, then any convex combination $\xi_{\alpha} = \alpha \xi^{i1} + (1 - \alpha) \xi^{i2}$ satisfies $\lambda^{\ast} \|\xi_{\alpha} - \hxi^{i}\| - \Psi(\xi_{\alpha}) \le \lambda^{\ast} \|\xi_{\alpha} - \hxi^{i}\| - \Psi^{j}(\xi_{\alpha}) \le \alpha \left[\lambda^{\ast} \|\xi^{i1} - \hxi^{i}\| - \Psi^{j}(\xi^{i1})\right] + (1 - \alpha) \left[\lambda^{\ast} \|\xi^{i2} - \hxi^{i}\| - \Psi^{j}(\xi^{i2})\right] = \min_{\xi \in \Xi} \{\lambda^{\ast} \|\xi - \hxi^{i}\| - \Psi(\xi)\}$, and thus $\xi_{\alpha} \in \argmin_{\xi \in \Xi} \{\lambda^{\ast} \|\xi - \hxi^{i}\| - \Psi(\xi)\}$.
Therefore, we can assume without loss of generality that if $\Psi(\xi^{i1}) = \Psi^{j_{1}}(\xi^{i1})$ and $\Psi(\xi^{i2}) = \Psi^{j_{2}}(\xi^{i2})$, then $j_{1} \neq j_{2}$.
Thus, we can consider distributions of the form
\[
\frac{1}{N} \sum_{i=1}^{N} \sum_{j=1}^{J} p_{ij} \delta_{\xi^{ij}}, \ \ \textrm{ where } \sum_{j=1}^{J} p_{ij} = 1, \ \Psi(\xi^{ij}) = \Psi^{j}(\xi^{ij}), \ \forall \ i=1,\ldots,N.
\]
If $\lambda^{\ast} = \kappa$, then by Lemma~\ref{lem:e-optimal lambda = kappa} and the concavity of $\Psi^{j}$, we can also consider distributions of the form above.
Therefore, the original primal problem is equivalent to
\[
\sup_{p_{ij} \geq 0, \; \xi^{ij} \in \Xi} \left\{\frac{1}{N} \sum_{i=1}^{N} \sum_{j=1}^{J} p_{ij} \Psi^{j}(\xi^{ij}) \; : \; \frac{1}{N} \sum_{i=1}^{N} \sum_{j=1}^{J} p_{ij} \|\xi^{ij} - \hxi^{i}\| \leq \theta, \; \sum_{j=1}^{J} p_{ij} = 1, \; \forall \ i = 1,\ldots,N\right\}.
\]
Next, let $\tilde{\xi}^{ij} \defi \hxi^{i} + p_{ij} (\xi^{ij} - \hxi^{i})$.
Then, by the positive homogeneity of norms and the convexity-preserving property of perspective functions (cf. Section 2.3.3 in \citet{boyd2004convex}), the primal problem can be reformulated as the following convex optimization problem:
\[
\sup_{\substack{p_{ij} \geq 0, \; \sum_{j} p_{ij} = 1 \\ \tilde{\xi}^{ij} \in \R^{K}}} \left\{\frac{1}{N} \sum_{i=1}^{N} \sum_{j=1}^{J} p_{ij} \Psi^{j}\left(\hxi^{i} + \frac{\tilde{\xi}^{ij} - \hxi^{i}}{p_{ij}}\right) \; : \; \frac{1}{N} \sum_{i=1}^{N} \sum_{j=1}^{J} \|\tilde{\xi}^{ij} - \hxi^{i}\| \leq \theta, \; \hxi^{i} + \frac{\tilde{\xi}^{ij} - \hxi^{i}}{p_{ij}} \in \Xi, \; \forall \ i,j\right\}.
\]
This establishes Theorem~4.4 in \citet{esfahani2015data}, which was obtained therein by a separate procedure of dualizing twice the reformulation~\eqref{eqn:dual_finite}.
\end{example}

\begin{example}[\textbf{Uncertainty Quantification}]
\label{eg:UQ}
Consider a metric space $(\Xi,d)$ such that every bounded subset in $(\Xi,d)$ is totally bounded, and any open proper Borel subset $C \subsetneq \Xi$.
Let $\Psi = -\mathds{1}_{C}$, and consider the uncertainty quantification problem
\begin{equation}
\sup_{\mu \in \frakM} \; \E_{\mu}[\Psi(\xi)] \ \ = \ \ \sup_{\mu \in \frakM} \; \E_{\mu}[-\mathds{1}_{C}(\xi)] \ \ = \ \ - \min_{\mu \in \frakM} \mu(C)
\label{eqn:uncertainty quantification}
\end{equation}
Next we show that it follows from Corollary~\ref{cor:existence}\ref{itm:iff} that problem~(\ref{eqn:uncertainty quantification}) has an optimal distribution.

It follows from the definition of $\kappa$ that $\kappa = 0$.
Since $C$ is open, $\Psi = -\mathds{1}_{C}$ is upper-semi-continuous.
It follows from Lemma~\ref{lemma:dual objective}\ref{itm:dual minimizer} that the dual objective function $h$ has a minimizer $\lambda^{\ast} \in [0,\infty)$.
Next we consider two cases: \\
Case~1:  There exists a dual minimizer $\lambda^{\ast} > 0$:
Then it follows from Corollary~\ref{cor:existence}\ref{itm:iff}\ref{itm:existenceCond1} that a worst-case distribution exists. \\
Case~2:  $\lambda^{\ast} = 0$ is the unique dual minimizer:
Since $C$ is a proper subset of $\Xi$, it follows that $\argmax_{\xi \in \Xi} \Psi(\xi) = \argmin_{\xi \in \Xi} \mathds{1}_{C}(\xi)$ is nonempty, and that
\[
\Phi(0,\zeta) \ \ = \ \ \inf_{\xi \in \Xi} \{-\Psi(\xi)\}
\ \ = \ \ \inf_{\xi \in \Xi} \mathds{1}_{C}(\xi) \ \ = \ \ 0.
\]
Also, for any $\lambda > 0$ it holds that
\begin{align*}
0 \ \ & = \ \ -\int_{\Xi} \Phi(0,\zeta) \nu(d\zeta)
\ \ = \ \ h(0)
\ \ < \ \ h(\lambda) \\
& = \ \ \lambda \theta^{p} - \int_{\Xi} \inf_{\xi \in \Xi} \{\lambda d^{p}(\xi,\zeta) + \mathds{1}_{C}(\xi)\} \nu(d\zeta) \\
& = \ \ \lambda \theta^{p} - \int_{\Xi} \min\left\{1, \ \lambda \inf_{\xi \in \Xi \setminus C} d^{p}(\xi,\zeta)\right\} \nu(d\zeta).
\end{align*}
Thus, for any $\lambda > 0$ it holds that
\[
\int_{\Xi} \min\left\{\frac{1}{\lambda}, \ \inf_{\xi \in \Xi \setminus C} d^{p}(\xi,\zeta)\right\} \nu(d\zeta) \ \ < \ \ \theta^{p}.
\]
Next, let $\lambda \to 0$, then it follows from the monotone convergence theorem that
\[
\int_{\Xi} \underline{D}_{0}(0,\zeta) \nu(d\zeta) \ \ = \ \ \int_{\Xi} \inf_{\xi \in \Xi \setminus C} d^{p}(\xi,\zeta) \ \ \leq \ \ \theta^{p}.
\]
Therefore, it follows from Corollary~\ref{cor:existence}\ref{itm:iff}\ref{itm:existenceCond3} that a worst-case distribution exists.

Depending on the metric~$d$, for $\zeta \in C$, $\argmin_{\xi \in \Xi \setminus C} d^{p}(\xi,\zeta)$ may or may not be on the boundary of $C$.
Next, suppose that $d$ is an intrinsic metric (so that, for $\zeta \in C$, it holds that $\argmin_{\xi \in \Xi \setminus C} d^{p}(\xi,\zeta) \subset \partial C$), and that $\nu$ has finite support, say $\nu = \frac{1}{N} \sum_{i=1}^{N} \delta_{\hxi^{i}}$.
Then the worst-case distribution $\mu^{\ast}$ of the problem~(\ref{eqn:uncertainty quantification})
has an intuitive form.
The worst-case distribution perturbs $\nu$ such that the set $C$ contains as little probability mass as possible, which can be achieved in a \emph{greedy} fashion as follows.
Suppose that $\{\hxi^{i}\}_{i=1}^{N}$ are sorted such that $\hxi^{1},\ldots,\hxi^{I} \in C$, $\hxi^{I+1},\ldots,\hxi^{N} \notin C$, and $d^{p}(\hxi^{1}, \Xi \setminus C) \leq \cdots \leq d^{p}(\hxi^{I}, \Xi \setminus C)$.
Then to save the total budget of perturbation, $\hxi^{I+1},\ldots,\hxi^{N}$ stay at the same place, and the points $\hxi^{i}$ with smaller indices~$i$ have priority to be transported to their closest points in $\partial C$.
Proceeding in this greedy way, points $\hxi^{1},\ldots,\hxi^{i_{0}-1}$ ($i_{0} \leq I$) are transported in full (with full mass $1 / N$) to their closest points in $\partial C$.
The next point $\hxi^{i_{0}}$ is transported in part or in full to its closest point in $\partial C$ --- it may be the case that transporting $\hxi^{i_{0}}$ in full to its closest point in $\partial C$ would violate the Wasserstein distance constraint.
Thus, only part $p_{0} / N$, with $0 < p_{0} \le 1$, of its mass is transported, and the remaining mass $(1 - p_{0}) / N$ stays at point $\hxi^{i_{0}}$ (see Figure \ref{fig:cc}).
Points $\hxi^{i_{0}+1},\ldots,\hxi^{N}$ stay at the same place.
Therefore the worst-case distribution has the form
\[
\label{eqn:mu_UQ}
\mu^{\ast} \ \ = \ \ \frac{1}{N} \sum_{i=1}^{i_{0}-1} \delta_{\xi^{i}_{\ast}} + \frac{p_{0}}{N} \delta_{\xi^{i_{0}}_{\ast}} + \frac{1 - p_{0}}{N} \delta_{\hxi^{i_{0}}} + \frac{1}{N} \sum_{i = i_{0}+1}^{N} \delta_{\hxi^{i}}.
\]
In fact, if $d^{p}(\hxi^{i_{0}-1}, \Xi \setminus C) < d^{p}(\hxi^{i_{0}}, \Xi \setminus C)$, then the dual optimizer $\lambda^{\ast}$ is such that
\[
\xi^{i}_{\ast} \ \ \in \ \ \argmin_{\xi \in \Xi} \{\lambda^{\ast} d^{p}(\xi,\hxi^{i}) + \mathds{1}_{C}(\xi)\}
\ \ = \ \ \argmin_{\xi \in \partial C} d^{p}(\xi,\hxi^{i})
\]
for all $i = 1,\ldots,i_{0}-1$, and
\[
\xi^{i_{0}}_{\ast} \ \ \in \ \ \argmin_{\xi \in \Xi} \{\lambda^{\ast} d^{p}(\xi,\hxi^{i_{0}}) + \mathds{1}_{C}(\xi)\}
\ \ = \ \ \begin{cases}
\{\hxi^{i_{0}}\} \cup \argmin_{\xi \in \partial C} d^{p}(\xi,\hxi^{i_{0}}) & \mbox{if } p_{0} < 1, \\
\argmin_{\xi \in \partial C} d^{p}(\xi,\hxi^{i_{0}}) & \mbox{if } p_{0} = 1.
\end{cases}
\]

\begin{figure}[h]
\centering
\includegraphics[width=0.3\textwidth]{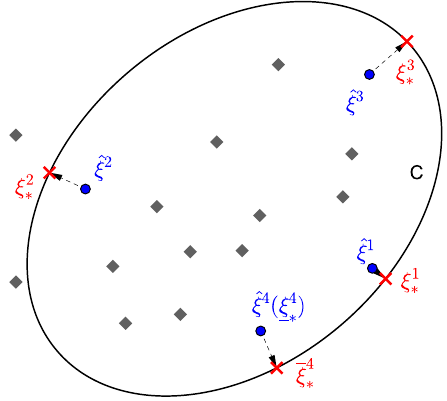}
\caption{When $\Psi = -\mathds{1}_C$, then the worst-case distribution perturbs the nominal distribution in a greedy fashion. The solid and diamond dots are the support of the nominal distribution $\nu$. $\hxi^{1},\hxi^{2},\hxi^{3}$ are the three interior points closest to $\partial C$ and thus are transported to $\xi_{\ast}^{1},\xi_{\ast}^{2},\xi_{\ast}^{3}$ respectively. $\hxi^{4}$ is the interior point fourth closest to $\partial C$, but its full mass cannot be transported to $\partial C$ due to the Wasserstein distance constraint, so it is split and the parts are moved to $\overline{\xi}_{\ast}^{4}$ and $\underline{\xi}_{\ast}^{4} = \hxi^{4}$.}
\label{fig:cc}
\end{figure}

% Indeed, using the notation in Corollary \ref{cor:finite}, we have
% \[
%   \argmin_{\xi \in \Xi} \{\lambda^{\ast} d^{p}(\xi,\hxi^{i_{0}}) - \Psi(\xi)\} =
% \]
%  $\over\in\{\zeta\}\cup\argmin_{\xi\in\partial C} d^{p}(\xi,\zeta)$, namely, $\mu^{\ast}$ either keeps $\zeta$ still, or perturbs it to the closest point on the boundary (so $\mathds{1}_C(\zeta)$ changes from 1 to 0). Since $\mu^{\ast}$ transports as much mass in $C$ outwards as possible, it transports mass in a \textsl{greedy} fashion.

In the discussion above we assumed that $C$ is open.
The next proposition shows that for the normed space $(\R^{K},\lVert\cdot\rVert)$, uncertainty quantification of any Borel set can be reduced to that of an open set.

\begin{proposition}[Continuity with respect to the boundary]
\label{prop:chance-constrained_C}
Consider a metric space $(\Xi,d)$ such that every bounded subset in $(\Xi,d)$ is totally bounded.
Let $\nu \in \cP(\Xi)$, $\theta > 0$, and $\frakM = \{\mu \in \cP(\Xi) \, : \, W_{p}(\mu,\nu) \leq \theta\}$.
Then for any proper Borel subset $C \subsetneq \Xi$, it holds that
\[
\inf_{\mu \in \frakM} \mu(C) \ \ = \ \ \min_{\mu \in \frakM} \mu(\interior(C)).
\]
\end{proposition}
% The result is quite intuitive. In fact, when $C$ is not open and $\partial C$ is non-empty, transporting mass to $\partial C$ may not change the objective from 1 to 0 as when $C$ is open. Instead, one can transport it to the point outside $C$ but arbitrarily close to $\partial C$. This explains why the worst-case probability is continuous with respect to $\partial C$.
\end{example}

% \begin{corollary}[Affinely-perturbed objective]
% \label{cor:robust_linear}
%     Suppose $\Psi(x,\xi) = a^{\top} x+ b$, where $\xi=[a;b]$. Assume the metric $d$ is induced by some norm $||\cdot||_q$ ($q\geq1$). Let $\nu=\frac{1}{N}\sum_{i=1}^{N}\delta_{\hxi^{i}}$ and $\hxi^{i}=[\hat{a}_{i};\hat{b}_{i}]$, $i=1,\ldots,N$. Then the DRSO problem (\ref{eqn:DRSO}) is equivalent to
%     \[
%       \min_{x\in X}~ \frac{1}{N}\sum_{i=1}^{N} (\hat{a}_{i}^{\top} x+\hat{b}_{i}) + \theta ||x||_{q^{\ast}},
%     \]
%     where $q^{\ast}$ is such that $1/q+1/q^{\ast}=1$.
% \end{corollary}

\begin{example}[\textbf{Finite Domain}]
\label{eg:finite domain}
Next consider the special case when $\Xi$ is finite, say $\Xi = \{\xi^{1},\ldots,\xi^{B}\}$ for some positive integer~$B$.
The nominal distribution~$\nu$ is given by $\nu = \sum_{i=1}^{B} \nu_{i} \delta_{\xi^{i}}$.
Then the DRSO becomes
\begin{equation}
\label{problem:primal bins}
\min_{x \in X} \max_{\mu \defi (\mu_{1},\ldots,\mu_{B}) \in \Delta_{B}} \left\{\sum_{i=1}^{B} \mu_{i} \Psi(x,\xi^{i}) \; : \; W_{p}(\mu,\nu) \leq \theta\right\}.
\end{equation}

\begin{corollary}
\label{cor:finite bins}
Problem~(\ref{problem:primal bins}) has a strong dual
\begin{equation}
\label{eqn:finite bins}
\min_{x \in X, \; \lambda \geq 0, \; y \defi (y_{1},\ldots,y_{B})} \left\{\lambda \theta^{p} + \sum_{i=0}^{B} \nu_{i} y_{i} \; : \; y_{i} \geq \Psi(x,\xi^{j}) - \lambda d^{p}(\xi^{i},\xi^{j}),\; \forall \; i,j = 1,\ldots,B\right\}.
\end{equation}
For any $\nu$ and $x$, the worst-case distribution can be computed by solving
\begin{equation}
\label{eqn:finite bins_assignment}
\max_{\mu \in \Delta_{B}, \; \gamma \in \R_{+}^{B \times B}} \left\{\sum_{i=1}^{B} \mu_{i} \Psi(x,\xi^{i}) \; : \; \sum_{i,j = 1}^{B} d^{p}(\xi^{i},\xi^{j}) \gamma_{ij} \leq \theta^{p}, \; \sum_{j=1}^{B} \gamma_{ij} = \mu_{i} \; \forall \; i, \; \sum_{i=1}^{B} \gamma_{ij} = \nu_{j} \; \forall \; j\right\}.
\end{equation}
\end{corollary}

\proof{Proof.}
Reformulation~(\ref{eqn:finite bins}) follows from Theorem~\ref{thm:strongDual}, and~(\ref{eqn:finite bins_assignment}) follows from the equivalent definition of Wasserstein distance in Example~\ref{eg:transportationLP}.
\hfillqed
\endproof
\end{example}

\section{Applications}
\label{sec:application}

In this section, we apply our results to various application problems, including on/off system control and intensity estimation for point processes.
In these problems, the space $\Xi$ is a space of counting measures (sample paths of a point process), which is non-convex and infinite dimensional, and the nominal distribution is a point process.
%In the third problem, the nominal distribution is an arbitrary probability distribution on a finite dimensional space.
Hence, the results in \citet{esfahani2015data} and \citet{zhao2018data} cannot be applied to these problems.

\subsection{On/Off System Control}
\label{sec:process}

In this problem, the decision maker faces a point process and controls a two-state (on/off) system.
The point process is assumed to be exogenous, that is, the arrival times are not affected by the on/off state of the system.
When the system is switched on, a cost of $c$ per unit time is incurred, and each arrival while the system is on contributes $1$~unit revenue.
When the system is off, no cost is incurred and no revenue is earned.
The decision maker wants to choose a control to maximize the total profit during a finite time horizon.
This problem is a prototype for problems in sensor networks and revenue management.

Let the finite time horizon be denoted with $[0,1]$, and let
\[
\Xi \ \ \defi \ \ \left\{\xi \ \defi \ \sum_{m=1}^{M} \delta_{\eta_{m}} \; : \; M \in \mathbb{Z}_{+}, \, \eta_{m} \in [0,1], \, m=1,\ldots,M\right\}
\]
be the set of sample paths of arrival times, which can also be described as the space of finite counting measures on $[0,1]$.
In many practical settings, the decision maker does not have a probability distribution for the point process.
Instead, the decision maker has observations of sample paths of the point process, which constitute an empirical point process as nominal distribution.
Specifically, suppose that we have data of $N$ sample paths $\hxi^{i} = \sum_{m=1}^{M_{i}} \delta_{\heta^{i}_{m}}$, $i = 1,\ldots,N$, where $M_{i}$ is a nonnegative integer, and $\heta^{i}_{m} \in [0,1]$ for all $i$ and $m$.
Then the nominal (empirical) distribution is $\nu = \frac{1}{N} \sum_{i=1}^{N} \delta_{\hxi^{i}}$.

Note that if one would maximize the expected profit with respect to the nominal distribution, then it would yield a degenerate control, in which the system is switched on only momentarily at the arrival time points of the observed sample paths.
Consequently, if future arrival times would differ from the historically observed arrival times by even a small amount, then the system would be switched off at future arrival times and no revenue would be earned.
Due to such degeneracy and instability of optimization with respect to the nominal distribution, we resort to the distributionally robust approach.

As pointed out before, the choice of ambiguity set $\frakM$ matters a great deal.
For example, suppose that one chose $\frakM$ to be the KL-divergence ball $\frakM_{\phi_{KL}}$, that is, the set of all distributions $\mu \in \cP(\Xi)$ such that $I_{\phi_{KL}}(\mu,\nu) \le \theta$ for some $\theta > 0$.
Then, for $I_{\phi_{KL}}(\mu,\nu) \le \theta$ to hold, $\mu$ can put positive probability on observed sample paths $\hxi^{i}$ only.
Thus, even if one would solve a DRSO problem with KL-divergence ball $\frakM_{\phi_{KL}}$, the resulting solution would be a degenerate control, in which the system is switched on only momentarily at the arrival time points of the historically observed sample paths, the same solution as when maximizing the expected profit with respect to the nominal distribution.
In this section we show that solving a DRSO problem with a Wasserstein ball results in sensible solutions.

%Next we define the metric $d$ on the space $\Xi$ of finite counting measures.
We assume that the metric $d$ on $\Xi$ satisfies the following conditions.
(Note that in this section, when we write the $W_{1}$ distance between two measures, we use the extended definition mentioned in Section~\ref{sec:wasserstein}.)
\begin{enumerate}[label=(\roman*)]
\item
The metric space $(\Xi,d)$ is a Polish space.
\item
For any nonnegative integer $M$, and any $\xi = \sum_{m=1}^{M} \delta_{\eta_{m}}$ and $\hxi = \sum_{m=1}^{M} \delta_{\heta_{m}}$, where $\{\eta_{m}\}_{m=1}^{M},\{\heta_{m}\}_{m=1}^{M} \subset [0,1]$, it holds that
\[
\label{itm:process_dassumption2}
d(\xi,\hxi) \ \ = \ \ W_{1}(\xi,\hxi) \ \ = \ \ \sum_{m=1}^{M} |\eta_{(m)} - \heta_{(m)}|,
\]
where $\{\eta_{(m)}\}$ and $\{\heta_{(m)}\}$ are the order statistics of $\{\eta_{m}\}$ and $\{\heta_{m}\}$ respectively.
\item
\label{itm:process_dassumption1}
For any Borel set $C \subset [0,1]$, $\theta \geq 0$, positive integer $M$, and $\hxi = \sum_{m=1}^{M} \delta_{\heta_{m}}$, where $\{\heta_{m}\}_{m=1}^{M} \subset [0,1]$, it holds that
\[
\inf_{\xi \in \Xi} \Big\{\xi(C): d(\xi,\hxi) \leq \theta \Big\} \ \ = \ \ \inf_{\xi \in \Xi} \Big\{\xi(C):  W_{1}(\xi,\hxi) \leq \theta \Big\}.
\]
\end{enumerate}
Note that condition~\ref{itm:process_dassumption2} is imposed only on $\xi,\hxi \in \Xi$ such that $\xi([0,1]) = \hxi([0,1])$, and that conditions~\ref{itm:process_dassumption2} and~\ref{itm:process_dassumption1} do not imply that $d = W_{1}$.
Examples of metrics~$d$ that satisfy the conditions above are
\[
d\left(\sum_{m=1}^{M} \delta_{\eta_{m}}, \sum_{l=1}^{L} \delta_{\heta_{l}}\right) \ \ = \ \ \sum_{m=1}^{\min\{M,L\}} \left|\eta_{(m)} - \heta_{(m)}\right| + \left|M - L\right|,
\]
\[
d\left(\sum_{m=1}^{M} \delta_{\eta_{m}}, \sum_{l=1}^{L} \delta_{\heta_{l}}\right) \ \ = \ \ \left\{\begin{array}{ll}
\max\{M,L\}, & M \neq L, \\
\sum_{m=1}^{M} \left|\eta_{(m)} - \heta_{(m)}\right|, & M = L,
\end{array}\right.
\]
and
\begin{equation}
\label{eqn:d_point_process}
d\left(\sum_{m=1}^{M} \delta_{\eta_{m}}, \sum_{l=1}^{L} \delta_{\heta_{l}}\right) \ \ = \ \ \left\{\begin{array}{ll}
+\infty, & M \neq L, \\
\sum_{m=1}^{M} \left|\eta_{(m)} - \heta_{(m)}\right|, & M = L.
\end{array}\right.
\end{equation}
These metrics are similar to the ones in \citet{barbour1992stein} and \citet{chen2004stein}.
% Given the metric $d$, the point processes on $[0,1]$ are then defined by the set $\cP(\Xi)$ of Borel probability measures on $\Xi$.

The set of point processes on $[0,1]$ is defined by the set $\cP(\Xi)$ of Borel probability measures on $\Xi$.
Given the metric~$d$, we choose the distance between two point processes $\mu,\nu \in \cP(\Xi)$ to be $W_{1}(\mu,\nu)$ as defined in~(\ref{eqn:def_wasserstein}).
Then the ambiguity set $\frakM = \{\mu \in \cP(\Xi) \, : \, W_{1}(\mu,\nu) \leq \theta\}$.
Let $X$ denote the set of all functions $x : [0,1] \mapsto \{0,1\}$ such that $x^{-1}(1) \defi \{t \in [0,1] \, : \, x(t) = 1\}$ is a Borel set.
The decision maker is looking for a control $x \in X$ that maximizes the total profit, by solving the problem
\begin{equation}
\label{eqn:problem_process}
v^{\ast} \ \ \defi \ \ \sup_{x \in X} \left\{v(x) \ \defi \ -c \int_{0}^{1} x(t) dt + \inf_{\mu \in \frakM} \E_{\xi \sim \mu}\big[\xi(x^{-1}(1))\big]\right\}.
\end{equation}

Next we investigate the structure of the optimal control.
Let $\interior(x^{-1}(1))$ denote the interior of the set $x^{-1}(1)$ on the space $[0,1]$ with canonical topology (and thus $0,1 \in \interior([0,1])$).

\begin{proposition}
\label{prop:process_policy}
Suppose that $\nu = \frac{1}{N} \sum_{i=1}^{N} \delta_{\hxi^{i}}$ with $\hxi^{i} = \sum_{m=1}^{M_{i}} \delta_{\heta^{i}_{m}}$.
For any control $x$, it holds that
\begin{align}
\inf_{\mu \in \frakM} \E_{\xi \sim \mu}[\xi(x^{-1}(1))] 
\ \ & = \ \ \inf_{\rho \in \cP(\Xi^2)} \Big\{\E_{(\xi,\hxi) \sim \rho}[\xi(x^{-1}(1))] \; : \;  \E_{(\xi,\hxi) \sim \rho}\big[W_{1}(\xi,\hxi)\big] \leq \theta, \ \pi^{2}_{\#} \rho = \nu\Big\} \label{eqn:process_tower} \\
& = \ \ \sup_{\lambda \geq 0} \left\{- \lambda \theta + \frac{1}{N} \sum_{i=1}^{N} \sum_{m=1}^{M_{i}} \min_{\eta \in [0,1]} \Big\{\mathds{1}_{\{\eta \in \interior(x^{-1}(1))\}} + \lambda \left|\eta - \heta^{i}_{m}\right|\Big\}\right\}.
\label{eq:process_dual}
\end{align}
Moreover, there exists a non-negative integer $J$ such that
\begin{equation}
\label{eqn:policy_unionOfIntervals}
v^{\ast} \ \ = \ \ \sup_{\substack{\underline{x}_{j}, \overline{x}_{j} \in [0,1], \\ \underline{x}_{j} < \overline{x}_{j} < \underline{x}_{j'} < \overline{x}_{j'} \; \forall \; 1 \leq j < j' \leq J}} \left\{v\left(\sum_{j=1}^{J} \mathds{1}_{[\underline{x}_{j},\overline{x}_{j}]}\right) \ \defi \ - c \sum_{j=1}^{J} (\overline{x}_{j} - \underline{x}_{j}) + \inf_{\mu \in \frakM} \E_{\xi \sim \mu} \big[\xi\big(\cup_{j=1}^{J} [\underline{x}_{j},\overline{x}_{j}]\big)\big]\right\}.
\end{equation}
\end{proposition}

Note that
\[
\inf_{\mu \in \frakM} \E_{\xi \sim \mu}[\xi(x^{-1}(1))] \ \ = \ \ \inf_{\gamma \in \cP(\Xi^{2})} \Big\{\E_{(\xi,\hxi) \sim \gamma}[\xi(x^{-1}(1))] \; : \; \E_{\gamma}[d(\xi,\hxi)] \leq \theta, \  \pi^{2}_{\#} \gamma = \nu\Big\}.
\]
Hence, (\ref{eqn:process_tower}) shows that without changing the optimal value, we can replace $d$ by $W_{1}$ in the constraint,
% , and enlarge the set of joint distributions from $\cP(\Xi^{2})$ to $\cP(\cB([0,1]) \times \Xi)$. 
and \eqref{eq:process_dual} provides an equivalent dual reformulation that can be interpreted as the dual problem of the uncertainty quantification problem in which the nominal distribution is a Borel measure $\frac{1}{n}\sum_{i=1}^n \sum_{m=1}^{M_i} \heta_m^i$ on $[0,1]$.
Then by Example \ref{eg:UQ}, it can be solved by a greedy algorithm.
Moreover, (\ref{eqn:policy_unionOfIntervals}) shows that if $\nu$ is an empirical point process, then it suffices to consider the set of controls such that the system is on during a finite disjoint union of intervals with positive length.

\begin{figure}[h]
\centering
\includegraphics[width=\textwidth]{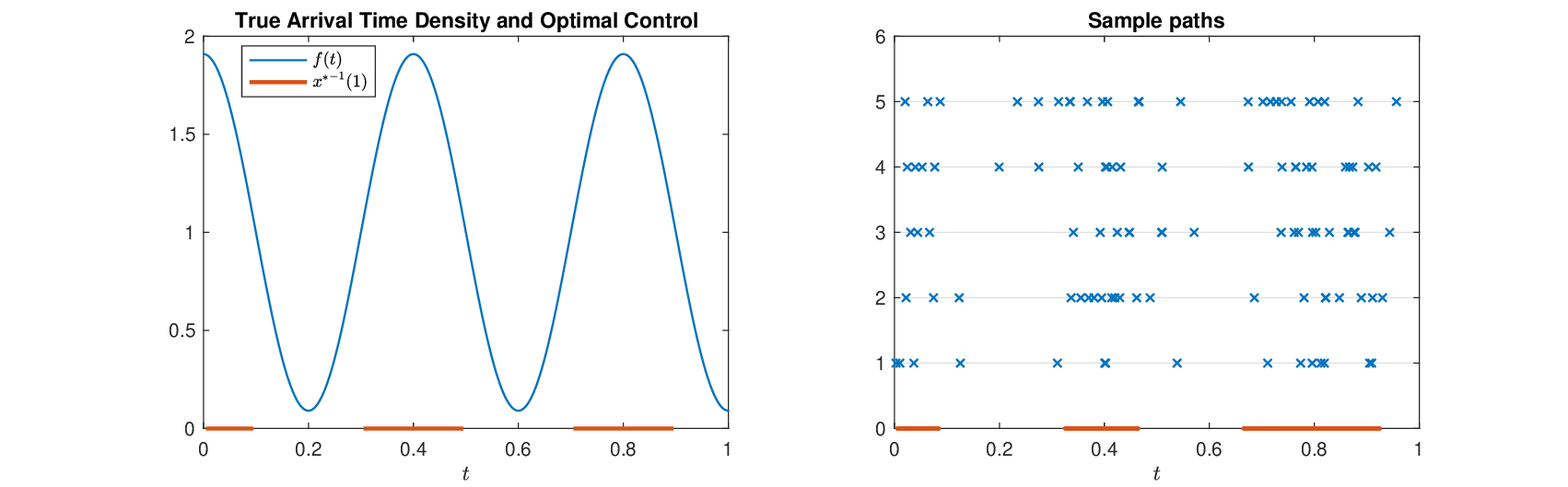}
\caption{Optimal on/off system control for the true process and the DRSO.}
\label{fig:process}
\end{figure}

\begin{example}
Figure~\ref{fig:process} shows results for an instance of the on/off system control problem.
Suppose that the number of arrivals is Poisson distributed with mean $\lambda$, and given the number of arrivals, the arrival points are i.i.d.\ with density $f(t)$, $t \in [0,1]$.
For example, $f \equiv 1$ corresponds to the Poisson process with rate $\lambda$.
The optimization problem based on knowing the process is $\max_{x} \int_{x^{-1}(1)} [-c + \lambda f(t)] dt$, with optimal control $x^{\ast}(t) = \mathds{1}_{\{\lambda f(t) > c\}}$.

Let $f(t) = k [a + \sin(w t + s)]$, with $a > 1$ and $k = 1 / [a + (\cos(s) - \cos(w+s)) / w]$.
Particularly, let $w = 5 \pi$, $s = \pi / 2$, $a = 1.1$, and $c = \lambda = 20$.
Thus ${x^{\ast}}^{-1}(1) = [0,0.1] \cup [0.3,0.5] \cup [0.7,0.9]$.
The left of Figure~\ref{fig:process} shows $f(t)$ and $x^{\ast}$.
Suppose we have $N = 5$ sample paths $\hxi^{1},\ldots,\hxi^{5}$, each of which contains $M_{i} \sim Poisson(\lambda)$ i.i.d.\ arrival points.
The right of Figure~\ref{fig:process} shows $5$~sample paths and the resulting DRSO solution.
Even with a relatively small number of sample paths, the two controls are close to each other, and the DRSO provides not only a sensible solution (unlike the approach of optimizing the expected profit with respect to the empirical distribution), but even a good solution for the true process control problem.
\end{example}

\subsection{Intensity Estimation for Non-homogeneous Poisson Processes}
\label{sec:NHPP}

Consider the problem of estimating the intensity function $a(t)$ of a non-homogeneous Poisson process $A(t)$ on $[0,T]$ using the maximum likelihood method.
Given $N$ i.i.d.~sample paths $\hxi^{i} = \sum_{m=1}^{M_{i}} \delta_{\heta^{i}_{m}}$, $i = 1,\ldots,N$, the log-likelihood function (see, e.g.~\citet{daley2003}) is written as
\[
\mathscr{L}(a) \ \ = \ \ -\int_{0}^{T} a(t) \, dt + \frac{1}{N} \sum_{i=1}^{N} \sum_{m=1}^{M_{i}} \ln\big(a(\heta^{i}_{m})\big).
\]
A common practice is to partition the time horizon $[0,T]$ into several intervals, and assume that $a(t)$ is piecewise constant with constant value on each of the chosen intervals.
Then the maximum likelihood estimator reduces to the average arrival rate on each interval.
Such an approach suffers from the drawback that the estimator is sensitive to the partition of the time horizon into intervals.
If the partition is very coarse, then the estimator remains constant during long intervals, and fails to capture time-varying arrival rates over shorter time periods.
On the other hand, if the partition is very fine, then many intervals have zero observations, and the estimator varies too erratically over short intervals of time.
%It appears that there is no systematic way to adaptively choose the partition for this sample average method.
A DRSO formulation with $\frakM$ chosen to be the KL-divergence ball $\frakM_{\phi_{KL}}$ has a similar problem, since the resulting estimator vanishes on intervals with no observations.

Consider the DRSO formulation with Wasserstein distance
\begin{equation}
\label{problem:intensityMLE}
\min_{a(t), \; 0 \le t \le T} \ \left\{v_{\theta}(a) \ \defi \ \int_{0}^{T} a(t) \, dt + \max_{\mu \in \frakM} \; \E_{\eta \sim \mu} \left[- \int_{0}^{T} \ln\big(a(t)\big) \, \eta(dt)\right]\right\},
\end{equation}
where $\frakM$ is the Wasserstein ball centered at the empirical process.
To facilitate further analysis, we choose~$d$ defined in~\eqref{eqn:d_point_process} as the distance between two counting measures.
Our strong duality results imply that the dual reformulation of~\eqref{problem:intensityMLE} is given by
\[
\min_{\substack{a(t) \\ \lambda \geq 0}} \ \left\{\int_{0}^{T} a(t) dt + \lambda \theta + \frac{1}{N} \sum_{i=1}^{N} \sum_{m=1}^{M_{i}} \sup_{\eta \in [0,T]} \Big\{- \ln\big(a(\xi)\big) - \lambda |\eta - \heta^{i}_{m}|\Big\}\right\}.
\]
% The following proposition gives the result that the optimal estimator is constant if the radius of the Wasserstein ball is sufficiently large.

% \begin{proposition}
% \label{prop:NHPP}
% For sufficiently large $\theta$, the optimal value $a_{\ast}(t)$ is constant.
% \end{proposition}

\begin{figure}
\centering
\includegraphics[width=0.3\textwidth]{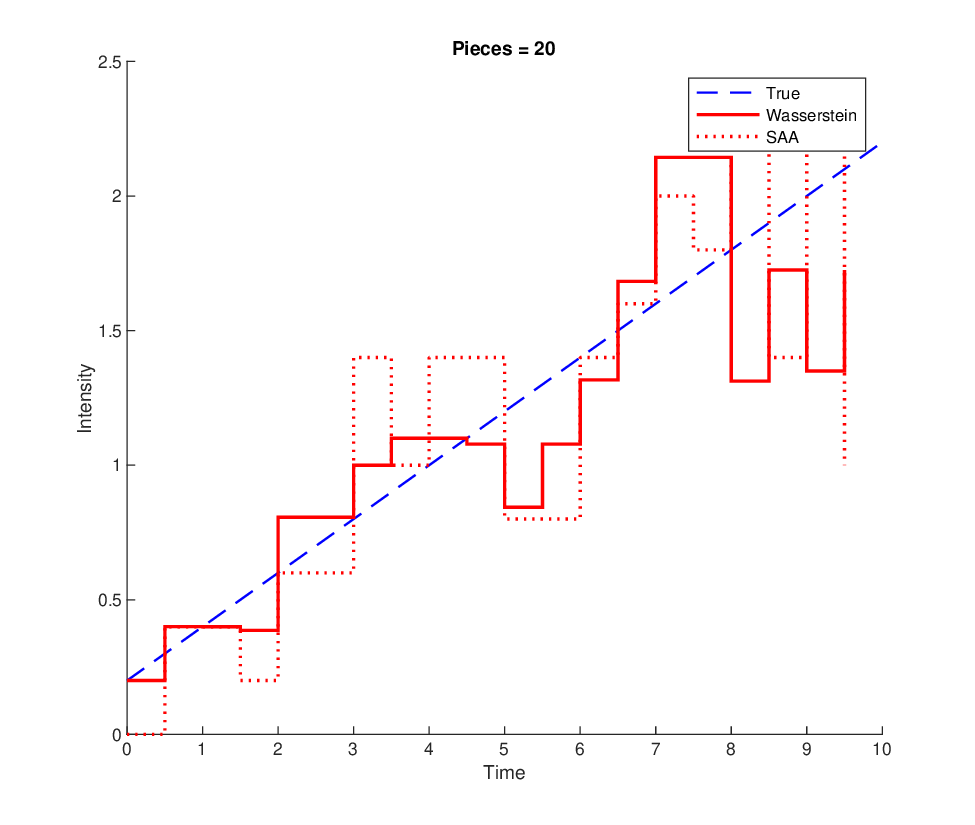}
\includegraphics[width=0.3\textwidth]{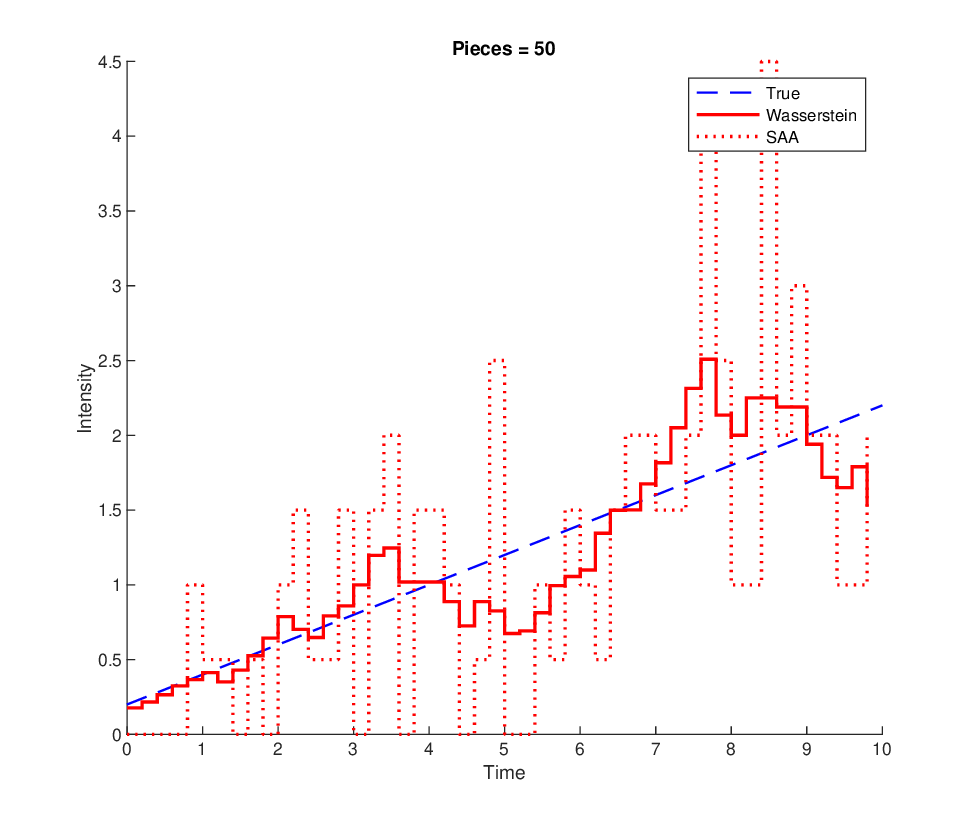}
\includegraphics[width=0.3\textwidth]{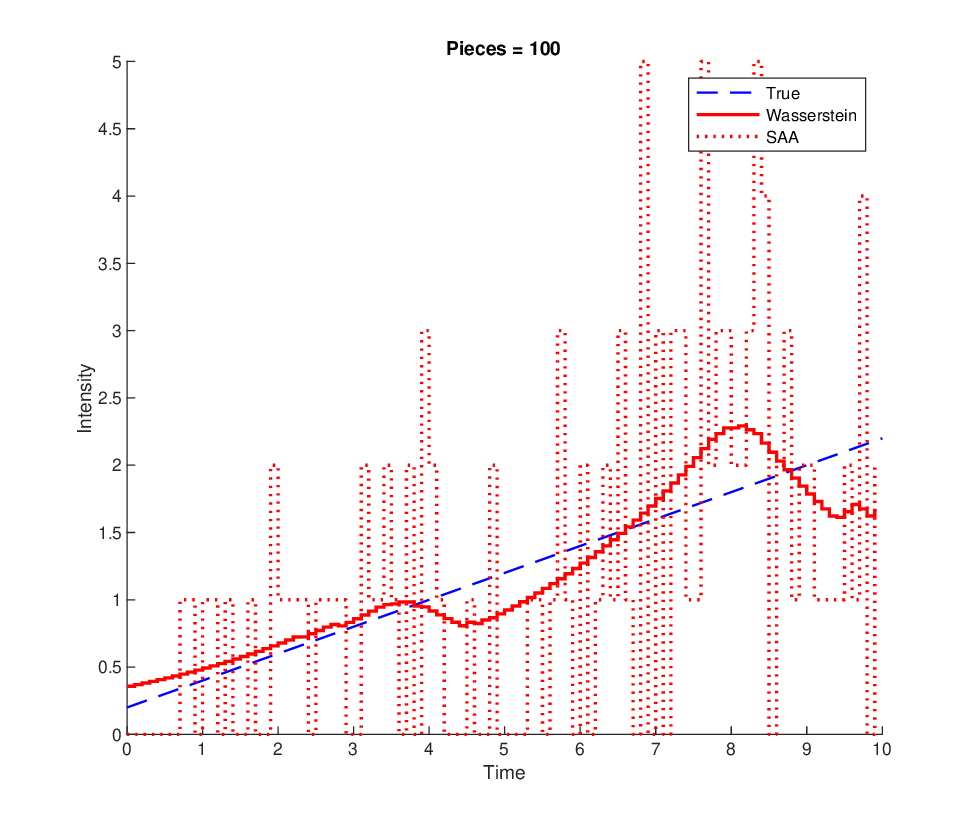}\\
\includegraphics[width=0.3\textwidth]{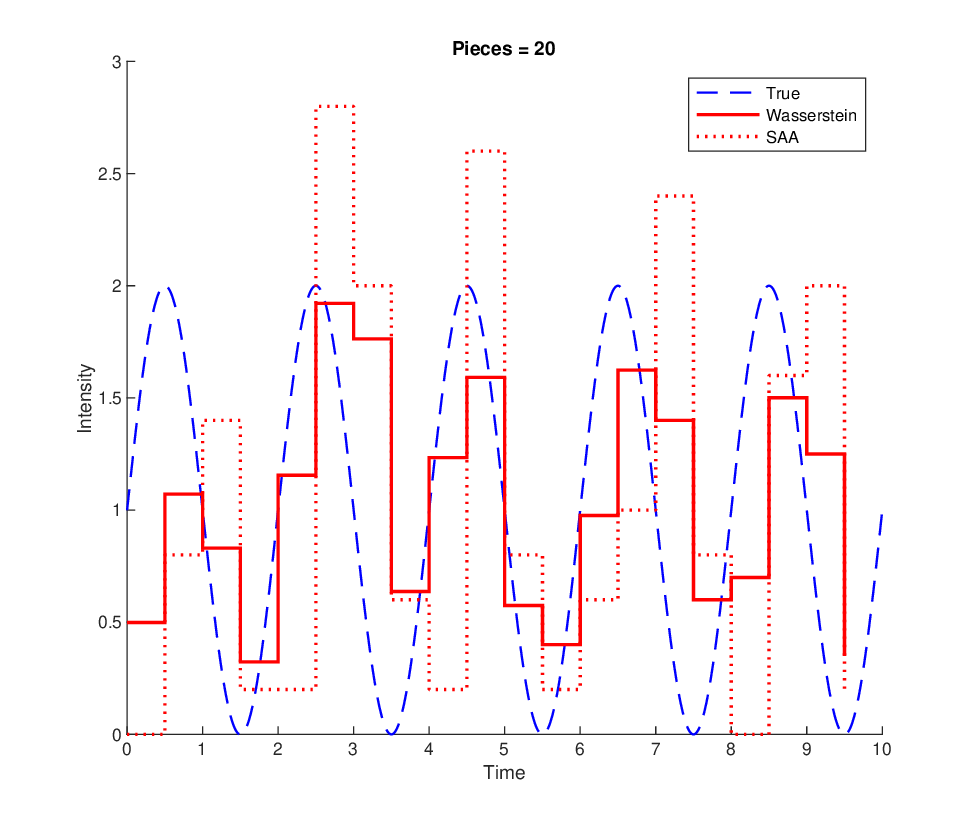}
\includegraphics[width=0.3\textwidth]{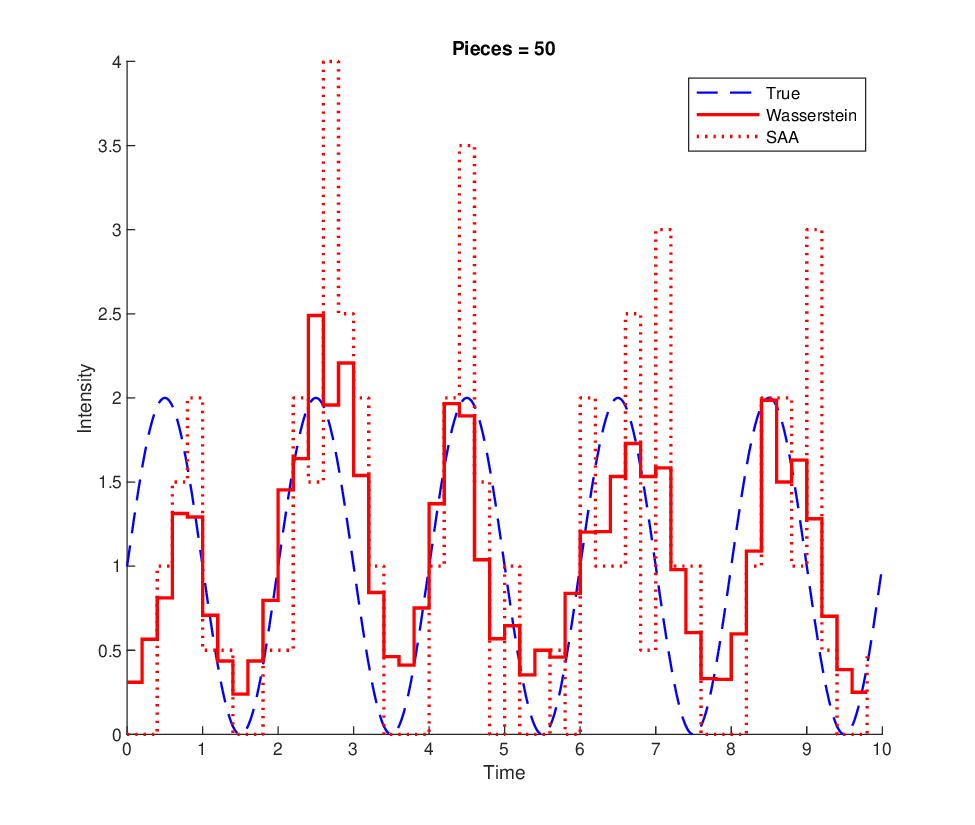}
\includegraphics[width=0.3\textwidth]{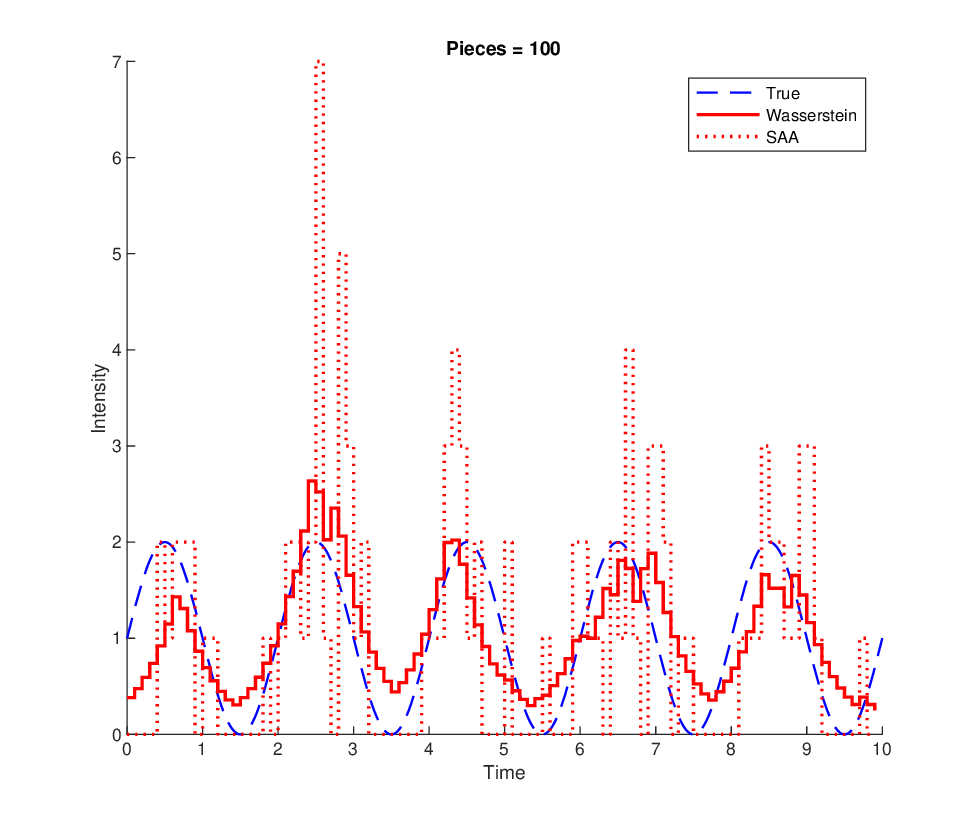}
\caption{Estimated intensity functions using Wasserstein DRSO and MLE}
\label{fig:NHPP}
\end{figure}

Next we present numerical results for underlying true intensity functions given by $a(t) = 0.2 + 0.2t$ and $a(t) = 1 + \sin(\pi t)$, $t \in [0,10]$.
The sample size (number of sample paths) is $N = 20$.
For both the maximum likelihood estimator and the DRSO estimator we optimized over piecewise constant $a(t)$ with number of pieces in $\{20,50,100\}$.
% Note that here piecewise constant assumption is only on the decision $a(t)$, but not on the data.
% This is different from the common practice, which first partitions the time horizon into bins and then count the empirical frequency falling in each bin.
The Wasserstein radius~$\theta$ was chosen via a cross-validation method, in which half of the sample paths were used for estimating $a_{\theta} \in \argmin_{a} v_{\theta}(a)$ for each value of $\theta$, and the other half of the sample paths were used for selecting the value of $\theta \in \argmax_{\theta} \mathscr{L}(a_{\theta})$ with the largest log-likelihood.
The estimation results for both the maximum likelihood estimator (MLE) and the Wasserstein DRSO estimator are shown in Figure~\ref{fig:NHPP}.
Table~\ref{tab:NHPP} shows the $L^{2}$ distances between the estimated intensity functions and the true intensity functions.
The Wasserstein DRSO estimator has superior performance in all cases.
Also, the performance of the DRSO estimator is less sensitive to the fineness of the partition for the piecewise constant function than the performance of the maximum likelihood estimator --- in these results the performance of the maximum likelihood estimator behaves poorly when the number of pieces is large.

\begin{table}[htbp]
\centering
\caption{$L^{2}$ distances between the true intensity functions and the estimated intensity functions of Wasserstein DRSO and MLE}
\label{tab:NHPP}
\begin{tabular}{l|lll|lll}
\toprule
& \multicolumn{3}{c|}{$0.2+0.2t$} & \multicolumn{3}{c}{$1+\sin(\pi t)$} \\
\cmidrule{2-7}    Pieces & 20    & 50    & 100   & 20    & 50    & 100 \\
\midrule
Wasserstein & 0.854 & 0.893 & 0.544 & 2.008 & 2.157 &  2.087 \\
MLE  & 1.510  & 6.536 & 11.906 & 6.160  & 6.591 & 11.766 \\
\bottomrule
\end{tabular}%
\label{tab:addlabel}%
\end{table}%

\subsection{Worst-case Value-at-Risk}
\label{sec:wcVaR}

Value-at-risk is a popular risk measure in finance.
Given a real-valued random variable~$Z$ that represents losses with measure $\nu$ that has a positive density on an interval, and $\alpha \in (0,1)$, the value-at-risk $\VaR^{\nu}_{\alpha}[Z]$ of $Z$ is defined by
\[
\label{eqn:VaR}
\VaR^{\nu}_{\alpha}[Z] \ \ \defi \ \ \inf\left\{t \; : \; \Pr_{\nu}\{Z < t\} \, > \, 1 - \alpha\right\}.
\]
% Note that the $\VaR_{\alpha}[Z]\leq0$ is equivalent to the chance constraint $\Pr\{Z\leq0\}\geq1-\alpha$.
In the spirit of \citet{ghaoui2003worst}, we consider the following worst-case $\VaR$ problem.
Consider a portfolio consisting of $K$ assets with allocation weights $w \defi (w_{1},\ldots,w_{K}) \geq 0$ such that $\sum_{k=1}^{K} w_{k} = 1$.
Let $\xi_{k}$ denote the (random) return of asset $k$, $k = 1,\ldots,K$.
% Assume that the metric~$d$ is induced by some norm $||\cdot||$ on $\R^{K}$.
Assume that $(\Xi,d) = (\R^{K},\|\cdot\|)$.
The worst-case $\VaR$ with respect to the set of probability distributions $\frakM$ is defined as
\[
\label{eqn:worst-caseVaR}
\VaR^{\frakM}_{\alpha}(w) \ \ \defi \ \ \sup_{\mu \in \frakM} \VaR^{\mu}_{\alpha}[-w^{\top} \xi]
\ \ = \ \ \min\left\{q \; : \; \inf_{\mu \in \frakM} \Pr_{\mu}\{-w^{\top} \xi < q\} \, > \, 1 - \alpha\right\}.
\]
Given $w$ and $q$, let $C \defi \{\xi \, : \, -w^{\top} \xi < q\}$.
Then $\inf_{\mu \in \frakM} \Pr_{\mu}\{-w^{\top} \xi < q\} = \inf_{\mu \in \frakM} \E_{\mu}[\mathds{1}_{C}(\xi)] \ \ = \ \ \min_{\mu \in \frakM} \mu(C)$, similar to the uncertainty quantification problem~\eqref{eqn:uncertainty quantification} considered in Example~\ref{eg:UQ}.
Using the rationale developed there, we obtain the following result.

\begin{proposition}
\label{prop:CVaR}
% Let $\alpha \in (0,1)$, $w \in \{w' \in \R^{K} \, : \, \sum_{k=1}^{K} w'_{k} = 1, \, w' \geq 0\}$, $q \geq \VaR_{\alpha}[-w^{\top} \xi]$, and $\theta > 0$.
% If $\E_{\nu}\big[(q + w^{\top} \xi)^{p} \mathds{1}_{\{\VaR_{\alpha}[-w^{\top}\xi] < -w^{\top}\xi < q\}}\big] < \theta^{p}$, then define
% \[
% \beta_{0} \ \ \defi \ \ \min\left\{1, \ \frac{\theta^{p} - \E_{\nu}\big[(q + w^{\top} \xi)^{p} \mathds{1}_{\{-w^{\top} \xi > \VaR_{\alpha}[-w^{\top} \xi]\}}\big]}{(\alpha - \nu\{\xi \, : \, -w^{\top} \xi > \VaR_{\alpha}[-w^{\top} \xi]\}) (q - \VaR_{\alpha}[-w^{\top} \xi])^{p}}\right\}.
% \]
% Then $\inf_{\mu \in \frakM} \Pr_{\mu}\{-w^{\top} \xi \leq q\} \geq 1 - \alpha$ is equivalent to
% \[
% \E_{\nu}\big[(q + w^{\top} \xi)^{p} \mathds{1}_{\{\VaR_{\alpha}[-w^{\top}\xi] < -w^{\top} \xi < q\}}\big] + \beta_{0} \E_{\nu}\big[(q + w^{\top} \xi)^{p} \mathds{1}_{\{-w^{\top} \xi = \VaR_{\alpha}[-w^{\top}\xi]\}}\big] \ \ \geq \ \ \theta^{p}.
% \]
% In particular, when $\nu$ has a continuous cumulative distribution function, then the condition above can be reduced to
% \[
% \E_{\nu}\big[(q + w^{\top} \xi)^{p} \mathds{1}_{\{\VaR_{\alpha}[-w^{\top}\xi] < -w^{\top}\xi < q\}}\big] \ \ \geq \ \ \theta^{p}.
% \]
Suppose that $\nu$ has a positive density on $(\Xi,d) = (\R^{K},\|\cdot\|)$, $\alpha \in (0,1)$, and $\theta > 0$.
Let $w$ be given, and let $\nu_{w}\{(-\infty,s)\} \defi \nu\{\xi \, : \, -w^{\top} \xi < s\}$ for all $s \in \R$.
Then, $\VaR^{\frakM}_{\alpha}(w)$ equals the unique solution $q$ of the following equation:
\[
\int_{\VaR^{\nu}_{\alpha}[-w^{\top} \xi]}^{q} (q-s)^{p} \nu_{w}(ds) \ \ = \ \ \theta^{p} \|w\|_{\ast}^{p}.
\]
\end{proposition}

% \begin{remark}
% It follows from Proposition~\ref{prop:CVaR} that the worst-case $\VaR$, $\VaR^{wc}_{\alpha}(w)$, is equal to the smallest value of $q$ such that
% \[
% \E_{\nu}\big[(q + w^{\top} \xi)^{p} \mathds{1}_{\{\VaR_{\alpha}[-w^{\top}\xi] < -w^{\top} \xi < q\}}\big] + \beta_{0} \E_{\nu}\big[(q + w^{\top} \xi)^{p} \mathds{1}_{\{-w^{\top} \xi = \VaR_{\alpha}[-w^{\top}\xi]\}}\big] \ \ \geq \ \ \theta^{p}
% \]
% and that $\VaR^{wc}_{\alpha}(w) \ge \VaR_{\alpha}[-w^{\top} \xi]$.
% \end{remark}

\begin{figure}[h]
\centering
\includegraphics[height=1.2in]{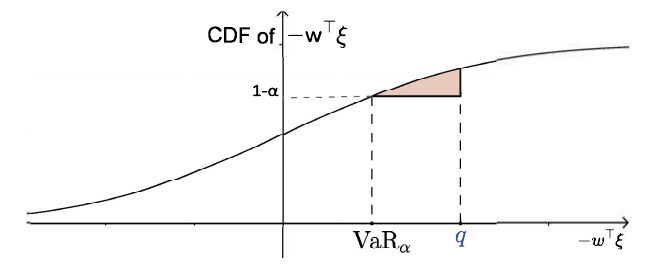}
\caption{Worst-case VaR. When $-w^{\top} \xi$ has a continuous cumulative distribution function, and $p = 1$, then $\VaR_{\alpha}^{\frakM}$ is equal to the value of $q$ such that the area of the shaded region is equal to $\theta \|w\|_{\ast}$.}
\label{fig:wcVaR}
\end{figure}

\begin{example}[\textbf{Worst-case $\VaR$ with Gaussian nominal distribution}]
Suppose that $\nu = N(\bar{\mu},\Sigma)$ and $p = 1$.
It follows that $-w^{\top} \xi \sim N(-w^{\top} \bar{\mu}, w^{\top} \Sigma w)$ and $\VaR^{\nu}_{\alpha}[-w^{\top} \xi] = -w^{\top} \bar{\mu} + \sqrt{w^{\top} \Sigma w} \Phi^{-1}(1-\alpha)$, where $\Phi = N(0,1)$.
By Proposition~\ref{prop:CVaR}, $\VaR_{\alpha}^{\frakM}(-w^{\top} \xi)$ is the value of $q$ such that (see Figure~\ref{fig:wcVaR})
\begin{equation}
\label{eqn:wcVaR_q}
f(q) \ \ \defi \ \ \frac{1}{\sqrt{2\pi w^{\top}\Sigma w}} \int_{\VaR^{\nu}_{\alpha}[-w^{\top}\xi]}^{q} (q-y) \, e^{-\frac{(y + w^{\top} \bar{\mu})^{2}}{2 w^{\top} \Sigma w}} \, dy \ \ = \ \ \theta \|w\|_{\ast}.
\end{equation}
% Let $q'$ be the solution of the equation
% \begin{equation}\label{eqn:q'}
%   f(q')\defiq'[\Phi(q')-(1-\alpha)]-\frac{1}{\sqrt{2\pi}}\left[e^{\frac{-\Phi^{-2}(1-\alpha)}{2}}-e^{\frac{-{q'}^{2}}{2}}\right]=\theta,
% \end{equation}
% where $\Phi$ is the cumulative distribution function of standard normal distribution.
%Then $\VaR_{\alpha}^{wc}(-w^{\top}\xi)=q'\sqrt{w^{\top}\Sigma w}-w^{\top}\bar{\mu}$.
Since $f(q)$ is increasing and continuous, \eqref{eqn:wcVaR_q} can be solved efficiently via a one-dimensional search algorithm.
\end{example}

Proposition~\ref{prop:CVaR} indicates that finding the worst-case $\VaR$ is tractable.
It should also be noted that finding the best allocation weight, i.e., optimizing over $w$, still may be hard, since the $\VaR$ constraint is essentially a chance constraint.

\section{Discussions}
\label{sec:discussion}

In this section, we discuss some of the advantages of using the Wasserstein ambiguity set.
In Section~\ref{sec:newsvendor}, we compare DRSO with the Wasserstein ambiguity set to DRSO with $\phi$-divergence ambiguity sets for the newsvendor problem.
In Section~\ref{sec:two-stage}, we illustrate how the close connection between DRSO and robust optimization (Corollary~\ref{cor:finite bins}\ref{itm:finiteDual_robustapprox}) can be used to solve DRSO problems.
% Finally in Section~\ref{sec:transportation}, we demonstrate the power of our constructive proof method by applying it to a class of distributionally robust transportation problems other than DRSOs.

\subsection{Newsvendor Problem: Comparison with \texorpdfstring{$\phi$}{phi}-divergence}
\label{sec:newsvendor}

In this section, we consider distributionally robust newsvendor problems, and compare results obtained with the Wasserstein ambiguity set and results obtained with $\phi$-divergence ambiguity sets.
In the newsvendor problem, the decision maker chooses the initial inventory level before the unknown demand is observed, facing both overage and underage costs.
If the demand distribution $\mu$ is known, the problem can be formulated as
\[
\min_{x \geq 0} \ \E_{\mu}[h [x - \xi]^+ + b [\xi - x]^+],
\]
where $x$ denotes the chosen initial inventory level, $\xi$ denotes the random demand, and $h$, $b$ represent respectively the overage and underage costs per unit.
Assume that $h,b > 0$, and $\mu$ is supported on $\{0,1,\ldots,B\}$ for some positive integer $B$.
Consider the usual metric $d(i,j) = |i - j|$ on $\{0,1,\ldots,B\}$.
% If the underlying distribution $\mu$ is known exactly, it is well-known that the set of optimal solutions is the whole interval of $\frac{b}{b+h}$-quantiles $\label{eqn:newsvendor_optimalQuantile}
%   [\inf\{t:F_{\mu}(t)\geq\frac{b}{b+h}\},\  \sup\{t:F_{\mu}(t)\leq\frac{b}{b+h}\}]
% $, where $F_{\mu}$ is the cumulative distribution function of $\mu$.
% Now suppose we are given $N_{i}$ observations of historical demand $i$, $i=0,\ldots,\bar{B}$, and let $N=\sum_{i=0}^{\bar{B}} N_{i}$ and $q_{i}\defiN_{i}/N$, $i=0,\ldots,\bar{B}$. The classical stochastic optimization theory would suggest using Sample Average Approximation (SAA) method to solve the SAA counterpart
% \[
%   \min_{x\geq0} \sum_{i=0}^{\bar{B}} q_{i}[h(x-i)_{+} + b(i-x)_{+}].
% \]
% Then the empirical $\frac{b}{b+h}$-quantile gives an optimal solution
% $
%   \hx\defi\inf_{i\geq0}\left\{i: \sum_{j=0}^{i} q_{j}\geq\frac{b}{b+h} \right\}.
% $
% Since the problem is merely one-dimensional, SAA solution $\hx$ should give an accurate approximation provided that the sample size $N$ is not too small (\citet{levi2015data}).
%On the other hand,

Usually the demand distribution $\mu$ is not known.
Then the decision maker may want to consider the DRSO problem
\[
\label{eqn:newsvendor}
\min_{x \geq 0} \ \sup_{\bp \in \Delta_{B}} \left\{\E_{\bp}[h [x - \xi]^+ + b [\xi - x]^+] \; : \; W_{p}(\bp,\bq) \leq \theta\right\}.
\]
Using Corollary~\ref{cor:finite bins}, we obtain a convex optimization reformulation
\[
\label{eqn:newsvendor_reform}
\min_{x \geq 0, \lambda \geq 0, y \defi (y_{0},\ldots,y_{B})} \left\{\lambda \theta^{p} + \sum_{i=0}^{B} \nu_{i} y_{i} \; : \; y_{i} \geq \max\big\{h (x - j), b (j - x)\big\} - \lambda |i - j|^{p}, \; \forall \; 0 \leq i,j \leq B\right\}.
\]
% In addition, given any solution $(x,\lambda)$, the worst-case distribution can then be computed by solving
% \[
%   \max_{p_{i}^h,p_{i}^b,\gamma_{ij}\geq 0}\left\{\sum_{i=0}^{\bar{B}} p_{i}^h h(x-i) + p_{i}^b b(i-x): \sum_{i,j}\gamma_{ij}|i-j|^{p}\leq\theta^{p},\ \sum_{j}\gamma_{ij}=p_{i}^h+p_{i}^b,\forall i, \ \sum_{i} \gamma_{ij}=q_{j},\forall j \right\}.
% \]
% Indeed, let $\phi^{\ast}(s)\defi\sup_{t\geq0}\{st-\phi(t)\}$ be its convex conjugate function.

\begin{table}
\caption{Examples of $\phi$-divergence}\label{tab:phi-divergence}
\small \setlength{\tabcolsep}{3pt}
\begin{tabular}{c|cccccc}\hline
Divergence & Kullback-Leibler & Burg entropy & $\chi^{2}$-distance & Modified $\chi^{2}$ & Hellinger & Total Variation\\ \hline
$\phi$ & $\phi_{kl}$ & $\phi_{b}$ & $\phi_{\chi^{2}}$ & $\phi_{m\chi^{2}}$ & $\phi_{h}$ & $\phi_{tv}$ \\ \hline
$\phi(t), t \geq 0$ & $t \log t$ & $-\log t$ & $\frac{1}{t}(t-1)^{2}$ & $(t-1)^{2}$ & $(\sqrt{t}-1)^{2}$ & $|t-1|$\\ \hline
$I_{\phi}(\bp,\bq)$ & $\sum p_{j} \log\left(\frac{p_{j}}{q_{j}}\right)$ & $\sum q_{j} \log\left(\frac{q_{j}}{p_{j}}\right)$ & $\sum \frac{(p_{j} - q_{j})^{2}}{p_{j}}$ & $\sum \frac{(p_{j} - q_{j})^{2}}{q_{j}}$ & $\sum (\sqrt{p_{j}} - \sqrt{q_{j}})^{2}$ & $\sum |p_{j} - q_{j}|$ \\ \hline
\end{tabular}
\end{table}

One may also consider $\phi$-divergence ambiguity sets (Table~\ref{tab:phi-divergence} shows some common $\phi$-divergences).
As mentioned in Section~\ref{sec:intro}, the worst-case distribution in a $\phi$-divergence ambiguity set may be problematic.
For example, when $\lim_{t \to \infty} \phi(t)/t = \infty$, such as $\phi_{kl}$ and $\phi_{m\chi^{2}}$, the $\phi$-divergence ambiguity set fails to include many relevant distributions.
Specifically, since $0 \phi(p_{j}/0) \defi p_{j} \lim_{t \to \infty} \phi(t)/t = \infty$ for all $p_{j} > 0$, the $\phi_{kl}$- and $\phi_{m\chi^{2}}$-divergence ambiguity sets do not include any distribution which is not absolutely continuous with respect to the nominal distribution $\bq$.
% For example, suppose $\Xi=\{0,1\}^{s}$ and the nominal distribution $\bq$ is supported on $\{\xi^{j}\}_{j=1}^{\bar{B}} \subset\Xi$. Then the ambiguity set only contains distributions that are supported on the given $\bar{B}$ data points $\{\xi^{j}\}_{i=1}^{\bar{B}}$. Unless $\bar{B}$ is large relative to $2^{s}$, the ambiguity set is far from rich enough to robustify the decision.

When $\lim_{t \to \infty} \phi(t)/t < \infty$, such as $\phi_{b}$, $\phi_{\chi^{2}}$, $\phi_h$, and $\phi_{tv}$, the situation is also bad.
Let $I_{0} \defi \{j \in \{0,\ldots,B\} \, : \, q_{j} > 0\}$.
Assume that $\Psi(\xi^{j})$ are different from each other, so that $j_{M} \in \argmax_{j \in \{0,\ldots,B\}} \{\Psi(\xi^{j}) \, : \, q_{j} = 0\}$ is unique.
Then according to \citet{ben2013robust} and \citet{doi:10.1287/educ.2015.0134}, the worst-case distribution $p^{\ast}$ satisfies
\begin{subequations}
\begin{eqnarray}
\frac{p_{j}^{\ast}}{q_{j}} & \ \ \in \ \ & \partial \phi^{\ast}\left(\frac{\Psi(\xi^{j}) - \beta^{\ast}}{\lambda^{\ast}}\right), \qquad \forall \ i \in I_{0}, \label{eqn:phi_I0} \\
p_{j}^{\ast} & \ \ = \ \ & 0, \qquad \forall \ j \notin I_{0} \cup \{j_{M}\}, \label{eqn:phi_{0}} \\
p_{j_{M}}^{\ast} & \ \ = \ \ & \left\{\begin{array}{ll}
1 - \sum_{i \in I_{0}} p_{j}^{\ast} & \textrm{ if } \beta^{\ast} = \Psi(\xi^{j_{M}}) - \lambda^{\ast} \lim_{t \to \infty} \phi(t)/t, \\
0 & \textrm{ if } \beta^{\ast} > \Psi(\xi^{j_{M}}) - \lambda^{\ast} \lim_{t \to \infty} \phi(t)/t, \end{array}\right. \label{eqn:phi_N}
\end{eqnarray}
\end{subequations}
for some $\lambda^{\ast} \geq 0$ and $\beta^{\ast} \geq \Psi(\xi^{j_{M}}) - \lambda^{\ast} \lim_{t \to \infty} \phi(t)/t$, where $\phi^{\ast}$ denotes the convex conjugate function of $\phi$.
% Let us consider two cases.
% When $\lim_{t\to\infty}\phi(t)/t=\infty$, from (\ref{eqn:phi_N}) we have $p_N=0$, which is also consistent with the previous discussion. In this case, note that $\phi^{\ast}$ is convex thus $\partial\phi^{\ast}$ is monotonically increasing. Then (\ref{eqn:phi_I0}) suggests $p_{i}^{\ast}/q_{i}$ is also monotonically increasing in $i$. Hence there exists a threshold $1<N_{0}<N$, such that $p_{i}^{\ast}/q_{i}\leq 1$ for $i\leq N_{0}$ and $p_{i}^{\ast}/q_{i}>1$ for $i>N_{0}$. Therefore the worst-case distribution has the tendency to put comparatively more weight on scenario $i$ with higher objective value $\Psi(\hxi^{i})$, which makes the worst-case distribution too conservative.
% When $\lim_{t\to\infty}\phi(t)/t<\infty$, although the ambiguity set is allowed to contain distributions that are not absolutely continuous with respect to the reference,
According to (\ref{eqn:phi_{0}}), the support of the worst-case distribution and that of the nominal distribution can differ by at most one point $\xi^{j_{M}}$.
According to (\ref{eqn:phi_N}), the probability mass $p^{\ast}_{j_{M}}$ is moved away from scenarios in $I_{0}$ to the worst scenario $\xi^{j_{M}}$.
% This is called ``popping'' behavior in \citet{doi:10.1287/educ.2015.0134}.
Note that in many applications where the support of $\xi$ is unknown, the choice of the underlying space $\Xi$ (such as $\{0,1,\ldots,B\}$) may be somewhat arbitrary, and therefore the worst scenario $\xi^{j_{M}}$ (such as $0$ or $B$) may be arbitrary.
Hence the worst-case behavior is too sensitive to arbitrary choices in the specification of $\Xi$.
% In the discussion above we assume that both the underlying distribution and the reference are discrete distributions. To avoid some drawbacks of $\phi$-divergence discussed previously, one may be willing to modify the discrete nominal distribution $\bq$ into continuous one via parametric or nonparametric density estimation (see \citet{ben2013robust,Jiang2015Data-driven} for more details). However, such manipulation unnecessarily introduces further ambiguity and the phenomenon in previous examples is not fully eliminated.

We present numerical examples of which the setup is similar to that in \citet{Wang2015Likelihood} and \citet{ben2013robust}.
Let $b = h = 1$, $B = 100$, and $N \in \{50,500\}$ representing small and large datasets.
Random data are generated from Binomial$(100,0.5)$ and Geometric$(0.1)$ truncated on $[0,100]$.
We estimate the radius of the ambiguity set such that it covers the underlying distribution with probability greater than $0.95$ (see Appendix C).
% For Burg entropy, we use the estimation developed in \citet{Wang2015Likelihood}, while for Wasserstein distance, the radius is determined by a classical concentration inequality for Wasserstein distance developed in \citet{bolley2007quantitative} (see Appendix \ref{sec:theta} for more details). Figure \ref{fig:numericalNewsvendor} shows the histograms of the worst-case distributions under the optimal initial inventory level determined by solving the DRSOs with Wasserstein and Burg entropy ambiguity set.

\begin{figure}
\centering
\subfloat[Binomial$(100,0.5)$, $N=500$]{
\label{fig:newsvendor500bino} %% label for first subfigure
\begin{minipage}[b]{0.36\textwidth}
\centering
\includegraphics[width=\textwidth]{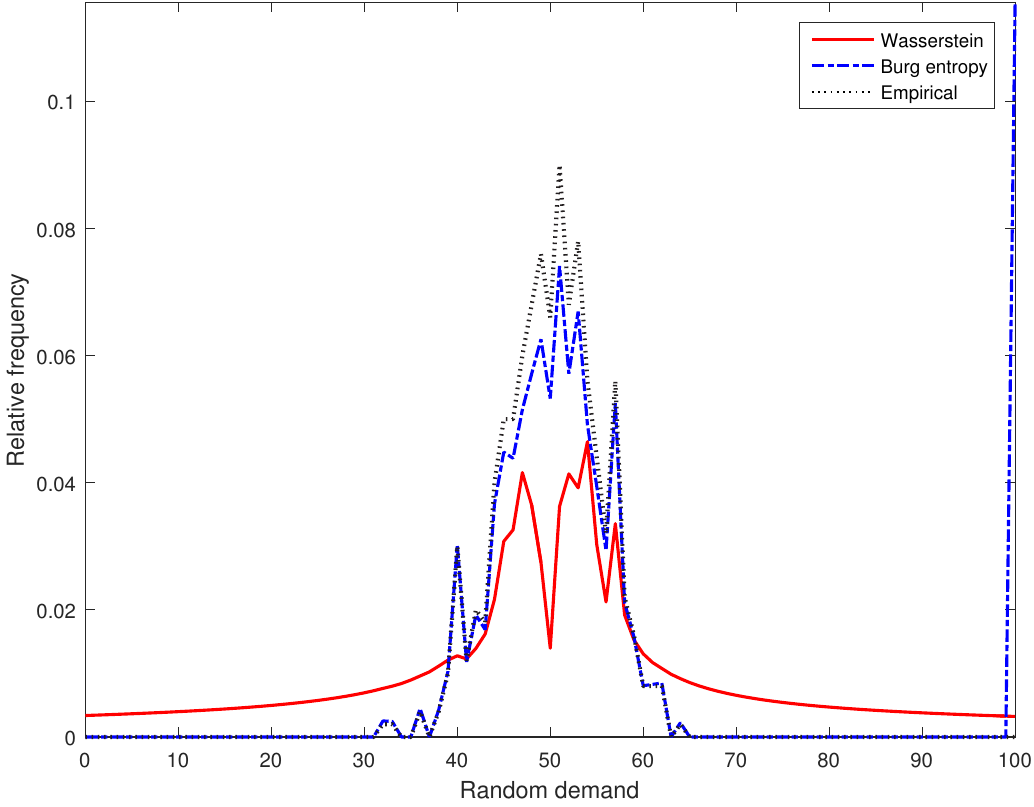}
\end{minipage}}%
\quad
\subfloat[Binomial$(100,0.5)$, $N=50$]{
\label{fig:newsvendor50bino} %% label for second subfigure
\begin{minipage}[b]{0.36\textwidth}
\centering
\includegraphics[width=\textwidth]{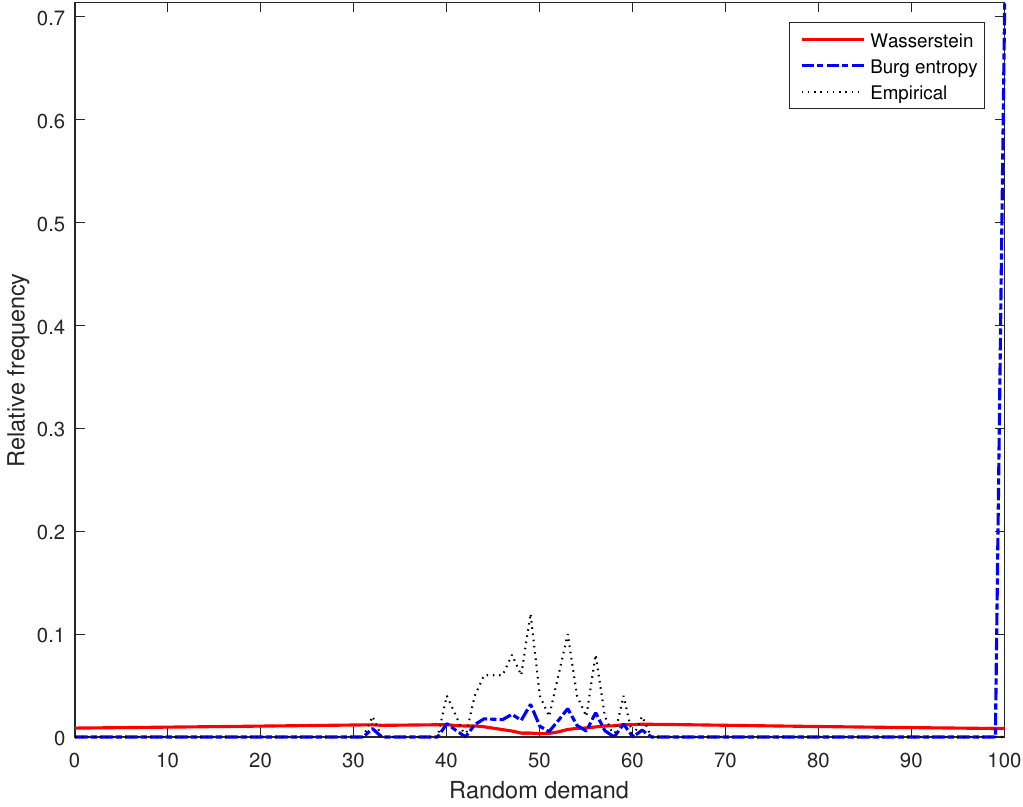}
\end{minipage}}\\
\subfloat[truncated Geometric$(0.1)$, $N=500$]{
\label{fig:newsvendor500geo} %% label for first subfigure
\begin{minipage}[b]{0.36\textwidth}
\centering
\includegraphics[width=\textwidth]{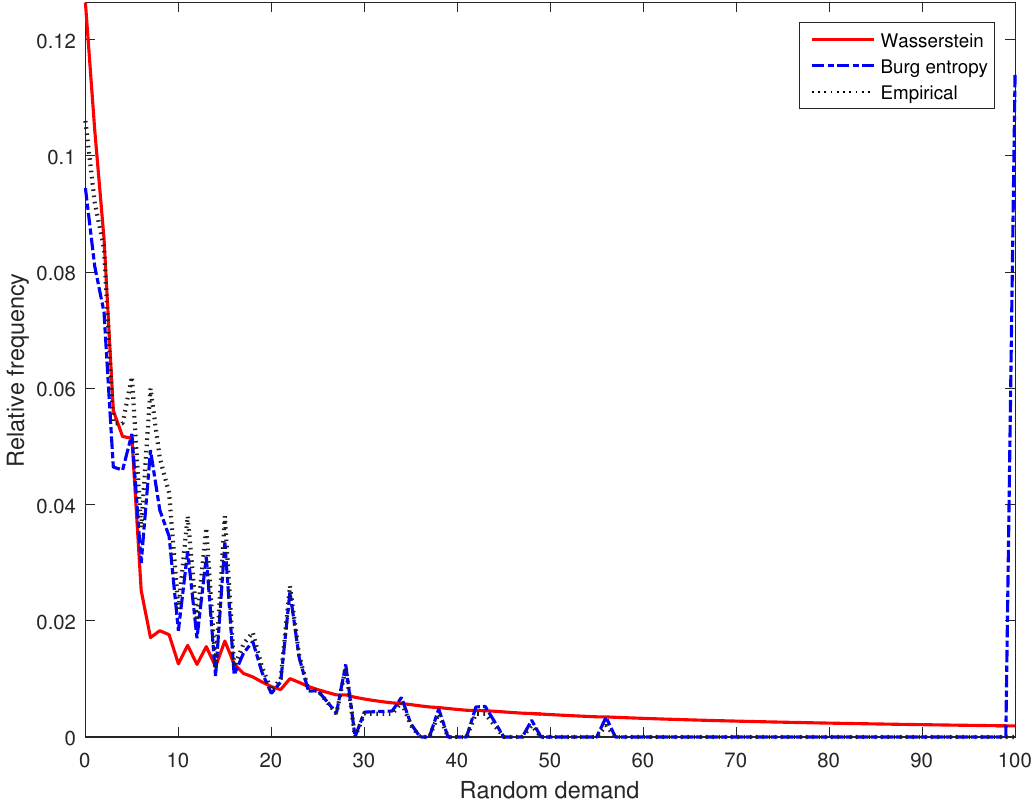}
\end{minipage}}%
\quad
\subfloat[truncated Geometric$(0.1)$, $N=50$]{
\label{fig:newsvendor50geo} %% label for second subfigure
\begin{minipage}[b]{0.36\textwidth}
\centering
\includegraphics[width=\textwidth]{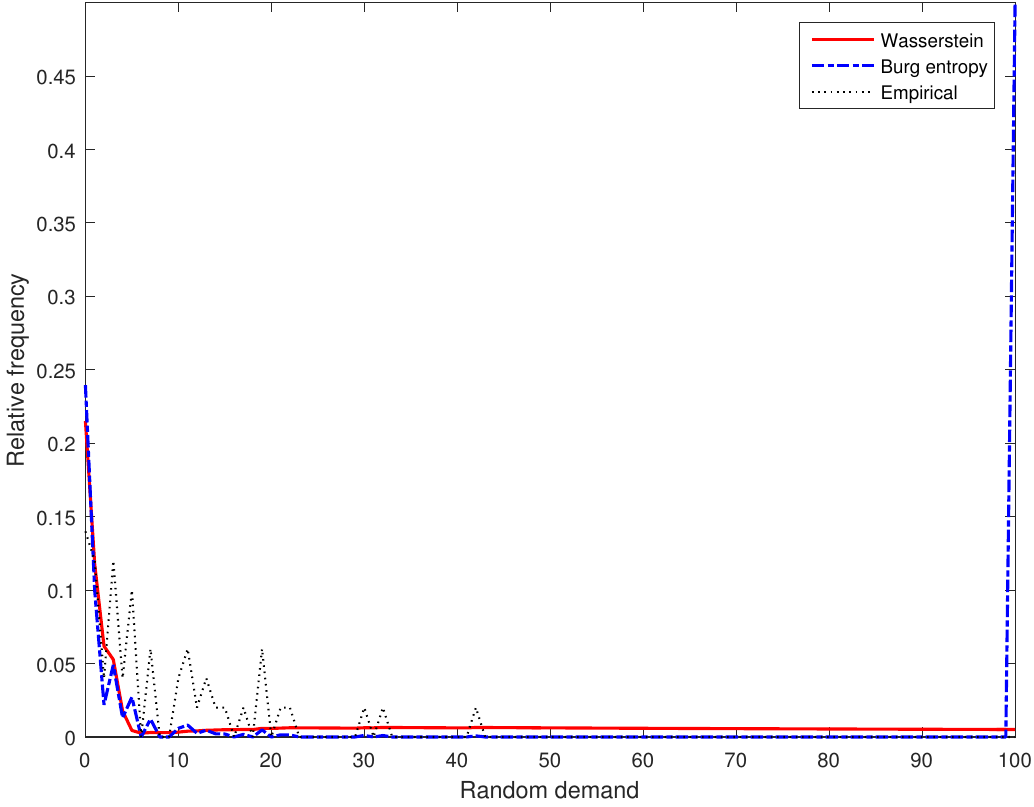}
\end{minipage}}
\caption{Histograms of worst-case distributions resulting from Wasserstein distance and Burg entropy}
\label{fig:numericalNewsvendor} %% label for entire figure
\end{figure}

When the underlying distribution is Binomial, the symmetry of the Binomial distribution and $b = h = 1$ implies that the optimal initial inventory level is close to $50$.
Intuitively, the corresponding worst-case distribution should make provision for both the possibility of high demand and the possibility of low demand.
This intuition is consistent with the worst-case distributions in the Wasserstein ambiguity set, shown by the solid curves in Figures~(\ref{fig:newsvendor500bino})--(\ref{fig:newsvendor50bino}).
The tails of these worst-case distributions are heavy on both the high side and the low side, and are quite smooth and reasonable for both small and large datasets.
In contrast, if Burg entropy is used to define the ambiguity set, then the worst-case distributions have disconnected support, as shown by the dashed curves in Figures~(\ref{fig:newsvendor500bino})--(\ref{fig:newsvendor50bino}).
There is a large spike on the boundary $B = 100$, displaying the ``popping'' behavior mentioned in \citet{doi:10.1287/educ.2015.0134}.
Especially when the dataset is small, the worst-case spike is huge, which makes the solution too conservative.

When the underlying distribution is Geometric, the worst-case distributions in the Wasserstein ambiguity set make provision for both the possibility of high demand and the possibility of low demand with fairly smooth variation of probability in between, as shown by solid curves in Figures~(\ref{fig:newsvendor500geo})--(\ref{fig:newsvendor50geo}).
As before, if Burg entropy is used, then the tail has unrealistic spikes on the boundary, and thus the worst-case distribution is very sensitive to the somewhat arbitrary choice of truncation value $B$.
In these examples the Wasserstein ambiguity set seems to yield a more reasonable worst-case distribution.

\subsection{Two-stage DRSO: Connection with Robust Optimization}
\label{sec:two-stage}

In Corollary~\ref{cor:finite}\ref{itm:finiteDual_robustapprox} we established the close connection between the DRSO problem and robust optimization.
More specifically, we showed that every DRSO problem can be approximated with high accuracy by robust optimization problems, which facilitates practical application of DRSO problems.
To illustrate this point, in this section we show the tractability of two-stage linear DRSO problems.

Consider the two-stage distributionally robust stochastic optimization problem
\begin{equation}
\label{eqn:two-stage}
\min_{x \in X} c^{\top} x + \sup_{\mu \in \frakM} \E_{\mu}[\Psi(x,\xi)],
\end{equation}
where $\Psi(x,\xi)$ is the optimal objective value of the second-stage problem
\[
\label{eqn:two-stage_second}
\Psi(x,\xi) \ \ \defi \ \ \min_{y \in \R^{m}} \left\{q(\xi)^{\top} y \; : \; T(\xi) x + W(\xi) y \leq h(\xi)\right\}
\]
and
\[
q(\xi) \ = \ q^{0} + \sum_{l=1}^{s} \xi_{l} q^{l}, \ \
T(\xi) \ = \ T^{0} + \sum_{l=1}^{s} \xi_{l} T^{l}, \ \
W(\xi) \ = \ W^{0} + \sum_{l=1}^{s} \xi_{l} W^{l}, \ \
h(\xi) \ = \ h^{0} + \sum_{l=1}^{s} \xi_{l} h^{l}.
\]
Assume that $p = 2$ and $\Xi = \R^{s}$ with Euclidean distance $d$.
In general, the two-stage problem (\ref{eqn:two-stage}) is NP-hard (see Section 3.3.2 in \cite{guslitser2002uncertainty}).
However, with tools from robust optimization, we are able to obtain a tractable approximation of (\ref{eqn:two-stage}).
Let $\frakM_{1} \defi \{(\xi^{1},\ldots,\xi^{N}) \in \Xi^{N} \, : \, \frac{1}{N} \sum_{i=1}^{N} \|\xi^{i} - \hxi^{i}\|_{2}^{2} \leq \theta^{2}\}$.
Using Corollary~\ref{cor:finite}\ref{itm:finiteDual_robustapprox} with $K=1$, we obtain an adjustable robust optimization approximation
\begin{equation}
\label{eqn:two-stage_approx}
\begin{aligned}
& \min_{x\in X} \left\{c^{\top} x + \sup_{(\xi^{i})_{i=1}^{N} \in \frakM_{1}} \frac{1}{N} \sum_{i=1}^{N} \Psi(x,\xi^{i})\right\} \\
= \ \ & \min_{\substack{x \in X, t \in \R, \\ y : \Xi \mapsto \R^m}} \left\{t \; : \;
\begin{array}{l}
t \ \geq \ c^{\top} x + \frac{1}{N} \sum_{i=1}^{N} q(\xi^{i})^{\top} y(\xi^{i}), \ \forall \ (\xi^{i})_{i=1}^{N} \in \frakM_{1}, \vspace{2mm} \\
T(\xi) x + W(\xi) y(\xi) \ \leq \ h(\xi), \ \forall \ \xi \in \bigcup_{i=1}^{N} \{\xi' \in \Xi \; : \; \|\xi' - \hxi^{i}\|_{2} \leq \theta \sqrt{N}\}
\end{array}\right\},
\end{aligned}
\end{equation}
where the second set of inequalities follows from the fact that $T(\xi)x+W(\xi)y(\xi)\leq h(\xi)$ should hold for any realization $\xi$ with positive probability for some distribution in $\frakM_{1}$.
Although problem~(\ref{eqn:two-stage_approx}) in general is still intractable, there has been a substantial literature on different approximations for problem~(\ref{eqn:two-stage_approx}).
One approach is to consider the affinely adjustable robust counterpart (AARC), as follows.
Consider $y$ that is an affine function of $\xi$:
\[
y(\xi) \ \ = \ \ y^{0} + \sum_{l=1}^{s} \xi_{l} y^{l}, \ \ \forall \ \xi \in \bigcup_{i=1}^{N} B^{i},
\]
for some chosen $y^{0},y^{l} \in \R^{m}$, where $B^{i} \defi \{\xi' \in \Xi \, : \, \|\xi' - \hxi^{i}\|_{2} \leq \theta \sqrt{N}\}$.
Then the AARC of (\ref{eqn:two-stage_approx}) is
\begin{equation}
\label{eqn:AARC}
\begin{aligned}
\min_{\substack{x \in X, t \in \R, \\ y^{l} \in \R^{m}, l = 0,\ldots,s}} & \Bigg\{t \; : \; c^{\top} x + \frac{1}{N} \sum_{i=1}^{N} \left(q^{0} + \sum_{l=1}^{s} \xi^{i}_{l} q^{l}\right)^{\top} \left(y^{0} + \sum_{l=1}^{s} \xi^{i}_{l} y^{l}\right) - t \ \leq \ 0, \ \forall \ (\xi^{i})_{i=1}^{N} \in \frakM_{1}, \\
& \left(T^{0} + \sum_{l=1}^{s} \xi_{l} T^{l}\right) x + \left(W^{0} + \sum_{l=1}^{s} \xi_{l} W^{l}\right) \left(y^{0} + \sum_{l=1}^{s} \xi_{l} y^{l}\right) - \left(h^{0} + \sum_{l=1}^{s} \xi_{l} h^{l}\right) \ \leq \ 0, \  \forall \ \xi \in \bigcup_{i=1}^{N} B^{i}\Bigg\}.
\end{aligned}
\end{equation}
Set $\zeta_{il} \defi \xi^{i}_{l} - \hxi^{i}_{l}$ for $i = 1,\ldots,N$ and $l = 1,\ldots,s$.
Note that $(\xi^{1},\ldots,\xi^{N}) \in \frakM_{1}$ if and only if
\[
\zeta \ \ \defi \ \ (\zeta_{il})_{i,l} \ \ \in \ \ \calU \ \ \defi \ \ \{(\zeta'_{il})_{i,l} \; : \; \sum_{i=1}^{N} \sum_{l=1}^{s} {\zeta'}_{il}^{2} \leq N \theta^{2}\}.
\]
Set $z \defi \left(x,t,\{y^{l}\}_{l=0}^{s}\right)$, and let
\begin{align*}
\alpha_{0}(z) \ \ & \defi \ \ - \left[c^{\top} x + \frac{1}{N} \sum_{i=1}^{N} \left(q^{0} + \sum_{l=1}^{s} \hxi^{i}_{l} q^{l}\right)^{\top} \left(y^{0} +\sum_{l=1}^{s} \hxi^{i}_{l} y^{l}\right) - t\right], \\
\beta_{0}^{il}(z) \ \ & \defi \ \ - \frac{\left[\left(q^{0} + \sum_{l'=1}^{s} \hxi^{i}_{l'} q^{l'}\right)^{\top} y^{l} + {q^{l}}^{\top} \left(y^{0} + \sum_{l'=1}^{s} \hxi^{i}_{l'} y^{l'}\right)\right]}{2N}, \qquad \forall \ i = 1,\ldots,N, \ l = 1,\ldots,s. \\
\Gamma_{0}^{(l,l')}(z) \ \ & \defi \ \ - \frac{{q^{l}}^{\top} {y^{l'}} + {q^{l'}}^{\top} {y^{l}}}{2N},
\qquad \forall \ l,l' = 1,\ldots,s.
\end{align*}
Then the first set of constraints in (\ref{eqn:AARC}) is equivalent to
\begin{equation}
\label{eqn:AARC_1stconstr}
\alpha_{0}(z) + 2 \sum_{i=1}^{N} \sum_{l=1}^{s} \beta_{0}^{il}(z) \zeta_{il} + \sum_{i=1}^{N} \sum_{l=1}^{s} \sum_{l'=1}^{s} \Gamma_{0}^{(l,l')}(z) \zeta_{il} \zeta_{il'} \ \ \geq \ \ 0, \qquad \forall \ \zeta \in \calU.
\end{equation}
It follows from Theorem~4.2 in \citet{ben2004adjustable} that (\ref{eqn:AARC_1stconstr}) holds if and only if there exists $\lambda_{0} \geq 0$ such that
\begin{align*}
(\alpha_{0}(z) - \lambda_{0}) v^{2} + 2 v \sum_{i=1}^{N} \sum_{l=1}^{s} \beta_{0}^{il}(z) w_{il} + \sum_{i=1}^{N} \sum_{l=1}^{s} \sum_{l'=1}^{s} \Gamma_{0}^{(l,l')}(z) w_{il} w_{il'} + \frac{\lambda_{0}}{N \theta^{2}} \sum_{i=1}^{N} \sum_{l=1}^{s} w_{il}^{2} \ \ \geq \ \ 0, \\
\forall \ v \in \R, \ w_{il} \in \R, \ i = 1,\ldots,N, \ l = 1,\ldots,s.
\end{align*}
In matrix form,
\begin{equation}
\label{eqn:AARC_1stconstr_SDP}
\exists \lambda_{0} \geq 0 \ : \ \left(\begin{array}{cc}
\Gamma_{0}(z) \otimes I_{N} + \frac{\lambda_{0}}{N \theta^{2}} I_{sN} \ \ & \ \ \mathrm{vec}(\beta_{0}(z)) \vspace{2mm} \\
\mathrm{vec}(\beta_{0}(z))^{\top} \ \ & \ \ \alpha_{0}(z) - \lambda_{0}
\end{array}\right) \ \ \succeq \ \ 0,
\end{equation}
where $I_{N}$ (resp. $I_{sN}$) is the $N \times N$ (resp. $sN \times sN$) identity matrix, $\otimes$ is the Kronecker product of matrices, $\mathrm{vec}(\beta_{0})$ is the vectorization of matrix $\beta_{0}(z)$, and $\Gamma_0(z)$ is a matrix whose $(l,l')$-element equal to $\Gamma_0^{(l,l')}(z)$.

Next we reformulate the second set of constraints in (\ref{eqn:AARC}).
Let $T_{j}^{l}$ and $W_{j}^{l}$ denote row~$j$ of $T^{l}$ and $W^{l}$, respectively, $j = 1,\ldots,J$.
For all $i = 1,\ldots,N$, $j = 1,\ldots,J$, and $l,l' = 1,\ldots,s$, set
\begin{align*}
\alpha_{ij}(z) \ \ & \defi \ \ -\left[\left(T_{j}^{0} + \sum_{l=1}^{s} \hxi^{i}_{l} T_{j}^{l}\right) x + \left(W_{j}^{0} + \sum_{l=1}^{s} \hxi^{i}_{l} W_{j}^{l}\right) \left(y^{0} + \sum_{l=1}^{s} \hxi^{i}_{l} y^{l}\right) - \left(h_{j}^{0} + \sum_{l=1}^{s} \hxi^{i}_{l} h_{j}^{l}\right)\right], \\
\beta_{ij}^{l}(z) \ \ & \defi \ \ -\frac{\left[T_{j}^{l} x + \left(W_{j}^{0} + \sum_{l=1}^{s} \hxi^{i}_{l} W_{j}^{l}\right) y^{l} + W_{j}^{l} \left(y^{0} + \sum_{l=1}^{s} \hxi^{i}_{l} y^{l}\right) - h_{j}^{l}\right]}{2}, \\
\Gamma_{j}^{(l,l')}(z) \ \ & \defi \ \ - \frac{W_{j}^{l} y^{l'} + W_{j}^{l'} y^{l}}{2}.
\end{align*}
Let $\eta^{i} \defi \xi - \hxi^{i}$, $\beta_{ij}(z) \defi (\beta_{ij}^{l}(z))_{l}$, and $\Gamma_{j}(z) \defi (\Gamma_{j}^{(l,l')}(z))_{l,l'}$.
Then the second set of constraints in~(\ref{eqn:AARC}) is equivalent to
\[
\alpha_{ij}(z) + 2 \beta_{ij}(z)^{\top} \eta^{i} + {\eta^{i}}^{\top} \Gamma_{j}(z) \eta^{i} \ \ \geq \ \ 0, \qquad \forall \ \eta^{i} \in \{\eta' \in \R^{s} \; : \; \|\eta'\|_{2} \leq \theta \sqrt{N}\}, \ i = 1,\ldots,N, \ j = 1,\ldots,J.
\]
Again by Theorem~4.2 in \citet{ben2004adjustable} the second set of constraints in~(\ref{eqn:AARC}) is equivalent to
\begin{equation}
\label{eqn:AARC_2ndconstr_SDP}
\exists \lambda_{ij} \geq 0 \ : \ \left(\begin{array}{cc}
\Gamma_{j}(z) + \frac{\lambda_{ij}}{N \theta^{2}} \; I_{s} \ \ & \ \ \beta_{ij}(z) \vspace{2mm} \\
\beta_{ij}(z)^{\top} \ \ & \ \ \alpha_{ij}(z) - \lambda_{ij}
\end{array}\right) \ \ \succeq \ \ 0, \qquad \forall \ i = 1,\ldots,N, \ j = 1,\ldots,J.
\end{equation}
Combining (\ref{eqn:AARC_1stconstr_SDP}) and (\ref{eqn:AARC_2ndconstr_SDP}) we obtain the following result.

\begin{proposition}
An exact semidefinite program reformulation of the AARC of (\ref{eqn:two-stage_approx}) is given by
\begin{equation}
\label{eqn:two-stage_SDP}
\min_{\substack{x \in X, t \in \R, y^{l} \in \R^{m}, l = 0,\ldots,s, \\ \lambda_{0}, \lambda_{ij} \geq 0, i = 1,\ldots,N, j = 1,\ldots,J}} \left\{t \ : \ (\ref{eqn:AARC_1stconstr_SDP}), (\ref{eqn:AARC_2ndconstr_SDP}) \textrm{ holds}\right\}.
\end{equation}
\end{proposition}

By Corollary~\ref{cor:finite}\ref{itm:finiteDual_robustapprox}, (\ref{eqn:two-stage_approx}) is a fairly good approximation of the original two-stage DRSO problem (\ref{eqn:two-stage}).
Hence, as long as the AARC of (\ref{eqn:two-stage_approx}) is reasonably good, its semidefinite program reformulation (\ref{eqn:two-stage_SDP}) provides a good tractable approximation of the two-stage linear DRSO (\ref{eqn:two-stage}).

\section{Conclusions}
\label{sec:conclusion}

In this paper, we developed a constructive proof method to derive the dual reformulation of distributionally robust stochastic optimization with Wasserstein distance for a general setting.
This approach allows us to obtain a precise structural description of the worst-case distribution.
It also facilitates a connection between distributionally robust stochastic optimization and robust optimization.
We showed how the results can be used to obtain theoretical and computational conclusions for a variety of problems.
%For future work, extensions to multi-stage distributionally robust stochastic optimization will be explored.

% Appendix here
% Options are (1) APPENDIX (with or without general title) or
%             (2) APPENDICES (if it has more than one unrelated sections)
% Outcomment the appropriate case if necessary
%
% \begin{APPENDIX}{<Title of the Appendix>}
% \end{APPENDIX}
%
%   or
%
\begin{APPENDICES}
% \section{<Title of Section A>}
% \section{<Title of Section B>}
% etc

\section{Proofs for Section~\ref{sec:wasserstein}}

\proof{Proof of Lemma~\ref{lemma:enlarge_{p}}.}
Let $(u_{0},v_{0})$ be any feasible solution for the maximization problem in~(\ref{eqn:Wp_dual}).
For any $t \in \R$ and any $\xi,\zeta \in \Xi$, let $u_{t}(\xi) \defi u_{0}(\xi) + t$ and $v_{t}(\zeta) \defi v_{0}(\zeta) - t$.
Then it follows that $u_{t}(\xi) + v_{t}(\zeta) \leq d^{p}(\xi,\zeta)$ for all $\xi,\zeta \in \Xi$, and
\[
\int_{\Xi} u_{t}(\xi) \mu(d\xi) + \int_{\Xi} v_{t}(\zeta) \nu(d\zeta)
\ \ = \ \ \int_{\Xi} u_{0}(\xi) \mu(d\xi) + \int_{\Xi} v_{0}(\zeta) \nu(d\zeta) + t [\mu(\Xi) - \nu(\Xi)].
\]
Since $\mu(\Xi) \neq \nu(\Xi)$,
\[
\sup_{t \in \R} \left\{\int_{\Xi} u_{t}(\xi) \mu(d\xi) + \int_{\Xi} v_{t}(\zeta) \nu(d\zeta)\right\} \ \ = \ \ \infty,
\]
and thus $W_{p}^{p}(\mu,\nu) = \infty$.
% \hfillqed
\endproof

% \proof{Proof of Lemma~\ref{lemma:Emu-Enu}.}
% Let $\gamma_{0}$ be a minimizer of the minimization problem in \eqref{eqn:def_wasserstein}. It follows that
% \begin{align*}
% \big|\E_{\mu}[\Psi(\xi)] - \E_{\nu}[\Psi(\xi)]\big|
% & \leq & \int_{\Xi \times \Xi} |\Psi(\xi) - \Psi(\zeta)| \gamma_{0}(d\xi,d\zeta) \\
% & \leq & \int_{\Xi \times \Xi} \left(L d^{p}(\xi,\zeta) + M\right) \gamma_{0}(d\xi,d\zeta) \\
% & = & L W_{p}^{p}(\mu,\nu) + M \\
% & \leq & L \theta^{p} + M.
% \end{align*}
% \hfillqed
% \endproof
% \vspace{1em}

\section{Proofs for Section \ref{sec:tractability}}

\subsection{Proofs for Section \ref{sec:dual_general}}

\subsubsection{Auxiliary results}
\label{sec:auxiliary results}

\begin{lemma}
\label{lemma:inequality}
Consider any $p \geq 1$ and any $\varepsilon > 0$.
Then there exists $C_{p}(\varepsilon) \geq 1$ such that
\[
(x+y)^{p} \ \ \leq \ \ (1 + \varepsilon) x^{p} + C_{p}(\varepsilon) y^{p}
\]
for all $x,y \geq 0$.
\end{lemma}

\proof{Proof of Lemma \ref{lemma:inequality}.}
Note that if $x = 0$, then the inequality holds for any $C_{p}(\varepsilon) \ge 1$.
Next we consider the case with $x > 0$, and we let $t \defi y/x$.
Let
\[
t_{0}(\varepsilon) \ \ \defi \ \ \sup\{t > 0 \; : \; 1 + \varepsilon \, \geq \, (1+t)^{p}\}.
\]
Note that $t_{0}(\varepsilon) > 0$.
Next let
\[
C_{p}(\varepsilon) \ \ \defi \ \ \max\left\{1, \ \sup_{t \geq t_{0}(\varepsilon)} \frac{(1+t)^{p-1}}{t^{p-1}}\right\}.
\]
Note that $C_{p}(\varepsilon) < \infty$ because $\lim_{t \to \infty} (1+t)^{p-1} / t^{p-1} = 1$.
Next, consider
\[
f(t) \ \ \defi \ \ 1 + \varepsilon + C_{p}(\varepsilon) t^{p} - (1+t)^{p}
\]
Note that $f(t) \ge 0$ for all $t \in [0,t_{0}(\varepsilon)]$.
Also, $f'(t) = C_{p}(\varepsilon) p t^{p-1} - p (1+t)^{p-1} \ge 0$ for all $t \in [t_{0}(\varepsilon),\infty)$.
Therefore $f(t) \geq 0$ for all $t \ge 0$, which establishes the inequality for $x > 0$.
\hfillqed
\endproof

\begin{lemma}
\label{lemma:D_bound}
Consider any $\zeta^{0} \in \Xi$. Then for any $\lambda > \lambda_{1} > \kappa$, there exists a constant $C > 0$ such that
\[
\frac{\lambda - \lambda_{1}}{2} \Dp(\lambda,\zeta) \ \ \leq \ \ \Phi(\lambda,\zeta) - \Phi(\lambda_{1},\zeta^{0}) + \lambda_{1} C d^{p}(\zeta,\zeta^{0})
\]
for all $\zeta \in \Xi$.
\end{lemma}

\proof{Proof of Lemma \ref{lemma:D_bound}.}
It follows from Lemma~\ref{lemma:inequality} with $\varepsilon \defi \frac{\lambda - \lambda_{1}}{2 \lambda_{1}}$ that
\[
\lambda_{1} d^{p}(\xi,\zeta^{0}) ~\leq ~\frac{\lambda + \lambda_{1}}{2} d^{p}(\xi,\zeta) + \lambda_{1} C_{p}(\varepsilon) d^{p}(\zeta,\zeta^{0})
\]
for all $\xi,\zeta,\zeta^{0} \in \Xi$.
Thus
\begin{eqnarray*}
\lambda d^{p}(\xi,\zeta) - \Psi(\xi) & = & \frac{\lambda - \lambda_{1}}{2} d^{p}(\xi,\zeta) - \Psi(\xi) + \frac{\lambda + \lambda_{1}}{2} d^{p}(\xi,\zeta) \\
& \geq & \frac{\lambda - \lambda_{1}}{2} d^{p}(\xi,\zeta) - \Psi(\xi) + \lambda_{1} d^{p}(\xi,\zeta^{0}) - \lambda_{1} C_{p}(\varepsilon) d^{p}(\zeta,\zeta^{0}) \\
& \geq & \frac{\lambda - \lambda_{1}}{2} d^{p}(\xi,\zeta) + \Phi(\lambda_{1},\zeta^{0}) - \lambda_{1} C_{p}(\varepsilon) d^{p}(\zeta,\zeta^{0}).
\end{eqnarray*}
Hence, for every $\xi \in \Xi$ that satisfies $\lambda d^{p}(\xi,\zeta) - \Psi(\xi) < \Phi(\lambda,\zeta) + \delta$ for some $\delta \geq 0$, it holds that
\begin{align*}
& \frac{\lambda - \lambda_{1}}{2} d^{p}(\xi,\zeta) \ \ < \ \ \Phi(\lambda,\zeta) - \Phi(\lambda_{1},\zeta^{0}) + \lambda_{1} C_{p}(\varepsilon) d^{p}(\zeta,\zeta^{0}) + \delta \\
%\Rightarrow \ & & \ \frac{\lambda - \lambda_{1}}{2} \sup_{\xi \in \Xi} \big\{d^{p}(\xi,\zeta) \; : \; \lambda d^{p}(\xi,\zeta) - \Psi(\xi) \, < \, \Phi(\lambda,\zeta) + \delta\big\} \ \ \leq \ \ \Phi(\lambda,\zeta) - \Phi(\lambda_{1},\zeta^{0}) + \lambda_{1} C_{p}(\varepsilon) d^{p}(\zeta,\zeta^{0}) + \delta \\
\Rightarrow \ \ \ & \frac{\lambda - \lambda_{1}}{2} \limsup_{\delta \downarrow 0} \Big\{\sup_{\xi \in \Xi} \big\{d^{p}(\xi,\zeta) \; : \; \lambda d^{p}(\xi,\zeta) - \Psi(\xi) \, < \, \Phi(\lambda,\zeta) + \delta\big\}\Big\} \\
& \hspace{70mm} \leq \ \ \limsup_{\delta \downarrow 0} \Big\{\Phi(\lambda,\zeta) - \Phi(\lambda_{1},\zeta^{0}) + \lambda_{1} C_{p}(\varepsilon) d^{p}(\zeta,\zeta^{0}) + \delta\Big\} \\
\Rightarrow \ \ \ & \frac{\lambda - \lambda_{1}}{2} \Dp(\lambda,\zeta) \ \ \leq \ \ \Phi(\lambda,\zeta) - \Phi(\lambda_{1},\zeta^{0}) + \lambda_{1} C_{p}(\varepsilon) d^{p}(\zeta,\zeta^{0})
\qedhere
\end{align*}
%Taking the supremum over $\xi\in\Xi$ on both sides and then the $\limsup$ with $\delta\downarrow0$, we obtain that
%\[
%\frac{\lambda - \lambda_{1}}{2} \Dp(\lambda,\zeta) \ \leq \ \Phi(\lambda,\zeta) - \Phi(\lambda_{1},\zeta^{0}) + \lambda_{1} C_{p}(\varepsilon) d^{p}(\zeta,\zeta^{0}).
%\qedhere
%\]
% \hfillqed
\endproof

\subsubsection{Proof of Proposition~\ref{prop:weakDuality}.}

Note that
\begin{align*}
v_{P} \ \ & = \ \ \sup_{\mu \in \cP(\Xi)} \left\{\int_{\Xi} \Psi(\xi) \mu(d\xi) \; : \; W_{p}^{p}(\mu,\nu) \leq \theta^{p}\right\} \\
& = \ \ \sup_{\gamma \in \cP(\Xi^{2})} \left\{\int_{\Xi^{2}} \Psi(\xi) \gamma(d\xi,d\zeta) \; : \; \int_{\Xi^{2}} d^{p}(\xi,\zeta) \gamma(d\xi,d\zeta) \leq \theta^{p}, \, \pi^{2}_{\#}\gamma = \nu\right\}.
\end{align*}
Thus, for any $\lambda \geq 0$ it holds that
\begin{align*}
v_{P} \ \ & \leq \ \ \sup_{\gamma \in \cP(\Xi^{2})} \left\{\int_{\Xi^{2}} \Psi(\xi) \gamma(d\xi,d\zeta) + \lambda \left[\theta^{p} - \int_{\Xi^{2}} d^{p}(\xi,\zeta) \gamma(d\xi,d\zeta)\right] \; : \; \pi^{2}_{\#}\gamma = \nu\right\}.
\end{align*}
Hence,
\begin{align*}
v_{P} \ \ & \leq \ \ \inf_{\lambda \geq 0} \left\{\lambda \theta^{p} + \sup_{\gamma \in \cP(\Xi^{2})} \int_{\Xi^{2}} \big[\Psi(\xi) - \lambda d^{p}(\xi,\zeta)\big] \gamma(d\xi,d\zeta) \; : \; \pi^{2}_{\#}\gamma = \nu\right\} \\
& \leq \ \ \inf_{\lambda \geq 0} \left\{\lambda \theta^{p} + \sup_{\gamma \in \cP(\Xi^{2})} \int_{\Xi^{2}} \sup_{\xi' \in \Xi} \big[\Psi(\xi') - \lambda d^{p}(\xi',\zeta)\big] \gamma(d\xi,d\zeta) \; : \; \pi^{2}_{\#}\gamma = \nu\right\} \\
& = \ \ \inf_{\lambda \geq 0} \left\{\lambda \theta^{p} + \int_{\Xi} \sup_{\xi \in \Xi} \big[\Psi(\xi) - \lambda d^{p}(\xi,\zeta)\big] \nu(d\zeta)\right\} \\
& = \ \ v_{D}.
\end{align*}

\hfillqed
% \endproof

\subsubsection{Proof of Lemma~\ref{lemma:kappa}}\leavevmode

\ref{lemma:kappa0}
We prove the result by contradiction.
Suppose that for some $\zeta^{0},\zeta^{1} \in \Xi$, it holds that
\[
\kappa^{0} \ \ \defi \ \ \limsup_{\xi \in \Xi \, : \, d(\xi,\zeta^{0}) \to \infty} \frac{\max\{0, \Psi(\xi) - \Psi(\zeta^{0})\}}{d^{p}(\xi,\zeta^{0})}
\ \ < \ \ \kappa^{1} \ \ \defi \ \ \limsup_{\xi \in \Xi \, : \, d(\xi,\zeta^{1}) \to \infty} \frac{\max\{0, \Psi(\xi) - \Psi(\zeta^{1})\}}{d^{p}(\xi,\zeta^{1})}
\]
($\kappa^{1} = \infty$ is allowed).
Choose any $\varepsilon \in (0,\kappa^{1} - \kappa^{0})$.
Then there exists an $R$ such that for all $\xi$ with $d(\xi,\zeta^{0}) > R$ it holds that
\begin{align*}
\Psi(\xi) - \Psi(\zeta^{1}) \ \ & = \ \ \Psi(\xi) - \Psi(\zeta^{0}) + \Psi(\zeta^{0}) - \Psi(\zeta^{1}) \\
& \le \ \ \max\{0, \Psi(\xi) - \Psi(\zeta^{0})\} + \Psi(\zeta^{0}) - \Psi(\zeta^{1}) \\
& < \ \ (\kappa^{0} + \varepsilon) d^{p}(\xi,\zeta^{0}) + [\Psi(\zeta^{0}) - \Psi(\zeta^{1})] \\
& \le \ \ (\kappa^{0} + \varepsilon) \left[d(\xi,\zeta^{1}) + d(\zeta^{1},\zeta^{0})\right]^{p} + [\Psi(\zeta^{0}) - \Psi(\zeta^{1})]
\end{align*}
Since $\kappa^{1} > 0$, it follows that
\begin{align*}
\kappa^{1} \ \ & = \ \ \limsup_{\xi \in \Xi \, : \, d(\xi,\zeta^{1}) \to \infty} \frac{\Psi(\xi) - \Psi(\zeta^{1})}{d^{p}(\xi,\zeta^{1})} \\
& \leq \ \ \limsup_{\xi \in \Xi \, : \, d(\xi,\zeta^{1}) \to \infty} \frac{(\kappa^{0} + \varepsilon) \left[d(\xi,\zeta^{1}) + d(\zeta^{1},\zeta^{0})\right]^{p} + [\Psi(\zeta^{0}) - \Psi(\zeta^{1})]}{d^{p}(\xi,\zeta^{1})} \\
& = \ \ \kappa^{0} + \varepsilon \ \ < \ \ \kappa^{1},
\end{align*}
which is a contradiction.

\ref{lemma:kappa_condition}
%(Sufficiency).
First we show that if there exists $\zeta^{0} \in \Xi$ and $L,M > 0$ such that $\Psi(\xi) - \Psi(\zeta^{0}) \leq L d^{p}(\xi,\zeta^{0}) + M$ for all $\xi \in \Xi$, then $\kappa < \infty$.
Let $\kappa^{0} \defi 0$ if $\Xi$ is bounded, and let
\[
\kappa^{0} \ \ \defi \ \ \limsup_{\xi \in \Xi \, : \, d(\xi,\zeta^{0}) \to \infty} \frac{\max\{0,\Psi(\xi) - \Psi(\zeta^{0})\}}{d^{p}(\xi,\zeta^{0})} \ \ \leq \ \ L \ \ < \ \ \infty
\]
if $\Xi$ is unbounded.
If $\Xi$ is unbounded, then it follows from~\ref{lemma:kappa0} that
\begin{equation}
\label{eqn:kappa_proof_0}
\kappa^{0} \ \ = \ \ \limsup_{\xi \in \Xi \, : \, d(\xi,\zeta) \to \infty} \frac{\max\{0, \Psi(\xi) - \Psi(\zeta)\}}{d^{p}(\xi,\zeta)} \ \ \ \ \forall \ \zeta \in \Xi.
\end{equation}
We are going to show that $\int_{\Xi} \Phi(\lambda,\zeta) \nu(d\zeta) > -\infty$ for all $\lambda > \kappa^{0}$, and therefore $\kappa \le \kappa^{0} < \infty$.

First we show that $\Phi(\lambda,\zeta) > -\infty$ for any $\lambda > \kappa^{0}$ and $\zeta \in \Xi$.
If $\Xi$ is bounded, then choose any $R(\zeta) > 0$ such that $d^{p}(\xi,\zeta) \le R(\zeta)$ for all $\xi \in \Xi$.
If $\Xi$ is unbounded, then it follows from~\eqref{eqn:kappa_proof_0} that for any $\zeta \in \Xi$, there is a $R(\zeta) > 0$ such that for all $\xi \in \Xi$ with $d^{p}(\xi,\zeta) > R(\zeta)$, it holds that
\[
\frac{\Psi(\xi) - \Psi(\zeta)}{d^{p}(\xi,\zeta)} \ \ < \ \ \frac{\lambda + \kappa^{0}}{2},
\]
that is, $(\frac{\lambda + \kappa^{0}}{2}) d^{p}(\xi,\zeta) - \Psi(\xi) > - \Psi(\zeta)$.
Thus, for all $\xi \in \Xi$ with $d^{p}(\xi,\zeta) > R(\zeta)$, it holds that
\begin{align*}
\lambda d^{p}(\xi,\zeta) - \Psi(\xi) \ \ & = \ \ \frac{\lambda + \kappa^{0}}{2} d^{p}(\xi,\zeta) - \Psi(\xi) + \frac{\lambda - \kappa^{0}}{2} d^{p}(\xi,\zeta) \\
& > \ \ - \Psi(\zeta) + \frac{\lambda - \kappa^{0}}{2} R(\zeta),
\end{align*}
and hence
\[
\inf_{\xi \in \Xi} \big\{\lambda d^{p}(\xi,\zeta) - \Psi(\xi) \, : \, d^{p}(\xi,\zeta) > R(\zeta)\big\} \ \ \ge \ \ - \Psi(\zeta) + \frac{\lambda - \kappa^{0}}{2} R(\zeta) \ \ > \ \ -\infty.
\]
Also, by assumption, for any $\xi \in \Xi$ it holds that
\begin{align*}
\Psi(\xi) - \Psi(\zeta^{0}) \ \ & \leq \ \ L d^{p}(\xi,\zeta^{0}) + M \\
& \leq \ \ L [d(\xi,\zeta) + d(\zeta,\zeta^{0})]^{p} + M \\
& \leq \ \ 2^{p-1} L [d^{p}(\xi,\zeta) + d^{p}(\zeta,\zeta^{0})] + M
\end{align*}
where the second inequality follows from the elementary inequality $(a+b)^{p} \leq 2^{p-1} (a^{p} + b^{p})$ for any $a,b \geq 0$ and $p \geq 1$.
Thus
\begin{align*}
\inf_{\xi \in \Xi} \big\{\lambda d^{p}(\xi,\zeta) - \Psi(\xi) \; : \; d^{p}(\xi,\zeta) \le R(\zeta)\big\}
\ \ & \ge \ \ \inf_{\xi \in \Xi} \big\{- \Psi(\xi) \; : \; d^{p}(\xi,\zeta) \le R(\zeta)\big\} \\
& \ge \ \ - \Psi(\zeta^{0}) - 2^{p-1} L R(\zeta) - 2^{p-1} L d^{p}(\zeta,\zeta^{0}) - M \ \ > \ \ -\infty.
\end{align*}
Therefore, $\Phi(\lambda,\zeta) > -\infty$ for all $\zeta \in \Xi$ and $\lambda > \kappa^{0}$.

Next we show that $\int_{\Xi} \Phi(\lambda,\zeta) \nu(d\zeta) > -\infty$ for any $\lambda > \kappa^{0}$.
Consider any $\lambda_{1} \in (\kappa^{0},\lambda)$ and any $\zeta^{0} \in \Xi$.
It follows from Lemma~\ref{lemma:D_bound} that there is a constant $C$ such that
\[
\Phi(\lambda,\zeta) \ \ \ge \ \ \frac{\lambda - \lambda_{1}}{2} \Dp(\lambda,\zeta) + \Phi(\lambda_{1},\zeta^{0}) - C d^{p}(\zeta,\zeta^{0})
\ \ \ge \ \ \Phi(\lambda_{1},\zeta^{0}) - C d^{p}(\zeta,\zeta^{0}).
\]
Thus
\[
\int_{\Xi} \Phi(\lambda,\zeta) \nu(d\zeta) \ \ \ge \ \ \Phi(\lambda_{1},\zeta^{0}) - C \int_{\Xi} d^{p}(\zeta,\zeta^{0}) \nu(d\zeta) \ \ > \ \ -\infty.
\]
Therefore $\kappa \le \kappa^{0} < \infty$.

%(Necessity).
Next we show that if there does not exist $\zeta^{0} \in \Xi$ and $L,M > 0$ such that $\Psi(\xi) - \Psi(\zeta^{0}) \leq L d^{p}(\xi,\zeta^{0}) + M$ for all $\xi \in \Xi$, then $\kappa = \infty$.
First, observe that if there exists $\zeta^{0} \in \Xi$ and $L,M > 0$ such that $\Psi(\xi) - \Psi(\zeta^{0}) \leq L d^{p}(\xi,\zeta^{0}) + M$ for all $\xi \in \Xi$, then for any $\xi,\zeta \in \Xi$ it holds that
\begin{align*}
\Psi(\xi) - \Psi(\zeta) \ \ & = \ \ \Psi(\xi) - \Psi(\zeta^{0}) + \Psi(\zeta^{0}) - \Psi(\zeta) \\
& \leq \ \ L d^{p}(\xi,\zeta^{0}) + M + \Psi(\zeta^{0}) - \Psi(\zeta) \\
& \leq \ \ L [d(\xi,\zeta) + d(\zeta,\zeta^{0})]^{p} + M + \Psi(\zeta^{0}) - \Psi(\zeta) \\
& \leq \ \ 2^{p-1} L [d^{p}(\xi,\zeta) + d^{p}(\zeta,\zeta^{0})] + M + \Psi(\zeta^{0}) - \Psi(\zeta)
\end{align*}
It follows that there exists $\zeta^{0} \in \Xi$ and $L,M > 0$ such that $\Psi(\xi) - \Psi(\zeta^{0}) \leq L d^{p}(\xi,\zeta^{0}) + M$ for all $\xi \in \Xi$ if and only if there exists $L' \defi 2^{p-1} L \geq 0$ and $M(\zeta) \defi 2^{p-1} L d^{p}(\zeta,\zeta^{0}) + M + \Psi(\zeta^{0}) - \Psi(\zeta) \in L^{1}(\nu)$ such that
$\Psi(\xi) - \Psi(\zeta) \leq L' d^{p}(\xi,\zeta) + M(\zeta)$ for all $\xi,\zeta \in \Xi$, that is, there exists $L' \geq 0$ and $M(\zeta) \in L^{1}(\nu)$ such that $- \Psi(\zeta) - M(\zeta) \leq \inf_{\xi \in \Xi} \left\{L' d^{p}(\xi,\zeta) - \Psi(\xi)\right\}$ for all $\zeta \in \Xi$.
Therefore, if there does not exist $\zeta^{0} \in \Xi$ and $L,M > 0$ such that $\Psi(\xi) - \Psi(\zeta^{0}) \leq L d^{p}(\xi,\zeta^{0}) + M$ for all $\xi \in \Xi$, then for any $\lambda \geq 0$ it holds that
\[
\inf_{\xi \in \Xi} \left\{\lambda d^{p}(\xi,\zeta) - \Psi(\xi)\right\} \ \ \notin \ \ L^{1}(\nu)
\]
which implies that $\kappa = \infty$.

\ref{lemma:kappa_expression}
It was established in the proof of \ref{lemma:kappa_condition} that if $\kappa < \infty$ then there exists $\zeta^{0} \in \Xi$ and $L,M > 0$ such that $\Psi(\xi) - \Psi(\zeta^{0}) \leq L d^{p}(\xi,\zeta^{0}) + M$ for all $\xi \in \Xi$, and then
\[
\kappa \ \ \le \ \ \kappa^{0} \ \ \defi \ \ \limsup_{\xi \in \Xi \, : \, d(\xi,\zeta^{0}) \to \infty} \frac{\max\{0,\Psi(\xi) - \Psi(\zeta^{0})\}}{d^{p}(\xi,\zeta^{0})}
\]
Next we show that $\kappa \ge \kappa^{0}$.
If $\kappa^{0} = 0$, then it follows from the definition of $\kappa$ that $\kappa \ge \kappa^{0}$.
Next, suppose that $\kappa^{0} > 0$, and consider any $\lambda \in [0,\kappa^{0})$.
We will show that $\inf_{\xi \in \Xi} \left\{\lambda d^{p}(\xi,\zeta) - \Psi(\xi)\right\} = -\infty$ for all $\zeta \in \Xi$.
If $\lambda = 0$, then it follows from $\kappa^{0} > 0$ that $\inf_{\xi \in \Xi} \left\{\lambda d^{p}(\xi,\zeta) - \Psi(\xi)\right\} = -\infty$ for all $\zeta$.
Next, consider any $\lambda \in (0,\kappa^{0})$, any $\zeta \in \Xi$, any $M > 0$, any $\lambda_{2} \in (\lambda,\kappa^{0})$, and any $\varepsilon \in (0, (\lambda_{2} - \lambda) / \lambda)$.
Since $[d(\xi,\zeta^{0}) + d(\zeta^{0},\zeta)]^{p} / d^{p}(\xi,\zeta^{0}) \to 1$ as $d^{p}(\xi,\zeta^{0}) \to \infty$, it follows that there exists $R_{1} > 0$ such that $d^{p}(\xi,\zeta) / d^{p}(\xi,\zeta^{0}) \le [d(\xi,\zeta^{0}) + d(\zeta^{0},\zeta)]^{p} / d^{p}(\xi,\zeta^{0}) \le 1 + \varepsilon$ for all $\xi \in \Xi$ such that $d^{p}(\xi,\zeta^{0}) > R_{1}$.
Choose any $R > \max\{R_{1}, [M - \Psi(\zeta^{0})] / (\lambda_{2} - \lambda - \lambda \varepsilon)\}$.
It follows from the definition of $\kappa^{0}$ that there exists $\xi \in \Xi$ such that $d^{p}(\xi,\zeta^{0}) > R$ and
\begin{align*}
\Psi(\xi) - \Psi(\zeta^{0}) \ \ & > \ \ \lambda_{2} d^{p}(\xi,\zeta^{0}) \\
& = \ \ \lambda d^{p}(\xi,\zeta^{0}) + (\lambda_{2} - \lambda) d^{p}(\xi,\zeta^{0}) \\
& \ge \ \ \lambda d^{p}(\xi,\zeta) + (\lambda_{2} - \lambda - \lambda \varepsilon) d^{p}(\xi,\zeta^{0}) \\
\Rightarrow \ \ \ \lambda d^{p}(\xi,\zeta) - \Psi(\xi) \ \ & < \ \ - \Psi(\zeta^{0}) - (\lambda_{2} - \lambda - \lambda \varepsilon) R \ \ < \ \ -M
\end{align*}
Thus, $\inf_{\xi \in \Xi} \left\{\lambda d^{p}(\xi,\zeta) - \Psi(\xi)\right\} = -\infty$ for all $\zeta$, and hence $\int_{\Xi} \inf_{\xi \in \Xi} \left\{\lambda d^{p}(\xi,\zeta) - \Psi(\xi)\right\} \nu(d\zeta) = -\infty$ for all $\lambda \in [0,\kappa^{0})$.
Therefore, $\kappa \ge \kappa^{0}$.
Next, recall that \ref{lemma:kappa0} established that $\kappa^{0}$ does not depend on the choice of $\zeta^{0}$, and therefore the result follows.
\hfillqed
% \endproof

\subsubsection{Proof of Lemma~\ref{lemma:measurable}.}\leavevmode

\ref{itm:measurable_phi}
By Definition~1.11 in \citet{ambrosio2000functions}, $\nu$ has an extension, still denoted by $\nu$, such that the measure space $(\Xi,\scrB_{\nu},\nu)$ is complete.
Note that for any $b \in \R$, it holds that
\begin{align*}
\{\zeta \in \Xi \; : \; \Phi(\lambda,\zeta) < b\}
\ \ & = \ \ \{\zeta \in \Xi \; : \; \exists \; \xi \in \Xi \mbox{ such that } \lambda d^{p}(\xi,\zeta) - \Psi(\xi) <  b\} \\
& = \ \ \pi^{2}\big(\{(\xi,\zeta) \in \Xi \times \Xi \; : \; \lambda d^{p}(\xi,\zeta) - \Psi(\xi) < b\}\big).
\end{align*}
Note that the set $\{(\xi,\zeta) \in \Xi \times \Xi \, : \, \lambda d^{p}(\xi,\zeta) - \Psi(\xi) < b\}$ on the right side is measurable.
Since $(\Xi,d)$ is Polish, it follows from the measurable projection theorem (cf. Theorem~8.3.2 in~\citet{aubin2009set}), that $\Phi(\lambda,\cdot)$ is $(\scrB_{\nu},\scrB(\R))$-measurable.

%Consider any $\lambda > \kappa$ and any $\delta > 0$.
Define functions $\overline{C},\underline{C}$ by
\begin{align*}
\overline{C}(\lambda,\zeta,\delta) \ \ & \defi \ \ \sup_{\xi \in \Xi} \{d^{p}(\xi,\zeta) \; : \; \lambda d^{p}(\xi,\zeta) - \Psi(\xi) \, < \, \Phi(\lambda,\zeta) + \delta\} \\
\underline{C}(\lambda,\zeta,\delta) \ \ & \defi \ \ \inf_{\xi \in \Xi} \{d^{p}(\xi,\zeta) \; : \; \lambda d^{p}(\xi,\zeta) - \Psi(\xi) \, < \, \Phi(\lambda,\zeta) + \delta\}.
\end{align*}
% Note that $\underline{C}(\lambda,\zeta,\delta) < \infty$ for all $\delta > 0$ and $\lambda \ge \kappa$, $\zeta \in \Xi$ such that $\Phi(\lambda,\zeta) \in \R$.
% It follows from~\ref{lemma:phi finite} that $\Phi(\lambda,\zeta) \in \R$ for all $\lambda > \kappa$ and $\zeta \in \Xi$.  ($\Phi(\kappa,\zeta)$ may or may not be finite.)
% Also, it follows from the proof of~\ref{lemma:phi_bound} that $\overline{C}(\lambda,\zeta,\delta) \in \R$ for all $\lambda > \kappa$, $\zeta \in \Xi$, and $\delta > 0$.
For any $b \in \R$ it holds that
\begin{align*}
\{\zeta \in \Xi \; : \; \overline{C}(\lambda,\zeta,\delta) > b\}
\ \ & = \ \ \{\zeta \in \Xi \; : \; \exists \; \xi \in \Xi \mbox{ such that } \lambda d^{p}(\xi,\zeta) - \Psi(\xi) \, < \, \Phi(\lambda,\zeta) + \delta, \; d^{p}(\xi,\zeta) \, > \, b\} \\
& = \ \ \pi^{2}\big(\{(\xi,\zeta) \in \Xi \times \Xi \; : \; \lambda d^{p}(\xi,\zeta) - \Psi(\xi) \, < \, \Phi(\lambda,\zeta) + \delta, \; d^{p}(\xi,\zeta) \, > \, b\}\big)
\end{align*}
and thus it follows from the measurable projection theorem that $\overline{C}(\lambda,\cdot,\delta)$ is $(\scrB_{\nu},\scrB(\R))$-measurable.
Similarly,
\begin{align*}
\{\zeta \in \Xi \; : \; \underline{C}(\lambda,\zeta,\delta) < b\}
& = \ \ \{\zeta \in \Xi \; : \; \exists \; \xi \in \Xi \mbox{ such that } \lambda d^{p}(\xi,\zeta) - \Psi(\xi) \, < \, \Phi(\lambda,\zeta) + \delta, \; d^{p}(\xi,\zeta) < b\} \\
& = \ \ \pi^{2}\big(\{(\xi,\zeta) \in \Xi \times \Xi \; : \; \lambda d^{p}(\xi,\zeta) - \Psi(\xi) \, < \, \Phi(\lambda,\zeta) + \delta, \; d^{p}(\xi,\zeta) < b\}\big)
\end{align*}
and thus $\underline{C}(\lambda,\cdot,\delta)$ is $(\scrB_{\nu},\scrB(\R))$-measurable.

Next, note that $\Dp(\lambda,\cdot) = \limsup_{\delta \downarrow 0} \overline{C}(\lambda,\cdot,\delta)$ and $\Dm(\lambda,\cdot) = \liminf_{\delta \downarrow 0} \underline{C}(\lambda,\cdot,\delta)$ are also $(\scrB_{\nu},\scrB(\R))$-measurable, because measurability is preserved under $\limsup$ and $\liminf$.
% Note that $\Dm(\lambda,\zeta) < \infty$ for all $\lambda \ge \kappa$, $\zeta \in \Xi$ such that $\Phi(\lambda,\zeta) \in \R$.
% Also, it follows from~\ref{lemma:phi_bound} that $\Dp(\lambda,\zeta) \in \R$ for all $\lambda > \kappa$ and $\zeta \in \Xi$.

For any $b \in \R$ it holds that
\begin{align*}
\{\zeta \in \Xi \; : \; \overline{D}_{0}(\lambda,\zeta) > b\}
\ \ & = \ \ \{\zeta \in \Xi \; : \; \exists \; \xi \in \Xi \mbox{ such that } \lambda d^{p}(\xi,\zeta) - \Psi(\xi) \, = \, \Phi(\lambda,\zeta), \; d^{p}(\xi,\zeta) \, > \, b\} \\
& = \ \ \pi^{2}\big(\{(\xi,\zeta) \in \Xi \times \Xi \; : \; \lambda d^{p}(\xi,\zeta) - \Psi(\xi) \, = \, \Phi(\lambda,\zeta), \; d^{p}(\xi,\zeta) \, > \, b\}\big)
\end{align*}
and thus it follows from the measurable projection theorem that $\overline{D}_{0}(\lambda,\cdot)$ is $(\scrB_{\nu},\scrB(\R))$-measurable.
Similarly,
\begin{align*}
\{\zeta \in \Xi \; : \; \underline{D}_{0}(\lambda,\zeta) < b\}
\ \ & = \ \ \{\zeta \in \Xi \; : \; \exists \; \xi \in \Xi \mbox{ such that } \lambda d^{p}(\xi,\zeta) - \Psi(\xi) \, = \, \Phi(\lambda,\zeta), \; d^{p}(\xi,\zeta) \, < \, b\} \\
& = \ \ \pi^{2}\big(\{(\xi,\zeta) \in \Xi \times \Xi \; : \; \lambda d^{p}(\xi,\zeta) - \Psi(\xi) \, = \, \Phi(\lambda,\zeta), \; d^{p}(\xi,\zeta) \, < \, b\}\big)
\end{align*}
and thus $\underline{D}_{0}(\lambda,\cdot)$ is $(\scrB_{\nu},\scrB(\R))$-measurable.

\ref{itm:measurable_selection}
%Consider any $\delta,\varepsilon > 0$.
%Define the multi-valued mappings $\overline{S},\underline{S} : \R_{+} \times \Xi \rightrightarrows \Xi$ by
%\begin{eqnarray*}
%\overline{S}(\lambda,\zeta) & \defi & \{\xi \in \Xi \; : \; \lambda d^{p}(\xi,\zeta) - \Psi(\xi) \, \leq \, \Phi(\lambda,\zeta) + \delta, \; d^{p}(\xi,\zeta) \, \geq \, \Dp(\lambda,\zeta) - \varepsilon\}, \\
%\underline{S}(\lambda,\zeta) & \defi \{\xi \in \Xi \; : \; \lambda d^{p}(\xi,\zeta) - \Psi(\xi) \, \leq \, \Phi(\lambda,\zeta) + \delta, \; d^{p}(\xi,\zeta) \, \leq \, \Dm(\lambda,\zeta) + \varepsilon\}.
%\end{eqnarray*}
% Note that $\underline{S}(\lambda,\zeta) \neq \varnothing$ if $\Dm(\lambda,\zeta) < \infty$, which holds for all $\lambda \ge \kappa$, $\zeta \in \Xi$ such that $\Phi(\lambda,\zeta) \in \R$, including all $\lambda > \kappa$ and $\zeta \in \Xi$.
% Similarly, $\overline{S}(\lambda,\zeta) \neq \varnothing$ if $\Dp(\lambda,\zeta) \in \R$, which holds for all $\lambda > \kappa$ and $\zeta \in \Xi$.
For each $\zeta \in \Xi$, it follows from the measurability of $\Psi$ and $d^{p}(\cdot,\zeta)$ that $\overline{F}^{\varepsilon}_{\delta}(\lambda,\zeta)$ and $\underline{F}^{\varepsilon}_{\delta}(\lambda,\zeta)$ are in $\scrB(\Xi)$.
Since $(\Xi,d)$ is Polish and $\nu$ is a complete finite measure, it follows from Aumann's measurable selection theorem (see, e.g. Theorem~18.26 in \citet{aliprantis2006infinite}) that $\nu$-measurable selections $\overline{T}^{\varepsilon}_{\delta}(\lambda,\cdot),\underline{T}^{\varepsilon}_{\delta}(\lambda,\cdot) : \Xi \mapsto \Xi$ exist such that $\overline{T}^{\varepsilon}_{\delta}(\lambda,\zeta) \in \overline{F}^{\varepsilon}_{\delta}(\lambda,\zeta)$ and $\underline{T}^{\varepsilon}_{\delta}(\lambda,\zeta) \in \underline{F}^{\varepsilon}_{\delta}(\lambda,\zeta)$ for $\nu$-almost all $\zeta \in \Xi$.
%When $\Phi(\kappa,\zeta)>-\infty$, using the similar argument as above, we can shown the existence of $\overline{T}^{\delta}(\zeta)$. When $\Xi=\R^{K}$, the above-mentioned measurable selections are well-defined even for $\delta,\varepsilon=0$.

\ref{itm:measurable_selection0}
The proof is the same as the proof of~\ref{itm:measurable_selection}.

\ref{itm:measurable_n}
For each $\zeta \in E$, it follows from the measurability of $\Psi$ and $d^{p}(\cdot,\zeta)$ that $F(\zeta) \in \scrB(\Xi)$.
Then using the same argument as in~\ref{itm:measurable_selection}, there exists a $\nu$-measurable selection $T : E \mapsto \Xi$ such that $T(\zeta) \in F(\zeta)$ for $\nu$-almost all $\zeta \in E$.

\ref{itm:measurable_M}
%Define a multi-valued mapping $S : \R_{+} \times \Xi \rightrightarrows \Xi$ by
%\[
%S(\varepsilon,\zeta) \ \ \defi \ \ \{\xi \in \Xi \; : \; \Psi(\xi) - \Psi(\zeta) \geq (\kappa - \varepsilon) d^{p}(\xi,\zeta), \ d^{p}(\xi,\zeta) \geq M(\zeta)\}.
%\]
For each $\zeta \in \Xi$, it follows from the measurability of $\Psi$, $M$, and $d^{p}(\cdot,\zeta)$ that $F(\zeta) \in \scrB(\Xi)$.
Then using the same argument as in~\ref{itm:measurable_selection}, there exists a $\nu$-measurable selection $T : \Xi \mapsto \Xi$ such that $T(\zeta) \in F(\zeta)$ for $\nu$-almost all $\zeta \in \Xi$.
\hfillqed
% \endproof

\subsubsection{Proof of Proposition~\ref{prop:kappa_infty}.}
If $\kappa = \infty$, then for any $n > 0$ it holds that
\[
\phi^{n}(\zeta) \ \ \defi \ \ \inf_{\xi \in \Xi} \{n d^{p}(\xi,\zeta) - \Psi(\xi) + \Psi(\zeta)\} \ \ \notin \ \ L^{1}(\nu).
\]
Observe that $\phi^{n}(\zeta) = \Phi(n,\zeta) + \Psi(\zeta)$.
Hence, for any $n > 0$, there exists $E^{n} \in \scrB_{\nu}$ with $\nu(E^{n}) > 0$, such that $\phi^{n}(\zeta) < -n$ for $\nu$-almost all $\zeta \in E^{n}$.
By Lemma~\ref{lemma:measurable}\ref{itm:measurable_n}, there exists a $\nu$-measurable mapping $T^{n} : E^{n} \mapsto \Xi$ such that
\[
T^{n}(\zeta) \ \ \in \ \ F^{n}(\zeta) \ \ \defi \ \ \{\xi \in \Xi \; : \; \Psi(\xi) - \Psi(\zeta) > n d^{p}(\xi,\zeta) + n\}
\]
for $\nu$-almost all $\zeta \in E^{n}$.
For $m = 1,2,\ldots$, consider the set
\[
E^{n}_{m} \ \ \defi \ \ \left\{\zeta \in E^{n} \; : \; d^{p}(T^{n}(\zeta),\zeta) \leq m\right\}.
\]
Note that $E^{n}_{m} \in \scrB_{\nu}$ for all $m$.
Then $\lim_{m \to \infty} E^{n}_{m} = E^{n}$ and thus $\lim_{m \to \infty} \nu(E^{n}_{m}) = \nu(E^{n}) > 0$.
Hence, for each $n$, there exists a $m^{n}$ such that $\nu(E^{n}_{m^{n}}) > 0$, and $\int_{E^{n}_{m^{n}}} d^{p}(T^{n}(\zeta),\zeta) nu(d\zeta) \leq m^{n} < \infty$.
Let $\overline{T}^{n}$ be the restriction of $T^{n}$ on $E^{n}_{m^{n}}$.
Note that $\overline{T}^{n}$ is also $\nu$-measurable.

For each $n$, let
\[
p^{n} \ \ \defi \ \ \min\left\{1, \ \frac{\theta^{p}}{\int_{E^{n}_{m^{n}}} d^{p}(\overline{T}^{n}(\zeta),\zeta) \nu(d\zeta)}\right\}
\]
and define the distribution
\[
\mu^{n} \ \ \defi \ \ p^{n} \overline{T}^{n}_{\#}\nu + (1 - p^{n}) \nu
\]
Then $\mu^{n}$ is a primal feasible solution, and
\[
\int_{\Xi} \Psi d\mu^{n} - \int_{\Xi} \Psi d\nu \ \ \geq \ \ p^{n} \left(n \int_{E^{n}_{m^{n}}} d^{p}(\overline{T}^{n}(\zeta),\zeta) \nu(d\zeta) + n\right) \ \ \geq \ \ \min\{n \theta^{p}, n\}.
\]
Since $n$ can be chosen arbitrarily large, it follows that $v_{P} = \infty = v_{D}$.
\hfillqed
% \endproof

% \proof{Proof of Proposition~\ref{prop:kappa_infty}.}
% Choose any compact set $E\subset\Xi$ with $\nu(E)>0$. Then for any $M,R>0$, there exists $\zeta_{M}$ such that
% \[
% \Psi(\zeta_{M})-\Psi(\xi)> d^{p}(\zeta_{M},\xi)>R,\quad \forall \xi\in E.
% \]
% Thereby we define a map $T_{M}:\Xi\mapsto\Xi$ by
% \[
% T_{M}(\zeta)\defi\zeta_{M}\mathds{1}_{\{\zeta\in E\}}+\zeta\mathds{1}_{\{\zeta\in\Xi\setminus E\}},\ \forall \zeta\in\Xi,
% \]
% and let
% \[
% \mu_{M}\defi p_{M}{T_{M}}_{\#}\nu+ \Big(1- p_{M}\Big)\nu,
% \]
% where $p_{M}\defi\frac{\theta^{p}}{\int_E d^{p}(\zeta_{M},\zeta)\nu(d\zeta)}$. Note that we can always make $p_{M}\leq 1$ by choosing sufficiently large $R$.
% Then $\mu_{M}$ is primal feasible and
% \[
% v_{P}-\int_{\Xi}\Psi(\xi)\nu(d\xi)\geq p_{M}\int_{E}[\Psi(T_{M}(\zeta))-\Psi(\zeta)]\nu(d\zeta)> p_{M} M\int_E d^{p}(T_{M}(\zeta),\zeta)\nu(d\zeta)= M \theta^{p},
% \]
% which goes to $\infty$ as $M\to\infty$. On the other hand, for any $\lambda\geq0$ and $\zeta\in E$, $\Phi(\lambda,\zeta)=-\infty$, whence $v_{D}=\infty$. Therefore we prove that $v_{P}=v_{D}$.
% \hfillqed
% \begin{figure}[h]
% \centering
% \includegraphics[width=.4\textwidth]{kappa-infty.pdf}
% \caption{Illustration of duality $v_{p}=\infty$ when $\kappa=\infty$.} \label{fig:kappaInfty}
% \end{figure}
% \endproof

\subsubsection{Proof of Lemma~\ref{lemma:phi}}\leavevmode

%\ref{lemma:phi finite}
%By definition, for all $\lambda \ge 0$ and all $\zeta \in \Xi$, it holds that $\Phi(\lambda,\zeta) \leq -\Psi(\zeta)$.
%Next, note that if $\Xi$ is bounded and $\Psi$ is bounded above, then $\kappa = 0$ and $\Phi(\lambda,\zeta) > -\infty$ for all $\lambda \ge 0$ and all $\zeta \in \Xi$, and if $\Xi$ is bounded and $\Psi$ is not bounded above, then $\kappa = \infty$ and $\Phi(\lambda,\zeta) = -\infty$ for all $\lambda \ge 0$ and all $\zeta \in \Xi$.
%Otherwise, if $\Xi$ is not bounded, then it follows from Lemma~\ref{lemma:kappa}\ref{lemma:kappa_expression} that if $\kappa < \infty$, then $\kappa = \kappa^{0}$.  Also, it follows from the proof of Lemma~\ref{lemma:kappa}\ref{lemma:kappa_condition} that $\Phi(\lambda,\zeta) > -\infty$ for all $\lambda > \kappa^{0}$ and all $\zeta \in \Xi$, and it follows from the proof of Lemma~\ref{lemma:kappa}\ref{lemma:kappa_expression} that $\Phi(\lambda,\zeta) = -\infty$ for all $\lambda < \kappa^{0}$ and all $\zeta \in \Xi$.

\ref{lemma:phi_monotone}
For any $\zeta \in \Xi$, $\Phi(\cdot,\zeta)$ is the infimum of nondecreasing functions.
Thus $\Phi(\cdot,\zeta)$ is nondecreasing for all $\zeta \in \Xi$.
Also, for any $\zeta \in \Xi$, $\Phi(\cdot,\zeta)$ is the infimum of continuous functions.
Thus $\Phi(\cdot,\zeta)$ is upper-semi-continuous for all $\zeta \in \Xi$.
Consider any sequence $\{\lambda_{n}\}_{n}$ such that $\lambda_{n} \downarrow \kappa$ as $n \to \infty$.
Since $\int_{\Xi} \Phi(\lambda_{n},\zeta) \nu(d\zeta) > -\infty$, it holds that there is a set $B_{n} \in \scrB_{\nu}(\Xi)$ such that $\nu(B_{n}) = 1$ and $\Phi(\lambda_{n},\zeta) > -\infty$ for all $\zeta \in B_{n}$.
Then it follows from $\Phi(\cdot,\zeta)$ being nondecreasing that $\Phi(\lambda,\zeta) > -\infty$ for all $\lambda \geq \lambda_{n}$ and all $\zeta \in B_{n}$.
Let $B \defi \cap_{n} B_{n}$.
Then $B \in \scrB_{\nu}(\Xi)$, and $\nu(B) = 1$, and $\Phi(\lambda,\zeta) > -\infty$ for all $\lambda > \kappa$ and all $\zeta \in B$.
Since $\Phi(\cdot,\zeta)$ is the infimum of affine functions of $\lambda$, and $\Phi(\lambda,\zeta) < \infty$ for all $\lambda \ge 0$ and all $\zeta \in \Xi$, and $\Phi(\lambda,\zeta) > -\infty$ for all $\lambda > \kappa$ and all $\zeta \in B$, it follows that $\Phi(\cdot,\zeta)$ is concave on $[0,\infty)$ for all $\zeta \in B$.

For the second part, consider any $\varepsilon \ge 0$, any $\lambda_{2} > \lambda_{1}$, and any $\zeta \in \Xi$, such that $\Phi(\lambda_{1},\zeta) > -\infty$.
Consider any $\xi' \in \Xi$ such that $d^{p}(\xi',\zeta) > \sup_{\xi \in \Xi} \left\{d^{p}(\xi,\zeta) \, : \, \lambda_{1} d^{p}(\xi,\zeta) - \Psi(\xi) \, \leq \, \Phi(\lambda_{1}, \zeta) + \varepsilon\right\}$.
If no such $\xi' \in \Xi$ exists, then it follows immediately that $\sup_{\xi \in \Xi} \left\{d^{p}(\xi,\zeta) \, : \, \lambda_{2} d^{p}(\xi,\zeta) - \Psi(\xi) \leq \Phi(\lambda_{2}, \zeta) + \varepsilon\right\} \le \sup_{\xi \in \Xi} \left\{d^{p}(\xi,\zeta) \, : \, \lambda_{1} d^{p}(\xi,\zeta) - \Psi(\xi) \leq \Phi(\lambda_{1}, \zeta) + \varepsilon\right\}$. \\
Otherwise, consider any $\xi''$ such that $\lambda_{1} d^{p}(\xi'',\zeta) - \Psi(\xi'') < \Phi(\lambda_{1}, \zeta) + \min\left\{\varepsilon, \ \left(\lambda_{2} - \lambda_{1}\right) \left[d^{p}(\xi',\zeta) - \sup_{\xi \in \Xi} \left\{d^{p}(\xi,\zeta) \, : \, \lambda_{1} d^{p}(\xi,\zeta) - \Psi(\xi) \, \leq \, \Phi(\lambda_{1}, \zeta) + \varepsilon\right\}\right]\right\}$.
Since $\lambda_{1} d^{p}(\xi'',\zeta) - \Psi(\xi'') < \Phi(\lambda_{1}, \zeta) + \varepsilon$ it follows that $d^{p}(\xi'',\zeta) \le \sup_{\xi \in \Xi} \left\{d^{p}(\xi,\zeta) \, : \, \lambda_{1} d^{p}(\xi,\zeta) - \Psi(\xi) \, \leq \, \Phi(\lambda_{1}, \zeta) + \varepsilon\right\} < d^{p}(\xi',\zeta)$.
Thus, $\lambda_{1} d^{p}(\xi',\zeta) - \Psi(\xi') > \Phi(\lambda_{1}, \zeta) + \varepsilon$ and
\begin{align*}
\lambda_{1} d^{p}(\xi'',\zeta) - \Psi(\xi'') \ \ & < \ \ \Phi(\lambda_{1}, \zeta) \\
& \qquad + \left(\lambda_{2} - \lambda_{1}\right) \left[d^{p}(\xi',\zeta) - \sup_{\xi \in \Xi} \left\{d^{p}(\xi,\zeta) \, : \, \lambda_{1} d^{p}(\xi,\zeta) - \Psi(\xi) \, \leq \, \Phi(\lambda_{1}, \zeta) + \varepsilon\right\}\right] \\
& < \ \ \lambda_{1} d^{p}(\xi',\zeta) - \Psi(\xi') - \varepsilon \\
& \qquad + \left(\lambda_{2} - \lambda_{1}\right) \left[d^{p}(\xi',\zeta) - \sup_{\xi \in \Xi} \left\{d^{p}(\xi,\zeta) \, : \, \lambda_{1} d^{p}(\xi,\zeta) - \Psi(\xi) \, \leq \, \Phi(\lambda_{1}, \zeta) + \varepsilon\right\}\right]
\end{align*}
\begin{align*}
\Rightarrow \ \ \ \lambda_{1} \left[d^{p}(\xi',\zeta) - d^{p}(\xi'',\zeta)\right] \ \ & > \ \ \Psi(\xi') - \Psi(\xi'') + \varepsilon \\
& - \left(\lambda_{2} - \lambda_{1}\right) \left[d^{p}(\xi',\zeta) - \sup_{\xi \in \Xi} \left\{d^{p}(\xi,\zeta) \, : \, \lambda_{1} d^{p}(\xi,\zeta) - \Psi(\xi) \, \leq \, \Phi(\lambda_{1}, \zeta) + \varepsilon\right\}\right] \\
\Rightarrow \ \ \ \lambda_{2} \left[d^{p}(\xi',\zeta) - d^{p}(\xi'',\zeta)\right] \ \ & > \ \ \Psi(\xi') - \Psi(\xi'') + \varepsilon + \left(\lambda_{2} - \lambda_{1}\right) \left[d^{p}(\xi',\zeta) - d^{p}(\xi'',\zeta)\right] \\
& - \left(\lambda_{2} - \lambda_{1}\right) \left[d^{p}(\xi',\zeta) - \sup_{\xi \in \Xi} \left\{d^{p}(\xi,\zeta) \, : \, \lambda_{1} d^{p}(\xi,\zeta) - \Psi(\xi) \, \leq \, \Phi(\lambda_{1}, \zeta) + \varepsilon\right\}\right] \\
& \ge \ \ \Psi(\xi') - \Psi(\xi'') + \varepsilon \\
\Rightarrow \ \ \ \lambda_{2} d^{p}(\xi',\zeta) - \Psi(\xi') \ \ & > \ \ \lambda_{2} d^{p}(\xi'',\zeta) - \Psi(\xi'') + \varepsilon
\ \ \ge \ \ \Phi(\lambda_{2}, \zeta) + \varepsilon.
\end{align*}
That is, every $\xi' \in \Xi$ such that $d^{p}(\xi',\zeta) > \sup_{\xi \in \Xi} \left\{d^{p}(\xi,\zeta) \, : \, \lambda_{1} d^{p}(\xi,\zeta) - \Psi(\xi) \, \leq \, \Phi(\lambda_{1}, \zeta) + \varepsilon\right\}$ satisfies $\lambda_{2} d^{p}(\xi',\zeta) - \Psi(\xi') > \Phi(\lambda_{2}, \zeta) + \varepsilon$, and thus $\sup_{\xi \in \Xi} \left\{d^{p}(\xi,\zeta) \, : \, \lambda_{2} d^{p}(\xi,\zeta) - \Psi(\xi) \leq \Phi(\lambda_{2}, \zeta) + \varepsilon\right\} \le \sup_{\xi \in \Xi} \left\{d^{p}(\xi,\zeta) \, : \, \lambda_{1} d^{p}(\xi,\zeta) - \Psi(\xi) \leq \Phi(\lambda_{1}, \zeta) + \varepsilon\right\}$.

For the third part, consider any $\lambda_{2} > \lambda_{1}$ and any $\zeta \in \Xi$ such that $\Phi(\lambda_{1},\zeta) > -\infty$.
For any $\delta_{i} > 0$, consider any $\xi_{i}^{\delta_{i}} \in \Xi$ such that $\lambda_{i} d^{p}(\xi_{i}^{\delta_{i}},\zeta) - \Psi(\xi_{i}^{\delta_{i}}) \leq \Phi(\lambda_{i},\zeta) + \delta_{i}$ for $i=1,2$.
It follows that
\begin{align*}
\lambda_{2} d^{p}(\xi_{2}^{\delta_{2}},\zeta) - \Psi(\xi_{2}^{\delta_{2}})
\ \ & \leq \ \ \Phi(\lambda_{2},\zeta) + \delta_{2} \\
& \leq \ \ \lambda_{2} d^{p}(\xi_{1}^{\delta_{1}},\zeta) - \Psi(\xi_{1}^{\delta_{1}}) + \delta_{2} \\
& = \ \ (\lambda_{2} - \lambda_{1}) d^{p}(\xi_{1}^{\delta_{1}},\zeta) + \lambda_{1} d^{p}(\xi_{1}^{\delta_{1}},\zeta) - \Psi(\xi_{1}^{\delta_{1}}) + \delta_{2} \\
& \leq \ \ (\lambda_{2} - \lambda_{1}) d^{p}(\xi_{1}^{\delta_{1}},\zeta) + \Phi(\lambda_{1},\zeta) + \delta_{1} + \delta_{2} \\
& \leq \ \ (\lambda_{2} - \lambda_{1}) d^{p}(\xi_{1}^{\delta_{1}},\zeta) + \lambda_{1} d^{p}(\xi_{2}^{\delta_{2}},\zeta) - \Psi(\xi_{2}^{\delta_{2}}) + \delta_{1} + \delta_{2} \\
\Rightarrow \ \ \ d^{p}(\xi_{2}^{\delta_{2}},\zeta) - \frac{\delta_{2}}{\lambda_{2} - \lambda_{1}}
\ \ & \leq \ \ d^{p}(\xi_{1}^{\delta_{1}},\zeta) + \frac{\delta_{1}}{\lambda_{2} - \lambda_{1}}
\end{align*}
\begin{align*}
\Rightarrow \ \ \ & \sup_{\xi \in \Xi} \left\{d^{p}(\xi,\zeta) - \frac{\delta_{2}}{\lambda_{2} - \lambda_{1}} \; : \; \lambda_{2} d^{p}(\xi,\zeta) - \Psi(\xi) \leq \Phi(\lambda_{2},\zeta) + \delta_{2}\right\} \\
& \le \ \ \inf_{\xi \in \Xi} \left\{d^{p}(\xi,\zeta) + \frac{\delta_{1}}{\lambda_{2} - \lambda_{1}} \; : \; \lambda_{1} d^{p}(\xi,\zeta) - \Psi(\xi) \leq \Phi(\lambda_{1},\zeta) + \delta_{1}\right\} \\
\Rightarrow \ \ \ & \limsup_{\delta_{2} \downarrow 0} \left\{\sup_{\xi \in \Xi} \left\{d^{p}(\xi,\zeta) - \frac{\delta_{2}}{\lambda_{2} - \lambda_{1}} \; : \; \lambda_{2} d^{p}(\xi,\zeta) - \Psi(\xi) \leq \Phi(\lambda_{2},\zeta) + \delta_{2}\right\}\right\} \\
& \le \ \ \liminf_{\delta_{1} \downarrow 0} \left\{\inf_{\xi \in \Xi} \left\{d^{p}(\xi,\zeta) + \frac{\delta_{1}}{\lambda_{2} - \lambda_{1}} \; : \; \lambda_{1} d^{p}(\xi,\zeta) - \Psi(\xi) \leq \Phi(\lambda_{1},\zeta) + \delta_{1}\right\}\right\} \\
\Rightarrow \ \ \ & \Dp(\lambda_{2},\zeta) \ \ \leq \ \ \Dm(\lambda_{1},\zeta).
\end{align*}
Also, it follows from the definition of $\Dp$ and $\Dm$ that $\Dm(\lambda_{1},\zeta) \leq \Dp(\lambda_{1},\zeta)$.

\ref{lemma:phi_bound}
It follows from the definition of $\Phi$ that for all $\xi,\zeta \in \Xi$ it holds that
\[
\Psi(\xi) \ \ \leq \ \ \lambda_{1} d^{p}(\xi,\zeta) - \Phi(\lambda_{1},\zeta).
\]
Also, for every $\xi \in \Xi$ that satisfies $\lambda_{2} d^{p}(\xi,\zeta) - \Psi(\xi) \leq \Phi(\lambda_{2},\zeta) + \delta$ for some $\delta \geq 0$, it holds that
\[
\lambda_{2} d^{p}(\xi,\zeta) - \Psi(\xi) - \delta \ \ \leq \ \ -\Psi(\zeta).
\]
Combining the two inequalities above yields that
\begin{eqnarray*}
& & \lambda_{2} d^{p}(\xi,\zeta) + \Psi(\zeta) - \delta \ \ \leq \ \ \lambda_{1} d^{p}(\xi,\zeta) - \Phi(\lambda_{1},\zeta) \\
\Rightarrow \ \ \ & & (\lambda_{2} - \lambda_{1}) d^{p}(\xi,\zeta) - \delta \ \ \leq \ \ -\Psi(\zeta) - \Phi(\lambda_{1},\zeta) \\
\Rightarrow \ \ \ & & \limsup_{\delta \downarrow 0} \left\{\sup_{\xi \in \Xi} \big\{(\lambda_{2} - \lambda_{1}) d^{p}(\xi,\zeta) - \delta \; : \; \lambda_{2} d^{p}(\xi,\zeta) - \Psi(\xi) \leq \Phi(\lambda_{2},\zeta) + \delta\big\}\right\} \ \ \leq \ \ -\Psi(\zeta) - \Phi(\lambda_{1},\zeta) \\
\Rightarrow \ \ \ & & (\lambda_{2} - \lambda_{1}) \Dp(\lambda_{2},\zeta) \ \ \leq \ \ -\Psi(\zeta) - \Phi(\lambda_{1},\zeta)
\end{eqnarray*}

\ref{lemma:phi_derivative}
%As in the proof of~\ref{lemma:phi_monotone}, there is a set $B \in \scrB_{\nu}(\Xi)$ such that $\nu(B) = 1$, and $\Phi(\lambda,\zeta) > -\infty$ for all $\lambda > \kappa$ and all $\zeta \in B$.
Consider any $\zeta \in B$ and any $\lambda_{2} > \lambda_{1} > \kappa$.
For any $\delta > 0$, choose any $\xi_{i}^{\delta} \in \Xi$ such that $\lambda_{i} d^{p}(\xi_{i}^{\delta},\zeta) - \Psi(\xi_{i}^{\delta}) \leq \Phi(\lambda_{i},\zeta) + \delta$ for $i=1,2$.
Then
\[
\Phi(\lambda_{1},\zeta) - \Phi(\lambda_{2},\zeta)
\ \ \leq \ \ \lambda_{1} d^{p}(\xi_{2}^{\delta},\zeta) - \Psi(\xi_{2}^{\delta}) - \Big[\lambda_{2} d^{p}(\xi_{2}^{\delta},\zeta) - \Psi(\xi_{2}^{\delta})\Big] + \delta
\ \ = \ \ (\lambda_{1} - \lambda_{2}) d^{p}(\xi_{2}^{\delta},\zeta) + \delta.
\]
Similarly,
$\Phi(\lambda_{2},\zeta) - \Phi(\lambda_{1},\zeta) \leq (\lambda_{2} - \lambda_{1}) d^{p}(\xi_{1}^{\delta},\zeta) + \delta$.
It follows that
\[
d^{p}(\xi_{2}^{\delta},\zeta) - \frac{\delta}{\lambda_{2} - \lambda_{1}}
\ \ \leq \ \ \frac{\Phi(\lambda_{2},\zeta) - \Phi(\lambda_{1},\zeta)}{\lambda_{2} - \lambda_{1}}
\ \ \leq \ \ d^{p}(\xi_{1}^{\delta},\zeta) + \frac{\delta}{\lambda_{2} - \lambda_{1}}.
\]
Then it follows from the definitions of $\Dp$ and $\Dm$ that
\[
\Dp(\lambda_{2},\zeta) \ \ \leq \ \ \frac{\Phi(\lambda_{2},\zeta) - \Phi(\lambda_{1},\zeta)}{\lambda_{2} - \lambda_{1}} \ \ \leq \ \ \Dm(\lambda_{1},\zeta).
\]
Since $\lambda_{1} \in \interior(\Dom(\Phi(\cdot,\zeta)))$, there is a $\lambda_{0} < \lambda_{1}$ such that $\Phi(\cdot,\zeta)$ is finite-valued and concave on $(\lambda_{0},\infty)$, and the left and right derivatives $\partial \Phi(\lambda,\zeta) / \partial_{\lambda\pm}$ exist for all $\lambda \in (\lambda_{0},\infty)$.
Setting $\lambda_{2} = \lambda$ and letting $\lambda_{1} \uparrow \lambda$, it follows that
\[
\Dp(\lambda,\zeta) \ \ \leq \ \ \frac{\partial \Phi(\lambda,\zeta)}{\partial \lambda-} \ \ \leq \ \ \lim_{\lambda_{1} \uparrow \lambda} \Dm(\lambda_{1},\zeta).
\]
Similarly, setting $\lambda_{1} = \lambda$ and letting $\lambda_{2} \downarrow \lambda$ in the inequality above, it follows that
\[
\lim_{\lambda_{2} \downarrow \lambda} \Dp(\lambda_{2},\zeta) \ \ \leq \ \ \frac{\partial \Phi(\lambda,\zeta)}{\partial \lambda+} \ \ \leq \ \ \Dm(\lambda,\zeta).
\]
% \endproof
\hfillqed

\subsubsection{Proof of Lemma~\ref{lemma:dual objective}.}\leavevmode

\ref{itm:dual objective finite}
It follows from Definition~\ref{def:kappa} of $\kappa$ that $\int_{\Xi} \Phi(\lambda,\zeta) \nu(d\zeta) = -\infty$ for all $\lambda < \kappa$ and $\int_{\Xi} \Phi(\lambda,\zeta) \nu(d\zeta)$ is finite for all $\lambda > \kappa$, and thus $h(\lambda) = \infty$ for all $\lambda < \kappa$ and $h(\lambda)$ is finite for all $\lambda > \kappa$.

\ref{itm:dual objective convex}
It follows from Lemma~\ref{lemma:phi}\ref{lemma:phi_monotone} that $h(\lambda)$ is the sum of a linear function $\lambda \theta^{p}$ and an (extended real-valued) convex function $-\int_{\Xi} \Phi(\lambda,\zeta) \nu(d\zeta)$ on $[0,\infty)$.
Thus $h$ is convex.

\ref{itm:dual objective lower semi-continuous}
Note that \ref{itm:dual objective finite} and \ref{itm:dual objective convex} imply that $h$ is continuous everywhere except possibly at $\kappa$.
We show that $h$ is lower semi-continuous at $\kappa$.
Consider any sequence $\{\lambda_{n}\}_{n}$ such that $\lambda_{n} \downarrow \kappa$.
Since $\Phi(\cdot,\zeta)$ is upper-semi-continuous for all $\zeta \in \Xi$, it follows that $\Phi(\kappa,\zeta) \geq \limsup_{n \to \infty} \Phi(\lambda_{n},\zeta)$.
Also, $\Phi(\lambda,\zeta) \leq -\Psi(\zeta)$ for all $\lambda$ and $\zeta$.
Thus $\liminf_{n \to \infty} [- \Phi(\lambda_{n},\zeta)] - \Psi(\zeta) \geq - \Phi(\kappa,\zeta) - \Psi(\zeta) \geq 0$ for all $\zeta$.
Hence it follows from Fatou's lemma that
\begin{align*}
\liminf_{n \to \infty} h(\lambda_{n}) - \int_{\Xi} \Psi(\zeta) \nu(d\zeta)
\ \ & = \ \ \liminf_{n \to \infty} \left\{\lambda_{n} \theta^{p} + \int_{\Xi} [- \Phi(\lambda_{n},\zeta) - \Psi(\zeta)] \nu(d\zeta)\right\} \\
& \geq \ \ \kappa \theta^{p} + \int_{\Xi} \liminf_{n \to \infty} [- \Phi(\lambda_{n},\zeta) - \Psi(\zeta)] \nu(d\zeta) \\
& \geq \ \ \kappa \theta^{p} + \int_{\Xi} [- \Phi(\kappa,\zeta) - \Psi(\zeta)] \nu(d\zeta) \\
& = \ \ h(\kappa) - \int_{\Xi} \Psi(\zeta) \nu(d\zeta)
\end{align*}
Since $\left|\int_{\Xi} \Psi(\zeta) \nu(d\zeta)\right| < \infty$, it follows that $\liminf_{n \to \infty} h(\lambda_{n}) \geq h(\kappa)$, and thus $h$ is lower semi-continuous.

\ref{itm:dual objective infinite}
Since $\Phi(\lambda,\zeta) \leq -\Psi(\zeta)$, it follows that $h(\lambda) \geq \lambda \theta^{p} + \int_{\Xi} \Psi(\zeta) \nu(d\zeta) \to \infty$ as $\lambda \to \infty$.

\ref{itm:dual minimizer}
The result follows from~\ref{itm:dual objective finite}--\ref{itm:dual objective infinite}.
\hfillqed
% \endproof

\subsubsection{Proof of Corollary~\ref{cor:existence}.} \leavevmode

\ref{itm:iff}
Recall from Lemma~\ref{lemma:phi}\ref{lemma:phi_monotone} that there is a set $B \in \scrB_{\nu}(\Xi)$ such that $\nu(B) = 1$, and $\Phi(\lambda,\zeta) > -\infty$ for all $\lambda > \kappa$ and all $\zeta \in B$.
Note that if $\Psi$ is upper-semi-continuous, and bounded subsets of $(\Xi,d)$ are totally bounded, then for $\delta = 0$, it holds that $\underline{F}(\lambda,\zeta)$ and $\overline{F}(\lambda,\zeta)$ in Lemma~\ref{lemma:measurable}\ref{itm:measurable_selection0} are nonempty for all $\lambda > \kappa$ and all $\zeta \in B$.
Next we show that $\overline{D}_{0}(\cdot,\zeta)$ and $\underline{D}_{0}(\cdot,\zeta)$ are nonincreasing.
Consider any $\lambda_{2} > \lambda_{1}$ and any $\zeta \in \Xi$ such that $\Phi(\lambda_{1},\zeta) > -\infty$.
Consider any $\xi_{i} \in \Xi$ such that $\lambda_{i} d^{p}(\xi_{i},\zeta) - \Psi(\xi_{i}) = \Phi(\lambda_{i},\zeta)$ for $i=1,2$.
Then it follows as in the proof of Lemma~\ref{lemma:phi}\ref{lemma:phi_monotone} that $d^{p}(\xi_{2},\zeta) \leq d^{p}(\xi_{1},\zeta)$.
Therefore $\overline{D}_{0}(\lambda_{2},\zeta) \leq \underline{D}_{0}(\lambda_{1},\zeta) \leq \overline{D}_{0}(\lambda_{1},\zeta)$.

Next we show that, for all $\zeta \in B$, it holds that $\overline{D}_{0}(\cdot,\zeta)$ is upper-semi-continuous and $\underline{D}_{0}(\cdot,\zeta)$ is lower semi-continuous at all $\lambda > \kappa$.
Consider any $\lambda > \kappa$ and any sequence $\{\lambda_{n}\}_{n}$ such that $\lambda_{n} \to \lambda$ as $n \to \infty$ and $\lambda_{n} \in \big((\lambda + \kappa)/2,\lambda + \delta)$ for all $n$, for some $\delta > 0$.
For each $n$ and each $\zeta \in B$, consider any $\xi^{n} \in \argmin_{\xi \in \Xi} \big\{\lambda_{n} d^{p}(\xi,\zeta) - \Psi(\xi)\big\}$.
Note that $d^{p}(\xi^{n},\zeta) \in \big[\overline{D}_{0}(\lambda + \delta,\zeta), \underline{D}_{0}\big((\lambda + \kappa)/2,\zeta\big)\big]$ for all $n$.
Since bounded subsets of $(\Xi,d)$ are totally bounded, it is sufficient to consider subsequences of $\{\xi^{n}\}_{n}$ that converge to some $\xi^{\ast} \in \Xi$.
It follows from the upper-semicontinuity of $\Psi$ and the continuity of $\Phi(\cdot,\zeta)$ at all $\lambda > \kappa$ that
\[
\lambda d^{p}(\xi^{\ast},\zeta) - \Psi(\xi^{\ast}) \ \ \leq \ \ \liminf_{n \to \infty} \{\lambda_{n} d^{p}(\xi^{n},\zeta) - \Psi(\xi^{n})\} \ \ = \ \ \liminf_{n \to \infty} \Phi(\lambda_{n},\zeta) \ \ = \ \ \Phi(\lambda,\zeta)
\]
and thus $\xi^{\ast} \in \argmin_{\xi \in \Xi} \big\{\lambda d^{p}(\xi,\zeta) - \Psi(\xi)\big\}$.
Since $\underline{D}_{0}(\lambda,\zeta) \leq d^{p}(\xi^{\ast},\zeta) = \lim_{n \to \infty} d^{p}(\xi^{n},\zeta) \leq \overline{D}_{0}(\lambda,\zeta)$, it follows that $\overline{D}_{0}(\cdot,\zeta)$ is upper-semi-continuous and $\underline{D}_{0}(\cdot,\zeta)$ is lower semi-continuous at all $\lambda > \kappa$ for all $\zeta \in B$.

Next we show that for all $\lambda > \kappa$ and all $\zeta \in B$, it holds that $\partial \Phi(\lambda,\zeta) / \partial \lambda- = \overline{D}_{0}(\lambda,\zeta)$ and $\partial \Phi(\lambda,\zeta) / \partial \lambda+ = \underline{D}_{0}(\lambda,\zeta)$.
Consider any $\zeta \in B$ and any $\lambda_{2} > \lambda_{1} > \kappa$.
Consider any $\xi^{i} \in \argmin_{\xi \in \Xi} \big\{\lambda_{i} d^{p}(\xi,\zeta) - \Psi(\xi)\big\}$ for $i=1,2$.
Then it follows as in the proof of Lemma~\ref{lemma:phi}\ref{lemma:phi_derivative} that
\[
d^{p}(\xi^{2},\zeta)
\ \ \leq \ \ \frac{\Phi(\lambda_{2},\zeta) - \Phi(\lambda_{1},\zeta)}{\lambda_{2} - \lambda_{1}}
\ \ \leq \ \ d^{p}(\xi^{1},\zeta).
\]
Then it follows from the definitions of $\overline{D}_{0}$ and $\underline{D}_{0}$ that
\[
\overline{D}_{0}(\lambda_{2},\zeta)\ \leq \ \frac{\Phi(\lambda_{2},\zeta) - \Phi(\lambda_{1},\zeta)}{\lambda_{2} - \lambda_{1}}\ \leq \ \underline{D}_{0}(\lambda_{1},\zeta).
\]
Setting $\lambda_{2} = \lambda$ and letting $\lambda_{1} \uparrow \lambda$, it follows from the upper-semicontinuity of $\overline{D}_{0}(\cdot,\zeta)$ that
\[
\overline{D}_{0}(\lambda,\zeta)\ \leq \ \frac{\partial}{\partial \lambda-} \Phi(\lambda,\zeta)\ \leq \ \lim_{\lambda_{1} \uparrow \lambda} \underline{D}_{0}(\lambda_{1},\zeta)\ \leq \ \lim_{\lambda_{1} \uparrow \lambda} \overline{D}_{0}(\lambda_{1},\zeta)\ \leq \ \overline{D}_{0}(\lambda,\zeta)
\]
and hence
\[
\frac{\partial}{\partial \lambda-} \Phi(\lambda,\zeta) \ \ = \ \ \overline{D}_{0}(\lambda,\zeta)
\]
Similarly, setting $\lambda_{1} = \lambda$ and letting $\lambda_{2} \downarrow \lambda$, it follows from the lower-semicontinuity of $\underline{D}_{0}(\cdot,\zeta)$ that
\[
\underline{D}_{0}(\lambda,\zeta)\ \leq \ \lim_{\lambda_{2} \downarrow \lambda} \underline{D}_{0}(\lambda_{2},\zeta)\ \leq \ \lim_{\lambda_{2} \downarrow \lambda} \overline{D}_{0}(\lambda_{2},\zeta)\ \leq \ \frac{\partial}{\partial \lambda+} \Phi(\lambda,\zeta)\ \leq \ \underline{D}_{0}(\lambda,\zeta)
\]
and hence
\[
\frac{\partial}{\partial \lambda+} \Phi(\lambda,\zeta) \ \ = \ \ \underline{D}_{0}(\lambda,\zeta)
\]

Next we show that if condition~\ref{itm:existenceCond1} or~\ref{itm:existenceCond2} or~\ref{itm:existenceCond3} holds, then there exists a primal optimal distribution.
First suppose that condition~\ref{itm:existenceCond1} holds: there exists a dual minimizer $\lambda^{\ast} > \kappa$.
Since for $\delta = 0$, it holds that $\underline{F}(\lambda^{\ast},\zeta)$ and $\overline{F}(\lambda^{\ast},\zeta)$ in Lemma~\ref{lemma:measurable}\ref{itm:measurable_selection0} are nonempty for all $\zeta \in B$, it follows that there exists $\nu$-measurable mappings $\overline{T},\underline{T} : \Xi \mapsto \Xi$ such that
\begin{align*}
\overline{T}(\zeta) \ \ & \in \ \ \Big\{\xi \in \Xi \; : \; \lambda^{\ast} d^{p}(\xi,\zeta) - \Psi(\xi) \, = \, \Phi(\lambda^{\ast},\zeta), \ d^{p}(\xi,\zeta) \, = \, \overline{D}_{0}(\lambda^{\ast},\zeta)\Big\}, \\
\underline{T}(\zeta) \ \ & \in \ \ \Big\{\xi \in \Xi \; : \; \lambda^{\ast} d^{p}(\xi,\zeta) - \Psi(\xi) \, = \, \Phi(\lambda^{\ast},\zeta) , \ d^{p}(\xi,\zeta) \, = \, \underline{D}_{0}(\lambda^{\ast},\zeta)\Big\}
\end{align*}
for $\nu$-almost all $\zeta \in \Xi$.
As in the proof of Lemma~\ref{lem:e-optimal lambda > kappa}, it follows from the first-order optimality conditions $\frac{\partial}{\partial \lambda-} h(\lambda^{\ast}) \leq 0$ and $\frac{\partial}{\partial \lambda+} h(\lambda^{\ast}) \geq 0$ that
\begin{align}
\label{eqn:dual right derivative}
\theta^{p} \ \geq \ & \frac{\partial}{\partial \lambda+} \left(\int_{\Xi} \Phi(\lambda^{\ast},\zeta) \nu(d\zeta)\right)
\ = \ \int_{\Xi} \frac{\partial}{\partial \lambda+} \Phi(\lambda^{\ast},\zeta) \nu(d\zeta)
\ = \ \int_{\Xi} \underline{D}_{0}(\lambda^{\ast},\zeta) \nu(d\zeta)
\ = \ \int_{\Xi} d^{p}(\underline{T}(\zeta),\zeta) \nu(d\zeta) \\
\theta^{p} \ \leq \ & \frac{\partial}{\partial \lambda-} \left(\int_{\Xi} \Phi(\lambda^{\ast},\zeta) \nu(d\zeta)\right)
\ = \ \int_{\Xi} \frac{\partial}{\partial \lambda-} \Phi(\lambda^{\ast},\zeta) \nu(d\zeta)
\ = \ \int_{\Xi} \overline{D}_{0}(\lambda^{\ast},\zeta) \nu(d\zeta)
\ = \ \int_{\Xi} d^{p}(\overline{T}(\zeta),\zeta) \nu(d\zeta). \nonumber
\end{align}
Let $q \in [0,1]$ be such that
\[
q \int_{\Xi} d^{p}(\underline{T}(\zeta),\zeta) \nu(d\zeta) + (1 - q) \int_{\Xi} d^{p}(\overline{T}(\zeta),\zeta) \nu(d\zeta) \ \ = \ \ \theta^{p}.
\]
Let
\begin{equation}
\label{eqn:cor_mu_ast}
\mu^{\ast} \ \ \defi \ \ q \underline{T}_{\#} \nu + (1 - q) \overline{T}_{\#} \nu.
\end{equation}
Then
\[
W_{p}^{p}(\mu^{\ast},\nu) \ \ \leq \ \ q \int_{\Xi} d^{p}(\underline{T}(\zeta),\zeta) \nu(d\zeta) + (1 - q) \int_{\Xi} d^{p}(\overline{T}(\zeta),\zeta) \nu(d\zeta)
\ \ = \ \ \theta^{p}
\]
and thus $\mu^{\ast}$ is primal feasible.
Also,
\begin{align*}
\int_{\Xi} \Psi(\xi) \mu^{\ast}(d\xi)
\ \ & = \ \ q \int_{\Xi} \Psi(\underline{T}(\zeta)) \nu(d\zeta) + (1 - q) \int_{\Xi} \Psi(\overline{T}(\zeta)) \nu(d\zeta) \\
& = \ \ q \int_{\Xi} \big[\lambda^{\ast} d^{p}(\underline{T}(\zeta),\zeta) - \Phi(\lambda^{\ast},\zeta)\big] \nu(d\zeta) + (1 - q) \int_{\Xi} \big[\lambda^{\ast} d^{p}(\overline{T}(\zeta),\zeta) - \Phi(\lambda^{\ast},\zeta)\big] \nu(d\zeta) \\
& = \ \ \lambda^{\ast} \theta^{p} - \int_{\Xi} \Phi(\lambda^{\ast},\zeta) \nu(d\zeta) \ \ = \ \ v_{D}.
\end{align*}
Therefore $\mu^{\ast}$ is primal optimal.

Suppose that condition~\ref{itm:existenceCond2} holds: $\lambda^{\ast} = \kappa > 0$ is the unique dual minimizer, $\nu\big(\{\zeta \in \Xi \, : \, \argmin_{\xi \in \Xi} \{\kappa d^{p}(\xi,\zeta) - \Psi(\xi)\} = \varnothing\}\big) = 0$, and
\[
\int_{\Xi} \underline{D}_{0}(\kappa,\zeta) \nu(d\zeta) \ \ \leq \ \ \theta^{p} \ \ \leq \ \ \int_{\Xi} \overline{D}_{0}(\kappa,\zeta) \nu(d\zeta).
\]
Then it follows in the same way as in the proof for condition~\ref{itm:existenceCond1} that there exists a primal optimal distribution.

Suppose that condition~\ref{itm:existenceCond3} holds: $\lambda^{\ast} = \kappa = 0$ is the unique dual minimizer, $\argmax_{\xi \in \Xi} \{\Psi(\xi)\}$ is nonempty, and
\[
\int_{\Xi} \underline{D}_{0}(0,\zeta) \nu(d\zeta) \ \ \leq \ \ \theta^{p}.
\]
Then, for $\delta = 0$, the sets $\underline{F}(\lambda^{\ast},\zeta)$ in Lemma~\ref{lemma:measurable}\ref{itm:measurable_selection0} are given by
\begin{align*}
\underline{F}(\lambda^{\ast},\zeta) \ \ & = \ \ \Big\{\xi \in \Xi \; : \; - \Psi(\xi) \, = \, \Phi(\lambda^{\ast},\zeta), \ d^{p}(\xi,\zeta) \, \leq \, \underline{D}_{0}(\lambda^{\ast},\zeta)\Big\} \\
& = \ \ \Big\{\xi \in \argmax_{\xi \in \Xi} \{\Psi(\xi)\} \; : \; d^{p}(\xi,\zeta) = \underline{D}_{0}(0,\zeta)\Big\}
\end{align*}
and are non-empty for $\nu$-almost all $\zeta \in \Xi$.
Thus there exists a $\nu$-measurable mapping $\underline{T} : \Xi \mapsto \Xi$ such that $\underline{T}(\zeta) \in \underline{F}(\lambda^{\ast},\zeta)$ for $\nu$-almost all $\zeta \in \Xi$.
Let $\mu^{\ast} \defi \underline{T}_{\#}\nu$.
Then
\[
W_{p}^{p}(\mu^{\ast},\nu) \ \ \leq \ \ \int_{\Xi} d^{p}(\underline{T}(\zeta),\zeta) \nu(d\zeta)
\ \ = \ \ \int_{\Xi} \underline{D}_{0}(0,\zeta) \nu(d\zeta)
\ \ \leq \ \ \theta^{p}
\]
and thus $\mu^{\ast}$ is primal feasible.
Furthermore,
\[
\int_{\Xi} \Psi(\xi) \mu^{\ast}(d\xi)
\ \ = \ \ \int_{\Xi} \Psi(\underline{T}(\zeta)) \nu(d\zeta)
\ \ = \ \ \max_{\xi \in \Xi} \Psi(\xi)
\ \ = \ \ v_{D}
\]
and thus $\mu^{\ast}$ is primal optimal.
Therefore we have shown that if condition~\ref{itm:existenceCond1} or~\ref{itm:existenceCond2} or~\ref{itm:existenceCond3} holds, then there exists a primal optimal distribution.

Next we show that if there exists a primal optimal distribution, then condition~\ref{itm:existenceCond1} or~\ref{itm:existenceCond2} or~\ref{itm:existenceCond3} holds.
Consider any primal feasible distribution $\mu$.
Let $\gamma \in \cP(\Xi \times \Xi)$ denote the corresponding optimal solution in definition~(\ref{eqn:def_wasserstein}) of Wasserstein distance $W_{p}(\mu,\nu)$, and let $\gamma_{\zeta}$ denote the corresponding conditional distribution of $\xi$ given $\zeta$.
Since $\mu$ is feasible, it holds that $\int_{\Xi} \int_{\Xi} d^{p}(\xi,\zeta) \gamma_{\zeta}(d\xi) \nu(d\zeta) \leq \theta^{p}$.
Lemma~\ref{lemma:dual objective}\ref{itm:dual minimizer} established existence of a dual minimizer $\lambda^{\ast} \in [\kappa,\infty)$.
Note that
\begin{align*}
v_{D} - \int_{\Xi} \Psi(\xi) \mu(d\xi)
\ \ & = \ \ \lambda^{\ast} \theta^{p} - \int_{\Xi} \Phi(\lambda^{\ast},\zeta) \nu(d\zeta) \\
& \qquad - \left[\int_{\Xi^2} [\Psi(\xi) - \lambda^{\ast} d^{p}(\xi,\zeta)] \gamma_{\zeta}(d\xi) \nu(d\zeta) + \int_{\Xi^{2}} \lambda^{\ast} d^{p}(\xi,\zeta) \gamma_{\zeta}(d\xi) \nu(d\zeta)\right] \\
& = \ \ \lambda^{\ast} \left[\theta^{p} - \int_{\Xi} \int_{\Xi} d^{p}(\xi,\zeta) \gamma_{\zeta}(d\xi) \nu(d\zeta)\right] \\
& \qquad + \int_{\Xi} \int_{\Xi} [\lambda^{\ast} d^{p}(\xi,\zeta) - \Psi(\xi) - \Phi(\lambda^{\ast},\zeta)] \gamma_{\zeta}(d\xi) \nu(d\zeta)
\end{align*}
For $\mu$ to be primal optimal, it must hold that $v_{D} - \int_{\Xi} \Psi(\xi) \mu(d\xi) = 0$.
Since $\lambda^{\ast} \geq 0$, $\theta^{p} - \int_{\Xi} \int_{\Xi} d^{p}(\xi,\zeta) \gamma_{\zeta}(d\xi) \nu(d\zeta) \geq 0$, and $\lambda^{\ast} d^{p}(\xi,\zeta) - \Psi(\xi) - \Phi(\lambda^{\ast},\zeta) \geq 0$ for all $(\xi,\zeta)$, it follows that all of the following must hold for $\mu$ to be primal optimal:
\begin{enumerate}[label=(\Alph*)]
\item
$\lambda^{\ast} \left[\theta^{p} - \int_{\Xi} \int_{\Xi} d^{p}(\xi,\zeta) \gamma_{\zeta}(d\xi) \nu(d\zeta)\right] = 0$.
\item
\label{itm:everywhere optimal}
$\int_{\Xi} [\lambda^{\ast} d^{p}(\xi,\zeta) - \Psi(\xi) - \Phi(\lambda^{\ast},\zeta)] \gamma_{\zeta}(d\xi) = 0$ for $\nu$-almost all $\zeta$, which in turn implies that $\argmin_{\xi \in \Xi} \{\lambda^{\ast} d^{p}(\xi,\zeta) - \Psi(\xi)\} \neq \varnothing$ for $\nu$-almost all $\zeta$, and the conditional distribution $\gamma_{\zeta}$ should be supported on $\argmin_{\xi \in \Xi} \{\lambda^{\ast} d^{p}(\xi,\zeta) - \Psi(\xi)\}$ for $\nu$-almost all $\zeta$.
\end{enumerate}

Next we show that these conditions imply that condition~\ref{itm:existenceCond1} or~\ref{itm:existenceCond2} or~\ref{itm:existenceCond3} holds.
Since there is a dual minimizer $\lambda^{\ast} \in [\kappa,\infty)$, one of the following conditions must hold:
\begin{enumerate}[label=\arabic*$^\circ$]
\item
There is a dual minimizer $\lambda^{\ast} > \kappa$.
\item
The unique dual minimizer satisfies $\lambda^{\ast} = \kappa > 0$.
\item
The unique dual minimizer satisfies $\lambda^{\ast} = \kappa = 0$.
\end{enumerate}
If~1$^\circ$ holds, then condition~\ref{itm:existenceCond1} holds, and the proof is complete. \\
Next suppose that~2$^\circ$ holds, and that $\mu$ is a primal optimal solution.
Condition \ref{itm:everywhere optimal} implies that $\nu\big(\{\zeta \in \Xi \, : \, \argmin_{\xi \in \Xi} \{\kappa d^{p}(\xi,\zeta) - \Psi(\xi)\} = \varnothing\}\big) = 0$.
Next we show that
\[
\int_{\Xi} \underline{D}_{0}(\kappa,\zeta) \nu(d\zeta) \ \ \leq \ \ \theta^{p} \ \ \leq \ \ \int_{\Xi} \overline{D}_{0}(\kappa,\zeta) \nu(d\zeta).
\]
It follows from~(A) that
\[
\theta^{p} \ \ = \ \ \int_{\Xi} \int_{\Xi} d^{p}(\xi,\zeta) \gamma_{\zeta}(d\xi) \nu(d\zeta) \ \ \leq \ \ \int_{\Xi} \overline{D}_{0}(\kappa,\zeta) \nu(d\zeta)
\]
If $\theta^{p} < \int_{\Xi} \underline{D}_{0}(\kappa,\zeta) \nu(d\zeta)$, then it follows as in~(\ref{eqn:dual right derivative}) that
\[
\frac{\partial}{\partial \lambda+} \int_{\Xi} \Phi(\kappa,\zeta) \nu(d\zeta) \ \ = \ \ \int_{\Xi} \underline{D}_{0}(\kappa,\zeta) \nu(d\zeta) \ \ > \ \ \theta^{p}.
\]
Then there exists a $\lambda > \kappa$ such that $\lambda \theta^{p} - \int_{\Xi} \Phi(\lambda,\zeta) \nu(d\zeta) < \kappa \theta^{p} - \int_{\Xi} \Phi(\kappa,\zeta) \nu(d\zeta)$, contradicting $\lambda^{\ast} = \kappa$ being a dual minimizer.
Therefore, if~2$^\circ$ holds, then condition~\ref{itm:existenceCond2} holds. \\
Next suppose that~3$^\circ$ holds.
Condition~\ref{itm:everywhere optimal} implies that $\argmax_{\xi \in \Xi} \{\Psi(\xi)\} \neq \varnothing$.
Suppose that $\mu$ is a primal optimal solution.
Then
\[
\int_{\Xi} \underline{D}_{0}(0,\zeta) \nu(d\zeta) \ \ \leq \ \ \int_{\Xi} \int_{\Xi} d^{p}(\xi,\zeta) \gamma_{\zeta}(d\xi) \nu(d\zeta) \ \ \leq \ \ \theta^{p}.
\]
Therefore, if~3$^\circ$ holds, then condition~\ref{itm:existenceCond3} holds.
% Let $\gamma_{0}$ be such that $\int_{\Xi^{2}}d(\xi,\zeta)\gamma_{0}(d\xi,d\zeta)=W_{p}(\mu,\nu)$, and let $\gamma_{\zeta}$ be the conditional distribution of $\xi$ given $\zeta$ under the joint distribution $\gamma_{0}$. Then
% \[
% \max_{\xi\in\Xi}\Psi(\xi)=\int_{\Xi} \Psi(\xi)\mu(d\xi) = \int_{\Xi^{2}}\Psi(\xi)\gamma_{0}(d\zeta,d\zeta)=\int_{\Xi}\int_{\Xi} \Psi(\xi)\gamma_{\zeta}(d\xi)\nu(d\zeta).
% \]
% It follows that
% \[
% \int_{\Xi}\Big[\int_{\Xi} \Psi(\xi)\gamma_{\zeta}(d\xi)- \max_{\xi\in\Xi}\Psi(\xi)\Big]\nu(d\zeta)=0.
% \]
% Then for $\nu$-almost every $\zeta$, the conditional distribution of $\xi$ should be supported on the set $\{\xi:\Psi(\xi)=\max_{\xi'\in\Xi}\Psi(\xi)\}$.
\vspace{0.5em}

%\ref{itm:small}
If $-\Psi(\zeta) \leq \inf_{\xi \in \Xi} \left\{\kappa d^{p}(\xi,\zeta) - \Psi(\xi)\right\}$ for $\nu$-almost $\zeta$, then $- \Psi(\zeta) \leq \kappa d^{p}(\xi,\zeta) - \Psi(\xi) \leq \lambda d^{p}(\xi,\zeta) - \Psi(\xi)$ for $\nu$-almost $\zeta$, for all $\xi \in \Xi$, and for all $\lambda \geq \kappa$.
Then $\Phi(\lambda,\zeta) = -\Psi(\zeta)$ and $h(\lambda) = \lambda \theta^{p} + \int_{\Xi} \Psi(\zeta) \nu(d\zeta)$ for all $\lambda \geq \kappa$, and hence the dual optimal solution $\lambda^{\ast} = \kappa$ for all $\theta > 0$.

Otherwise there exists a set $E \subset \Xi$ such that $\nu(E) > 0$ and $-\Psi(\zeta) > \Phi(\kappa,\zeta)$ for all $\zeta \in E$, and thus $- \int_{\Xi} \Psi(\zeta) \nu(d\zeta) > \int_{\Xi} \Phi(\kappa,\zeta) \nu(d\zeta)$.
Then by continuity (follows from concavity) of $\int_{\Xi} \Phi(\cdot,\zeta) \nu(d\zeta)$ on $[\kappa,\infty)$, there exists a $\lambda_{0} > \kappa$ such that $- \int_{\Xi} \Psi(\zeta) \nu(d\zeta) > \int_{\Xi} \Phi(\lambda_{0},\zeta) \nu(d\zeta)$.
Using the assumptions that $\Psi$ is upper-semi-continuous and that bounded subsets of $(\Xi,d)$ are totally bounded, as well as Lemma \ref{lemma:measurable}, it follows that there exists a $\nu$-measurable map $T(\lambda_{0},\cdot) : \Xi \mapsto \Xi$ such that $\lambda_{0} d^{p}(T(\lambda_{0},\zeta),\zeta) - \Psi(T(\lambda_{0},\zeta)) = \Phi(\lambda_{0},\zeta)$, and
\[
\varepsilon \ \ \defi \ \ \int_{\Xi} d^{p}(T(\lambda_{0},\zeta),\zeta) \nu(d\zeta) \ \ > \ \ 0
\]
since otherwise $- \int_{\Xi} \Psi(\zeta) \nu(d\zeta) = \int_{\Xi} \Phi(\lambda_{0},\zeta) \nu(d\zeta)$.
Then, since $\Phi(\kappa,\zeta) \leq \kappa d^{p}(T(\lambda_{0},\zeta),\zeta) - \Psi(T(\lambda_{0},\zeta))$, it follows that for any $\theta < \varepsilon^{1/p}$ it holds that
\begin{align*}
h(\kappa) \ \ & = \ \ \kappa \theta^{p} - \int_{\Xi} \Phi(\kappa,\zeta) \nu(d\zeta) \\
& \geq \ \ \kappa \theta^{p} - \int_{\Xi} \kappa d^{p}(T(\lambda_{0},\zeta),\zeta) \nu(d\zeta) + \int_{X} \Psi(T(\lambda_{0},\zeta)) \nu(d\zeta) \\
& = \ \ \kappa \theta^{p} - \kappa \varepsilon + \int_{\Xi} \lambda_{0} d^{p}(T(\lambda_{0},\zeta),\zeta) \nu(d\zeta) - \int_{\Xi} \Phi(\lambda_{0},\zeta) \nu(d\zeta)\\
& = \ \ \kappa \theta^{p} + (\lambda_{0} - \kappa) \varepsilon - \int_{\Xi} \Phi(\lambda_{0},\zeta) \nu(d\zeta)\\
& > \ \ \kappa \theta^{p} + (\lambda_{0} - \kappa) \theta^{p} - \int_{\Xi} \Phi(\lambda_{0},\zeta) \nu(d\zeta)
\ \ = \ \ h(\lambda_{0})
\end{align*}
and therefore $\lambda = \kappa$ cannot be dual optimal for $\theta < \varepsilon^{1/p}$.
\vspace{0.5em}

\ref{itm:form}
It follows from the proof of \ref{itm:iff} that whenever there is a worst-case distribution, then condition~\ref{itm:existenceCond1} or~\ref{itm:existenceCond2} or~\ref{itm:existenceCond3} holds.
As shown in the proof of \ref{itm:iff}, in each case there are maps $\overline{T}^{\ast}, \underline{T}^{\ast} : \Xi \mapsto \Xi$ such that \eqref{eqn:transportMap} holds, and there is a $p^{\ast} \in [0,1]$ such that
$\mu^{\ast} \defi p^{\ast} \overline{T}^{\ast}_{\#} \nu + (1 - p^{\ast}) \underline{T}^{\ast}_{\#} \nu$ is a worst-case distribution.

Next, consider $\gamma^{T} \in \cP(\Xi \times \Xi)$.
Note that for any measurable set $B \subset \Xi$ it holds that
\[
\pi^{1}_{\#}\gamma^{T}(B) \ \ = \ \ \gamma^{T}(B \times \Xi)
\ \ = \ \ p^{\ast} \nu(B) + \left(1 - p^{\ast}\right) \nu(B)
\ \ = \ \ \nu(B)
\]
and
\begin{align*}
\pi^{2}_{\#}\gamma^{T}(B) \ \ & = \ \ \gamma^{T}(\Xi \times B)
\ \ = \ \ p^{\ast} \nu\left(\left\{\zeta \; : \; \overline{T}^{\ast}(\zeta) \in B\right\}\right) + \left(1 - p^{\ast}\right) \nu\left(\left\{\zeta \; : \; \underline{T}^{\ast}(\zeta) \in B\right\}\right) \\
& = \ \ p^{\ast} \overline{T}^{\ast}_{\#}\nu(B) + \left(1 - p^{\ast}\right) \underline{T}^{\ast}_{\#}\nu(B)
\ \ = \ \ \mu^{\ast}(B)
\end{align*}
Thus $\gamma^{T}$ is a feasible joint distribution in the definition~(\ref{eqn:def_wasserstein}) of $W_{p}^{p}(\mu^{\ast},\nu)$.
Next, consider any $\gamma \in \cP(\Xi \times \Xi)$ such that $\pi^{1}_{\#}\gamma = \nu$ and $\pi^{2}_{\#}\gamma = \mu^{\ast}$, and let $\gamma_{\zeta}$ denote the corresponding conditional distribution of $\xi$ given $\zeta$.
Then
\begin{align*}
& \int_{\Xi \times \Xi} \left(\lambda^{\ast} d^{p}(\xi,\zeta) - \Psi(\xi)\right) \gamma(d\xi,d\zeta) \ \ = \ \ \int_{\Xi} \int_{\Xi} \left(\lambda^{\ast} d^{p}(\xi,\zeta) - \Psi(\xi)\right) \gamma_{\zeta}(d\xi) \nu(d\zeta) \\
& \ \ \ge \ \ \int_{\Xi} \min_{\xi \in \Xi} \left\{\lambda^{\ast} d^{p}(\xi,\zeta) - \Psi(\xi)\right\} \nu(d\zeta) \\
& \ \ = \ \ p^{\ast} \int_{\Xi} \left(\lambda^{\ast} d^{p}\left(\overline{T}^{\ast}(\zeta),\zeta\right) - \Psi\left(\overline{T}^{\ast}(\zeta)\right)\right\} \nu(d\zeta) + \left(1 - p^{\ast}\right) \int_{\Xi} \left(\lambda^{\ast} d^{p}\left(\underline{T}^{\ast}(\zeta),\zeta\right) - \Psi\left(\underline{T}^{\ast}(\zeta)\right)\right\} \nu(d\zeta) \\
& \ \ = \ \ \int_{\Xi \times \Xi} \left(\lambda^{\ast} d^{p}(\xi,\zeta) - \Psi(\xi)\right) \gamma^{T}(d\xi,d\zeta) \\
\Rightarrow \ \ \ & \int_{\Xi \times \Xi} \lambda^{\ast} d^{p}(\xi,\zeta) \gamma(d\xi,d\zeta) - \int_{\Xi \times \Xi} \Psi(\xi) \gamma(d\xi,d\zeta)
\ \ \ge \ \ \int_{\Xi \times \Xi} \lambda^{\ast} d^{p}(\xi,\zeta) \gamma^{T}(d\xi,d\zeta) - \int_{\Xi \times \Xi} \Psi(\xi) \gamma^{T}(d\xi,d\zeta) \\
\Rightarrow \ \ \ & \int_{\Xi \times \Xi} \lambda^{\ast} d^{p}(\xi,\zeta) \gamma(d\xi,d\zeta) - \int_{\Xi } \Psi(\xi) \mu^{\ast}(d\xi)
\ \ \ge \ \ \int_{\Xi \times \Xi} \lambda^{\ast} d^{p}(\xi,\zeta) \gamma^{T}(d\xi,d\zeta) - \int_{\Xi} \Psi(\xi) \mu^{\ast}(d\xi) \\
\Rightarrow \ \ \ & \int_{\Xi \times \Xi} d^{p}(\xi,\zeta) \gamma(d\xi,d\zeta)
\ \ \ge \ \ \int_{\Xi \times \Xi} d^{p}(\xi,\zeta) \gamma^{T}(d\xi,d\zeta)
\end{align*}
Therefore, $\gamma^{T}$ is an optimal joint distribution in the definition~(\ref{eqn:def_wasserstein}) of $W_{p}^{p}(\mu^{\ast},\nu)$.

% Let $\gamma^{\ast} \in \cP(\Xi \times \Xi)$ denote the corresponding optimal solution in definition~(\ref{eqn:def_wasserstein}) of Wasserstein distance $W_{p}(\mu^{\ast},\nu)$, and let $\gamma^{\ast}_{\zeta} = p^{\ast} \delta_{\overline{T}^{\ast}(\zeta)} + (1 - p^{\ast}) \delta_{\underline{T}^{\ast}(\zeta)}$ denote the corresponding conditional distribution of $\xi$ given $\zeta$.
% It follows from condition~\ref{itm:everywhere optimal} in the proof of \ref{itm:iff} that $0 = \int_{\Xi} [\lambda^{\ast} d^{p}(\xi,\zeta) - \Psi(\xi) - \Phi(\lambda^{\ast},\zeta)] \gamma^{\ast}_{\zeta}(d\xi) = p^{\ast} [\lambda^{\ast} d^{p}(\overline{T}^{\ast}(\zeta),\zeta) - \Psi(\overline{T}^{\ast}(\zeta)) - \Phi(\lambda^{\ast},\zeta)] + (1 - p^{\ast}) [\lambda^{\ast} d^{p}(\underline{T}^{\ast}(\zeta),\zeta) - \Psi(\underline{T}^{\ast}(\zeta)) - \Phi(\lambda^{\ast},\zeta)]$ for $\nu$-almost all $\zeta$.
% Thus, if $p^* > 0$, then $\lambda^{\ast} d^{p}(\overline{T}^{\ast}(\zeta),\zeta) - \Psi(\overline{T}^{\ast}(\zeta)) = \Phi(\lambda^{\ast},\zeta)$, that is, $\overline{T}^{\ast}(\zeta) \in \argmin_{\xi \in \Xi}\{\lambda^{\ast} d^{p}(\xi,\zeta) - \Psi(\xi)\}$, and if $p^* < 1$, then $\lambda^{\ast} d^{p}(\underline{T}^{\ast}(\zeta),\zeta) - \Psi(\underline{T}^{\ast}(\zeta)) = \Phi(\lambda^{\ast},\zeta)$, that is, $\underline{T}^{\ast}(\zeta) \in \argmin_{\xi \in \Xi}\{\lambda^{\ast} d^{p}(\xi,\zeta) - \Psi(\xi)\}$.
\vspace{0.5em}

\ref{itm:strongDual_concave}
% For any primal feasible solution $\mu$, let $\gamma^{\mu}$ be a minimizer in the definition \eqref{eqn:def_wasserstein} of $W_{p}(\mu,\nu)$, and for any $\zeta \in \Xi$, let $\gamma_{\zeta}^{\mu}$ be the conditional distribution of $\xi$ given $\zeta$ when the joint distribution of $(\xi,\zeta)$ is $\gamma^{\mu}$.
% Consider $T^{\mu} : \Xi \mapsto \Xi$ given by
% \[
% T^{\mu}(\zeta) \ \ \defi \ \ \E_{\gamma_{\zeta}}[\xi].
% \]
% Then it follows from $\Xi$ being convex and complete that $T^{\mu}(\zeta) \in \Xi$ for all $\zeta \in \Xi$.
% It follows from $d^{p}(\cdot,\zeta)$ being convex for all $\zeta \in \Xi$ and Jensen's inequality that
% \[
% W_{p}^{p}(T^{\mu}_{\#}\nu,\nu) \ \ \leq \ \ \int_{\Xi} d^{p}(\E_{\gamma_\zeta}[\xi],\zeta) \nu(d\zeta)
% \ \ \leq \ \ \int_{\Xi} \E_{\gamma_{\zeta}}[d^{p}(\xi,\zeta)] \nu(d\zeta)
% \ \ = \ \ W_{p}^{p}(\mu,\nu) \ \ \leq \ \ \theta^{p}.
% \]
% Also, it follows from $\Psi$ being concave that
% \[
% \int_{\Xi} \Psi(T^{\mu}(\zeta)) \nu(d\zeta)
% \ \ = \ \ \int_{\Xi} \Psi(\E_{\gamma_{\zeta}}[\xi]) \nu(d\zeta)
% \ \ \ge \ \ \int_{\Xi} \E_{\gamma_{\zeta}}[\Psi(\xi)] \nu(d\zeta)
% \ \ = \ \ \E_{\mu}[\Psi].
% \]
% Hence $T^{\mu}_{\#}\nu$ is a feasible solution with objective value no worse than $\mu$, and thus the result follows.
% \hfillqed
Observe that
\[
\left\{T_{\#}\nu \; : \; W_{p}(T_{\#}\nu,\nu) \leq \theta,\; T : \Xi \mapsto \Xi \textrm{ is } \nu\textrm{-measurable}\right\} \ \ \subset \ \ \left\{\mu \in \cP(\Xi) \; : \; W_p(\mu,\nu) \leq \theta\right\}
\]
and thus
\[
v_{P} \ \ \geq \ \ \sup\left\{\E_{T_{\#}\nu}[\Psi(\xi)] \; : \; W_{p}(T_{\#}\nu,\nu) \leq \theta,\; T : \Xi \mapsto \Xi \textrm{ is } \nu\textrm{-measurable}\right\}.
\]
To show that equality holds, we consider the 2~cases in Lemmas~\ref{lem:e-optimal lambda > kappa} and~\ref{lem:e-optimal lambda = kappa}:
\begin{itemize}
\item
Case~1: $h$ has a minimizer $\lambda^{\ast} > \kappa$.
\end{itemize}
Consider any $\delta,\varepsilon > 0$ and any $\lambda_{1},\lambda_{2}$ such that $\kappa < \lambda_{1} < \lambda^{\ast} < \lambda_{2}$.
Let $q^{\varepsilon}_{\delta}(\lambda_{1},\lambda_{2}) \in [0,1]$, $q^{\delta} = \theta^{p} / (\theta^{p} + \delta)$, $\overline{T}^{\varepsilon}_{\delta}(\lambda_{1},\cdot), \underline{T}^{\varepsilon}_{\delta}(\lambda_{2},\cdot) : \Xi \mapsto \Xi$ be as in the proof of Lemma~\ref{lem:e-optimal lambda > kappa}.
Let $T^{\varepsilon}_{\delta}(\lambda_{1},\lambda_{2},\cdot) : \Xi \mapsto \Xi$ be given by
\[
T^{\varepsilon}_{\delta}(\lambda_{1},\lambda_{2},\zeta) \ \ \defi \ \ q^{\delta} q^{\varepsilon}_{\delta}(\lambda_{1},\lambda_{2}) \overline{T}^{\varepsilon}_{\delta}(\lambda_{1},\zeta) + q^{\delta} \big(1 - q^{\varepsilon}_{\delta}(\lambda_{1},\lambda_{2})\big) \underline{T}^{\varepsilon}_{\delta}(\lambda_{2},\zeta) + (1 - q^{\delta}) \zeta.
\]
Note that $T^{\varepsilon}_{\delta}(\lambda_{1},\lambda_{2},\zeta) \in \Xi$ for all $\zeta \in \Xi$ and that $T^{\varepsilon}_{\delta}(\lambda_{1},\lambda_{2},\cdot)$ is $\nu$-measurable.
It follows from $d^{p}(\cdot,\zeta)$ being convex that
\begin{align*}
W_{p}^{p}(T^{\varepsilon}_{\delta}(\lambda_{1},\lambda_{2},\cdot)_{\#}\nu,\nu) \ \ & \leq \ \ \int_{\Xi} d^{p}\left(q^{\delta} q^{\varepsilon}_{\delta}(\lambda_{1},\lambda_{2}) \overline{T}^{\varepsilon}_{\delta}(\lambda_{1},\zeta) + q^{\delta} \big(1 - q^{\varepsilon}_{\delta}(\lambda_{1},\lambda_{2})\big) \underline{T}^{\varepsilon}_{\delta}(\lambda_{2},\zeta) + (1 - q^{\delta}) \zeta, \ \zeta\right) \nu(d\zeta) \\
& \leq \ \
q^{\delta} q^{\varepsilon}_{\delta}(\lambda_{1},\lambda_{2}) \int_{\Xi} d^{p}(\overline{T}^{\varepsilon}_{\delta}(\lambda_{1},\zeta),\zeta) \nu(d\zeta) + q^{\delta} \big(1 - q^{\varepsilon}_{\delta}(\lambda_{1},\lambda_{2})\big) \int_{\Xi} d^{p}(\underline{T}^{\varepsilon}_{\delta}(\lambda_{2},\zeta),\zeta) \nu(d\zeta) \\
& = \ \ q^{\delta} \bigg(\theta^{p} + \big[1 - 2 q^{\varepsilon}_{\delta}(\lambda_{1},\lambda_{2})\big] \delta\bigg)
\ \ \leq \ \ \theta^{p}
\end{align*}
Also, it follows from $\Psi$ being concave that
\begin{align}
\E_{T^{\varepsilon}_{\delta}(\lambda_{1},\lambda_{2},\cdot)_{\#}\nu}[\Psi(\xi)] \ \ & = \ \ \int_{\Xi} \Psi\left(q^{\delta} q^{\varepsilon}_{\delta}(\lambda_{1},\lambda_{2}) \overline{T}^{\varepsilon}_{\delta}(\lambda_{1},\zeta) + q^{\delta} \big(1 - q^{\varepsilon}_{\delta}(\lambda_{1},\lambda_{2})\big) \underline{T}^{\varepsilon}_{\delta}(\lambda_{2},\zeta) + (1 - q^{\delta}) \zeta\right) \nu(d\zeta) \nonumber \\
& \geq \ \ q^{\delta} q^{\varepsilon}_{\delta}(\lambda_{1},\lambda_{2}) \int_{\Xi} \Psi(\overline{T}^{\varepsilon}_{\delta}(\lambda_{1},\zeta)) \nu(d\zeta) + q^{\delta} \big(1 - q^{\varepsilon}_{\delta}(\lambda_{1},\lambda_{2})\big) \int_{\Xi} \Psi(\underline{T}^{\varepsilon}_{\delta}(\lambda_{2},\zeta)) \nu(d\zeta) \nonumber \\
& \qquad + (1 - q^{\delta}) \int_{\Xi} \Psi(\zeta) \nu(d\zeta) \nonumber \\
& \geq \ \ q^{\delta} \lambda_{1} \bigg[\theta^{p} + \big(1 - 2 q^{\varepsilon}_{\delta}(\lambda_{1},\lambda_{2})\big) \delta\bigg] - q^{\delta} q^{\varepsilon}_{\delta}(\lambda_{1},\lambda_{2}) \int_{\Xi} \Phi(\lambda_{1},\zeta) \nu(d\zeta) \nonumber \\
& \qquad - q^{\delta} \big(1 - q^{\varepsilon}_{\delta}(\lambda_{1},\lambda_{2})\big) \int_{\Xi} \Phi(\lambda_{2},\zeta) \nu(d\zeta) - q^{\delta} \varepsilon + (1 - q^{\delta}) \int_{\Xi} \Psi(\zeta) \nu(d\zeta).
\label{eqn:e-optimal perturbation 1}
\end{align}
Thus, given any $\varepsilon > 0$, choose $\lambda_{1}^{\varepsilon} \in \left(\max\{\kappa, \lambda^{\ast} - \varepsilon\}, \, \lambda^{\ast}\right)$ such that $\int_{\Xi} \Phi(\lambda_{1}^{\varepsilon},\zeta) \nu(d\zeta) \le \int_{\Xi} \Phi(\lambda^{\ast},\zeta) \nu(d\zeta) + \varepsilon$ and $\left(\lambda^{\ast} - \lambda_{1}^{\varepsilon}\right) \theta^{p} \le \varepsilon$, choose $\lambda_{2}^{\varepsilon} \in \left(\lambda^{\ast}, \, \lambda^{\ast} + \varepsilon\right)$ such that $\int_{\Xi} \Phi(\lambda_{2}^{\varepsilon},\zeta) \nu(d\zeta) \le \int_{\Xi} \Phi(\lambda^{\ast},\zeta) \nu(d\zeta) + \varepsilon$, and choose $\delta^{\varepsilon} \in (0,\varepsilon)$ such that $2 \lambda_{1}^{\varepsilon} \delta^{\varepsilon} \le \varepsilon$, $- \delta^{\varepsilon} \int_{\Xi} \Phi(\lambda^{\ast},\zeta) \nu(d\zeta) \le \theta^{p} \varepsilon$, and $- \delta^{\varepsilon} \int_{\Xi} \Psi(\zeta) \nu(d\zeta) \le \theta^{p} \varepsilon$.
Set $T^{\varepsilon} \defi T^{\varepsilon}_{\delta^{\varepsilon}}(\lambda_{1}^{\varepsilon},\lambda_{2}^{\varepsilon},\cdot)$.
Then it follows from~\eqref{eqn:e-optimal perturbation 1} that
\[
\E_{T^{\varepsilon}_{\#}\nu}[\Psi(\xi)]
\ \ \geq \ \  \lambda^{\ast} \theta^{p} - \int_{\Xi} \Phi(\lambda^{\ast},\zeta) \nu(d\zeta) - 6 \varepsilon
\ \ = \ \ v_{D} - 6 \varepsilon
\]
Since $\varepsilon > 0$ can be arbitrarily small, it follows that
\[
\sup\left\{\E_{T_{\#}\nu}[\Psi(\xi)] \; : \; W_{p}(T_{\#}\nu,\nu) \leq \theta,\; T : \Xi \mapsto \Xi \textrm{ is } \nu\textrm{-measurable}\right\} \ \ = \ \ v_{P}
\]

\begin{itemize}
\item
Case~2: $\kappa$ is the unique minimizer of $h$.
\end{itemize}
Consider any $\lambda > \kappa$.
Consider any $\varepsilon \in \left(0, (\lambda - \kappa) \theta^{p} - \int_{\Xi} \left[\Phi(\lambda,\zeta) - \Phi(\kappa,\zeta)\right] \nu(d\zeta)\right)$.
Let $\underline{T}_{\varepsilon}(\lambda,\cdot) : \Xi \rightarrow \Xi$ be as in the proof of Lemma~\ref{lem:e-optimal lambda = kappa}.
Then
\[
W_{p}^{p}(\underline{T}_{\varepsilon}(\lambda,\cdot)_{\#}\nu,\nu) \ \ \leq \ \ \int_{\Xi} d^{p}(\underline{T}_{\varepsilon}(\lambda,\zeta),\zeta) \nu(d\zeta)
\ \ < \ \ \theta^{p}
\]
\begin{itemize}
\item
Case~2.1: $\kappa = 0$.
\end{itemize}
Then
\begin{align*}
& \E_{\underline{T}_{\varepsilon}(\lambda,\cdot)_{\#}\nu}[\Psi(\xi)]
\ \ = \ \ \int_{\Xi} \Psi(\xi) \underline{T}_{\varepsilon}(\lambda,\cdot)_{\#}\nu(d\xi)
\ \ \geq \ \ - \int_{\Xi} \Phi(\lambda,\zeta) \nu(d\zeta) - \varepsilon \\
& \Rightarrow \ \ \ \lim_{\lambda \downarrow 0} \E_{\underline{T}_{\varepsilon}(\lambda,\cdot)_{\#}\nu}[\Psi(\xi)]
\ \ \geq \ \ \lim_{\lambda \downarrow 0} \left\{- \int_{\Xi} \Phi(\lambda,\zeta) \nu(d\zeta) - \varepsilon\right\}
\ \ = \ \ v_{D} - \varepsilon \\
& \Rightarrow \ \ \ \sup\left\{\E_{T_{\#}\nu}[\Psi(\xi)] \; : \; W_{p}(T_{\#}\nu,\nu) \leq \theta,\; T : \Xi \mapsto \Xi \textrm{ is } \nu\textrm{-measurable}\right\} \ \ = \ \ v_{P}
\end{align*}
\begin{itemize}
\item
Case~2.2: $\kappa > 0$.
\end{itemize}
Consider any $\lambda > \kappa$, $\kappa' \in (0,\kappa)$, and $R > \theta^{p}$.
Let $\overline{T}^{R}(\kappa',\cdot) : \Xi \mapsto \Xi$ and $q_{\varepsilon}^{R}(\kappa',\lambda) \in (0,1)$ be as in the proof of Lemma~\ref{lem:e-optimal lambda = kappa}.
Let $T_{\varepsilon}^{R}(\kappa',\lambda,\cdot) : \Xi \mapsto \Xi$ be given by
\[
T_{\varepsilon}^{R}(\kappa',\lambda,\zeta) \ \ \defi \ \ q_{\varepsilon}^{R}(\kappa',\lambda) \underline{T}_{\varepsilon}(\lambda,\zeta) + \left(1 - q_{\varepsilon}^{R}(\kappa',\lambda)\right) \overline{T}^{R}(\kappa',\zeta)
\]
Note that $T_{\varepsilon}^{R}(\kappa',\lambda,\zeta) \in \Xi$ for all $\zeta \in \Xi$ and that $T_{\varepsilon}^{R}(\kappa',\lambda,\cdot)$ is $\nu$-measurable.
It follows from $d^{p}(\cdot,\zeta)$ being convex that
\begin{align*}
W_{p}^{p}(T_{\varepsilon}^{R}(\kappa',\lambda,\cdot)_{\#}\nu,\nu) \ \ & \leq \ \ \int_{\Xi} d^{p}\left(q_{\varepsilon}^{R}(\kappa',\lambda) \underline{T}_{\varepsilon}(\lambda,\zeta) + \left(1 - q_{\varepsilon}^{R}(\kappa',\lambda)\right) \overline{T}^{R}(\kappa',\zeta), \ \zeta\right) \nu(d\zeta) \\
& \leq \ \ q_{\varepsilon}^{R}(\kappa',\lambda) \int_{\Xi} d^{p}(\underline{T}_{\varepsilon}(\lambda,\zeta),\zeta) \nu(d\zeta) + \left(1 - q_{\varepsilon}^{R}(\kappa',\lambda)\right) \int_{\Xi} d^{p}\left(\overline{T}^{R}(\kappa',\zeta),\zeta\right) \nu(d\zeta) \ \ = \ \ \theta^{p}.
\end{align*}
Also, it follows from $\Psi$ being concave that
\begin{align}
\E_{T_{\varepsilon}^{R}(\kappa',\lambda,\cdot)_{\#}\nu}[\Psi(\xi)] \ \ & = \ \ \int_{\Xi} \Psi\left(q_{\varepsilon}^{R}(\kappa',\lambda) \underline{T}_{\varepsilon}(\lambda,\zeta) + \left(1 - q_{\varepsilon}^{R}(\kappa',\lambda)\right) \overline{T}^{R}(\kappa',\zeta)\right) \nu(d\zeta) \nonumber \\
& \geq \ \ q_{\varepsilon}^{R}(\kappa',\lambda) \int_{\Xi} \Psi(\underline{T}_{\varepsilon}(\lambda,\zeta)) \nu(d\zeta) + \left(1 - q_{\varepsilon}^{R}(\kappa',\lambda)\right) \int_{\Xi} \Psi\left(\overline{T}^{R}(\kappa',\zeta)\right) \nu(d\zeta) \nonumber \\
& \geq \ \ \kappa' \theta^{p} - q_{\varepsilon}^{R}(\kappa',\lambda) \int_{\Xi} \Phi(\lambda,\zeta) \nu(d\zeta) - q_{\varepsilon}^{R}(\kappa',\lambda) \varepsilon + \left(1 - q_{\varepsilon}^{R}(\kappa',\lambda)\right) \int_{\Xi} \Psi(\zeta) \nu(d\zeta).
\label{eqn:e-optimal perturbation 2}
\end{align}
Thus, given any $\varepsilon > 0$, choose $\lambda_{1}^{\varepsilon} \in \left(\kappa - \varepsilon, \, \kappa\right)$ such that $(\kappa - \lambda_{1}^{\varepsilon}) \theta^{p} \le \varepsilon$, choose $\lambda_{2}^{\varepsilon} \in \left(\kappa, \, \kappa + \varepsilon\right)$ such that $(\lambda_{2}^{\varepsilon} - \kappa) \theta^{p} \le \varepsilon$, choose $\varepsilon' \in \left(0, (\lambda_{2}^{\varepsilon} - \kappa) \theta^{p} - \int_{\Xi} \left[\Phi(\lambda_{2}^{\varepsilon},\zeta) - \Phi(\kappa,\zeta)\right] \nu(d\zeta)\right)$ such that $\varepsilon' \left|\int_{\Xi} \Phi(\lambda_{2}^{\varepsilon},\zeta) \nu(d\zeta)\right| \le \varepsilon$ and $\varepsilon' \left|\int_{\Xi} \Psi(\zeta) \nu(d\zeta)\right| \le \varepsilon$ (note that $\varepsilon' \le \varepsilon$), and choose $R \ge \theta^{p} + \theta^{p} / \varepsilon'$.
Set $T^{\varepsilon} = T_{\varepsilon'}^{R}(\lambda_{1}^{\varepsilon},\lambda_{2}^{\varepsilon},\cdot)$.
Then it follows from~\eqref{eqn:e-optimal perturbation 2} that
\[
\E_{T^{\varepsilon}_{\#}\nu}[\Psi(\xi)] \ \ \geq \ \ \kappa \theta^{p} - \int_{\Xi} \Phi(\kappa,\zeta) \nu(d\zeta) - 5 \varepsilon
\ \ = \ \ v_{D} - 5 \varepsilon
\]
Since $\varepsilon > 0$ can be arbitrarily small, it follows that
\[
\sup\left\{\E_{T_{\#}\nu}[\Psi(\xi)] \; : \; W_{p}(T_{\#}\nu,\nu) \leq \theta,\; T : \Xi \mapsto \Xi \textrm{ is } \nu\textrm{-measurable}\right\} \ \ = \ \ v_{P}
\]

Next, suppose that a worst-case distribution exists.
Then it follows from \ref{itm:iff} that condition~\ref{itm:existenceCond1} or \ref{itm:existenceCond2} or \ref{itm:existenceCond3} holds.
\begin{itemize}
\item
Case~a/b: Condition~\ref{itm:existenceCond1} or \ref{itm:existenceCond2} holds.
\end{itemize}
Let $\underline{T}, \overline{T} : \Xi \mapsto \Xi$ and $q \in [0,1]$ be as in the proof of \ref{itm:iff}.
Let $T^{\ast} : \Xi \mapsto \Xi$ be given by $T^{\ast}(\zeta) \defi q \underline{T}(\zeta) + (1-q) \overline{T}(\zeta)$.
Note that $T^{\ast}(\zeta) \in \Xi$ for all $\zeta \in \Xi$ and that $T^{\ast}$ is $\nu$-measurable.
It follows from $d^{p}(\cdot,\zeta)$ being convex that
\begin{align*}
W_{p}^{p}(T^\ast_{\#}\nu,\nu) \ \ & \leq \ \ \int_{\Xi} d^{p}(q \underline{T}(\zeta) + (1-q) \overline{T}(\zeta),\zeta) \nu(d\zeta) \\
& \leq \ \ q \int_{\Xi} d^{p}(\underline{T}(\zeta),\zeta) \nu(d\zeta) + (1-q) \int_{\Xi} d^{p}(\overline{T}(\zeta),\zeta) \nu(d\zeta)
\ \ = \ \ \theta^{p}.
\end{align*}
Also, it follows from $\Psi$ being concave that
\begin{align*}
\E_{T^{\ast}_{\#}\nu}[\Psi(\xi)] \ \ & = \ \ \int_{\Xi} \Psi(q \underline{T}(\zeta) + (1-q) \overline{T}(\zeta)) \nu(d\zeta) \\
& \geq \ \ q \int_{\Xi} \Psi(\underline{T}(\zeta)) \nu(d\zeta) + (1-q) \int_{\Xi} \Psi(\overline{T}(\zeta)) \nu(d\zeta)
\ \ = \ \ v_{D}.
\end{align*}
This shows that $T^{\ast}_{\#}\nu$ is a primal optimal solution.
\begin{itemize}
\item
Case~c: Condition~\ref{itm:existenceCond3} holds.
\end{itemize}
Let $\underline{T} : \Xi \mapsto \Xi$ be as in the proof of \ref{itm:iff}.
It follows from the proof of \ref{itm:iff} that $\underline{T}_{\#} \nu$ is primal optimal.
\hfillqed

\endproof

\subsection{Proofs for Section~\ref{sec:dual_finite}}

% \begin{lemma}
% \label{lemma:T_epsilon}
% Let $C$ be a Borel set in $\Xi=\R^K$ with nonempty boundary $\partial C$. Then for any $\varepsilon > 0$, there exists a Borel map $T_{\varepsilon} : \partial C \mapsto \Xi \setminus \cl(C)$ such that $d(\xi,T_{\varepsilon}(\xi)) < \varepsilon$ for all $\xi \in \partial C$.
% \end{lemma}

% \proof{Proof of Lemma~\ref{lemma:T_epsilon}.}
%   Observe that the graph
%   \[
%     G_\epsilon = \{(\xi,\zeta)\in\partial C\times (\Xi\setminus\cl(C)):\; 0<d(\xi,\zeta)<\epsilon\}
%   \]
%   is Borel since the distance function is Borel. Thus the measurable selection theorem ((see, e.g. Theorem~18.26 in \citet{aliprantis2006infinite}) ensures the existence of $T_\epsilon$. 
% % Since $\Xi$ is separable, $\partial C$ has a countable dense subset $\{\xi^{i}\}_{i=1}^{\infty}$. For each $\xi^{i}$, there exists $\xi'^{i}\in\Xi\setminus\cl(\Xi)$ such that $0<2d(\xi^{i},\xi'^{i})<\varepsilon$. Thus $\partial C=\bigcup_{i=1}^{\infty} B_{\varepsilon_{i}}(\xi^{i})$, where $B_{\varepsilon}(\xi^{i})$ is the open ball centered at $\xi^{i}$ with radius $\varepsilon$. Define
% % \[
% %   i^{\ast}(\xi)\defi \min_{i\geq0}\{i:\xi\in B_{\varepsilon}(\xi^{i})\},\quad\xi\in\partial C,
% % \]
% % and
% % \[
% %   T_{\varepsilon}(\xi)\defi\xi'_{i^{\ast}(\xi)},\quad\xi\in\partial C.
% % \]
% % Then $T_{\varepsilon}$ satisfies the requirements in the lemma.
% \hfillqed
%   \endproof

\proof{Proof of Proposition~\ref{prop:chance-constrained_C}.}

Consider $\min_{\mu \in \frakM} \mu(\interior(C))$.
It was shown in Example~\ref{eg:UQ} that Corollary~\ref{cor:existence} applies, and thus a worst-case distribution $\mu^{\ast} \in \argmin_{\mu \in \frakM} \mu(\interior(C))$ exists.

First, suppose that $\mu^{\ast}(\interior(C)) = 1$.
Since $1 = \mu^{\ast}(\interior(C)) = \min_{\mu \in \frakM} \mu(\interior(C)) \le \inf_{\mu \in \frakM} \mu(C) \le \mu^{\ast}(C) \le 1$, it follows that $\inf_{\mu \in \frakM} \mu(C) = \min_{\mu \in \frakM} \mu(\interior(C))$ (and $\mu^{\ast} \in \argmin_{\mu \in \frakM} \mu(C)$).

Next, suppose that $\mu^{\ast}(\interior(C)) < 1$, and consider any $\varepsilon > 0$.
Note that for every $\xi \in (C \setminus \interior(C)) \subset \partial C$ it holds that $\{\zeta \in \Xi \setminus C \, : \, 0 < d(\xi,\zeta) < \varepsilon\} \neq \varnothing$, and that the graph
\[
G_{\varepsilon} \ \ \defi \ \ \left\{(\xi,\zeta) \in (C \setminus \interior(C)) \times (\Xi \setminus C) \; : \; 0 < d(\xi,\zeta) < \varepsilon\right\}
\]
is Borel in $\Xi \times \Xi$.
It follows from the measurable selection theorem (see, e.g. Theorem~18.26 in \citet{aliprantis2006infinite}) that there exists a Borel-measurable map $T^{\varepsilon} : \Xi \mapsto \Xi$ such that $T^{\varepsilon}(\xi) = \xi$ for all $\xi \in \interior(C) \cup (\Xi \setminus C)$, and $T^{\varepsilon}(\xi) \in \Xi \setminus C$ such that $0 < d(\xi,T^{\varepsilon}(\xi)) < \varepsilon$ for all $\xi \in C \setminus \interior(C)$.
For any $q^{\varepsilon} \in (0,1)$, let
\[
\mu^{\varepsilon} \ \ \defi \ \ (1 - q^{\varepsilon}) T^{\varepsilon}_{\#}\mu^{\ast} + q^{\varepsilon} \nu.
\]
Let $\gamma^{T\varepsilon} \in \argmin_{\gamma \in \cP(\Xi \times \Xi)} \left\{\int_{\Xi \times \Xi} d^{p}(\xi,\zeta) \gamma(d\xi,d\zeta) \, : \, \pi^{1}_{\#}\gamma = T^{\varepsilon}_{\#}\mu^{\ast}, \pi^{2}_{\#}\gamma = \nu\right\}$.
Let $\nu^{2} \in \cP(\Xi \times \Xi)$ be given by $\nu^{2}(B) \defi \nu(\{\zeta \, : \, (\zeta,\zeta) \in B\})$.
Note that $\pi^{1}_{\#}\nu^{2} = \pi^{2}_{\#}\nu^{2} = \nu$.
Consider $\gamma^{\varepsilon} \defi (1 - q^{\varepsilon}) \gamma^{T\varepsilon} + q^{\varepsilon} \nu^{2}$.
Note that
\begin{align*}
\pi^{1}_{\#}\gamma^{\varepsilon} \ \ & = \ \ (1 - q^{\varepsilon}) \pi^{1}_{\#}\gamma^{T\varepsilon} + q^{\varepsilon} \pi^{1}_{\#}\nu^{2}
\ \ = \ \ (1 - q^{\varepsilon}) T^{\varepsilon}_{\#}\mu^{\ast} + q^{\varepsilon} \nu
\ \ = \ \ \mu^{\varepsilon} \\
\pi^{2}_{\#}\gamma^{\varepsilon} \ \ & = \ \ (1 - q^{\varepsilon}) \pi^{2}_{\#}\gamma^{T\varepsilon} + q^{\varepsilon} \pi^{2}_{\#}\nu^{2}
\ \ = \ \ (1 - q^{\varepsilon}) \nu + q^{\varepsilon} \nu
\ \ = \ \ \nu
\end{align*}
Thus
\begin{align*}
W^{p}_{p}(\mu^{\varepsilon}, \nu) \ \ & = \ \ \min_{\gamma \in \cP(\Xi \times \Xi)} \left\{\int_{\Xi \times \Xi} d^{p}(\xi,\zeta) \gamma(d\xi,d\zeta) \, : \, \pi^{1}_{\#}\gamma = \mu^{\varepsilon}, \pi^{2}_{\#}\gamma = \nu\right\} \\
& \le \ \ \int_{\Xi \times \Xi} d^{p}(\xi,\zeta) \gamma^{\varepsilon}(d\xi,d\zeta) \\
& = \ \ (1 - q^{\varepsilon}) \int_{\Xi \times \Xi} d^{p}(\xi,\zeta) \gamma^{T\varepsilon}(d\xi,d\zeta) + q^{\varepsilon} \int_{\Xi \times \Xi} d^{p}(\xi,\zeta) \nu^{2}(d\xi,d\zeta) \\
& = \ \ (1 - q^{\varepsilon}) \int_{\Xi \times \Xi} d^{p}(\xi,\zeta) \gamma^{T\varepsilon}(d\xi,d\zeta) \\
& = \ \ (1 - q^{\varepsilon}) W^{p}_{p}(T^{\varepsilon}_{\#}\mu^{\ast}, \nu)
\end{align*}
Next, let $\gamma^{\ast} \in \argmin_{\gamma \in \cP(\Xi \times \Xi)} \left\{\int_{\Xi \times \Xi} d^{p}(\xi,\zeta) \gamma(d\xi,d\zeta) \, : \, \pi^{1}_{\#}\gamma = \mu^{\ast}, \pi^{2}_{\#}\gamma = \nu\right\}$.
Consider $\gamma^{T} \in \cP(\Xi \times \Xi)$ given by $\gamma^{T}(B) \defi \gamma^{\ast}\big(\{(\xi,\zeta) \, : \, (T^{\varepsilon}(\xi),\zeta) \in B\}\big)$.
Note that
\begin{align*}
\pi^{1}_{\#}\gamma^{T}(A) \ \ & = \ \ \gamma^{T}(A \times \Xi)
\ \ = \ \ \gamma^{\ast}\big(\{(\xi,\zeta) \; : \; T^{\varepsilon}(\xi) \in A, \zeta \in \Xi\}\big)
\ \ = \ \ \gamma^{\ast}\big((T^{\varepsilon})^{-1}(A) \times \Xi\big) \\
& = \ \ \pi^{1}_{\#}\gamma^{\ast}\big((T^{\varepsilon})^{-1}(A)\big)
\ \ = \ \ \mu^{\ast}\big((T^{\varepsilon})^{-1}(A)\big)
\ \ = \ \ T^{\varepsilon}_{\#}\mu^{\ast}(A) \\
\pi^{2}_{\#}\gamma^{T}(A) \ \ & = \ \ \gamma^{T}(\Xi \times A)
\ \ = \ \ \gamma^{\ast}\big(\{(\xi,\zeta) \; : \; T^{\varepsilon}(\xi) \in \Xi, \zeta \in A\}\big)
\ \ = \ \ \gamma^{\ast}\big(\Xi \times A\big) \\
& = \ \ \pi^{2}_{\#}\gamma^{\ast}(A)
\ \ = \ \ \nu(A)
\end{align*}
Thus $\pi^{1}_{\#}\gamma^{T} = T^{\varepsilon}_{\#}\mu^{\ast}$ and $\pi^{2}_{\#}\gamma^{T} = \nu$.
Hence
\begin{align*}
W^{p}_{p}(T^{\varepsilon}_{\#}\mu^{\ast}, \nu) \ \ & = \ \ \min_{\gamma \in \cP(\Xi \times \Xi)} \left\{\int_{\Xi \times \Xi} d^{p}(\xi,\zeta) \gamma(d\xi,d\zeta) \, : \, \pi^{1}_{\#}\gamma = T^{\varepsilon}_{\#}\mu^{\ast}, \pi^{2}_{\#}\gamma = \nu\right\} \\
& \le \ \ \int_{\Xi \times \Xi} d^{p}(\xi,\zeta) \gamma^{T}(d\xi,d\zeta) \\
& = \ \ \int_{\Xi \times \Xi} d^{p}(T^{\varepsilon}(\xi),\zeta) \gamma^{\ast}(d\xi,d\zeta) \\
& \le \ \ \int_{\Xi \times \Xi} \big[d(T^{\varepsilon}(\xi),\xi) + d(\xi,\zeta)\big]^{p} \gamma^{\ast}(d\xi,d\zeta) \\
& \le \ \ \int_{\Xi \times \Xi} \big[\varepsilon + d(\xi,\zeta)\big]^{p} \gamma^{\ast}(d\xi,d\zeta) \\
& \le \ \ \left(\varepsilon + \left(\int_{\Xi \times \Xi} d^{p}(\xi,\zeta) \gamma^{\ast}(d\xi,d\zeta)\right)^{1/p}\right)^{p} \\
& = \ \ \left(\varepsilon + W_{p}(\mu^{\ast}, \nu)\right)^{p},
\end{align*}
where the last inequality follows from the triangle inequality for the $L^{p}$-norm.
Therefore
\[
W^{p}_{p}(\mu^{\varepsilon}, \nu) \ \ \leq \ \ (1 - q^\varepsilon) \left(\varepsilon + W_{p}(\mu^{\ast}, \nu)\right)^{p}
\ \ \leq \ \ (1 - q^\varepsilon) \left(\varepsilon + \theta\right)^{p}.
\]
Choose $q^{\varepsilon} = 1 - \theta^{p} / (\varepsilon + \theta)^{p}$.
Then $W_{p}(\mu^{\varepsilon}, \nu) \leq \theta$, and thus $\mu^{\varepsilon} \in \frakM$ for all $\varepsilon > 0$.
Also,
\begin{align*}
\min_{\mu \in \frakM} \mu(\interior(C)) \ \ \le \ \ \inf_{\mu \in \frakM} \mu(C) \ \ \le \ \ \mu^{\varepsilon}(C) \ \ & = \ \ (1 - q^{\varepsilon}) T^{\varepsilon}_{\#} \mu^{\ast}(C) + q^{\varepsilon} \nu(C) \\
& = \ \ (1 - q^{\varepsilon}) \mu^{\ast}(\interior(C)) + q^{\varepsilon} \nu(C) \\[1mm]
& \le \ \ (1 - q^{\varepsilon}) \min_{\mu \in \frakM} \mu(\interior(C)) + q^{\varepsilon}
\end{align*}
and therefore it follows that $\inf_{\mu \in \frakM} \mu(C) = \min_{\mu \in \frakM} \mu(\interior(C))$.
\hfillqed

\section{Proofs for Section~\ref{sec:application}}

\subsection{Proofs for Section~\ref{sec:process}}

\proof{Proof of Proposition~\ref{prop:process_policy}.}

It follows from Lemmas~\ref{lem:e-optimal lambda > kappa} and~\ref{lem:e-optimal lambda = kappa} that
\begin{align*}
& \inf_{\mu \in \frakM} \E_{\xi \sim \mu}[\xi(x^{-1}(1))] \\
= \ \ & \inf_{\mu \in \cP(\Xi)} \left\{\E_{\xi \sim \mu}[\xi(x^{-1}(1))] \; : \; \min_{\gamma \in \cP(\Xi^2)} \left\{\int_{\Xi^2} d(\xi,\zeta) \gamma(d\xi,d\zeta) \; : \; \pi^{1}_{\#}\gamma = \mu, \pi^{2}_{\#}\gamma = \nu\right\} \le \theta\right\} \\
= \ \ & \inf_{\substack{\{\xi^{i}_{\pm}\}_{i=1}^{n} \subset \Xi \\ q_{1},q_{2} \in [0,1], \; \{\theta^{i}_{\pm}\}_{i=1}^{n} \subset \R_+}} \begin{multlined}[t] \left\{\frac{1}{n} \sum_{i=1}^{n} \left(q_{1} \xi^{i}_+(x^{-1}(1)) + q_{2} \xi^{i}_-(x^{-1}(1)) + (1 - q_{1} - q_{2}) \hxi^{i}(x^{-1}(1))\right) \; : \right. \\
\left. q_{1} + q_{2} \leq 1, \; d(\xi^{i}_{\pm},\hxi^{i}) \le \theta^{i}_{\pm} \; \forall \; 1 \leq i \leq n, \; \frac{1}{n} \sum_{i=1}^{n} \left(q_{1} \theta^{i}_+ + q_{2} \theta^{i}_-\right) \leq \theta\right\} \end{multlined} \\
= \ \ & \inf_{\substack{q_{1},q_{2} \in [0,1], \\ \{\theta^{i}_{\pm}\}_{i=1}^{n} \subset \R_+}} \begin{multlined}[t] \left\{\frac{1}{n} \sum_{i=1}^{n} q_{1} \bigg(\inf_{\xi_+^{i} \in \Xi} \left\{\xi^{i}_+(x^{-1}(1)) \; : \; d(\xi^{i}_+,\hxi^{i}) \le \theta_+^{i}\right\}
+ q_{2} \inf_{\xi_-^{i} \in \Xi} \left\{\xi^{i}_-(x^{-1}(1)) \; : \; d(\xi^{i}_-,\hxi^{i}) \le \theta_-^{i}\right\}\right. \\
\left. + (1 - q_{1} - q_{2}) \hxi^{i}(x^{-1}(1))\bigg) \ : \ q_{1} + q_{2} \leq 1, \; \frac{1}{n} \sum_{i=1}^{n} \left(q_{1} \theta^{i}_+ + q_{2} \theta^{i}_-\right) \leq \theta\right\}. \end{multlined}
\end{align*}
It follows from condition~\ref{itm:process_dassumption1} that
\[
\inf_{\xi^{i}_{\pm} \in \Xi} \left\{\xi^{i}_{\pm}(x^{-1}(1)) \; : \; d(\xi^{i}_{\pm},\hxi^{i}) \le \theta^{i}_{\pm}\right\}
\ \ = \ \ \inf_{\xi^{i}_{\pm} \in \Xi} \left\{\xi^{i}_{\pm}(x^{-1}(1)) \; : \; W_{1}(\xi^{i}_{\pm},\hxi^{i}) \le \theta^{i}_{\pm}\right\}.
\]
Thus
\begin{align*}
& \inf_{\mu \in \frakM} \E_{\xi \sim \mu}[\xi(x^{-1}(1))] \\
= \ \ & \inf_{\substack{q_{1},q_{2} \in [0,1], \\ \{\theta^{i}_{\pm}\}_{i=1}^{n} \subset \R_+}} \begin{multlined}[t] \left\{\frac{1}{n} \sum_{i=1}^{n} q_{1} \bigg(\inf_{\xi_+^{i} \in \Xi} \left\{\xi^{i}_+(x^{-1}(1)) \; : \; W_{1}(\xi^{i}_+,\hxi^{i}) \le \theta_+^{i}\right\}
+ q_{2} \inf_{\xi_-^{i} \in \Xi} \left\{\xi^{i}_-(x^{-1}(1)) \; : \; W_{1}(\xi^{i}_-,\hxi^{i}) \le \theta_-^{i}\right\}\right. \\
\left. + (1 - q_{1} - q_{2}) \hxi^{i}(x^{-1}(1))\bigg) \ : \ q_{1} + q_{2} \leq 1, \; \frac{1}{n} \sum_{i=1}^{n} \left(q_{1} \theta^{i}_+ + q_{2} \theta^{i}_-\right) \leq \theta\right\} \end{multlined} \\
%= \ \ & \inf_{\substack{\{\xi^{i}_{\pm}\}_{i=1}^{n} \subset \Xi \\ q_{1},q_{2} \in [0,1], \; \{\theta^{i}_{\pm}\}_{i=1}^{n} \subset \R_+}} \begin{multlined}[t] \left\{\frac{1}{n} \sum_{i=1}^{n} \left(q_{1} \xi^{i}_+(x^{-1}(1)) + q_{2} \xi^{i}_-(x^{-1}(1)) + (1 - q_{1} - q_{2}) \hxi^{i}(x^{-1}(1))\right) \; : \right. \\
%\left. q_{1} + q_{2} \leq 1, \; W_{1}(\xi^{i}_{\pm},\hxi^{i}) \le \theta^{i}_{\pm} \; \forall \; 1 \leq i \leq n, \; \frac{1}{n} \sum_{i=1}^{n} \left(q_{1} \theta^{i}_+ + q_{2} \theta^{i}_-\right) \leq \theta\right\} \end{multlined} \\
= \ \ & \inf_{\mu \in \cP(\Xi)} \left\{\E_{\xi \sim \mu}[\xi(x^{-1}(1))] \; : \; \min_{\gamma \in \cP(\Xi^2)} \left\{\int_{\Xi^2} W_{1}(\xi,\zeta) \gamma(d\xi,d\zeta) \; : \; \pi^{1}_{\#}\gamma = \mu, \pi^{2}_{\#}\gamma = \nu\right\} \le \theta\right\} \\
= \ \ & \sup_{\lambda \geq 0} \left\{- \lambda \theta + \frac{1}{n} \sum_{i=1}^{n} \inf_{\xi \in \Xi} \left\{\xi(x^{-1}(1)) + \lambda W_{1}(\xi,\hxi^{i})\right\}\right\} \\
\ge \ \ & \sup_{\lambda \geq 0} \left\{- \lambda \theta + \frac{1}{n} \sum_{i=1}^{n} \inf_{\mu \in \cB([0,1])} \left\{\mu(x^{-1}(1)) + \lambda W_{1}(\mu,\hxi^{i})\right\}\right\},
\end{align*}
where the last equality follows from Theorem~\ref{thm:strongDual}, and the inequality holds because $\Xi \subset \cB([0,1])$.
Note that the inner problem $\inf_{\mu \in \cB([0,1])} \left\{\mu(x^{-1}(1)) + \lambda W_{1}(\mu,\hxi^{i})\right\}$ is similar to the problem considered in Example~\ref{eg:UQ}.
It follows from Example~\ref{eg:UQ} and Proposition~\ref{prop:chance-constrained_C} that
\begin{align*}
\inf_{\mu \in \cB([0,1])} \left\{\mu(x^{-1}(1)) + \lambda W_{1}(\xi,\hxi^{i})\right\}
\ \ & = \ \ \inf_{\mu \in \cB([0,1])} \left\{\mu\big(\interior(x^{-1}(1))\big) + \lambda W_{1}(\mu,\hxi^{i})\right\} \\
& = \ \ \sum_{m=1}^{M_{i}} \min_{\eta \in [0,1]} \left\{\mathds{1}\{\eta \in \interior(x^{-1}(1))\} + \lambda \left|\eta - \heta^{i}_{m}\right|\right\}.
\end{align*}
Hence
\begin{align}
& \inf_{\mu \in \frakM} \E_{\xi \sim \mu}[\xi(x^{-1}(1))] \nonumber \\
\label{eqn:counting vs Borel}
\ge \ \ & \sup_{\lambda \geq 0} \left\{- \lambda \theta + \frac{1}{n} \sum_{i=1}^{n} \sum_{m=1}^{M_{i}} \min_{\eta \in [0,1]} \left\{\mathds{1}\{\eta \in \interior(x^{-1}(1))\} + \lambda \left|\eta - \heta^{i}_{m}\right|\right\}\right\} \\
= \ \ & \inf_{\mu \in \cB([0,1])} \left\{\mu\big(\interior(x^{-1}(1))\big) \; : \; \min_{\gamma \in \cB([0,1]^2)} \left\{\int_{[0,1]^2} \left|\eta - \heta\right| \gamma(d\eta,d\heta) \; : \; \pi^{1}_{\#}\gamma = \mu, \pi^{2}_{\#}\gamma = \hat{\nu}\right\} \le \theta\right\} \nonumber,
\end{align}
where $\hat{\nu} \defi \frac{1}{n} \sum_{i=1}^{n} \sum_{m=1}^{M_{i}} \delta_{\heta^{i}_{m}} \in \cB([0,1])$, and the last equality follows from Theorem~\ref{thm:strongDual}.
This last problem is the same as the problem considered in Example~\ref{eg:UQ}, and it follows that it has an optimal solution $\hat{\mu} \in \cB([0,1])$ of the form
\[
\hat{\mu} \ \ = \ \ \frac{1}{n} \sum\limits_{\substack{1 \leq i \le n, \; 1 \le m \le M_{i} \\ (i,m) \neq (i_{0},m_{0})}} \delta_{\eta^{i}_{m}} + \frac{1 - p_{0}}{n} \delta_{\eta^{i_{0}}_{m_{0}-}} + \frac{p_{0}}{n} \delta_{\eta^{i_{0}}_{m_{0}+}},
\]
and that
\[
W_{1}(\hat{\mu},\hat{\nu})
\ \ = \ \ \frac{1}{n} \sum\limits_{\substack{1 \leq i \le n, \; 1 \le m \le M_{i} \\ (i,m) \neq (i_{0},m_{0})}} \left|\eta^{i}_{m} - \heta^{i}_{m}\right| + \frac{1 - p_{0}}{n} \left|\eta^{i_{0}}_{m_{0}-} - \heta^{i_{0}}_{m_{0}}\right| + \frac{p_{0}}{n} \left|\eta^{i_{0}}_{m_{0}+} - \heta^{i_{0}}_{m_{0}}\right| \ \ \le \ \ \theta.
\]
Consider the $2n$ sample paths $\{\xi^{i}_{\pm}\}_{i=1}^{n} \subset \Xi$ given by
\[
\xi^{i}_{\pm} \ \ \defi \ \ \sum\limits_{\substack{1 \le m \le M_{i} \\ (i,m) \neq (i_{0},m_{0})}} \delta_{\eta^{i}_{m}} + \delta_{\eta^{i_{0}}_{m_{0}\pm}} \mathds{1}\{i = i_{0}\},
\]
and consider the point process $\mu^* \in \cP(\Xi)$ with support on the $2n$ sample paths, given by
\[
\mu^* \ \ = \ \ \frac{1 - p_{0}}{n} \sum_{i=1}^{n} \delta_{\xi^{i}_{-}} + \frac{p_{0}}{n} \sum_{i=1}^{n} \delta_{\xi^{i}_{+}}.
\]
Note that
\begin{align*}
& \E_{\xi \sim \mu^*}\left[\xi\big(\interior(x^{-1}(1))\big)\right] \\
= \ \ & \frac{1 - p_{0}}{n} \left(\sum_{i=1}^{n} \sum\limits_{\substack{1 \le m \le M_{i} \\ (i,m) \neq (i_{0},m_{0})}} \mathds{1}\{\eta^{i}_{m} \in \interior(x^{-1}(1))\} + \mathds{1}\{\eta^{i_{0}}_{m_{0}-} \in \interior(x^{-1}(1))\}\right) \\
& + \frac{p_{0}}{n} \left(\sum_{i=1}^{n} \sum\limits_{\substack{1 \le m \le M_{i} \\ (i,m) \neq (i_{0},m_{0})}} \mathds{1}\{\eta^{i}_{m} \in \interior(x^{-1}(1))\} + \mathds{1}\{\eta^{i_{0}}_{m_{0}+} \in \interior(x^{-1}(1))\}\right) \\
= \ \ & \frac{1}{n} \sum\limits_{\substack{1 \leq i \le n, \; 1 \le m \le M_{i} \\ (i,m) \neq (i_{0},m_{0})}} \mathds{1}\{\eta^{i}_{m} \in \interior(x^{-1}(1))\} + \frac{1 - p_{0}}{n} \mathds{1}\{\eta^{i_{0}}_{m_{0}-} \in \interior(x^{-1}(1))\} + \frac{p_{0}}{n} \mathds{1}\{\eta^{i_{0}}_{m_{0}+} \in \interior(x^{-1}(1))\} \\
= \ \ & \hat{\mu}\big(\interior(x^{-1}(1))\big),
\end{align*}
and that
\begin{align*}
& W_{1}(\mu^*,\nu) \\
\ \ \le \ \ & \frac{1 - p_{0}}{n} \sum_{i=1}^{n} d(\xi^{i}_{-},\hxi^{i}) + \frac{p_{0}}{n} \sum_{i=1}^{n} d(\xi^{i}_{+},\hxi^{i}) \\
= \ \ & \frac{1 - p_{0}}{n} \left(\sum_{i=1}^{n} \sum\limits_{\substack{1 \le m \le M_{i} \\ (i,m) \neq (i_{0},m_{0})}} \left|\eta^{i}_{m} - \heta^{i}_{m}\right| + \left|\eta^{i_{0}}_{m_{0}-} - \heta^{i_{0}}_{m_{0}}\right|\right) + \frac{p_{0}}{n} \left(\sum_{i=1}^{n} \sum\limits_{\substack{1 \le m \le M_{i} \\ (i,m) \neq (i_{0},m_{0})}} \left|\eta^{i}_{m} - \heta^{i}_{m}\right| + \left|\eta^{i_{0}}_{m_{0}+} - \heta^{i_{0}}_{m_{0}}\right|\right) \\
= \ \ & \frac{1}{n} \sum\limits_{\substack{1 \leq i \le n, \; 1 \le m \le M_{i} \\ (i,m) \neq (i_{0},m_{0})}} \left|\eta^{i}_{m} - \heta^{i}_{m}\right| + \frac{1 - p_{0}}{n} \left|\eta^{i_{0}}_{m_{0}-} - \heta^{i_{0}}_{m_{0}}\right| + \frac{p_{0}}{n} \left|\eta^{i_{0}}_{m_{0}+} - \heta^{i_{0}}_{m_{0}}\right| \ \ \le \ \ \theta,
\end{align*}
where the first equality follows from condition~\ref{itm:process_dassumption2}.
That is, $\mu^* \in \frakM$, and $\mu^*$ is an optimal solution for $\inf_{\mu \in \frakM} \E_{\xi \sim \mu}\left[\xi\big(\interior(x^{-1}(1))\big)\right]$.
It follows as in Proposition~\ref{prop:chance-constrained_C} that $\inf_{\mu \in \frakM} \E_{\xi \sim \mu}\left[\xi(x^{-1}(1))\right] = \inf_{\mu \in \frakM} \E_{\xi \sim \mu}\left[\xi\big(\interior(x^{-1}(1))\big)\right]$, and hence inequality~\eqref{eqn:counting vs Borel} holds as an equality. 

Let $\hat{\Xi} \defi \{\heta^{i}_{m} \, : \, i=1,\ldots,N, m=1,\ldots,M_{i}\}$.
Next we show that it suffices to optimize over the set of functions $x : [0,1] \mapsto \{0,1\}$ such that each connected component of $x^{-1}(1)$ contains at least one point in $\hat{\Xi}$.
Consider any $x : [0,1] \mapsto \{0,1\}$ and any connected component $C_{0}$ of $x^{-1}(1)$ such that $C_{0} \cap \hat{\Xi} = \varnothing$.
Then for every $\heta^{i}_{m}$ and every $\lambda \ge 0$, it holds that \[\min_{\eta \in [0,1]} \left\{\mathds{1}\{\eta \in \interior(x^{-1}(1))\} + \lambda \left|\eta - \heta^{i}_{m}\right|\right\} = \min_{\eta \in [0,1] \setminus C_{0}} \left\{\mathds{1}\{\eta \in \interior(x^{-1}(1))\} + \lambda \left|\eta - \heta^{i}_{m}\right|\right\},
\]
and thus 
\[\inf_{\mu \in \frakM} \E_{\xi \sim \mu}\left[\xi\big(\interior(x^{-1}(1))\big)\right] = \inf_{\mu \in \frakM} \E_{\xi \sim \mu}\left[\xi\big(\interior(x^{-1}(1)) \setminus C_{0}\big)\right].
\]
Hence, $v\left(\mathds{1}\{\interior\big(x^{-1}(1)\big) \setminus C_{0}\}\right) \ge v\left(\mathds{1}\{\interior\big(x^{-1}(1)\big)\}\right)$.
Similarly, $v\left(\mathds{1}\{x^{-1}(1) \setminus C_{0}\}\right) \ge v(x)$.
Since $\hat{\Xi}$ is finite, there exists $\{\underline{x}_{j},\overline{x}_{j}\}_{j=1}^{J}$, where $J \leq \textrm{card}(\hat{\Xi})$, such that \eqref{eqn:policy_unionOfIntervals} holds.
\hfillqed
\endproof

\subsection{Proofs for Section~\ref{sec:wcVaR}}

\proof{Proof of Proposition~\ref{prop:CVaR}.}
For each $q$, let $C_{q} \defi \{\zeta \, : \, -w^{\top} \zeta < q\}$.
For any $q > \VaR^{\nu}_{\alpha}[-w^{\top} \xi]$, a necessary condition for $\nu$ to be transported to $\mu \in \cP(\Xi)$ such that $\VaR^{\mu}_{\alpha}[-w^{\top} \xi] = q$ is that $\nu_{w}\{(-\infty,q)\} - (1 - \alpha)$ of $\nu$-probability is transported from $C_{q}$ to $\Xi \setminus C_{q}$.
Note that the metric~$d$ given by $d(\zeta,\xi) = \|\xi - \zeta\|$ is an intrinsic metric, and thus the least cost way to transport probability from $C_{q}$ to $\Xi \setminus C_{q}$ is to transport it to $\partial C_{q}$.
It follows from Example~\ref{eg:UQ} that given $\theta$ and $q$, there exists a worst-case distribution $\mu^{\ast}_{q}$ such that $\mu^{\ast}_{q}(C_{q}) = \min_{\mu \in \frakM} \mu(C_{q})$, and that $\mu^{\ast}_{q}$ is obtained from $\nu$ by greedily transporting probability from points in $C_{q}$ closest to $\partial C_{q}$ to $\partial C_{q}$.
%Next we provide more details.

Specifically, there exists an $\xi^*_{w}$ such that $\|\xi^*_{w}\| = 1$ and $\|w\|_{\ast} \defi \sup\left\{\left|w^{\top} \xi\right| \, : \, \|\xi\| = 1\right\} = w^{\top} \xi^*_{w}$.
Consider any $\zeta \in C_{q}$ and let $s = -w^{\top} \zeta < q$.
Let $\xi_{\zeta} \defi \zeta - (q - s) \xi^*_{w} / \|w\|_{\ast}$.
Note that $-w^{\top} \xi_{\zeta} = -w^{\top} \zeta + (q - s) w^{\top} \xi^*_{w} / \|w\|_{\ast} = s + (q - s) \|w\|_{\ast} / \|w\|_{\ast} = q$, and thus $\xi_{\zeta} \in \Xi \setminus C_{q}$.
Consider any $\xi \in \Xi \setminus C_{q}$, that is, $-w^{\top} \xi \ge q$.
It follows from the definition of $\|w\|_{\ast}$ that $\|\xi - \zeta\| \ge \left|w^{\top} (\xi - \zeta)\right| / \|w\|_{\ast} \ge -w^{\top} (\xi - \zeta) / \|w\|_{\ast} \ge (q - s) / \|w\|_{\ast}$.
Note that $\|\xi_{\zeta} - \zeta\| = (q - s) \|\xi^*_{w}\| / \|w\|_{\ast} = (q - s) / \|w\|_{\ast}$.
Thus $\xi_{\zeta} \in \argmin\left\{\|\xi - \zeta\| \, : \, \xi \in \Xi \setminus C_{q}\right\}$.

Note that for $\zeta \in C_{q}$ it holds that $\min\left\{\|\xi - \zeta\| \, : \, \xi \in \Xi \setminus C_{q}\right\} = \left(q + w^{\top} \zeta\right) / \|w\|_{\ast}$ is increasing in $w^{\top} \zeta$, and thus a worst-case distribution $\mu^{\ast}_{q} \in \argmin_{\mu \in \frakM} \mu(C_{q})$ is obtained from $\nu$ by transporting probability from points $\zeta \in C_{q}$ with the smallest values of $w^{\top} \zeta$ to $\partial C_{q}$.
For any $t < q$, let $\mu_{t}^{q}$ denote the distribution obtained from $\nu$ by transporting probability from points in $\{\zeta \, : \, t \le -w^{\top} \zeta < q\}$ to $\partial C_{q}$, that is, for any $A \in \scrB_{\nu}(\Xi)$, $\mu_{t}^{q}(A) = \nu\left([A \cap C_{t}] \cup \{\zeta \, : \, t \le -w^{\top} \zeta < q, \, \xi_{\zeta} \in A\} \cup [A \cap (\Xi \setminus C_{q})]\right)$.
Note that
\[
\VaR^{\mu_{t}^{q}}_{\alpha}[-w^{\top} \xi] \ \ = \ \ \left\{\begin{array}{ll}
q & \mbox{ if } t \le \VaR^{\nu}_{\alpha}[-w^{\top} \xi] \\ \VaR^{\nu}_{\alpha}[-w^{\top} \xi] & \mbox{ if } t > \VaR^{\nu}_{\alpha}[-w^{\top} \xi]
\end{array}\right..
\]
Also note that $W_{p}^{p}(\mu_{t}^{q},\nu)$ is nonincreasing in~$t$.
Therefore
\[
\VaR^{\frakM}_{\alpha}[-w^{\top} \xi] \ \ = \ \ \sup\left\{q \; : \; t = \VaR^{\nu}_{\alpha}[-w^{\top} \xi], \, W_{p}^{p}(\mu_{t}^{q},\nu) \le \theta^{p}\right\}.
\]
Note that, with $t = \VaR^{\nu}_{\alpha}[-w^{\top} \xi]$, it holds that
\[
W_{p}^{p}(\mu_{t}^{q},\nu) \ \ = \ \ \int_{\VaR^{\nu}_{\alpha}[-w^{\top}\xi]}^{q} \frac{(q - s)^{p}}{\|w\|_{\ast}^{p}} \, \nu_{w}(ds),
\]
and thus
\[
\VaR^{\frakM}_{\alpha}[-w^{\top} \xi] \ \ = \ \ \sup\left\{q \; : \; \int_{\VaR^{\nu}_{\alpha}[-w^{\top}\xi]}^{q} \frac{(q - s)^{p}}{\|w\|_{\ast}^{p}} \, \nu_{w}(ds) \le \theta^{p}\right\}.
\]
It follows from the definition of $\VaR^{\nu}_{\alpha}[-w^{\top} \xi]$ that $\nu_{w}\{(-\infty,q)\} > (1 - \alpha)$ for all $q > \VaR^{\nu}_{\alpha}[-w^{\top} \xi]$, and thus $W_{p}^{p}(\mu_{t}^{q},\nu)$ is increasing in~$q$ on $[\VaR^{\nu}_{\alpha}[-w^{\top} \xi],\infty)$.
In addition, $W_{p}^{p}(\mu_{t}^{q},\nu)$ is continuous in~$q$.
Therefore, $\VaR^{\frakM}_{\alpha}[-w^{\top} \xi]$ is equal to the unique solution~$q$ of
\[
\int_{\VaR^{\nu}_{\alpha}[-w^{\top}\xi]}^{q} (q - s)^{p} \, \nu_{w}(ds) \ \ = \ \ \theta^{p} \|w\|_{\ast}^{p},
\]
and $\mu_{t}^{q}$ with $t = \VaR^{\nu}_{\alpha}[-w^{\top} \xi]$ and the resulting~$q$ is a worst-case distribution.
\hfillqed
\endproof

\section{Selecting Radius \texorpdfstring{$\theta$}{theta}}
\label{sec:theta}

The following approach for selecting the radius of the Wasserstein ball uses a result for Wasserstein distance from \citet{bolley2007quantitative}.
Let $\nu_{N}$ denote the empirical distribution given by $N$ i.i.d.\ observations from the underlying distribution $\nu_{0}$.
In Theorem 1.1 (see also Remark 1.4) of \citet{bolley2007quantitative}, it is shown that
$\Pr\left[W_{1}(\nu_{N},\nu_{0}) > \theta\right] \leq C(\theta) e^{-\lambda N \theta^{2} / 2}$ for some constant $\lambda$ dependent on $\nu_{0}$, and $C$ dependent on $\theta$.
Since their result holds for general distributions, here we simplify it for our purpose by explicitly computing the constants $\lambda$ and $C$.
For a more general analysis, we refer the reader to Section 2.1 in \citet{bolley2007quantitative}.

By assumption $\supp \nu_{0} \subset [0,\bar{B}]$, and thus the truncation step in \citet{bolley2007quantitative} is no longer needed.
Hence the probability bound (2.12) (see also (2.15)) of \citet{bolley2007quantitative} is reduced to
\[
\Pr\left[W_{1}(\nu_N,\nu_{0}) > \theta\right] \ \ \leq \ \ \max\left\{8 e \frac{\bar{B}}{\delta}, \; 1\right\}^{\mathcal{N}(\delta / 2)} e^{-\lambda N (\theta - \delta)^{2} / 8},
\]
for some constants $\lambda > 0$, any $\delta \in (0,\theta)$, where $e$ is the natural logarithm, and $\mathcal{N}(\delta / 2)$ is the minimal number of balls of radius $\delta / 2$ needed to cover the support of $\nu_{0}$.
In this case, $\mathcal{N}(\delta / 2) = \bar{B} / \delta$.
Next we compute $\lambda$.
By Theorem 1.1 of \citet{bolley2007quantitative}, $\lambda$ is the constant in the Talagrand inequality
\[
W_{1}(\mu,\nu_{0}) \ \ \leq \ \ \sqrt{\frac{2}{\lambda} I_{\phi_{kl}}(\mu,\nu_{0})},
\]
where the Kullback-Leibler divergence of $\mu$ with respect to $\nu$ is defined by $I_{\phi_{kl}}(\mu,\nu_{0}) = \infty$ if $\mu$ is not absolutely continuous with respect to $\nu_{0}$, otherwise $I_{\phi_{kl}}(\mu,\nu_{0}) = \int f \log f \, d\nu_{0}$, where $f$ is the Radon-Nikodym derivative $d\mu / d\nu_{0}$.
Corollary~4 in \citet{bolley2005weighted} shows that $\lambda$ can be chosen as
\[
\lambda \ \ = \ \ \left[\inf_{\zeta^{0} \in \Xi, \; \alpha > 0} \frac{1}{\alpha}\left(1 + \log\int e^{\alpha d^{2}(\xi,\zeta^{0})} \nu(d\xi)\right)\right]^{-1},
\]
which can be estimated from data.
Finally, we obtain a concentration inequality
\begin{equation}
\label{eqn:Wassserstein_concentration}
\Pr\left[W_{1}(\nu_N,\nu_{0})>\theta\right] \ \ \leq \ \ \max\left\{8 e \frac{\bar{B}}{\delta}, \; 1\right\}^{\bar{B} / \delta} e^{-\lambda N (\theta - \delta)^{2} / 8}.
\end{equation}
In the numerical experiment, we chose $\delta$ to make the right side of (\ref{eqn:Wassserstein_concentration}) as small as possible, and $\theta$ was chosen such that the right side of (\ref{eqn:Wassserstein_concentration}) was equal to $0.05$.

\end{APPENDICES}

% Acknowledgments here
\section*{Acknowledgments.}
% Enter the text of acknowledgments here
The authors would like to thank Wilfrid Gangbo, David Goldberg, Alex Shapiro, and Weijun Xie for several stimulating discussions, and Yuan Li for providing the image (Figure~(\ref{fig:bird_true})).
% References here (outcomment the appropriate case)

% CASE 1: BiBTeX used to constantly update the references
%   (while the paper is being written).
\bibliographystyle{informs2014} % outcomment this and next line in Case 1
\bibliography{ref_wassersteinMOR} % if more than one, comma separated

% CASE 2: BiBTeX used to generate mypaper.bbl (to be further fine tuned)
%\input{mypaper.bbl} % outcomment this line in Case 2

\end{document}